\documentclass[opre,nonblindrev]{informs3itai}
\pdfoutput=1
\usepackage[english]{babel}
\usepackage{nicefrac}

\usepackage{amsmath}
\usepackage{graphicx}
\usepackage{amssymb}
\usepackage{tabularx}
\usepackage{verbatim}
\usepackage{bbm}
\usepackage{booktabs}
\newcommand{\wbar}{\overline}

\newcommand{\eProof}{\hfill \qed}
\newcommand{\bsq}{\vrule height .9ex width .8ex depth -.1ex}

\usepackage[title]{appendix}

\usepackage{xcolor}

\newcommand{\calS}{\mathcal{S}}
\newcommand{\calC}{\mathcal{C}}
\newtheorem{theorem}{Theorem}
\newtheorem{lemma}{Lemma}
\newtheorem{corollary}{Corollary}

\newtheorem{remark}{Remark}

\newtheorem{definition}{Definition}

\begin{document}

\TITLE{A Hierarchical Approach to Robust Stability of Multiclass Queueing Networks}

\RUNTITLE{Robust Stability of Multiclass Queueing Networks}

\ARTICLEAUTHORS{%
    \AUTHOR{Feiyang Zhao\thanks{University of Texas at Austin, Department of Mechanical Engineering}, Itai Gurvich\thanks{Northwestern University, Kellogg School of Management}, and John Hasenbein\footnotemark[1]}}

\ABSTRACT{We re-visit the global---relative to control policies---stability of multiclass queueing networks. In these, as is known, it is generally insufficient that the nominal utilization at each server is below 100\%. Certain policies, although work conserving, may destabilize a network that satisfies the nominal-load conditions; additional conditions on the primitives are needed for global stability (stability under any work-conserving policy). The global-stability region was fully characterized for two-station networks in \cite{dai1996global}, but a general framework for networks with more than two stations remains elusive. In this paper, we offer progress on this front by considering a subset of non-idling control policies, namely queue-ratio (QR) policies. These include as special cases all static-priority policies. With this restriction, we are able to introduce a complete framework that applies to networks of any size. Our framework breaks the analysis of robust QR stability (stability under any QR policy) into (i) robust state-space collapse and (ii) robust stability of the Skorohod problem (SP) representing the fluid workload. Sufficient conditions for both are specified in terms of simple optimization problems. We use these optimization problems to prove that the family of QR policies satisfies a weak form of convexity relative to policies. A direct implication of this convexity is that: if the SP is stable for all static-priority policies (the ``extreme'' QR policies), then it is also stable under any QR policy. While robust QR stability is weaker than global stability, our framework recovers necessary and sufficient conditions for global stability in specific networks.}
\maketitle
\vspace*{-0.8cm} 
\section{Introduction\label{Section: Introduction}}
\setcounter{equation}{0}
\renewcommand{\theequation}{\thesection.\arabic{equation}}

Stochastic processing networks are often used to model telecommunication networks, as well as production and service processes. In multiclass queueing networks, inputs flow through multiple processing stations, each serving jobs of possibly different types. These jobs, in turn, may differ in their processing requirement (e.g., the processing-time distribution). A control policy specifies to the server which job to serve at any given time. The control policy is chosen so as to satisfy certain criteria or meet prescribed goals. The most basic performance criterion for a policy is stability, which is closely related to throughput maximality and the efficient use of resources. A policy is throughput optimal if---when being used---no queue ``explodes’’ if the nominal utilization of each server is less than 100\%, a requirement that is known as the ``nominal-load condition.''  

\begin{figure}[h!]
\centering
\includegraphics[scale=0.23]{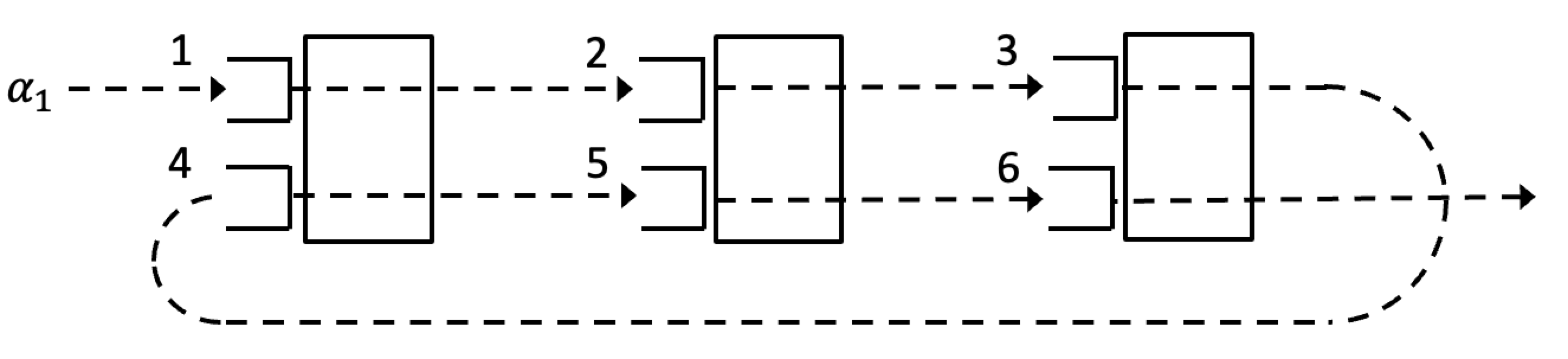}
\caption{DHV network}
\label{fig:DHV Network}
\end{figure}

A network is said to be globally stable if it is stable under {\em any} non-idling policy as long as the nominal-load conditions are met. Global stability is a practically appealing property. It means that decentralized decisions---as are common in large telecommunication networks or in service networks (an emergency room, for example)---are ``safe.'' All single-station networks are globally stable, as are all multiclass queueing networks that are feedforward. The nominal-load conditions, while necessary, are generally insufficient for multiclass queueing networks with feedback. Three well-known examples of this effect occur in the Lu-Kumar and Rybko-Stolyar (also known as Kumar-Seidman) two-station networks \cite{dai1996stability,kumar1990dynamic,rybko1992ergodicity} and the three-station network studied by Dai, Hasenbein, and Vande Vate (DHV) in \cite{DHV} and depicted in Figure \ref{fig:DHV Network}; see \cite{bramson2008stability} for a survey of stability results and methods. 

\begin{figure}[h!]
\centering
\includegraphics[scale=0.06]{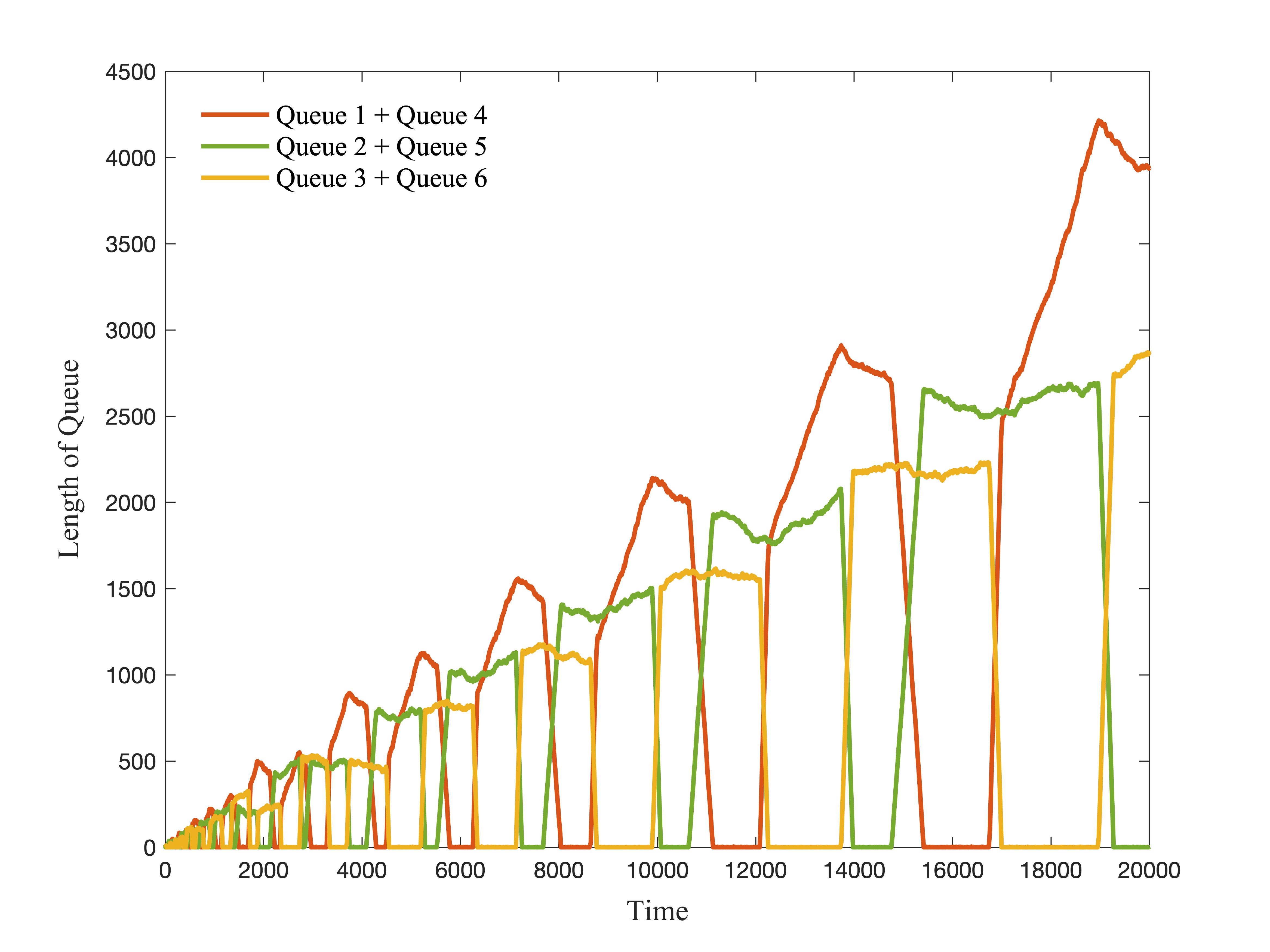}
\caption{A simulation of a DHV network with class 2, 4, and 6 having high priority, and with $\rho<1$.}
\label{fig:Simulation of DHV Network}
\end{figure}

In Figure \ref{fig:Simulation of DHV Network} we plot a single sample path in the simulation of the DHV network in Figure \ref{fig:DHV Network}. In this simulation, classes 2, 4, and 6 have priority at their respective stations. We observe that, as time progresses, not only do the individual queues oscillate but, also, the total queue length grows and ``blows up'' as the simulation progresses: the queueing network is unstable. This is so even though the network parameters satisfy the {\em nominal load} conditions, i.e., $\rho_1 = \alpha_1(m_1+m_4)<1$, $\rho_2 = \alpha_1(m_2+m_5)<1$, $\rho_3 = \alpha_1(m_3+m_6)<1$, where $\alpha_1$ is the exogenous arrival rate and $m_k$ is the mean service time of class-$k$ jobs. Under the nominal-load conditions, there {\em exists} a policy that makes the network stable---for instance, the policy that gives priority to classes $1,2$ and $3$---but the network is not {\em globally} stable: there exist non-idling policies, one of which is policy giving high priority to classes 2, 4, and 6, that destabilize the otherwise stabilizable network. 

For the network to be globally stable, additional conditions must be imposed on the parameters. In the DHV network with the 2-4-6 priority policy, these classes form a so-called \textit{pseudostation}: when server $2$ is working on queue $2$, queue $5$ is not being served and there is no input to a (possibly empty) queue $6$. These processing conflicts between classes 2, 4, and 6 suggest that adding the requirement
$\alpha_1(m_2+m_4+m_6)<2$ to the nominal-load conditions suffices for the stability of the high-priority classes of the network under this policy when the network is in heavy traffic; this is proved in \cite{chen2001existence}. \footnote{Generally, this phenomenon occurs whenever two or more classes of jobs block each other when given high priority in their respective stations. The set of these classes are said to constitute \textit{virtual station} if only one job in the set can be processed at any given time. The set is called a \textit{pseudostation} if the number of classes that can be processed simultaneously is less than the number of stations (but not necessarily equal to one); see \cite{bramson2008stability} and \cite{hasenbein1997necessary}. }

Global stability is one notion of robustness. A network that is globally stable is robust to the choice of customer prioritization by individual servers. Global fluid stability was fully characterized for two-station networks in \cite{dai1996global} and studied in an ad-hoc manner for specific three-station networks; see, e.g., \cite{DHV,dumas1997multiclass}. There is, as of now, no framework completely characterizing global stability for multiclass queueing networks with more than two stations. 

This paper describes progress in this direction. We introduce a useful, but more restricted, notion of global stability, and a novel framework for the derivation of sufficient conditions for this stability.
Instead of allowing for arbitrary non-idling policies, we restrict attention to the family of so-called {\em queue-ratio} (QR) policies; see \cite{gurvich2010service}. This is a large family that includes static-priority policies and max-weight as special instances, but is more flexible. The main restriction is that the scheduling decisions at a given station/server are made based on the queues at that station and not the queues at other stations. 

A queue-ratio policy is characterized by its {\em ratio matrix} which specifies a (targeted) proportional allocation of a station's workload to each of the classes. Consider, as an example, the DHV network in Figure \ref{fig:DHV Network} and station 1 there; suppose that $m_1=m_4$. The first column of the ratio matrix specifies the targeted proportional allocation of station 1's workload to classes 1 and 4. If, for instance, the entries in that column are $1/2$ for both classes $1$ and $4$, the policy (localized to station 1) seeks to equalize queue 1 and 4; it is the longest queue policy. The space of queue-ratio policies is represented by the polyhedron that defines feasible ratio matrices. The corners of this polyhedron are the static-priority policies, where all the jobs of a station are concentrated in one class: the lowest priority class; the formal description appears in \S \ref{Sec:Mathematical Background}. 

We study the fluid model of queue-ratio policies. The stability of the fluid model guarantees that of the queueing network. Intuitively, the fluid model of queue-ratio policies should exhibit a form of {\em state-space collapse} (SSC). After some time (that may depend on the initial state), the workload in each station should be distributed to this station's queues following the targeted proportions encoded in the ratio matrix. At that point in time and onward, the fluid workload (whose dimension is the number of stations) has all the information needed to infer the queue trajectories. 

This supports a hierarchical view of robust stability. A network is robustly queue-ratio stable if {\em for any ratio matrix} (i) the fluid model exhibits SSC, and (ii) the (collapsed) fluid workload, whose trajectories are characterized by a {\em Skorohod Problem} (SP), empties in finite time. If both are achieved, then---for any ratio matrix---the fluid queues empty in finite time, which subsequently guarantees the stability of the network. Because queue-ratio policies are specified by the ratio matrix, robust queue-ratio stability is conceptually reduced to robustness relative to that matrix. 

We draw on the literature to formulate algebraic sufficient conditions for both requirements: SSC and SP stability {\em for a given ratio matrix}. These algebraic conditions correspond to the drift of suitable linear Lyapunov functions. For the SP, the algebraic condition is one we label as Chen-$\calS$ inspired by \cite{chen1996sufficient}; for SSC, we require a linear Lyapunov function reminiscent of the $L^2$-condition in \cite{chen2000sufficient}.

Robust queue-ratio stability then corresponds to the ``robustness'' of these algebraic conditions to the ratio matrix. This connection reformulates robust stability in the language of robust optimization: given a network and its parameters, it is possible to test for the robust stability of the SP by solving a robust optimization problem in which the uncertainty set is the family of ratio matrices. 

Our main result is a weak form of convexity for SP stability: we show that if the SP is stable under two ratio matrices that differ in only one column (corresponding to the class proportions in one station), it is also stable under any convex combination of these two matrices. The important implication is that if the SP is stable under static priority policies, it is stable under {\em all} queue-ratio policies. 

We establish a conceptually similar result for robust SSC. We show that a linear Lyapunov function exists for SSC under {\em any} queue-ratio policy, provided that one exists for each of the static priority policies. 

In \cite{dai_vate_twostation}, the authors established that---for two-station networks---the fluid model's global stability region is the intersection of the stability regions of all static priority policies. While not fully characterizing the stability region, our results do expand the sufficiency of static priority policies---as corners of the space of queue-ratio policies---for the stability of all queue-ratio policies.

Theorem \ref{thm: full_robust_stability} provides sufficient conditions for robust QR stability. It combines the two key ingredients---stated in Theorems \ref{Thm: full statement for corner->inter} and \ref{thm: LEGO paper (Theorem 2)}---of robust SP stability and robust SSC. The SP, for a static priority policy, is characterized by its drift vector $\theta_{\Delta}$ and its so-called reflection matrix $R_{\Delta}$. The formal definition of Chen-$\calS$ is in Definition \ref{def:Chen-S}; and that of linear attraction for SSC appears in Definition \ref{def: linear attraction}. Our early statement of the main theorem serves as a roadmap for what follows.

\vspace{-0.2cm}
\begin{theorem}
\label{thm: full_robust_stability}
A network is robustly queue-ratio stable if,
\begin{itemize} 
\item[(i)] {\bf Robust SP Stability.} For any static priority policy, the reflection matrix $R_{\Delta}$ of its SP is invertible and the SP with data $(R_\Delta,\theta_\Delta)$ is Chen-$\calS$. As a collection, these reflection matrices all have positive determinants or all have negative determinants.
\item[(ii)] {\bf Robust SSC.} Linear attraction holds with  a common test vector.
\end{itemize} 
\end{theorem}

\begin{figure}[h!]
\centering
\includegraphics[scale=0.6]{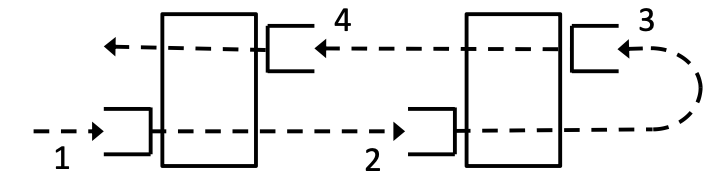}
\caption{The two-station and four-class Lu-Kumar network}
\label{fig:Lu-Kumar Network}
\end{figure}

We use our framework to revisit three well-known examples from the literature. The Lu-Kumar network in Figure \ref{fig:Lu-Kumar Network} is a classical illustration of the challenges in global stability. Here, our framework recovers the sufficient condition $\alpha_1(m_2+m_4)<1$ known from the literature \cite{dai1996stability} in the case that both stations are bottlenecks. The same is true for the two-station network in \cite{dai2004stability}. 

For the DHV network in Figure \ref{fig:DHV Network} our framework shows that if all stations are bottlenecks (a balanced network), it suffices to augment the nominal-load conditions ($\alpha_1(m_1+m_4)<1$, $\alpha_1(m_2+m_5)<1$, and $\alpha_1(m_3+m_6)<1$) with the requirement $\alpha_1(m_2+m_4+m_6)<2$ to have robust queue-ratio stability. Interestingly, no requirements on the mean service times are needed for robust SP stability. Rather, the condition $\alpha_1(m_2+m_4+m_6)<2$ arises entirely from the analysis of robust SSC. The gap between our condition and the necessary and sufficient conditions in \cite{DHV} is consistent with our focus on decentralized policies. The policy used in \cite{DHV} to establish their necessary condition is a centralized one, in which the priority in a station depends on the queues in other stations. From a practical perspective, if centralized policies are implementable, it is best to use those that are known to stabilize the network under the nominal-load conditions. 
State space collapse and the SP are central constructs in the study of diffusion approximations \cite{bramson1998SSCheavylimit,williams1998diffusion}. Through the examples, our framework casts existing results in new light and helps underscore interesting connections between robust stability, fluid models, and diffusion approximations. 

\subsection{Literature Review}

In addition to the literature surveyed in our introduction, several papers are important for broader context. That the stability of a network's fluid model---its attraction to the origin in finite time---implies positive (Harris) recurrence (stability) of the underlying stochastic queueing network, has been known since the seminal work in {\cite{dai1995positive,rybko1992ergodicity,stolyar1995stability}.

The observation that SSC in the fluid model and SP stability, combined, are sufficient for fluid model stability was made first in \cite{delgado2010state}. This observation creates a bridge between questions of stability and those pertaining to central-limit-theorem (or diffusion) approximations for networks in heavy traffic. In \cite{bramson1998state,bramson1998SSCheavylimit}, Bramson proved that uniform attraction in the fluid model implies SSC in the diffusion model. Williams subsequently proved in \cite{williams1998diffusion} heavy-traffic limit theorems for multiclass queueing networks. Two key requirements for these limit theorems are (i) SSC and (ii) that the reflection matrix of the SP is completely-$\calS$. 

State-space collapse (SSC) under static-priority policies plays a key role in our results. The papers \cite{chen2001existence,chen2000sufficient} study the existence of diffusion approximations in multiclass queueing networks with priority service disciplines. The sufficient condition we derive for the robust queue-ratio stability of the DHV network is the same as the one identified in \cite{chen2001existence} as sufficient for SSC under the 2-4-6 priority policy. SSC and SP stability are central to the so-called limit-interchange problem; see \cite{gurvich2014validity} under queue-ratio policies, and \cite{ye2018justifying} for a more general framework.

The stability of the Skorohod problem is the other key ingredient of our work. \cite{bernard34kharroubi} showed that having a reflection matrix that is completely-$\calS$ is necessary and sufficient for the existence of solutions of the SP. In \cite{chen1996sufficient}, Chen provides sufficient conditions for the stability of the SP. We use his linear Lyapunov functions to study the SP's {\em robust} stability. 

The emergence of fluid models as a tool to study stability gave a boost to the study of global stability. \cite{dai1996stability} studies the global stability of the two-station Lu-Kumar network.
\cite{down1997piecewise} introduces a linear programming problem whose solution, if positive, guarantees that the network is stable under any non-idling policy. In \cite{dai_vate_twostation}, a complete characterization of global (fluid) stability for 2-station networks is provided, and it is shown that the global stability region is the intersection of all static priority stability regions in 2-station networks. A specific three-station network, the one in Figure \ref{fig:DHV Network}, is studied in  \cite{chen2001existence}, \cite{chen2002piecewise}, and \cite{DHV}. See the surveys \cite{bramson2008stability} and \cite{dai2020processing} for further coverage of global stability.

\vspace{0.2cm}

{\bf Organization and notation.} In \S 2, we describe the mathematical primitives: multiclass queueing networks, queue-ratio (QR) policies, fluid models, state-space collapse (SSC), and the Skorohod problem (SP). \S 3 analyzes robust stability of the SP. In \S 4, we study robust SSC. Throughout we revisit networks that were studied in the literature. All theorems are proved in \S \ref{sec:mainproofs} and all supporting lemmas are proved in the appendix. 

\vspace{0.2cm}

Throughout the paper, vectors are interpreted as column vectors. The transpose of a vector or matrix $v$ is $v'$ or $v^\top$. Given a subset $a \subseteq \{1,..., d\}$, $u_a$ stands for the sub-vector of $u$ corresponding to indices in the set $a$. Given a matrix $R$, $R_{\cdot j}$ stands for the $j^{th}$ column of $R$ and $R_{j\cdot}$ stands for the $j^{th}$ row of $R$. The sub-matrix $R_{ab}$ of $R$ is obtained by taking the rows in the index set $a$ and the columns in the index set $b$; we write $R_a$ for the square sub-matrix $R_{aa}$. We use $det(R)$ to denote the determinant of the matrix $R$. The vector $e$ has all elements equal to 1 and $I$ is the identity matrix. For a set $b$, $|b|$ stands for the cardinality of the set. Finally, $\mathbbm{1}\{A\} = 1$ if $A$ is true, and $\mathbbm{1}\{A\} = 0$ otherwise.

\section{Mathematical Background}
\label{Sec:Mathematical Background}

\setcounter{equation}{0}
\renewcommand{\theequation}{\thesection.\arabic{equation}}

\noindent {\bf Network primitives.} We adopt standard queueing-network notation; see, e.g., \cite{williams1998diffusion}. There are $J$ servers (or stations) and $K$ classes; in multiclass queueing networks $K>J$. The many-to-one mapping from classes to stations is described by a $J\times K$ \textit{constituency matrix} $C$ where $C_{jk} = 1$ if class $k$ is served at station $j$ (in which case we write $s(k) = j$), and $C_{jk} = 0$ otherwise. The matrix $C$ indicates which station serves which class. For $j=1,\ldots,J$, we denote by  $\mathcal{C}(j)$ the set of classes served by server $j$. We denote by $M$ the $K\times K$ diagonal matrix with the mean service time of class $k$, $m_k >0$, as the $k^{th}$ diagonal element; $\mu_k = 1/m_k$ is the service rate at which class $k$ jobs are processed when server $s(k)$ is working on this class’s jobs. Customer routing is Bernoulli and captured by a \textit{routing matrix} $P$ in which $P_{kl}$ is the probability that a class $k$ job becomes a class $l$ job upon completion of service. We assume that the network is open, namely that $I+P+P^2+\ldots$ is convergent and the matrix $(I - P')$ is invertible; in this case, $(I - P')^{-1} = (I+P+P^2+\ldots)'$. We let $Q = (I - P')^{-1}$, so that  $Q_{kl}$ is the expected number of visits to class $k$ starting from class $l$. Exogenous arrivals to class $k$ follow a (possibly delayed) renewal process with rate  $\alpha_k \geq 0$;  $\alpha = (\alpha_1,\ldots,\alpha_K)'$. The effective arrival rate to class $k$ is denoted by $\lambda_k$; $\lambda = Q\alpha = (\lambda_1,\ldots,\lambda_K)'$. The traffic intensities are expressed as $\rho = CM\lambda$, where $\rho_j$ is the intensity of traffic for station $j$. The {\em nominal-load conditions} are $\rho_j <1$ for all $j \in \{1,\ldots,J\}$. We refer to $(C, M, P, \alpha)$ as the network primitives. These describe the structure of the multiclass queueing network.

On the stochastic primitives---service times and renewal arrival processes---we impose the independence and distributional assumptions in \cite{dai1995positive} that allow deducing network stability from the fluid-model stability; see (1.2)-(1.5) there. It is a standard argument to show that subsequential limits of the fluid-scaled processes converge to the fluid model equations specified in full in \eqref{equ:begin fluid model}-\eqref{equ:end fluid model} below. These facts allow us to apply the stability results in \cite{dai_vate_twostation} and deduce stability of the queueing network from that of the fluid model; see the monograph \cite{bramson2008stability} for a thorough review of the fluid limits and their connection to stability.

\vspace{0.2cm}

\paragraph{\bf Queue-Ratio Policies.} Let $Z_k(t)$ be the number of class-$k$ jobs in the system at time $t$, and $W_j(t)$ be the workload in station $j$ at time $t$, defined by $W_j(t) = \sum_{k\in \mathcal{C}(j)}m_kZ_k(t)$. Under a preemptive queue-ratio policy, server $j$ serves at time $t$ the customer at the head of a class-$k$ queue where 
$$k\in \{\ell \in \mathcal{C}(j):Z_l(t)-\delta_lW_j(t) >0 \}.$$
Here $\delta_{\ell}\times m_{\ell}$ is the proportion of station-$j$'s workload ``targeted'' for queue $\ell$. The weights $\delta_k,k=1,\ldots,K$ are collected in a \textit{ratio matrix} $\Delta$ that has $\Delta_{kj} = \delta_k$ if class $k$ is served in station $j$ and $\Delta_{kj} = 0$ otherwise; $\Delta$ has $J$ columns and $K$ rows. The proportions $\delta_{\ell}\times m_{\ell},\ell\in\calC(j)$ satisfy $\sum_{\ell\in\calC(j)}  \delta_{\ell} m_{\ell} = 1$; in vector notation $CM\Delta=I$.  

To complete the description of the policy, one must specify the choice of the class to be served when the set  $\{\ell \in \mathcal{C}(j):Z_l(t)-\delta_lW_j(t) >0 \}$ contains more than one class. One natural choice is the class with the maximal imbalance $Z_l(t)-\delta_lW_j(t) >0$, but one could also define a fixed priority order within this set; a queue-ratio policy is fully specified by the matrix $\Delta$ and the ``tie-breaking'' rule. Because $CM\Delta=I$, $0 = CM(Z-\Delta W) = \sum_{\ell\in\mathcal{C}(j)}m_\ell(Z_\ell(t)-\delta_\ell W_j(t))$, not all the terms $(Z_\ell (t)-\delta_\ell W_j(t), \ell\in\mathcal{C}(j))$ can be positive, therefore $|\{\ell \in \mathcal{C}(j):Z_l(t)-\delta_lW_j(t) >0 \}|<|\calC(j)|$. In the DHV network, where each station has only two classes, this set can contain at most one class at each time $t$. 

The QR family of policies contains static priority policies as special instances; for these, the matrix $\Delta$ takes on a special form.

\begin{definition}[Static-priority $\Delta$]
\label{defin:staticpriorityDelta}
  The ratio matrix is said to be a static-priority $\Delta$ if $\Delta_{kj} = 1/m_k$ for some $k\in\calC(j)$ (the lowest priority class) and $\Delta_{lj} = 0$ for all other $l \in \mathcal{C}(j)$.
\end{definition}

Thus, any static-priority policy can be implemented via a suitable choice of a ratio matrix $\Delta$. The max-weight policy with quadratic exponents is, as well, a queue-ratio policy. In this policy, the class in service in station $j$ at time $t$ is chosen from the set $\argmax_{i\in\calC(j)} c_i\mu_iZ_i^2(t)$. This policy is equivalent to the QR policy with $$\delta_i=\frac{\nicefrac{m_i}{c_i}}{\sum_{l\in\calC(j)}\nicefrac{m_l^2}{c_l}}.$$ 
This equivalence is demonstrated in the proof of \cite[Theorem 3.2]{gurvich2009scheduling}.

\vspace{0.2cm}

\paragraph{\bf Fluid model: State-Space Collapse and the Skorohod Problem.} We overline a symbol if it is related to the fluid model (vs.\ the queueing network). For instance, $Z_k(t)$ represents the number of class-$k$ jobs in the system at time $t$, and $\wbar{Z}_k(t)$ represents the fluid level of class $k$ in the corresponding fluid model at time $t$.

Let $\wbar{Y}_j(t)$ be the cumulative amount of time (in the fluid model) that the server in station $j$ is idle by time $t$, and let $\wbar{T}_k(t)$ be the cumulative time station $s(k)$ has allocated to class-$k$ fluid by time $t$, so that $\wbar{Y}_j(t)=t-\sum_{k\in\calC(j)}\wbar{T}_k(t)$. $\wbar{T}(\cdot)$ depends on the the policy followed by the network; through $\wbar{T}(\cdot)$, the rest of the network dynamics are determined. For station $j$ and class $k\in\calC(j)$, $\wbar{\epsilon}_k(t)=\wbar{Z}_k(t)-\delta_k \wbar{W}_j(t)$ captures the distance of the fluid level of class $k$ from the targeted  proportion of workload at time $t$. Server $j$ prioritizes queues $k$ with  $\wbar{\epsilon}_k(t) >0$ over those with $\wbar{\epsilon}_k(t) <0$. Under a static priority policy, for all but the lowest priority class in a station, $\delta_{k}=0$ so that $\wbar{\epsilon}_k(t) = \wbar{Z}_k(t) \geq 0$ at all times $t$. 

The fluid model equations for a queue-ratio policy are given by 
\begin{align}
     &\wbar{Z}(t) = \wbar{Z}(0)+\alpha t - (I-P')M^{-1}\wbar{T}(t), \label{equ:begin fluid model}\\
     & \wbar{W}(t) = CM\wbar{Z}(t), \\
     & \wbar{\epsilon}(t) = \wbar{\epsilon}(0) + (I-\Delta CM)(\alpha t - (I-P')M^{-1}\wbar{T}(t)), \\
     & \wbar{T}(0) = 0 \text{ and } \wbar{T}(t) \text{ is nondecreasing}, \\
     & \wbar{Y}(t) = et-C\wbar{T}(t) \text{ is nondecreasing}, \\
     & \int_0^\infty \wbar{W}(s)d\wbar{Y}_j(s) = 0, j\in \mathcal{J}, \\
     & \int_0^t\wbar{\epsilon}_k^{+}(s)d(s-\wbar{T}_k^{+}(s)) = 0, k\in \mathcal{K}, \label{equ:end fluid model}
\end{align} where $\wbar{T}_k^{+} = \sum_{i\leq k:s(i)=s(k), \wbar{\epsilon}_i(t) >0}\wbar{T}_i(t)$, $\wbar{\epsilon}_k^+(t) = \sum_{i\leq k:s(i) = s(k), \wbar{\epsilon}_i(t) >0}\wbar{\epsilon}_i(t)$, when assuming, without
loss of generality, a tie-breaking rule that a server chooses first the class with the lowest index among the ones with $\wbar{\epsilon} >0$.

\begin{definition}[State-Space Collapse (SSC)]
\label{def:SSC}
The network satisfies state-space collapse if there exists $t_0\geq 0$ (possibly dependent on the initial condition $\wbar{Z}(0)$) such that $\wbar{Z}(t) = \Delta \wbar{W}(t)$ for all $t \geq t_0$. 
\end{definition}

Delgado \cite{delgado2010state} made a useful observation relating the stability of a network's fluid model to state-space collapse and the stability of a Skorohod problem, which we review and use here. Suppose that $CMQ\Delta$ is an invertible matrix and let $R = (CMQ\Delta)^{-1}$. Consider a time $t\geq t_0$ (where $t_0$ is as in Definition \ref{def:SSC}). Using the fluid model equations we obtain
\begin{align*}
         & \wbar{Z}(t) = \wbar{Z}(t_0)+\alpha (t-t_0) - (I-P')M^{-1}(\wbar{T}(t)-\wbar{T}(t_0)),
         \end{align*} so that, because $\wbar{Z}=\Delta\wbar{W}(t)$ for all $t\geq t_0$, \begin{align*}  & CMQ\wbar{Z}(t) = CMQ\wbar{Z}(t_0)+CMQ\alpha (t-t_0) - CMQ(I-P')M^{-1}(\wbar{T}(t)-\wbar{T}(t_0))\\
     \stackrel{\wbar{Z}(t)=\Delta \wbar{W}(t)}{\Leftrightarrow} \;\;&  CMQ\Delta\wbar{W}(t) = CMQ\Delta\wbar{W}(t_0)+CMQ\alpha (t-t_0) - CMQQ^{-1}M^{-1}(\wbar{T}(t)-\wbar{T}(t_0))\\
    \Leftrightarrow \;\;&  R^{-1}\wbar{W}(t)= R^{-1}\wbar{W}(t_0)+CMQ\alpha (t-t_0) - C(\wbar{T}(t)-\wbar{T}(t_0))\\
    \Leftrightarrow \;\;&  R^{-1}\wbar{W}(t) = R^{-1}\wbar{W}(t_0)+CM\lambda (t-t_0) - e(t-t_0) +(\wbar{Y}(t)-\wbar{Y}(t_0))\\
    \Leftrightarrow \;\;&  \wbar{W}(t) = \wbar{W}(t_0)+R(\rho- e)(t-t_0) +R(\wbar{Y}(t)-\wbar{Y}(t_0)),
\end{align*}
where, recall, $\wbar{Y}_j$ can increase only at times $t$ where $\wbar{W}_j(t)=0$. 
The final row in this derivation means that, from time $t_0$ and onward, $\wbar{W}(t)$ is a solution to the so-called Skorohod problem as formally defined below, with $\theta = R(\rho-e)$.  
\vspace{-0.1cm}
\begin{definition}[Skorohod Problem (SP)]
\label{def:SP}
Given a vector $\theta \in \Bbb{R}^J$ and $\wbar{W}(0) \geq 0$, the pair $(\wbar{Y}, \wbar{W}) \in \calC^{J\times J}$ solves the linear Skorohod problem if 
$$\wbar{W}(t) = \wbar{W}(0) + \theta t + R\wbar{Y}(t),$$
where $\wbar{Y}(0) = 0$ and $\wbar{Y}(t)$ is non-decreasing and satisfies the complementarity condition
\begin{equation*}
    \int_0^\infty \wbar{W}_j(t)d\wbar{Y}_j(t) = 0 \text{, for all } j \in \mathcal{J}.
\end{equation*}
This Skorohod problem is stable if for any solution $(\wbar{Y}, \wbar{W})$ there exists a $T \in [0,\infty)$ such that
$$\wbar{W}(t) = 0\;\;\;\forall t\geq T. $$
\end{definition}
\vspace{-0.1cm}

It is evident from the above derivations that if (i) the fluid model admits state-space collapse, and (ii) the Skorohod problem is stable, then the fluid model reaches $0$ in finite time.

To discuss stability of the Skorohod problem we require additional ingredients. 
\vspace{-0.1cm}
\begin{definition}[Completely-$\calS$]
\label{def:completelyS}
 A square matrix $R$ is an S-matrix if there exists a positive vector u such that $Ru > 0$. It is a completely-$\calS$ matrix if all of its principal sub-matrices are $\calS$-matrices.
\end{definition}
\vspace{-0.1cm}

That the reflection $R$ is completely-$\calS$ is a necessary and sufficient condition for the existence of solutions to the Skorohod problem (see \cite{bernard34kharroubi}). Chen \cite{chen1996sufficient} establishes sufficient conditions on the  matrix $R$ and the vector $\theta$ that guarantee the stability of the SP. Theorem \ref{thm: Chen-S (Theorem 2.5)} below is a re-statement of Theorem 2.5 in \cite{chen1996sufficient}.
\vspace{-0.1cm}
\begin{definition}[Chen-$\calS$]
\label{def:Chen-S}
Let $R=(CMQ\Delta)^{-1}$ and  $\theta = R(\rho-e)$. We say that the pair $(R, \theta)$ is Chen-$\calS$ if $R$ is completely-$\calS$ and there exists a positive vector $h \in{\mathbb{R}^J}$  such that 
\begin{align}
h_a^{'}\left[\theta_a + R_{ab}u\right] <0,
\label{equ:ChenS_def}
\end{align} for any partition $(a,b)$ of $\mathcal{J}$, and all $u \in \left\{ v\in R_{+}^{|b|}: \theta_b + R_bv=0 \right\}$.
\end{definition}
\vspace{-0.1cm}
\begin{theorem}[Theorem 2.5 in \cite{chen1996sufficient}]
\label{thm: Chen-S (Theorem 2.5)}
The Skorohod problem is stable if $(R, \theta)$ is Chen-$\calS$.
\end{theorem}

Given $\Delta$ we write $R_{\Delta}$ and $\theta_{\Delta}$ to make explicit their dependence on the ratio matrix. 
When $\Delta$ is clear from the context and fixed, $R_{ab}=(R_{\Delta})_{ab}$ is the sub-matrix of $R$ that has the rows in the set $a$ and the columns in the set $b$. The Chen-$\calS$ requirement, which provides a sufficient condition for the stability of the SP, is a central building block in our analysis: robust SP stability depends on whether (or not) the pair $(R_\Delta,\theta_\Delta)$ is Chen-$\calS$ for all possible ratio matrices $\Delta$. The algebraic formulation of Chen-$\calS$ allows us to represent robust SP stability (more precisely robust Chen-$\calS$) through a robust optimization problem, in which the uncertainty set is the set of ratio matrices $\Delta$.
\begin{lemma}
\label{lem: RO_ChenS}
The SP is robustly queue-ratio stable if $\rho <e $ and the min-max problem
\begin{equation}
\label{RO:ChenS}
\begin{split}
{\min_{\Delta \in \Bbb{R}_{+}^{K\times J}: CM\Delta = I} \;\; \max_{h,\alpha^b,\epsilon, \eta^b}}\quad\epsilon\\
s.t.\quad &h_a^\top \theta_a - \alpha^{b\top} \theta_b \leq -\epsilon \qquad \forall (a,b)\\
& h_a^\top R_{ab} \leq \alpha^{b\top}R_b  \qquad \forall (a,b) \\
& \epsilon e\le R_b\eta^b \qquad \forall  b \subseteq J\\
& \alpha^b \in \mathbb{R}^{|b|}, h\geq e, \epsilon \leq M , \eta^b \geq e^b
\end{split}
\end{equation}
has a strictly positive solution, where $M$ is a large enough positive number.
\end{lemma}


Before proceeding, it is worthwhile relating Chen-$\calS$ to another well-known, but stronger, condition: Schur-$\calS$. 
\vspace{-0.2cm}
\begin{definition}[Schur-$\calS$]
\label{def:Schur-S}
A $J\times J$ matrix R is Schur-$\calS$ if all its principal sub-matrices are non-singular and there exists a positive vector $h$ such that $h_a^{'}\left[R_a - R_{ab} R_b^{-1}R_{ba}\right] > 0$ for any partition $(a,b)$ of $\mathcal{J}$.
\end{definition}
\vspace{-0.2cm}

It is proved in Corollary 2.8 of \cite{chen1996sufficient} that, if $R$ is both completely-$\calS$ and Schur-$\calS$, then the Skorohod problem is stable if $R^{-1}\theta < 0$. Evidently, Schur-$\calS$ is an algebraically easier condition to work with, but it is more restrictive as the next lemma shows. While Lemma \ref{lem:Schur-S is stronger than Chen-S} is not formally stated in \cite{chen1996sufficient}, it follows from the proof of Corollary 2.6 there.

\vspace{-0.2cm}

\begin{lemma}
\label{lem:Schur-S is stronger than Chen-S}
If a reflection matrix $R$ is both completely-$\calS$ and Schur-$\calS$, then the nominal-load condition, $\rho <e$, implies that $(R,\theta)$ is Chen-$\calS$.
\end{lemma}

\section{Robust Stability of the Skorohod Problem}
\label{sec:Robust Stability of the Skorohod Problem}

\setcounter{equation}{0}
\renewcommand{\theequation}{\thesection.\arabic{equation}}

The main result of this section is that SP stability {\em for all queue-ratio policies} is inherited from SP stability of static priority policies---the ``corners'' of the queue-ratio policy space. Theorem \ref{Thm: completely-S for corner->inter} is a pre-requisite. Static-priority $\Delta$ are as defined in Definition \ref{defin:staticpriorityDelta}.
\vspace{-0.2cm}
\begin{theorem}
\label{Thm: completely-S for corner->inter}
Suppose that, for all static priority $\Delta$: ($i$) $R_{\Delta}^{-1}=CMQ\Delta$ is invertible and $R_\Delta$ is completely-$\calS$;  and ($ii$) all $det(R_{\Delta}^{-1})$ have the same sign. Then, $R_\Delta^{-1}$ is invertible and $R_\Delta$ is completely-$\calS$ for any ratio matrix $\Delta$.
\end{theorem} 
\vspace{-0.2cm}

Recall that the existence of solutions to the SP is guaranteed if $R$ is completely-$\calS$; see \cite{bernard34kharroubi, mandelbaum1987complementarity}. Theorem \ref{Thm: completely-S for corner->inter} guarantees that if the SP is well defined at the corners, it is well defined also in the interior of the space of queue-ratio policies. A square matrix is invertible if and only if its determinant is non-zero. Imposing the same-determinant-sign requirement on the corners, allows us to prove that the reflection matrices for ``interior'' $\Delta$ have a non-zero determinant; the details appear in Lemma \ref{lem: ensure R_Delta invertible}.

The next theorem is an expansion of Theorem \ref{Thm: completely-S for corner->inter}, and shows that the Chen-$\calS$ property (a stronger property than Completely-$\calS$)  has the same characteristic: it is inherited from the static priority ``corners'' to all ratio matrices. In turn, SP stability, for {\em any} queue-ratio policy is inherited from the static-priority corners. 
\vspace{-0.2cm}
\begin{theorem}
\label{Thm: full statement for corner->inter}
Suppose that for any static-priority ratio matrix $\Delta$: ($i$) $R_{\Delta}^{-1}=CMQ\Delta$ is invertible and $(R_\Delta,\theta_\Delta)$ is Chen-$\calS$;  and ($ii$) all $det(R_{\Delta}^{-1})$ have the same sign. Then, $R_\Delta^{-1}$ is invertible and $(R_\Delta, \theta_\Delta)$ is Chen-$\calS$ for any ratio matrix $\Delta$.
\end{theorem} 
\vspace{-0.2cm}

\begin{remark}[Numerical Certification]{\em For a given network and network parameters, it is not necessary to check each of the corners to conclude that the SP is Chen-$\calS$ for any ratio matrix $\Delta$. Instead, one can solve the min-max problem \eqref{RO:ChenS}. A strictly positive objective function value certifies that the SP for the network is robustly QR stable. If the optimal value of \eqref{RO:ChenS} is not strictly positive, then the optimization problem will return at least one ``culprit'' for the violation of our sufficient condition. That is, it will return a $\Delta$ for which the SP's sufficient condition for stability is violated. In the absence of an optimization engine that can solve the problem \eqref{RO:ChenS} efficiently, Theorem \ref{Thm: full statement for corner->inter} immediately suggests an ``enumeration plus optimization'' approach: for each static-priority $\Delta$ solve the interior problem in \eqref{RO:ChenS} until either all such $\Delta$ have returned a positive objective value or a $\Delta$ was found for which that value is not strictly positive; in Figure \ref{fig:Chen-S_DHV_linear_example} this is the corner corresponding to the 2-4-6 static priority policy.} \hfill \bsq
\label{rem:algorithm} 
\end{remark} 

\begin{remark}[Beyond Queue-ratio Policies] {\em The details of queue-ratio policies are not directly used in our analysis of the Skorohod problem. What plays a central role, instead, is that the fluid model reaches state-space collapse. This is true, for example, for Kelly networks operating under First-In-First-Out (FIFO); see \cite{bramson1998state}. In this case, there also exists a matrix $\Delta_{FIFO}$ and a time $t_0$ such that for all $t\geq t_0$, $\wbar{Z}(t)=\Delta_{FIFO} \wbar{W}(t)$. In particular, if the SP is Chen-$\calS$ for all static-priority $\Delta$, it is Chen-$\calS$ for $\Delta_{FIFO}$ which lies in the interior of the space of ratio matrices. In this way, our analysis of SP stability has implications for other policies as long as these have an affine form of state-space collapse; i.e., there exists a matrix $\Delta_{\pi}$ for the policy $\pi$ such that for some $t_0$, $\wbar{Z}(t)=\Delta_{\pi}\wbar{W}(t)$ for all $t\geq t_0$.} We restrict our attention to queue-ratio policies for concreteness as this allows us to develop a complete framework that includes the derivation of state-space collapse as an ingredient. \hfill \bsq
\end{remark}

The interpretation of Theorems \ref{Thm: completely-S for corner->inter} and \ref{Thm: full statement for corner->inter} as a ``convexity'' result is valid only in a restricted sense. Consider two networks, network $1$ and network $2$, that have identical primitives $(C, M, P, \alpha)$, but are  operated under different queue ratio policies with ratio matrices $\Delta_1$ and $\Delta_2$. Suppose further that $(R_{\Delta_1},\theta_{\Delta_1})$ and 
$(R_{\Delta_2},\theta_{\Delta_2})$
are Chen-$\calS$. This does not imply that $(R_{\Delta},\theta_{\Delta})$ is Chen-$\calS$ for any $\Delta = \lambda\Delta_1+(1-\lambda)\Delta_2$ (for $\lambda\in (0,1)$). However, convexity {\em does} hold if the network ratio matrices are  neighbors in the following sense. 

\vspace{-0.2cm}

\begin{definition}[Ratio-matrix Neighbors]
\label{def: policy neighbors}
We say that a collection of ratio matrices $\Delta_i,i=1,2,...$ are ratio-matrix neighbors if they differ in at most one column.
\end{definition}

\vspace{-0.4cm}

\begin{theorem}[Orthogonal Convexity]
\label{thm:Chen-S for a line} Fix network primitives $(C,M,P,\alpha)$ and two  ratio-matrix neighbors $\Delta_1, \Delta_2$ and a third ratio matrix $\Delta_3 = \lambda\Delta_1+(1-\lambda)\Delta_2$, $\lambda\in(0,1)$. Then, $\Delta_1,\Delta_2,\Delta_3$ are ratio-matrix neighbors. Assume that the reflection matrices $R_{\Delta_i}=(CMQ\Delta_i)^{-1},i=1,2,3$ are invertible. Then, 
\begin{itemize} 
\item[(i)] if $R_{\Delta_1}$ and $R_{\Delta_2}$ are both completely-$\calS$, so is $R_{\Delta_3}$;
\item[(ii)] 
if $(R_{\Delta_1},\theta_{\Delta_1})$ and $(R_{\Delta_2},\theta_{\Delta_2})$ are both Chen-$\calS$, so is $(R_{\Delta_3},\theta_{\Delta_3})$.\end{itemize} 
\end{theorem}

\vspace{0.1cm}

\paragraph{\bf Example: The DHV Network.} The DHV network has 3 stations (or servers) and 6 customers classes; see Figure \ref{fig:DHV Network}. Jobs enter the system at rate $\alpha_1$ and are served in the three stations in sequence. After being served at station 3, jobs feed back to station 1 and are again processed by all stations in sequence before departing. 

A ratio matrix $\Delta$ must satisfy $CM\Delta =I$ which, in this network, translates to 
\begin{equation*}
\begin{cases}
         m_1\delta_1 + m_4\delta_4 = 1,\\
         m_2\delta_2 + m_5\delta_5 = 1,\\
         m_3\delta_3 + m_6\delta_6 = 1.
\end{cases}
\end{equation*}
The admissible region of $\Delta \in \mathbb{R}^6_{+}$ can be represented by the 3-dimensional space with $\delta_2\in \left[0, 1/m_2\right]$, $\delta_4\in \left[0, 1/m_4\right]$, $\delta_6\in \left[0, 1/m_6\right]$, since once $\delta_4$ is known, we can calculate $\delta_1$, etc.

We plot the admissible region in Figure \ref{fig:Chen-S_DHV_linear_example} for specific parameters. In the admissible region, a point corresponds to a queue-ratio matrix, thus a queue-ratio policy. We create a dense grid and for each ratio matrix $\Delta$ in the grid we test if $(R_{\Delta},\theta_{\Delta})$ is Chen-$\calS$. The sub-region in green is its Chen-$\calS$ region: for any $\Delta$ in this region $(R_{\Delta},\theta_{\Delta})$ is Chen-$\calS$. The Chen-$\calS$ region is evidently non-convex. Nevertheless, as stated in Theorem \ref{thm:Chen-S for a line}, any line segment defined by neighboring $\Delta$ matrix---a line segment that is parallel to an axis by definition---is in the Chen-$\calS$ region. In this example, one of the static-priority policies (in which classes 2, 4, and 6 have high priority at their corresponding stations) is not Chen-$\calS$ and, in particular, the Chen-$\calS$ property does not hold for all ratio matrices $\Delta$. 

\begin{figure}[h!]
\centering
\includegraphics[scale=0.2]{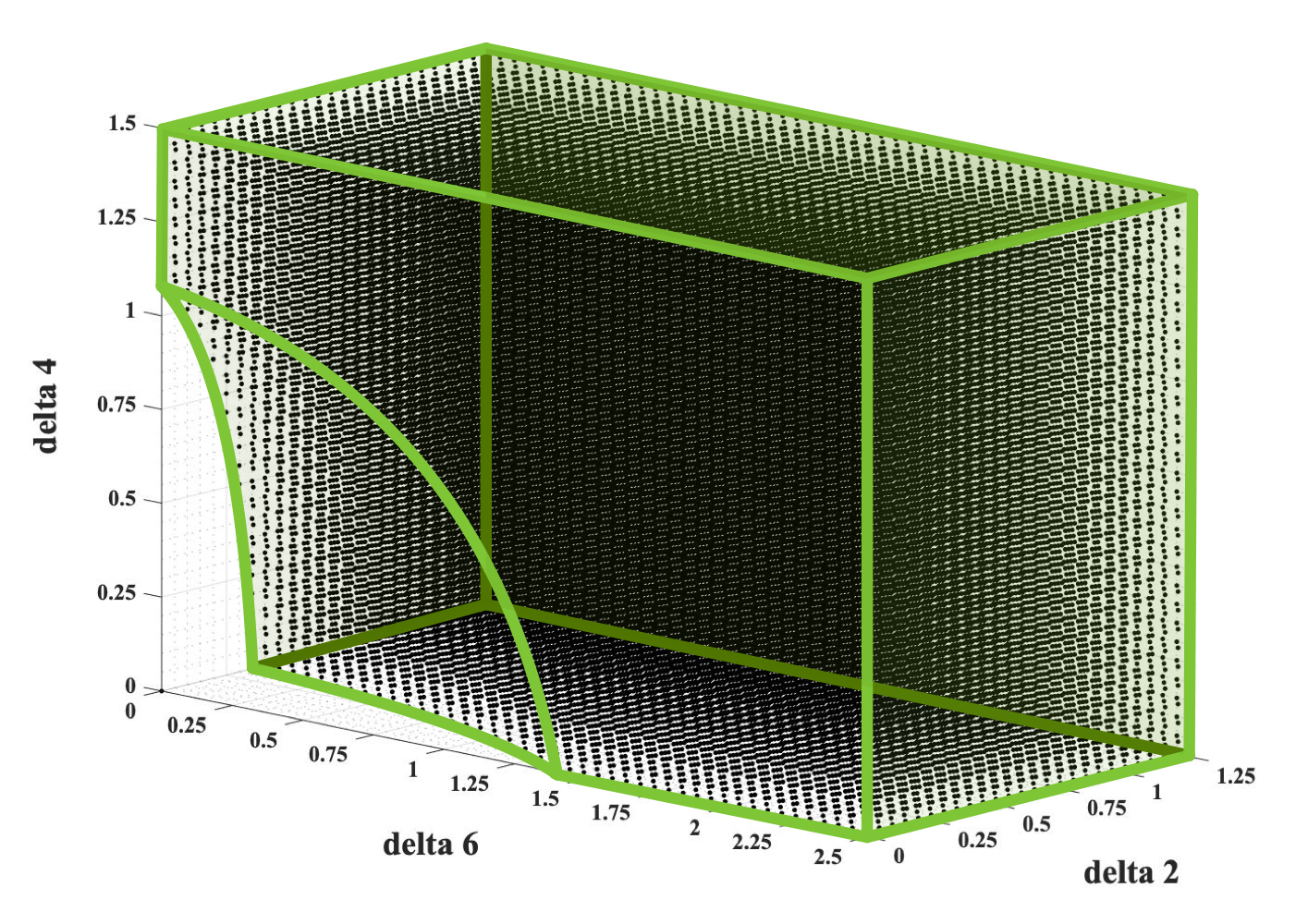}
\caption{The Chen-$\calS$ region in the DHV network with parameters $m_1 = 0.1$, $m_2 = 0.8$, $m_3 = 0.1$, $m_4=0.65$, $m_5 =0.1$, $m_6 = 0.4$ and $\alpha_1 = 0.811$.}
\label{fig:Chen-S_DHV_linear_example}
\end{figure}

To prove Theorems \ref{Thm: completely-S for corner->inter} and \ref{Thm: full statement for corner->inter}, we first establish a linear relationship among $R_{\Delta_1}$, $R_{\Delta_2}$, and $R_{\Delta_3}$ in Theorem \ref{thm: linear combination for R}. This linear relationship is then used to prove Theorem \ref{thm:Chen-S for a line}. Lemma \ref{lem: ensure R_Delta invertible} concerns the preservation of invertibility of the reflection matrix. The detailed proofs for the theorems are in \S \ref{sec:mainproofs} and the proofs for supporting lemmas can be found in Appendix \ref{Appendix: Proofs of Supporting Lemmas}.

\vspace{-0.2cm}

\begin{theorem}
\label{thm: linear combination for R}
Fix network primitives  $(C,M,P,\alpha)$ and two ratio-matrix neighbors $\Delta_1, \Delta_2$ and a third ratio matrix $\Delta_3 = \lambda\Delta_1+(1-\lambda)\Delta_2$, $\lambda\in(0,1)$. Then, $\Delta_1, \Delta_2, \Delta_3$  are ratio-matrix neighbors. Assume that the reflection matrices $R_{\Delta_i} = (CMQ\Delta_i)^{-1}, i = 1,2,3$ are invertible. Then, $R_{\Delta_3}$ is a linear combination of $R_{\Delta_1}$ and $R_{\Delta_2}$. In particular, $R_{\Delta_3} = \beta R_{\Delta_1} + (1-\beta)R_{\Delta_2}$, where 
\begin{equation*}
    \beta = \frac{\lambda det(R_{\Delta_2})}{\lambda det(R_{\Delta_2})+ (1-\lambda) det(R_{\Delta_1})} \in (0,1),
\end{equation*}
or equivalently,
 \begin{equation*}
     \frac{R_{\Delta_3}}{det(R_{\Delta_3})} = \lambda \frac{R_{\Delta_1}}{det(R_{\Delta_1})} + (1-\lambda) \frac{R_{\Delta_2}}{det(R_{\Delta_2})}. 
\end{equation*}  
\end{theorem} 

\vspace{0.2cm}

Theorems \ref{thm:Chen-S for a line} and \ref{thm: linear combination for R} are based on the assumption that the reflection matrices are invertible. Lemma \ref{lem: ensure R_Delta invertible} deduces the invertibility of $R_{\Delta_3}$ from that of $R_{\Delta_1}$ and $R_{\Delta_2}$. 

\vspace{-0.2cm}

\begin{lemma}
\label{lem: ensure R_Delta invertible}
Fix network primitives  $(C,M,P,\alpha)$ and two ratio-matrix neighbors $\Delta_1, \Delta_2$ and a third ratio matrix $\Delta_3 = \lambda\Delta_1+(1-\lambda)\Delta_2$, $\lambda\in(0,1)$. Then, $\Delta_1, \Delta_2, \Delta_3$ are ratio-matrix neighbors.  
Suppose $R_{\Delta_1}^{-1} = CMQ\Delta_1$ and $R_{\Delta_2}^{-1} = CMQ\Delta_2$ are invertible. If $det(R_{\Delta_1}^{-1})$ and $det(R_{\Delta_2}^{-1})$ have the same sign, then $R_{\Delta_3}^{-1} = CMQ\Delta_3$ is invertible.
\end{lemma}


Therefore, Lemma \ref{lem: ensure R_Delta invertible} delineates the invertibility relationship between a reflection matrix corresponding to a ratio matrix $\Delta$ and the reflection matrices corresponding to $\Delta$'s ratio-matrix neighbors. Theorem \ref{thm:Chen-S for a line} proves that, assuming all the reflection matrices are invertible, the completely-$\calS$ and Chen-$\calS$ property of the reflection matrix corresponding to a ratio matrix $\Delta$ can be inherited from the reflection matrices corresponding to $\Delta$'s ratio-matrix neighbors. Combined, these results imply that the completely-$\calS$ and Chen-$\calS$ properties of the reflection matrix corresponding to a ratio matrix can be inherited from the reflection matrices corresponding to extreme ratio matrices, which is precisely the statement of Theorems \ref{Thm: completely-S for corner->inter} and \ref{Thm: full statement for corner->inter}. For example, as shown in Figure \ref{fig:DHV_process_new}, in the DHV network, suppose the admissible region of $\Delta$ is the cuboid, and each point inside the admissible region corresponds to a queue-ratio matrix, and thus to a queue-ratio policy. For an arbitrary policy inside the cuboid, denoted as $A$, the Chen-$\calS$ (or completely-$\calS$) property of $A$ aligns with those of its ratio-matrix neighbors $B$ and $C$ situated on the faces of the cuboid. Similarly, the Chen-$\calS$ (or completely-$\calS$) property of $B$ (or $C$) corresponds with its ratio-matrix neighbors $D$ and $E$ (or $F$ and $G$) located along the edges of the cuboid, the Chen-$\calS$ (or completely-$\calS$) property of which, in turn, are consistent with those of the corner points. Consequently, if all the reflection matrices corresponding to extreme ratio matrices are Chen-$\calS$ (or completely-$\calS$), the reflection matrix corresponding to an arbitrary ratio matrix is Chen-$\calS$ (or completely-$\calS$) as well. 
\begin{figure}[h!]
\centering
\includegraphics[scale=0.15]{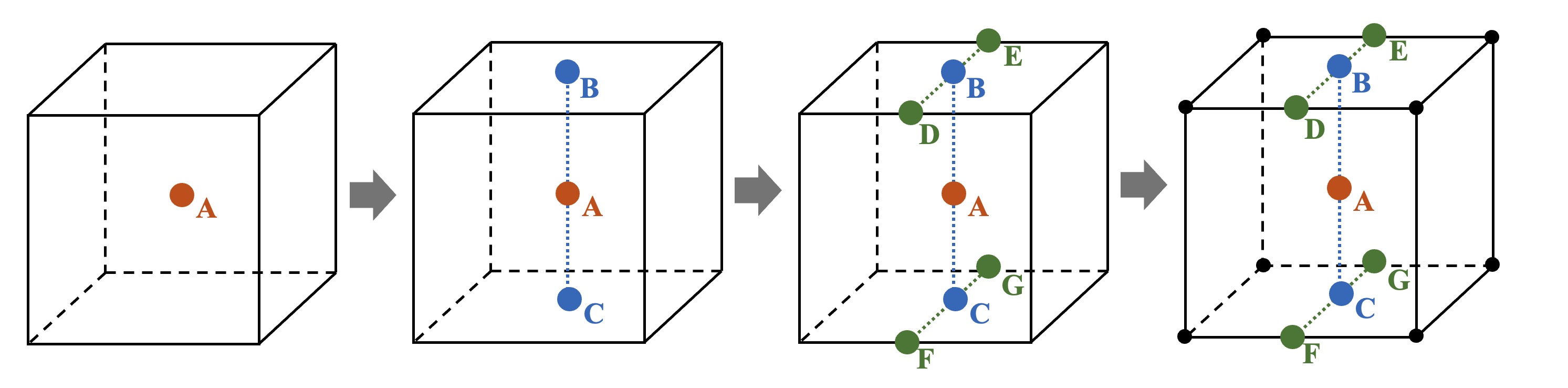}
\caption{The Chen-$\calS$ and Completely-$\calS$ properties at the interior are inherited from the corners.}
\label{fig:DHV_process_new}
\end{figure}
\vspace{-0.5cm}
We now apply our framework to three well-known multiclass queueing networks. 
In \S \ref{DHV Network Example(SP)}, \S \ref{DHV2 Example(SP)} and \S \ref{Lu-Kumar Example(SP)}, we use our results to obtain sufficient conditions for robust SP stability of the three example networks. In \S \ref{Balanced DHV Network Example (SSC)}, \S \ref{Balanced DHV2 Network Example (SSC)} and \S \ref{Balanced LuKumar Network Example (SSC)}, we derive sufficient conditions for robust SSC of the three example networks. There, we also combine these to state sufficient conditions for robust QR-stability in the three networks.

\subsection{Example 1: The DHV Network}
\label{DHV Network Example(SP)}

\vspace{0.2cm}

We say that a network is {\em balanced} when all servers are bottlenecks, that is $\rho_i=(CM\lambda)_i$ is identical for all stations $i$. For the DHV network in Figure \ref{fig:DHV Network}, we assume that $m_1+m_4=m_2+m_5=m_3+m_6=1$\footnote{\label{footnote:barm} Balance means for this network $m_1+m_4=m_2+m_5=m_3+m_6=\bar{m}$ for some $\bar{m}>0$. We can take $\bar{m}=1$ without loss of generality because the fluid model is equivalent to the one with $\bar{m}\neq 1$ via suitable change of the time units.}. The nominal-load condition $\rho_i<1$ for $i = 1,2,3$ then reduces to the requirement that $\alpha_1 <1$. Table \ref{table: Balanced DHV Network} shows that, for any $\alpha_1<1$ and each of the eight static priority policies, $(R_\Delta,\theta_\Delta)$ is Chen-$\calS$ and $R_{\Delta}$ is invertible and has a positive determinant. It then follows from Theorem \ref{Thm: full statement for corner->inter} that the refection matrix $R_{\Delta}$ is invertible and $(R_\Delta,\theta_\Delta)$ is Chen-$\calS$ {\em for any} ratio matrix $\Delta$. Corollary \ref{Corollary: Balanced DHV SP Stable} is now a special case of Theorem \ref{thm: Chen-S (Theorem 2.5)}. The details underlying this table appear in Appendix \ref{Appendix: Chen-S in Balanced DHV Networks}.

\begin{table}[!htbp]
\centering
\begin{tabular}{c c c}
\toprule
High-priority Classes & Chen-$\calS$ & $det(R_\Delta)$\\
\midrule
1, 2, 3 & Yes & $1>0$\\
1, 2, 6 & Yes & $1>0$\\
1, 5, 3 & Yes & $1>0$\\
1, 5, 6 & Yes & $1>0$\\
4, 2, 3 & Yes & $m_1>0$ \vspace{1.3mm}\\
4, 2, 6 & Yes & $\displaystyle\frac{m_1m_3m_5}{m_1m_3m_5+m_2m_4m_6}>0$ \vspace{1.3mm}  \\
4, 5, 3 & Yes & $m_2>0$\\
4, 5, 6 & Yes & $m_3>0$\\
\bottomrule
\end{tabular}
\caption{Properties of the 8 static-priority policies in balanced DHV networks}
\label{table: Balanced DHV Network}
\end{table}

\begin{corollary}
\label{Corollary: Balanced DHV SP Stable}
        In a balanced DHV network that satisfies the nominal-load conditions, the Skorohod Problem is stable for any ratio matrix $\Delta$.
\end{corollary}
\vspace{0.2cm}

The admissible region of $\Delta$ for one specific balanced DHV network appears in Figure \ref{fig:Chen-S_Balanced_DHV_example}. The corner points correspond to static-priority policies. The green region is the Chen-$\calS$ region, where all $(R_\Delta,\theta_\Delta)$ are Chen-$\calS$; unlike the case in Figure \ref{fig:Chen-S_DHV_linear_example}, here it fully covers the space of ratio matrices. 

\begin{figure}[h!]
\centering
\includegraphics[scale=0.21]{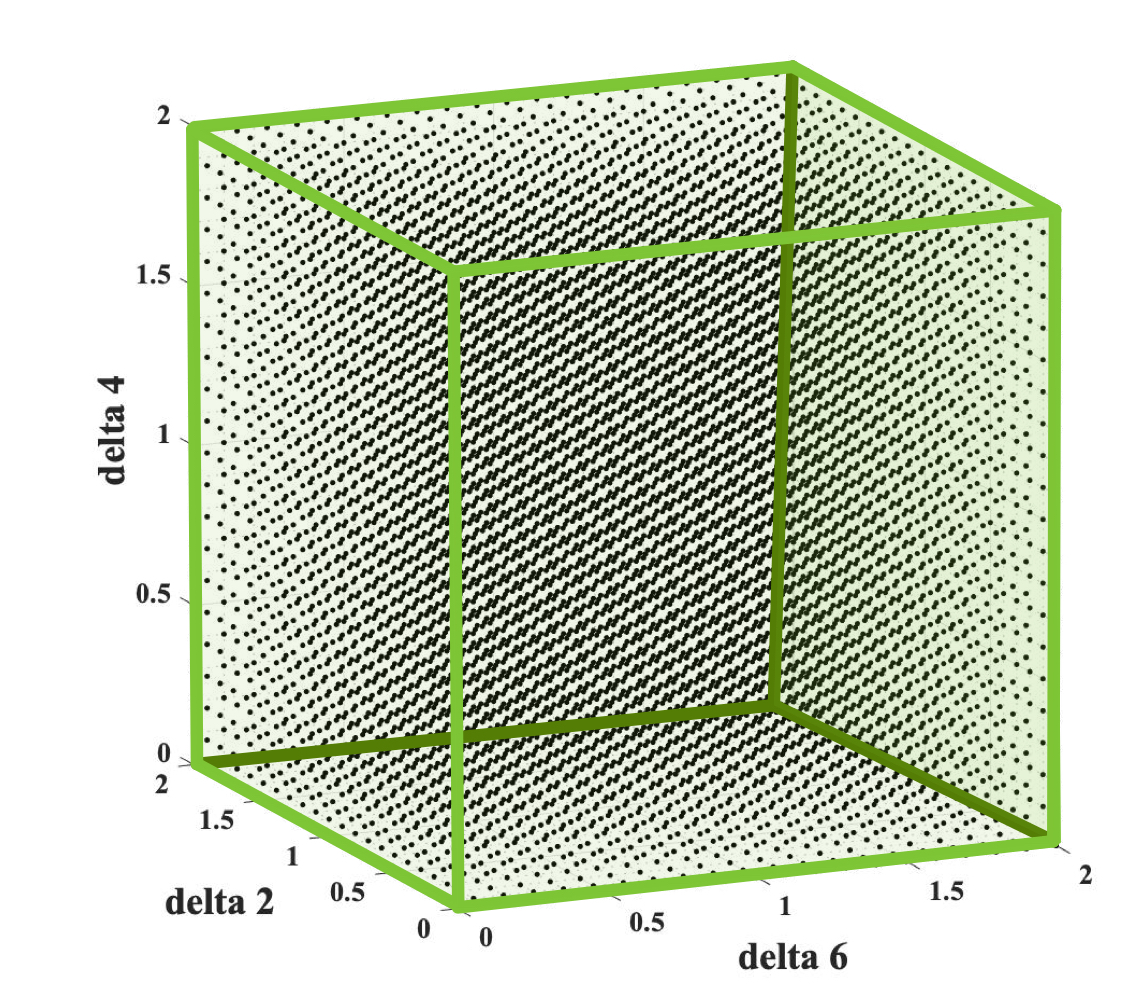}
\caption{For the balanced DHV network with $m_1 = m_2 = m_3 = m_4 = m_5 = m_6 = 0.5$ and $\alpha_1 = 0.811$, the Chen-$\calS$ region covers the whole space of ratio matrices.}
\label{fig:Chen-S_Balanced_DHV_example}
\end{figure}

\subsection{Example 2: The Push Started Lu-Kumar Network}
\label{DHV2 Example(SP)}

The push started Lu-Kumar network \cite{dai2004stability} is the 2-station-5-class network in Figure \ref{fig:DHV2 Network}. Jobs visit the two stations following the route $1 \rightarrow 2  \rightarrow 1 \rightarrow 1 \rightarrow 2$. After completing activity $5$ in station 2, a job leaves the system.

\begin{figure}[h!]
\centering
\includegraphics[scale=0.23]{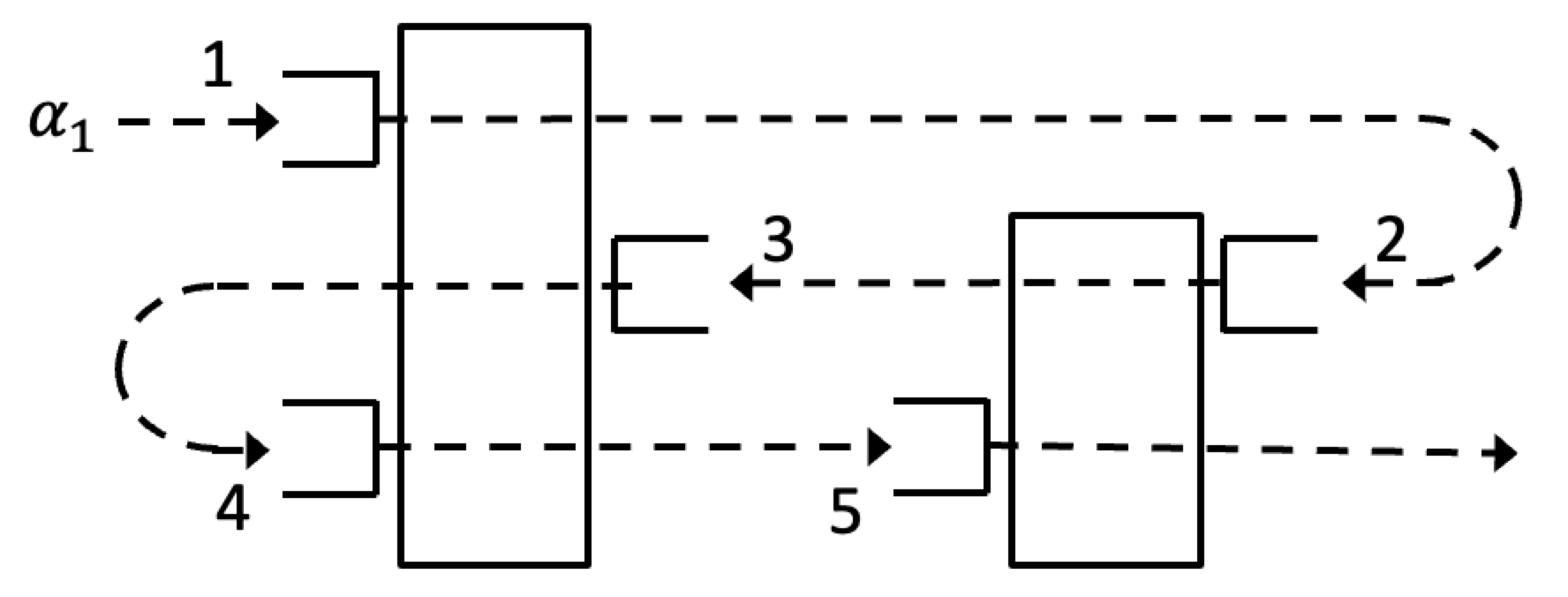}
\caption{A push started Lu-Kumar network}
\label{fig:DHV2 Network}
\end{figure}

\begin{figure}[h!]
\centering
\includegraphics[scale=0.19]{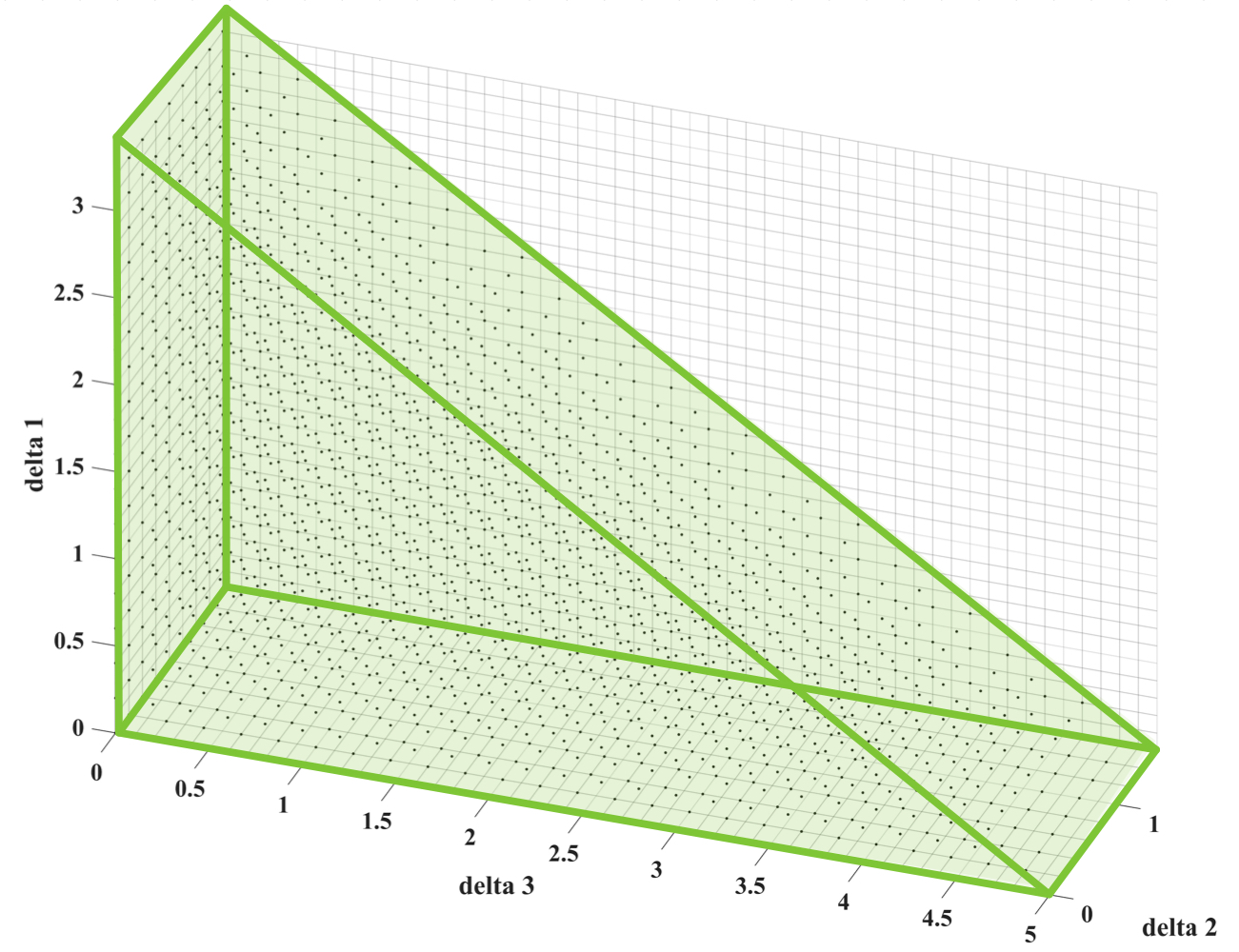}
\caption{The Chen-$\calS$ region in a balanced push started Lu-Kumar network with $m_1 = 0.3$, $m_2 = 0.6$, $m_3 = 0.2$, $m_4=0.5$, $m_5 =0.4$ and $\alpha_1 = 0.8$. These parameters satisfy $m_2m_4 > m_3m_5$.}
\label{fig:Chen-S_DHV2_full_example}
\end{figure}

In a balanced push started Lu-Kumar network $m_1+m_3+m_4=1$ and $m_2+m_5=1$; see Footnote \ref{footnote:barm}. The nominal-load condition then requires that $\alpha_1<1$. We prove that if, in addition to the nominal-load condition, $m_2m_4 > m_3m_5$, then all static-priority $\Delta$ are (i) invertible, (ii) have the same sign of the determinant, and (iii) the pair  $(R_\Delta,\theta_\Delta)$ is Chen-$\calS$; this evidence is collected in Table \ref{table: Balanced DHV2 Network} and is proved in Appendix \ref{Appendix: Chen-S in Balanced DHV2 Networks}. Theorem \ref{Thm: full statement for corner->inter} then guarantees that {\em for any} ratio matrix $\Delta$, $(R_\Delta,\theta_\Delta)$ is Chen-$\calS$, provided that $m_2m_4 > m_3m_5$. This, in turn, implies Corollary 
\ref{Corollary: Balanced DHV2 SP Stable}. 

\begin{table}[!htbp]
\centering
\begin{tabular}{c c c}
\toprule
Lowest-priority Classes in Each Station & Chen-$\calS$ & $det(R_\Delta)$\\
\midrule \vspace{1.5mm}
1, 2 & Yes & $m_2$ \vspace{1.5mm}\\
 1, 5 & Yes & $m_1$ \vspace{1.5mm} \\
3, 2 & Yes & $\displaystyle\frac{m_3}{1-m_1}$ \vspace{1.5mm}\\
3, 5 & Yes & $\displaystyle\frac{m_3}{1-m_1}$ \vspace{1.5mm} \\
4, 2 & Yes, if $m_2m_4>m_3m_5$ & $\displaystyle\frac{m_2m_4}{m_2m_4-m_3m_5}$ \vspace{1.5mm}\\
4, 5 & Yes & $1$ \\
\bottomrule
\end{tabular}
\caption{Properties of the 6 static-priority policies in balanced push started Lu-Kumar networks}
\label{table: Balanced DHV2 Network}
\end{table}
\vspace{-0.2cm}
\begin{corollary}
\label{Corollary: Balanced DHV2 SP Stable}
In a balanced push-started Lu-Kumar network that satisfies the nominal-load conditions the Skorohod problem is stable for any ratio matrix $\Delta$, provided that $m_2m_4 > m_3m_5$.
\end{corollary}
\vspace{-0.2cm}

Recall that a ratio matrix $\Delta$ must satisfy $CM\Delta =I$. In the balanced push started Lu-Kumar network, this specializes to 
\begin{equation*}
\begin{cases}
         m_1\delta_1 + m_3\delta_3 + m_4\delta_4 = 1,\\
         m_2\delta_2 + m_5\delta_5 = 1.
\end{cases}
\end{equation*}
We represent the admissible region of $\Delta \in \mathbb{R}^5_{+}$ by the 3-dimensional space of $(\delta_1,\delta_2,\delta_3)$ with the constraints $m_1\delta_1+m_3\delta_3\leq 1$, and $m_2\delta_2\leq 1$. The values of  $\delta_4$ and $\delta_5$ are then calculated from the others via the identity $CM\Delta= I$.  

Figure \ref{fig:Chen-S_DHV2_full_example} depicts the Chen-$\calS$ region for a balanced push started Lu-Kumar network example and parameters that satisfy $m_2m_4 > m_3m_5$. The 6 corners correspond to the 6 possible static priority policies. In turn, the Chen-$\calS$ region (in green) is equal to the full space of ratio matrices. 

\begin{figure}[h!]
\centering
\includegraphics[scale=0.23]{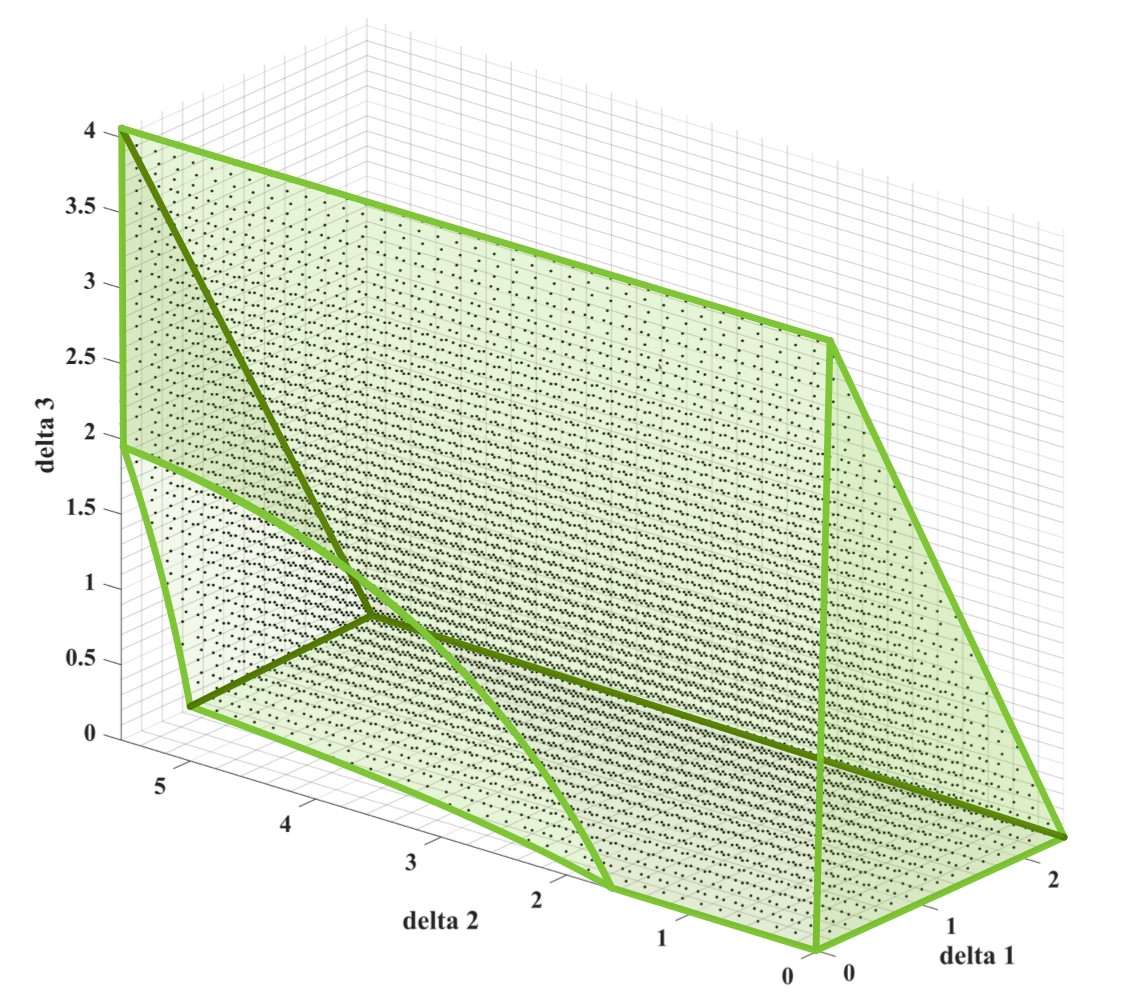}
\caption{The Chen-$\calS$ region in a balanced push started Lu-Kumar network example with $m_1 = 0.41$, $m_2 = 0.18$, $m_3 = 0.24$, $m_4=0.35$, $m_5 =0.82$ and $\alpha_1 = 0.8$, which does not satisfy $m_2m_4 > m_3m_5$.}
\label{fig:DHV2_NotFull_ChenS}
\end{figure}

Figure \ref{fig:DHV2_NotFull_ChenS} shows the space of ratio matrices and the Chen-$\calS$ region for a balanced push-started Lu-Kumar network with $m_1 = 0.41$, $m_2 = 0.18$, $m_3 = 0.24$, $m_4=0.35$, $m_5 =0.82$ and $\alpha_1 = 0.8$, which does not satisfy $m_2m_4 > m_3m_5$. The static-priority ratio matrix, corresponding to classes 4 and 2 having the lowest priorities in their respective stations, is not Chen-$\calS$; there is, then,  a subset of $\Delta$ matrices for which $(R_{\Delta}, \theta_{\Delta})$ is not Chen-$\calS$.

\subsection{Example 3: The Lu-Kumar Network}
\label{Lu-Kumar Example(SP)}

The Lu-Kumar network is the 2-station and 4-class network in Figure \ref{fig:Lu-Kumar Network}. Jobs visit the two stations following the route $1 \rightarrow 2  \rightarrow 2 \rightarrow 1$. After completing activity $4$ in station 1, a job leaves the system.

\begin{figure}[h!]
\begin{center} 
~~~ \includegraphics[scale=0.18]{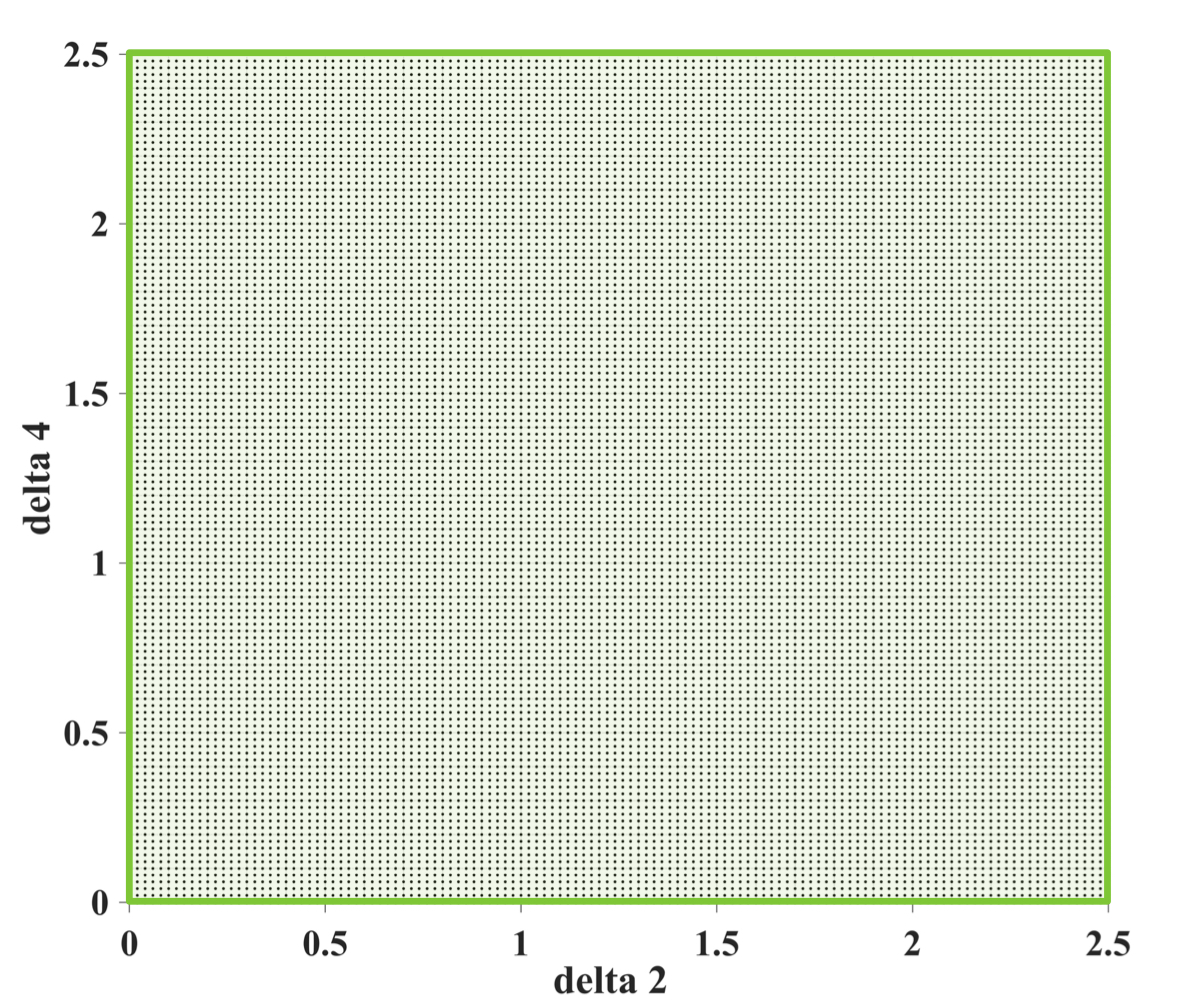}~~~\includegraphics[scale=0.18]{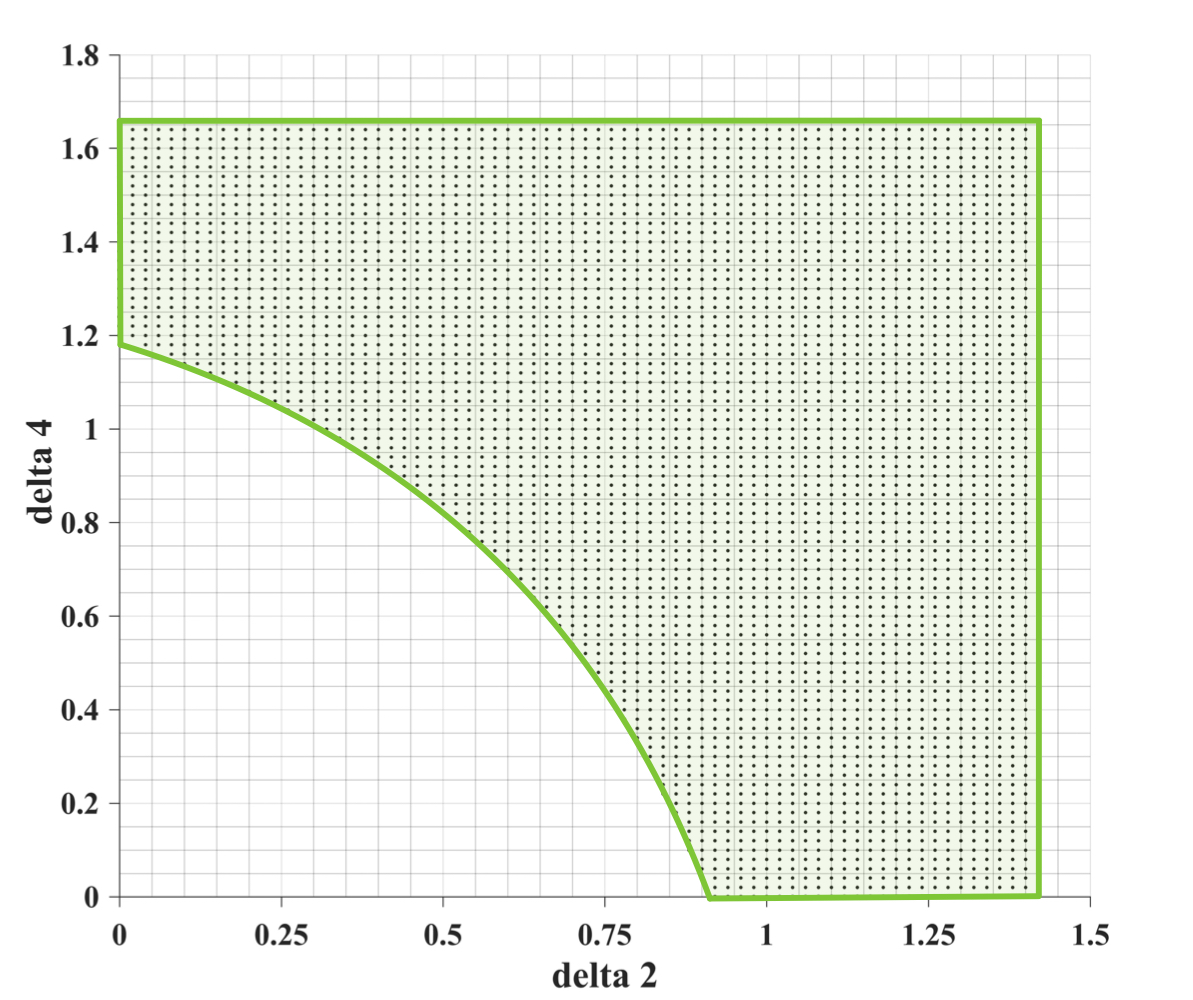}\end{center} 
\caption{(LEFT) The Chen-$\calS$ region in a balanced Lu-Kumar network with $m_1 = 0.6$, $m_2 = 0.4$, $m_3 = 0.6$, $m_4=0.4$ and $\alpha_1 = 0.9$. These parameters satisfy $m_2+m_4 < 1$. (RIGHT) The Chen-$\calS$ region in a balanced Lu-Kumar network example with $m_1 = 0.4$, $m_2 = 0.7$, $m_3 = 0.3$, $m_4=0.6$ and $\alpha_1 = 0.9$, which does not satisfy $m_2+m_4 < 1$. 
\label{fig:Chen-S_LuKumar_full_example}}
\end{figure}

In the balanced Lu-Kumar network $m_1+m_4 = m_2+m_3=1$; see Footnote \ref{footnote:barm}. The nominal-load condition requires that $\alpha_1<1$. We prove that if, in addition to the nominal-load condition, we have $m_2+m_4 < 1$, then all static-priority $\Delta$ are (i) invertible, (ii) have the same sign of the determinant, and (iii) the pair  $(R_\Delta,\theta_\Delta)$ is Chen-$\calS$. This evidence is collected in Table \ref{table: Balanced LuKumar Network} and is proved in Appendix \ref{Appendix: Chen-S in Balanced LuKumar Networks}. Theorem \ref{Thm: full statement for corner->inter} then guarantees that $(R_\Delta,\theta_\Delta)$ is Chen-$\calS$ {\em for any} ratio matrix $\Delta$  provided that $m_2+m_4 < 1$.

\begin{table}[!htbp]
\centering
\begin{tabular}{c c c}
\toprule
High-priority Classes & Chen-$\calS$ & $det(R_\Delta)$\\
\midrule \vspace{1.5mm}
1, 2 & Yes & $1$ \vspace{1.5mm}\\
1, 3 & Yes & $m_2$ \vspace{1.5mm} \\
4, 2 & Yes if $m_2+m_4 <1$ & $\displaystyle\frac{m_1m_3}{1-m_4-m_2}$ \vspace{1.5mm}\\
4, 3 & Yes & $m_2$\\
\bottomrule
\end{tabular}
\caption{Properties of the 4 static-priority policies in balanced Lu-Kumar networks}
\label{table: Balanced LuKumar Network}
\end{table}
\vspace{-0.2cm}
\begin{corollary}
\label{Corollary: Balanced LuKumar SP Stable} In a balanced Lu-Kumar network satisfying the nominal-load conditions, the Skorohod problem is stable for any ratio matrix $\Delta$ if $m_2+m_4<1$.
\end{corollary}
\vspace{-0.2cm}

Figure \ref{fig:Chen-S_LuKumar_full_example}(LEFT) depicts the Chen-$\calS$ region for a balanced Lu-Kumar network that satisfies $m_2+m_4 < 1$; the Chen-$\calS$ region (in green) is equal to the full space of ratio matrices. Figure \ref{fig:Chen-S_LuKumar_full_example}(RIGHT) displays, in contrast, an example that violates the requirement $m_2+m_4 <1$. For the $\Delta$ corresponding to a policy which gives high priority to classes 4 and 2, the pair $(R_{\Delta}, \theta_{\Delta})$ is not Chen-$\calS$. There is, in this example, a family of $\Delta$ matrices, in the vicinity of the corresponding corner, for which $(R_{\Delta}, \theta_{\Delta})$ is not Chen-$\calS$.

\section{Robust State-Space Collapse}
\setcounter{equation}{0}
\renewcommand{\theequation}{\thesection.\arabic{equation}}

Solutions to the fluid model equations \eqref{equ:begin fluid model}-\eqref{equ:end fluid model} are Lipschitz continuous, hence almost everywhere differentiable. Following standard terminology, we refer to times $t$ where the solution is differentiable as ``regular.'' For a static priority policy $\pi$, we use $\mathcal{L}_\pi$ to denote the lowest priority class in each station, and $\mathcal{H}_\pi = \mathcal{K} \backslash \mathcal{L}_\pi$. To make explicit the connection between a static-priority policy ratio matrix $\Delta$ and the high-priority classes, we sometimes write $\mathcal{H}_{\pi}^{\Delta}$. 

The main theorem of this section relates state-space collapse for arbitrary queue-ratio policies to the linear attraction of static-priority policies. In the fluid model of static-priority policies, only the identity of the lowest-priority class matters; the ordering among the high-priority classes does not.

\vspace{-0.3cm}

\begin{definition}[Linear Attraction]
\label{def: linear attraction}
    We say that linear attraction holds for the static-priority policy with ratio matrix $\Delta$ if there exists a test vector $h \in\mathbb{R}^K_{++}$ and a constant $r >0$ such that any solution to the fluid model equations \eqref{equ:begin fluid model}-\eqref{equ:end fluid model} satisfies 
    $$\sum_{k\in \mathcal{H}_\pi^{\Delta}:\wbar{Z}_k>0} h_k\dot{\wbar{Z}}_k(t) \leq -r \text{ for any regular time t with } ||\wbar{Z}_{\mathcal{H}_\pi^{\Delta}}(t)|| >0.$$
\end{definition}

\vspace{-0.2cm}

\begin{theorem}
\label{thm: LEGO paper (Theorem 2)}
Suppose that linear attraction holds for all static priority $\Delta$ with a common test vector $h\in \mathbb{R}^K_{++}$ and a constant $r >0$. Then, for any ratio matrix $\Delta$ and any initial condition $\wbar{Z}(0)$, there exists a $t_0$ such that $\wbar{\epsilon}(t) = \wbar{Z}(t) - \Delta CM \wbar{Z}(t)  =0$ for all $t\geq t_0$.
\end{theorem}


Theorem \ref{thm: LEGO paper (Theorem 2)} stipulates that queue-ratio SSC is inherited from the static priority policies: SSC of all static priority policies implies SSC for any queue-ratio policy. Intuitively, a queue-ratio policy moves between static priority policies, depending on the balance status of the queues: the classes in the set $\{\ell \in \mathcal{C}(j):Z_\ell(t)-\delta_\ell W_j(t) >0 \}$ are prioritized at time $t$ over all others. The proof of Theorem \ref{thm: LEGO paper (Theorem 2)} (in \S \ref{sec:mainproofs}) proceeds by showing that the set of possible QR actions at a given state of the queues is the same as that set for {\em some} priority permutation and then leveraging the assumed linear attraction of all priority policies. 

\subsection{Example 1: The DHV Network}
\label{Balanced DHV Network Example (SSC)}

 The detailed analysis of linear attraction for the 8 static-priority policies appears in Appendix \ref{Appendix: The SSC Inequalities for Balanced DHV Networks} and establishes that the requirements of Theorem \ref{thm: LEGO paper (Theorem 2)} hold for the balanced DHV network if and only if  there exists a vector $h$ that satisfies the following inequalities: 
 \begin{align}
    & h_1, h_2, h_3, h_4, h_5, h_6 >0 \label{DHV_eq_1}\\ 
    & h_2(\mu_1-\mu_2)-h_6\mu_6 <0  
    \label{DHV_eq_2}\\
    & h_4(\mu_3-\mu_4)-h_2\mu_2 <0  
    \label{DHV_eq_3}\\
    & h_6(\mu_5-\mu_6)-h_4\mu_4 <0  
    \label{DHV_eq_4}\\
    & h_1(\alpha_1-\mu_1)+h_2(\mu_1-\mu_2)+h_3(\mu_2-\mu_3) <0  
    \label{DHV_eq_5}\\
    & h_1(\alpha_1-\mu_1)+h_2(\mu_1-\mu_2) <0 
    \label{DHV_eq_6}\\
    & h_1(\alpha_1-\mu_1)+h_3(\mu_1-\mu_3) <0  
    \label{DHV_eq_7}\\
    & h_2(\alpha_1-\mu_2)+h_3(\mu_2-\mu_3) <0 
    \label{DHV_eq_8}\\
    & h_3(\mu_2-\mu_3)+h_4(\mu_3-\mu_4)-h_2\mu_2 <0  
    \label{DHV_eq_9}\\
    & h_2\left[\mu_1\left(1-\frac{\mu_3}{\mu_4}\right)-\mu_2\right]+h_3(\mu_2-\mu_3) <0  
    \label{DHV_eq_10}\\
    & h_4(\mu_2-\mu_4)-h_2\mu_2 <0  
    \label{DHV_eq_11}\\
    & h_4(\mu_3-\mu_4)-h_3\mu_3 <0  
    \label{DHV_eq_12}\\
    & h_4(\mu_3-\mu_4)+h_5(\mu_4-\mu_5)-h_3\mu_3 <0  
    \label{DHV_eq_13}\\
    & h_3\left[\mu_2\left(1-\frac{\mu_4}{\mu_5}\right)-\mu_3\right]+h_4(\mu_3-\mu_4) <0  
    \label{DHV_eq_14}\\
    & h_5(\mu_3-\mu_5)-h_3\mu_3 <0  
    \label{DHV_eq_15}\\
    & h_5(\mu_4-\mu_5)-h_4\mu_4 <0  
    \label{DHV_eq_16}\\
    & h_5(\mu_4-\mu_5)+h_6(\mu_5-\mu_6)-h_4\mu_4 <0  
    \label{DHV_eq_17}\\
    & h_4\left[\mu_3\left(1-\frac{\mu_5}{\mu_6}\right)-\mu_4\right]+h_5(\mu_4-\mu_5) <0  
    \label{DHV_eq_18}\\
    & h_6(\mu_4-\mu_6)-h_4\mu_4 <0  
    \label{DHV_eq_19}\\
    & h_6(\mu_5-\mu_6)-h_5\mu_5 <0  
    \label{DHV_eq_20}\\
    & h_1(\alpha_1-\mu_1)+h_6(\mu_5-\mu_6)-h_5\mu_5 <0  
    \label{DHV_eq_21}\\
    & h_5\left[\mu_4\left(1-\frac{\alpha_1}{\mu_1}\right)-\mu_5\right]+h_6(\mu_5-\mu_6) <0  
    \label{DHV_eq_22}\\
    & h_1(\alpha_1-\mu_1)+h_2(\mu_1-\mu_2)-h_6\mu_6  <0  
    \label{DHV_eq_23}\\
    & h_1(\alpha_1-\mu_1)+h_2(\mu_1-\mu_2) <0  
    \label{DHV_eq_24}\\
    & h_1(\alpha_1-\mu_1)+h_3(\mu_2-\mu_3) <0.  
    \label{DHV_eq_25}
\end{align}

Exactly how these inequalities arise from the eight priority policies is made clear in the appendices. Equation \eqref{SSC_proof_eq_2} in the proof of Theorem \ref{thm: LEGO paper (Theorem 2)} grounds the general ``recipe'' for producing these inequalities. We note that inequality \eqref{DHV_eq_23} is immediately satisfied if \eqref{DHV_eq_6} is satisfied; and inequality \eqref{DHV_eq_24} is the same as \eqref{DHV_eq_6}; 23 inequalities remain. 

The analysis of these inequalities in  Appendix \ref{Appendix: Sufficient Conditions for SSC of Balanced DHV Networks} establishes the following.

\vspace{-0.2cm}
\begin{corollary}
\label{Corollary: Balanced DHV SSC} Consider the balanced DHV network. There exists $\bar{\alpha}\in (0,1)$ (dependent on $M$) such that if $\alpha_1> \bar{\alpha}$ and 
\begin{itemize} 
\item[(i)] the nominal-load conditions are satisfied: $\rho_1=\rho_2=\rho_3<1$ (equivalently, $\alpha_1 <1$),
\item[(ii)] $m_2+ m_4+m_6<2$,\\
then the network satisfies SSC under any queue-ratio policy.
\end{itemize} 
\end{corollary}

The exact value of $\bar{\alpha}$ appears in Appendix \ref{Appendix: The SSC Inequalities for Balanced DHV Networks}. The requirement that $\alpha$ is ``large enough'' is for SSC only---it is not required for SP Stability---and arises from our insistence on using linear attraction for SSC. In the absence of this requirement, some high priority queues grow before decreasing, a behavior that is inconsistent with linear Lyapunov attraction. It is plausible that this requirement can be relaxed; see  further discussion in \S \ref{sec:conclusions}. Importantly, the requirement that $\alpha_1>\bar{\alpha}$, is immediately satisfied if the network is in heavy traffic. 

Combining Corollary \ref{Corollary: Balanced DHV SSC} with the robust SP stability in Corollary \ref{Corollary: Balanced DHV SP Stable}, we conclude the following.  
\vspace{-0.2cm}
\begin{corollary}
\label{Corollary: Balanced DHV global stability condition} Consider the balanced DHV network. There exists $\bar{\alpha}\in (0,1)$ (dependent on $M$) such that if $\alpha_1> \bar{\alpha}$ and 
\begin{itemize} 
\item[(i)] the nominal-load conditions are satisfied: $\rho_1=\rho_2=\rho_3<1$ (equivalently, $\alpha_1 <1$),
\item[(ii)] $m_2+ m_4+m_6<2$,\\
then the network 
is robustly queue-ratio stable.\end{itemize} 
\end{corollary}
\vspace{-0.2cm}

The condition $m_2+ m_4+m_6<2$ is weaker than the sufficient global-stability condition $m_2+m_4m_6<1$ in \cite{DHV}. The gap between our condition and the one in \cite{DHV} is consistent with our focus on queue-ratio policies. The policies considered in \cite{DHV} are general non-idling policies, of which the family of queue-ratio policies is a subset.

\subsection{Example 2: The Push-Started Lu-Kumar Network
\label{Balanced DHV2 Network Example (SSC)}}

In this network (see Figure \ref{fig:DHV2 Network}), three classes are served at station 1, so we must specify the choice of the class for service when the set $\{\ell \in \mathcal{C}(j):Z_l(t)-\delta_lW_j(t) >0 \}$ contains more than one class. We fix the priority order to $1>3>4$ as the ``tie-breaking'' rule: if both queue 1 and queue 3 are above their target, class 1 has priority over class 3.

 The detailed analysis of linear attraction for the 6 static-priority policies appears in Appendix \ref{Appendix: The SSC Inequalities for Balanced DHV2 Networks} and establishes that the requirements of Theorem \ref{thm: LEGO paper (Theorem 2)} hold for this network if and only if  there exists a vector $h$ that satisfies the following inequalities: 
\begin{align}
    & h_1, h_2, h_3, h_4, h_5 >0 
    \label{DHV2_eq_1}\\ 
    & -h_4\mu_4 + h_5(\mu_4-\mu_5) <0  
    \label{DHV2_eq_2}\\
    & h_3(\mu_2-\mu_3)+h_4\mu_3 <0  
    \label{DHV2_eq_3}\\
    & -h_2\mu_2+h_3(\mu_2-\mu_3) + h_4\mu_3 <0  
    \label{DHV2_eq_4}\\
    & -h_3\mu_3 + h_4\mu_3 <0  
    \label{DHV2_eq_5}\\
    & h_1(\alpha_1-\mu_1)+h_2(\mu_1-\mu_2) <0  
    \label{DHV2_eq_6}\\
    & h_1(\alpha_1-\mu_1)+h_3\mu_2 <0  
    \label{DHV2_eq_7}\\
    & m_2m_4>m_4m_5  
    \label{DHV2_eq_8}\\
    & h_1(\alpha_1-\mu_1) + h_2(\mu_1-\mu_2)+h_3\mu_2<0 
    \label{DHV2_eq_9}\\
    & h_1(\alpha_1-\mu_1)+h_3\mu_1 <0 \;\;\text{if}\;\; \mu_1 \leq \mu_2.  
    \label{DHV2_eq_10}
\end{align}

These arise from the different priority orders. Equation \eqref{DHV2_eq_8}, specifically, arises from the priority order $\pi\{1,3,5\}$. Inequalities \eqref{DHV2_eq_4} and \eqref{DHV2_eq_5} are immediately satisfied if \eqref{DHV2_eq_2} is. The inequality \eqref{DHV2_eq_6} is satisfied if \eqref{DHV2_eq_9} is, and \eqref{DHV2_eq_10} is satisfied if \eqref{DHV2_eq_7} is. Thus, we need only consider the inequality $m_2m_4>m_3m_5$ and the remaining five inequalities:
\begin{equation*}
\left\{\begin{array}{lr}
     h_1, h_2, h_3, h_4, h_5 >0 \\ 
     -h_4\mu_4 + h_5(\mu_4-\mu_5) <0\\
     h_3(\mu_2-\mu_3)+h_4\mu_3 <0 \\
    h_1(\alpha_1-\mu_1)+h_3\mu_2 <0  \\
    h_1(\alpha_1-\mu_1) + h_2(\mu_1-\mu_2)+h_3\mu_2<0. \\
        \end{array}
    \right.
\end{equation*}

We make two useful observations: (i) If $\mu_2 \geq \mu_3$ ($m_2 \leq m_3$), the inequality $ h_3(\mu_2-\mu_3)+h_4\mu_3 <0$ cannot be satisfied with positive $h_3$ and $h_4$. When  $m_2>m_3$, in contrast, there exists a positive vector $h$ with $h_1 \geq h_3 \geq h_4 \geq h_5$ and $h_1 \geq h_2$ satisfying the five inequalities. Thus, $m_2>m_3$ is a sufficient condition for the existence of an $h$ that satisfies the five inequalities; and (ii) for the balanced network $m_2m_4>m_3m_5$ implies $m_2<m_3$. This is because, with $m_1+m_3+m_4 =m_2+m_5 = 1$, $m_2m_4>m_3m_5 \Leftrightarrow m_2m_4>m_3(1-m_2) \Leftrightarrow m_2(m_3+m_4)>m_3 \Leftrightarrow
m_2(1-m_1)>m_3 
\Rightarrow m_2>m_3
$. 

Combining these observations, we conclude that $m_2m_4>m_3m_5$ suffices for the existence of a strictly positive vector $h$ that satisfies all the inequalities. Theorem \ref{thm: LEGO paper (Theorem 2)}, then yields the following corollary. Below, the tie-breaking rule is $1>3>4$.

\begin{corollary}
\label{Corollary: Balanced DHV2 SSC} Consider the balanced push-started Lu-Kumar network. There exists $\bar{\alpha}\in (0,1)$ (dependent on $M$) such that if $\alpha_1> \bar{\alpha}$, and 
\begin{itemize} 
\item[(i)] the nominal-load condition is satisfied: $\rho_1=\rho_2<1$ (equivalently, $\alpha_1 <1$),
\item[(ii)] $m_2m_4 > m_3m_5$,\\
then the network satisfies SSC under any queue-ratio policy.\end{itemize} 
\end{corollary}


Combining this with Corollary \ref{Corollary: Balanced DHV2 SP Stable}, we conclude the following.  
\begin{corollary}
\label{Corollary: Balanced DHV2 global stability condition}Consider the balanced push-started Lu-Kumar network. There exists $\bar{\alpha}\in (0,1)$ (dependent on $M$) such that if $\alpha_1> \bar{\alpha}$, and 
\begin{itemize} 
\item[(i)] the nominal-load condition is satisfied: $\rho_1=\rho_2<1$ (equivalently, $\alpha_1 <1$),
\item[(ii)] $m_2m_4 > m_3m_5$,\\
then the network is robustly queue-ratio stable.\end{itemize} 
\end{corollary}

With load balance, the requirement $m_2m_4>m_3m_5$ is equivalent to the requirement $\frac{m_3}{1-m_1}+m_5<1$. The latter is the sufficient condition for global stability in \cite{dai2004stability}. This is because, in this case, $m_2m_4>m_3m_5 \Leftrightarrow (1-m_5)(1-m_1-m_3) > m_3m_5 
\Leftrightarrow 
1-m_1>m_3+m_5-m_1m_5
\Leftrightarrow 
1-m_1>m_3+m_5(1-m_1)
\Leftrightarrow 
\frac{m_3}{1-m_1} +m_5 <1$. This last condition, known already from \cite{dai2004stability}, is intuitively understood by considering that, once server 1 allocates capacity to class 1, the residual network is a Lu-Kumar network where the global stability condition would be $m_3 +m_5 <1$. Since part of server 1's time is dedicated to class 1, we essentially ``rescale'' the average service time for class 3 by $1/({1-m_1)}$ to obtain a modified version of this inequality: $\frac{m_3}{1-m_1} +m_5 <1$.

\subsection{Example 3: The Lu-Kumar Network
\label{Balanced LuKumar Network Example (SSC)}}

The detailed analysis of linear attraction for the 4 static-priority policies appears in Appendix \ref{Appendix: The SSC Inequalities for Balanced LuKumar Networks} and establishes that the requirements of Theorem \ref{thm: LEGO paper (Theorem 2)} hold for the this network if and only if $m_2+m_4<1$. 
\begin{corollary}
\label{Corollary: Balanced LuKumar SSC} A balanced Lu-Kumar network satisfies SSC under any queue-ratio policy if
\begin{itemize} 
\item[(i)] the nominal-load condition is satisfied: $\rho_1=\rho_2<1$ (equivalently, $\alpha_1 <1$);
\item[(ii)] $m_2+m_4<1$.
\end{itemize} 
\end{corollary}
Combining this with Corollary \ref{Corollary: Balanced LuKumar SP Stable}, we conclude the following.  
\begin{corollary}
\label{Corollary: Balanced DHV2 global stability condition}
A balanced Lu-Kumar network is robustly queue-ratio stable if
\begin{itemize} 
\item[(i)] the nominal-load condition is satisfied: $\rho_1=\rho_2<1$ (equivalently, $\alpha_1 <1$);
\item[(ii)]  $m_2+m_4<1$.\end{itemize} 
\end{corollary}
The requirement $ m_2+m_4<1$ is necessary and sufficient for global stability; see \cite{dai1996stability}.

\vspace*{0.5cm} 
\section{Proofs of Theorems} 
\label{sec:mainproofs}
\setcounter{equation}{0}
\renewcommand{\theequation}{\thesection.\arabic{equation}}

For the proof of Theorem \ref{thm: linear combination for R}, we will need Lemma \ref{lem:decomposition with A} below. Suppose $\Delta_1$ and $\Delta_2$ are different only in the $k^{th}$ column, $R_{\Delta_3}^{-1} = CMQ(\lambda\Delta_1+(1-\lambda)\Delta_2) = CMQ\Delta_1 +(1-\lambda)CMQ(\Delta_2-\Delta_1) = R_{\Delta_1}^{-1} + (1-\lambda)CMQ(\Delta_2-\Delta_1)$, where $CMQ(\Delta_2-\Delta_1)$ is a matrix with only one non-zero column, say the $k^{th}$ column.  We can rewrite $CMQ(\Delta_2-\Delta_1)$ as 
\begin{align}
& uv^\top = CMQ(\Delta_2-\Delta_1) ,\text{where}  \label{3.1}\\
& \text{$u$ is a column vector whose elements are the same as the $k^{th}$ column in $CMQ(\Delta_2-\Delta_1)$,} \label{3.2}\\
& \text{$v$ is a column vector with the $k^{th}$ element 1 and all other elements 0.} \label{3.3}
\end{align} 
Thus, we have
\begin{equation}
\label{equ:R3_with_uv}
    R_{\Delta_3}^{-1} = R_{\Delta_1}^{-1} + (1-\lambda)CMQ(\Delta_2-\Delta_1)=R_{\Delta_1}^{-1} + (1-\lambda)uv^\top
\end{equation}
Similarly, 
\begin{equation}
\label{equ:R3_with_uv_2}
   R_{\Delta_3}^{-1} = CMQ\Delta_2 +\lambda CMQ(\Delta_1-\Delta_2) = R_{\Delta_2}^{-1} - \lambda CMQ(\Delta_2-\Delta_1) =R_{\Delta_2}^{-1} -\lambda uv^\top.
\end{equation}

\vspace{-0.2cm}

\begin{lemma}
\label{lem:decomposition with A}
Fix network primitives  $(C,M,P,\alpha)$ and two ratio-matrix neighbors $\Delta_1, \Delta_2$ that differ only in the $k^{th}$ column and a third ratio matrix $\Delta_3 = \lambda\Delta_1+(1-\lambda)\Delta_2$, $\lambda\in(0,1)$. Then, $\Delta_1, \Delta_2, \Delta_3$  are ratio-matrix neighbors. Assume that the reflection matrices $R_{\Delta_i} = (CMQ\Delta_i)^{-1}, i = 1,2,3$ are invertible, and let $u,v$ be the vectors that satisfy \eqref{3.1}-\eqref{3.3} for $\Delta_1$ and $\Delta_2$. Then,
$R_{\Delta_1}$ and $R_{\Delta_2}$ are proportional in the $k^{th}$ row with the ratio
\begin{equation}
\label{equ(mainBody):R1_R2_proportional_in_k_row}
(R_{\Delta_1})_{k.} = \frac{det(R_{\Delta_2}^{-1})}{det(R_{\Delta_1}^{-1})} (R_{\Delta_2})_{k.},
\end{equation}
and $R_{\Delta_3} = R_{\Delta_1}(I-(1-\lambda)A) = R_{\Delta_2}(I+\lambda A)$, where
\begin{equation*}
 A = \frac{uv^\top R_{\Delta_1}}{1+(1-\lambda)v^\top R_{\Delta_1} u} = \frac{uv^\top R_{\Delta_2}}{1-\lambda v^\top R_{\Delta_2} u}.
\end{equation*}
\end{lemma}

\vspace{2mm}

\begin{proof}{Proof of Theorem \ref{thm: linear combination for R}.} 
\label{proof: linear combination for R}
By Lemma \ref{lem:decomposition with A}, \begin{equation} R_{\Delta_3} = R_{\Delta_1}(I-(1-\lambda)A) = R_{\Delta_2}(I+\lambda A),\label{eq:itais2}\end{equation} where, recall, 
\begin{equation*}
    A = \frac{uv^\top R_{\Delta_1}}{1+(1-\lambda)v^\top R_{\Delta_1} u} = \frac{uv^\top R_{\Delta_2}}{1-\lambda v^\top R_{\Delta_2} u}.   
\end{equation*}
Defining
\begin{equation*}
    \beta_{\lambda} = 1-\frac{1-\lambda}{1-\lambda v^\top R_{\Delta_2} u},
\end{equation*}
we have 
\begin{equation*}
\begin{split}
  A & = \frac{1-\beta_\lambda}{1-\lambda}uv^\top R_{\Delta_2} = \frac{1-\beta_\lambda}{1-\lambda}(R_{\Delta_2}^{-1}-R_{\Delta_1}^{-1})R_{\Delta_2} =\frac{1-\beta_\lambda}{1-\lambda}(I-R_{\Delta_1}^{-1}R_{\Delta_2}),   
\end{split}
\end{equation*}
so that 
\begin{equation*}
\begin{split}
   R_{\Delta_3} & = R_{\Delta_1}(I-(1-\lambda)A) = R_{\Delta_1}(I-(1-\beta_\lambda)(I-R_{\Delta_1}^{-1}R_{\Delta_2}))\\
   & = R_{\Delta_1} - (1-\beta_\lambda)R_{\Delta_1}+(1-\beta_\lambda)R_{\Delta_1}R_{\Delta_1}^{-1}R_{\Delta_2}  = \beta_\lambda R_{\Delta_1} + (1-\beta_\lambda)R_{\Delta_2}.
\end{split}
\end{equation*}

Since 
\begin{equation*}
    1-\beta_\lambda = \frac{1-\lambda}{1-\lambda v^\top R_{\Delta_2} u} , \; A =\frac{1-\beta_\lambda}{1-\lambda}(I-R_{\Delta_1}^{-1}R_{\Delta_2}),
\end{equation*}
we have 
\begin{equation}
  R_{\Delta_1}A = \frac{1-\beta_\lambda}{1-\lambda}(R_{\Delta_1} - R_{\Delta_2}).   \label{equ:R1_A}
\end{equation}
By \eqref{eq:itais2} we have
\begin{equation*}
    \lambda R_{\Delta_2} A = R_{\Delta_3} - R_{\Delta_2} = \beta_\lambda R_{\Delta_1} + (1-\beta_\lambda)R_{\Delta_2} - R_{\Delta_2} = \beta_\lambda(R_{\Delta_1}- R_{\Delta_2}), 
\end{equation*}
so that 
\begin{equation}
  R_{\Delta_2} A = \frac{\beta_\lambda}{\lambda}(R_{\Delta_1} - R_{\Delta_2}).   
 \label{equ:R2_A}
\end{equation}
Consequently, from \eqref{equ:R1_A} and \eqref{equ:R2_A},
\begin{equation}
\label{equ:R1A_R2A}
    \frac{R_{\Delta_1} A }{R_{\Delta_2} A }= \frac{(1-\beta_\lambda)\lambda}{(1-\lambda)\beta_\lambda}.
\end{equation}
Since
\begin{equation*}
    \beta_\lambda=1-\frac{1-\lambda}{1-\lambda v^\top R_{\Delta_2} u} = \frac{\lambda(1-v^\top R_{\Delta_2} u )}{1-\lambda v^\top R_{\Delta_2} u},  
\end{equation*}
we have 
\begin{equation*}
\begin{split}
    \frac{R_{\Delta_1} A }{R_{\Delta_2} A } &= \frac{(1-\beta_\lambda)\lambda}{(1-\lambda)\beta_\lambda} 
     = \frac{\lambda}{1-\lambda}\frac{\displaystyle\frac{1-\lambda}{1-\lambda v^\top R_{\Delta_2} u}}{\displaystyle\frac{\lambda(1-v^\top R_{\Delta_2} u )}{1-\lambda v^\top R_{\Delta_2} u}} = \frac{1}{1-v^\top R_{\Delta_2} u}
     = \frac{1}{1-(R_{\Delta_2})_{k.}((R_{\Delta_2}^{-1})_{.k} - (R_{\Delta_1}^{-1})_{.k})} \\
    & = \frac{1}{(R_{\Delta_2})_{k.}(R_{\Delta_1}^{-1})_{.k}} 
     = \frac{det(R_{\Delta_2}^{-1})}{det(R_{\Delta_1}^{-1})(R_{\Delta_1})_{k.}(R_{\Delta_1}^{-1})_{.k}}  = \frac{det(R_{\Delta_2}^{-1})}{det(R_{\Delta_1}^{-1})} 
     = \frac{det(R_{\Delta_1})}{det(R_{\Delta_2})},
\end{split}
\end{equation*}
where the fourth equation follows from $v^\top R_{\Delta_2} = (R_{\Delta_2})_{k.}$, $u = (R_{\Delta_2}^{-1})_{.k} - (R_{\Delta_1}^{-1})_{.k}$, and the sixth equation follows \eqref{equ(mainBody):R1_R2_proportional_in_k_row}.

Thus, from \eqref{equ:R1A_R2A},
\begin{equation*}
    \frac{1-\beta_\lambda}{\beta_\lambda} = \frac{1-\lambda}{\lambda}\frac{det(R_{\Delta_1})}{det(R_{\Delta_2})} \;\Rightarrow \; \beta_\lambda = \frac{\lambda det(R_{\Delta_2})}{\lambda det(R_{\Delta_2})+ (1-\lambda) det(R_{\Delta_1})}.   
\end{equation*} 
Because 
\begin{equation*}
    \frac{\partial \beta_\lambda}{\partial \lambda}=\frac{det(R_{\Delta_2})det(R_{\Delta_1})}{(\lambda det(R_{\Delta_2})+ (1-\lambda) det(R_{\Delta_1}))^2},   
\end{equation*}
and its sign is determined by $det(R_{\Delta_2})det(R_{\Delta_1})$, $\beta_\lambda$ is monotone in $\lambda$. Since $\beta_{\lambda=0} = 0$, $\beta_{\lambda=1} = 1$, we conclude that $\beta_\lambda \in (0,1)$. 
In the analysis below $\lambda$ is fixed, so to simplify notation, we define $\beta =\beta_\lambda$.

\vspace{0.2cm}

Since $\Delta_1$ and $\Delta_2$ are ratio-matrix neighbors as defined in Definition \ref{def: policy neighbors} and differ only in the $k^{th}$ column, the matrices $R_{\Delta_1}^{-1} =CMQ\Delta_1$,  $R_{\Delta_2}^{-1}=CMQ\Delta_2$,   $R_{\Delta_3}^{-1}=\lambda R_{\Delta_1}^{-1} +(1-\lambda) R_{\Delta_2}^{-1}$ as well differ only in the $k^{th}$ column, and 
\begin{equation*}
\begin{split}
    & det(R_{\Delta_3}^{-1}) = \lambda det(R_{\Delta_1}^{-1}) + (1-\lambda) det(R_{\Delta_1}^{-1})
    \Leftrightarrow \;\;
    \frac{1}{det(R_{\Delta_3})} = \frac{\lambda}{det(R_{\Delta_1})} + \frac{1-\lambda}{det(R_{\Delta_2})} \\
     \Leftrightarrow \;\; 
    & \lambda det(R_{\Delta_2}) +(1-\lambda) det(R_{\Delta_1}) = \frac{det(R_{\Delta_2})det(R_{\Delta_1})}{det(R_{\Delta_3})}.  
\end{split}
\end{equation*}
Thus, we have
\begin{equation*}
      \beta = \displaystyle\frac{\lambda det(R_{\Delta_2})}{\lambda det(R_{\Delta_2})+ (1-\lambda) det(R_{\Delta_1})} = \frac{\lambda det(R_{\Delta_3})}{det(R_{\Delta_1})}, \text{ and }
      1-\beta = \displaystyle\frac{(1-\lambda) det(R_{\Delta_3})}{det(R_{\Delta_2})}.
\end{equation*}
Consequently, $R_{\Delta_3} = \beta R_{\Delta_1} + (1-\beta)R_{\Delta_2}$ can also be written as
 \begin{equation*}
     \frac{R_{\Delta_3}}{det(R_{\Delta_3})} = \lambda \frac{R_{\Delta_1}}{det(R_{\Delta_1})} + (1-\lambda) \frac{R_{\Delta_2}}{det(R_{\Delta_2})}. 
\end{equation*} \eProof
\end{proof}

\begin{proof}{Proof of Theorem \ref{thm:Chen-S for a line}.}
\label{proof:Chen-S for a line}
We first consider the completely-$\calS$ property. For fixed $\Delta$ (in turn fixed $R = (CMQ\Delta)^{-1}$), it follows from the definition of completely-$\calS$ that the property holds if the following LP has a strictly positive optimal value:
\begin{equation}
\begin{split}
P(\Delta) = 
{\max_{\epsilon, u^b}}\quad &\epsilon\\
s.t.\quad & R_bu^b\ge \epsilon e^b \qquad \forall  b \subseteq J\\
& u^b\ge e^b \qquad \forall  b \subseteq J\\
& \epsilon \le 1 \\
& \epsilon \in \mathbb{R}, u^b \in \mathbb{R}^{|b|}
\end{split}.
\label{LP:completelyS_prime}
\end{equation}

The constraint $\epsilon\leq 1$ is added here to guarantee that the problem is bounded. Otherwise, if there exists a feasible solution $(\epsilon,u^b)$ with $\epsilon>0$, then $(\kappa\epsilon,\kappa u^b)$ is feasible for any $\kappa>0$ so the problem would be unbounded.  

With this constraint, if $P(\Delta) >0$, then $P(\Delta) =1$. Also notice that the LP is always feasible. Indeed, $u^b = e^b$, $\forall b \in J$ and $\epsilon = \max\{x: \; x\leq 1,x \leq R_be^b,\; \forall b \subseteq J \}$ form one feasible solution. 

The dual of \eqref{LP:completelyS_prime} is 
\begin{equation}
\begin{aligned}
D(\Delta) = \;
& \underset{z^b}{\text{min}}
&&  1 - \sum_{b\subseteq J}e^{b\top}z^b + \sum_{b\subseteq J}e^{b\top}(R_\Delta)_b^\top z^b\\
& \text{s.t.}
&& z^b \in \mathbb{R}_{+}^{|b|}\\
\end{aligned}.
\label{LP:completelyS_dual}
\end{equation}
 An obvious feasible solution of \eqref{LP:completelyS_dual} is $z^b = 0$, and the corresponding objective value is 1. Thus, the dual problem is always feasible. As a result, strong duality holds. 
 
 Since $R_{\Delta_1}$ and $R_{\Delta_1}$ are both completely-$\calS$, $P(\Delta_1) = P(\Delta_2) = 1$. Strong duality then guarantees that $D(\Delta_1)=P(\Delta_1)=1$, $D(\Delta_2)=P(\Delta_2)=1$. Because $\Delta_1$ and $\Delta_2$ are ratio-matrix neighbors, we have by Theorem \ref{thm: linear combination for R}, that 
$R_{\Delta} = \beta R_{\Delta_1} +(1-\beta)R_{\Delta_2}$ for $\Delta = \lambda \Delta_1 +(1-\lambda) \Delta_2$. Then,
\begin{equation*}
\begin{aligned}
D(\Delta) = \;
& \underset{z^b}{\text{min}}
&&  \beta \left(1- \sum_{b\subseteq J}e^{b\top}z^b
+ \sum_{b\subseteq J}e^{b\top}(R_{\Delta_1})_b^\top z^b\right) + (1-\beta) \left(1 - \sum_{b\subseteq J}e^{b\top}z^b  +
\sum_{b\subseteq J}e^{b\top} (R_{\Delta_2})_b^\top z^b \right)
\\
& \text{s.t.}
&& z^b \in \mathbb{R}_{+}^{|b|}\\
\end{aligned},
\end{equation*}

\begin{equation*}
\begin{aligned}
D(\Delta_1) = \;
& \underset{z^b}{\text{min}}
&&  1 - \sum_{b\subseteq J}e^{b\top}z^b + \sum_{b\subseteq J}e^{b\top}(R_{\Delta_1})_b^\top z^b\\
& \text{s.t.}
&& z^b \in \mathbb{R}_{+}^{|b|}\\
\end{aligned},\quad
\begin{gathered}
\begin{aligned}
D(\Delta_2) = \;
& \underset{z^b}{\text{min}}
&&  1 - \sum_{b\subseteq J}e^{b\top}z^b + \sum_{b\subseteq J}e^{b\top}(R_{\Delta_2})_b^\top z^b\\
& \text{s.t.}
&& z^b \in \mathbb{R}_{+}^{|b|}\\
\end{aligned}.
\end{gathered}
\end{equation*}

Because the feasible regions of $D(\Delta_1)$, $D(\Delta_2)$ and $D(\Delta)$ are all $\{z^b: z^b \in \mathbb{R}_{+}^{|b|} \}$, $D(\Delta)  \geq \beta D(\Delta_1) + (1-\beta) D(\Delta_2) = \beta +(1-\beta) = 1$. 

Due to strong duality, we then have that $1\leq D(\Delta) = P(\Delta) \leq 1$. It must then hold that $D(\Delta) =P(\Delta)= 1$, and thus the reflection matrix for $\Delta$ is completely-$\calS$. 

\vspace{0.2cm}


According to Lemma \ref{lem: RO_ChenS}, for a fixed $\Delta$, Chen-$\calS$ holds if the following LP has strictly positive optimal value:
\begin{equation}
\label{LP:ChenS_prime}
\begin{split}
{\max_{h,\alpha^b,\epsilon, \eta^b}}\quad\epsilon\\
s.t.\quad &h_a^\top \theta_a - \alpha^{b\top} \theta_b \leq -\epsilon \qquad \forall (a,b)\\
& h_a^\top R_{ab} \leq \alpha^{b\top}R_b  \qquad \forall (a,b) \\
& \epsilon e\le R_b\eta^b \qquad \forall  b \subseteq J\\
& \alpha^b \in \mathbb{R}^{|b|}, h\geq e, \epsilon \leq M , \eta^b \geq e^b
\end{split}.
\end{equation}

\noindent The dual LP of \eqref{LP:ChenS_prime} is
\begin{equation*}
\begin{aligned}
& \underset{}{\text{min}}
&& e^{J\top} \sum_{(a,b)}N^{a\top}(R_{ab}R_b^{-1}\theta_b-\theta_a)x^a +M-\sum_{(a,b)}x^a - \sum_{b\subseteq J}e^{b\top}z^b + \sum_{b\subseteq J}e^{b\top}R_b^\top z^b\\
& \text{s.t.}
&& x^a \in \mathbb{R}_{+}, z^b \in \mathbb{R}_{+}^{|b|}\\
\end{aligned},
\end{equation*} 
where the matrix $N^a$ maps the column vector to its sub-vector with the indices in the set $a$.

From here, an identical argument to that for completely-$\calS$ proves that 
if $(R_{\Delta_1},\theta_{\Delta_1})$ and $(R_{\Delta_2},\theta_{\Delta_2})$ are both Chen-$\calS$, so is $(R_{\Delta_3},\theta_{\Delta_3})$; we omit the details.
\eProof
\end{proof}

\vspace{0.2cm}

\begin{proof}
{Proof of Theorem \ref{thm: LEGO paper (Theorem 2)}.}
\label{proof: LEGO paper (Theorem 2)}

Let $G(u) = \alpha - (I-P')M^{-1}u$, $F(u) = (I-\Delta CM)G(u)$; with this notation, $\dot{\wbar{Z}}(t)= G(\dot{\wbar{T}}(t))$. 
A solution $\wbar{X} = (\wbar{Z}(t),\wbar{W}(t),\wbar{\epsilon}(t),\wbar{T}(t))$ to the fluid model equations  \eqref{equ:begin fluid model}-\eqref{equ:end fluid model} satisfies, for any regular $t \geq 0$,  $$(\wbar{Z}(t),\wbar{\epsilon}(t)) \in \Bbb{X},~~~ \dot{\wbar{T}}(t) \in\Bbb{U}_\Delta(\wbar{Z}(t),\wbar{\epsilon}(t)), ~~~\dot{\wbar{Z}}(t) \in \{G(u): u \in \Bbb{U}_\Delta(\wbar{Z}(t),\wbar{\epsilon}(t))\},$$ where $\Bbb{X} = \{(z,\epsilon)\in \Bbb{R}_+^{K}\times \Bbb{R}^K: \epsilon = z-\Delta CMz\}$ is the family of possible states for $\wbar{Z}$ and $\wbar{\epsilon}$, and $$\Bbb{U}_\Delta(z,\epsilon) = \{u \in \Bbb{R}_+^K: e-Cu \geq 0, (CMz)'(e-Cu) = 0, \epsilon_k^+(1-u_k^+)=0, \mathbbm{1}\{z_k = 0\}G_k(u)\ = 0 \text{, for all } k \in \mathcal {K}\},$$
is the set of possible controls $u$ in state $(z,\epsilon)$. Define $\wbar{Z}_k^+(t) = \sum_{i\geq k:s(i) = s(k), \wbar{Z}_i(t) >0}\wbar{Z}_i(t)$. For a static-priority policy $\pi$, $\epsilon_k^+=z_k^+$ for all $k \in \mathcal{H}_\pi$, and 
$z_k^+=0$ implies $\epsilon_k^+=0$ for all $k \in \mathcal{L}_\pi$,
so that the set $\Bbb{U}_\Delta(z,\epsilon)$ is, for static-priority policies, replaced by $$\Bbb{U}_\pi(z):= \{u \in \Bbb{R}_+^K: e-Cu \geq 0, (CMz)'(e-Cu) = 0, z_k^+(1-u_k^+)=0, \mathbbm{1}\{z_k = 0\}G_k(u)\ = 0 \text{, for all } k \in \mathcal {K}\}.$$ We note that $\Bbb{U}_\pi(z)$ depends only on the sign of $z$, not its actual value. 

Recall that $\dot{\wbar{Z}}(t) \in \{G(u): u \in \Bbb{U}_\Delta(\wbar{Z}(t),\wbar{\epsilon}(t))\}$. To prove the theorem it then suffices to establish the existence of an $h\in \mathbb{R}^K_{++}$ and a constant $r >0$ such that, for any $\Delta$, 
\begin{align}
\sum_{k:\epsilon_k>0}h_kG_k(u) \leq - r, \text{ for any state } (z,\epsilon) \in \Bbb{X}, \text{ with } \lVert \epsilon \rVert >0 \text{ and all } u \in \Bbb{U}_\Delta(z,\epsilon)\text{.}  \label{SSC_proof_eq_1}
\end{align}
We then have for all regular time with $\lVert \wbar{\epsilon}(t)\rVert >0$, 
$$-r \geq \sum_{k:\wbar{\epsilon}_k(t)>0}h_k \dot{\wbar{Z}}_k(t) = \sum_{j\in \mathcal{J}} \sum_{k\in \mathcal{C}(j):\wbar{\epsilon}_k(t)>0} h_k(\dot{\wbar{\epsilon}}_k(t)+\delta_{kj}\dot{\wbar{W}}_j(t))\text{.}$$
Since $\wbar{W}(t)$ and $\sum_{k\in \mathcal{C}(j):\wbar{\epsilon}_k(t)>0} h_k\wbar{\epsilon}_k(t)$ are non-negative, there exists $t_0 < \infty$ such that for all $t\geq t_0$, either $\wbar{W}_j(t) = (CM\wbar{Z}(t))_j = 0$ or $\sum_{k\in \mathcal{C}(j):\wbar{\epsilon}_k(t)>0} \wbar{\epsilon}_k(t) = 0$. In the former case, since $\wbar{\epsilon}_k(t) \leq \wbar{Z}_k(t)$, this implies $\wbar{\epsilon}_k(t) = 0$. Thus, in either case, $\wbar{\epsilon}_k(t) = 0$ for all $t \geq t_0$. 

We prove next that $h$ and $r$ satisfying  \eqref{SSC_proof_eq_1} exist if linear attraction holds for all static priority policies with common $h$ and $r$. We divide the possible values of $(z,\epsilon)\in \mathbb{R}^K_{+}$ into three cases. In all cases, the argument proceeds by showing that the set of feasible controls $U_{\Delta}(z,\epsilon)$ is the same as $\Bbb{U}_{\pi}(z,\epsilon)$ for some priority permutation $\pi$ and then leveraging the assumed linear attraction of all priority policies. 

\vspace{0.2cm}

\noindent \textbf{Case 1:} For each station $j \in \mathcal{J}$, there exists a class $i\in \mathcal{C}(j)$ with $\epsilon_i >0$.

In this case no station is empty. Let $\mathcal{H}_\Delta(z,\epsilon) = \{i \in  \mathcal{K}: \epsilon_i >0\}$. We construct a policy priority permutation $_\epsilon\pi$ where the classes in $\mathcal{H}_\Delta(z,\epsilon)$ take the highest priorities in their stations. 
Then, $\delta_i = 0$ for all $i \in  \mathcal{K}$ with $\epsilon_i>0$, and $\epsilon_i = z_i - \delta_i(CMZ)_{s(i)} = z_i$. Thus, if $\epsilon_i>0$, then $z_i >0$.
Given $i \in \mathcal{H}_\Delta(z,\epsilon)$, there exists no $k \not\in \mathcal{H}_\Delta(z,\epsilon)$ with $s(k) = s(i) = j$ such that $_\epsilon\pi_k(j) \geq {_\epsilon\pi_i(j)}$. For each $j \in \mathcal{J}$, pick $i\in \mathcal{C}(j)$ with $\epsilon_i <0$ and assign it as the lowest priority class in station $j$. Such an $i$ exists by the assumption that there exists an $i$ with $\epsilon_i>0$ for each station and the fact that $CM\epsilon = 0$ for any $(z,\epsilon) \in \Bbb{X}$. Then we assign the remaining classes to priorities arbitrarily. Let $\mathcal{L}_{_\epsilon\pi}$ be the classes with the lowest priority in each station, and $\mathcal{H}_{_\epsilon\pi} = \mathcal{K} \backslash \mathcal{L}_{_\epsilon\pi}$. 

Since $z_i = \epsilon_i - \delta_i(CMZ)_{s(i)} >0$ for all $i$ with $\epsilon_i >0$, with the definition of $_\epsilon\pi$, we have $\Bbb{U}_{_\epsilon\pi}(z) = \Bbb{U}_\Delta(z,\epsilon)$. According to $_\epsilon\pi$, we have $u_k=1$ for $k \in \mathcal{H}_\Delta(z,\epsilon)$ and $u_k = 0$ otherwise. In particular, by the assumption of linear attraction under static-priority policies, we have that 
\begin{align}
\sum_{k \in \mathcal{H}_\Delta(z,\epsilon)} h_kG_k(u) \leq \sum_{k\in\mathcal{H}_{_\epsilon\pi}(z):z_k >0}h_kG_k(u) \leq -r \text{, for all } u\in  \Bbb{U}_{_\epsilon\pi}(z)= \Bbb{U}_\Delta(z,\epsilon)\text{.} \label{SSC_proof_eq_2}
\end{align}
The first inequality follows because $\mathcal{H}_\Delta(z,\epsilon) \subseteq \mathcal{H}_{_\epsilon\pi}(z)$ 
and $u_k = 0$ (so that $G_k(u) \geq 0$) for all $k \not\in \mathcal{H}_\Delta(z,\epsilon)$. The second inequality follows from the assumed linear attraction under static-priority policies. Thus, we have (\ref{SSC_proof_eq_1}) in this case.

\vspace{0.2cm}

The remaining cases are when $(z,\epsilon)$ are such that $j\in \mathcal{J}$ exists with $\epsilon_i = 0$ for all $i \in \mathcal{C}(j)$. We assume that there is a single such station, call it station $j$. The same argument extends easily to the case of multiple such stations.

\vspace{0.2cm}

\noindent \textbf{Case 2:} $\epsilon_i = 0$ for all $i \in \mathcal{C}(j)$, and $z_i = 0$ for all $i \in \mathcal{C}(j)$. 

In this case, since $\epsilon_i = z_i = 0, \mathcal{H}_\Delta(z,\epsilon) = \mathcal{H}_\Delta(z,\epsilon)/\mathcal{C}(j)$ and $\mathcal{H}_{_\epsilon\pi}(z) = \mathcal{H}_{_\epsilon\pi}(z)/\mathcal{C}(j)$, neither of the sums in (\ref{SSC_proof_eq_2}) include $\mathcal{C}(j)$. Also, we still have $\Bbb{U}_{_\epsilon\pi}(z) = \Bbb{U}_\Delta(z,\epsilon)$ since the constraint $ \mathbbm{1} \{z_k = 0\}G_k(u) = 0$ is active for these classes in both $\Bbb{U}_{_\epsilon\pi}(z)$ and $\Bbb{U}_\Delta(z,\epsilon)$. Thus, the argument in Case $1$ also applies.

\vspace{0.2cm}

\noindent \textbf{Case 3:} $\epsilon_i = 0$ for all $i \in \mathcal{C}(j)$, but $z_i >0$ for some $i \in \mathcal{C}(j)$.

Under static priority, we have in station $j$ that $\delta_{kj} = 0$ for all but one class, and we can use again the argument from Case 1. Thus, we assume that there are distinct classes $k,l \in \mathcal{C}(j)$ for which $\delta_k, \delta_l >0$. Since $\epsilon_i = 0$ for all $i \in \mathcal{C}(j)$, we have $z_k, z_l >0$. To complete this case, we need two auxiliary lemmas. 



\vspace{-0.2cm}
\begin{lemma}
\label{lemma: LEGO paper (Lemma 15)}
Consider the queue-ratio fluid network model (\ref{equ:begin fluid model})-(\ref{equ:end fluid model}). Let $t$ be a regular time. If $j$ is such that $\sum_{k\in\mathcal{C}(j)}\wbar{\epsilon}_k(t) = 0$, it must be the case that $\dot{\wbar{\epsilon}}_k(t) = 0$ for all $k \in \mathcal{C}(j)$. In particular, at every regular time $t$, we have $\dot{T}(t) \in \tilde{\Bbb{U}}_\Delta(\wbar{Z}(t),\wbar{\epsilon}(t))$ where
$$\tilde{\Bbb{U}}_\Delta(z,\epsilon) := \Bbb{U}_\Delta(z,\epsilon) \cap \{u \in \Bbb{R}_{+}^{K}: \mathbbm{1}\{\epsilon_k = 0, k\in \mathcal{C}(j)\}F_k(u)=0 \text{ for all }j\in \mathcal{J} \text{ and }k\in \mathcal{C}(j)\}\text{.} $$
\end{lemma}

\vspace{-0.2cm}

We let $z^{-j}$ be the entries of the vector $z\in\Bbb{R}_{+}^{K}$ corresponding to entries $k \not\in \mathcal{C}(j)$. Similarly, $_\epsilon\pi^{-j}$ is the restriction of $_\epsilon\pi$ to stations other than station $j$.


\begin{lemma}
\label{lemma: LEGO paper (Lemma 16)}
Fix $(z,\epsilon)$ and let $j$ be the only station such that $\epsilon_i = 0$ for all $i \in \mathcal{C}(j)$. Suppose that there exist $u \in \Bbb{U}_\Delta(z,\epsilon), j \in \mathcal{J}$ and $k \in \mathcal{C}(j)$ so that $G_k(u) \leq 0$, $\delta_k >0$ and $F_l(u) = 0$ for all $l\in \mathcal{C}(j)$. Then, $G_l(u) \leq 0$ for all $l \in \mathcal{C}(j)$. Moreover, $u$ is a convex combination of vectors in $\cap_{\mathcal{Q}}\Bbb{U}_{(_\epsilon\pi^{-j},\pi(j))}(z^{-j},x^j)$, where $\mathcal{Q} = \{ \pi(j),x^j:x_k^j = 0 \text{ for } k\in \mathcal{H}_{(_\epsilon\pi^{-j},\pi(j))}\}$.
\end{lemma}

The latter part of Lemma \ref{lemma: LEGO paper (Lemma 16)} means that $u_j$, the allocation in station $j$, is a convex combination of allocations under static priorities in station $j$ when the high-priority queues are empty.

\vspace{0.2cm}

Since $\epsilon_i = 0$ for all $i \in \mathcal{C}(j)$, Lemma \ref{lemma: LEGO paper (Lemma 15)} implies that $F_k(u) = 0$. Let $u \in \Bbb{U}_\Delta(z,\epsilon)$, then $u_k = 1$ for all $k \not\in \mathcal{C}(j)$ since $\epsilon_i >0$ in each such station. (1) If $G_k(u) \geq 0$ for all $k$ with $\delta_k >0$ and all $u$ in this station, then this is true also under static priority cases because $G_k(u) = 0$ implies $F_k(u) = 0$. Thus, from linear attraction for static priority policies, we know that the queues with $\epsilon_i >0$ will decrease. (2) If there exists a $u$ such that $G_k(u) \leq 0$, then for station $j$, from Lemma \ref{lemma: LEGO paper (Lemma 16)}, $u$ is a convex combination of static priorities in this station. Thus, $u_{z,\epsilon} = \sum_{\pi',z^j}\alpha_{\pi',z}u_{\pi',z^j,z^{-j}}$. Since $G_k(u)$ is a linear (in particular, convex) function of $u$, we have 
$$\sum_{k\in\mathcal{H}_\Delta(z,\epsilon)}h_kG_k(u) = \sum_{\pi',z}\alpha_{\pi',z}\sum_{k\in\mathcal{H}_{\pi'}:z_k >0}h_kG_k(u_{\pi',z}) \leq - r\text{.}$$

Note that the classes counted in $\{k\in\mathcal{H}_\Delta(z,\epsilon)\}$ and $\{k\in\mathcal{H}_{\pi'}:z_k >0\}$ are the same. The reason is that the classes in station $j$ all have $\epsilon_i = 0$ and do not appear in $\{k\in\mathcal{H}_\Delta(z,\epsilon)\}$; they also do not appear in $\{k\in\mathcal{H}_{\pi'}:z_k >0\}$ since we constructed the set so that $z_k = 0$ for all $k \in \mathcal{H}_{\pi'}$.

This completes the proof of \eqref{SSC_proof_eq_1} and hence of state-space collapse under any queue-ratio discipline. \eProof
\end{proof}

\vspace{0.2cm}

\section{Summary\label{sec:conclusions} }
\setcounter{equation}{0}
\renewcommand{\theequation}{\thesection.\arabic{equation}}

What we put forward in this paper is a framework to study policy robustness---with regards to stability---in multiclass queueing networks, with the ``uncertainty set'' containing all queue-ratio policies. 
The framework is a hierarchical one that breaks robust stability into two ingredients: (i) robust stability of the Skorohod problem (SP), and (ii) robust state-space collapse (SSC).

We prove that the stability under a policy within this uncertainty set is inherited from that of the ``vertices'' of the set; these ``vertices'' represent static-priority policies. In other words, in our framework, one must only verify the stability of the corner points of the uncertainty set to deduce policy robustness. Our methodology brings together techniques from optimization and applied probability to revisit and make progress on a fundamental question in stochastic modeling.

It is plausible that the specific implementation of each step in our framework can be strengthened. Our treatment of the Skorohod problem is relatively complete, and this is the focus of much of the theoretical development. There is likely room to strengthen the results pertaining to state-space collapse. Our requirement of linear attraction (with a common test vector $h$) is useful for tractability but is, admittedly, a strong one. It is this requirement that necessitates, in our example networks, the consideration of an arrival rate $\alpha$ that is ``large enough.'' Going beyond linear attraction is a natural next step. 

It is also plausible that one can go beyond queue-ratio policies to cover policies that have non-linear state-space collapse (i.e., where the matrix $\Delta$ depends itself on the workload). We leave that, too, for future research. 
\bibliographystyle{abbrv}
\bibliography{paper_references}
\newpage 
\begin{appendices}

\section{Proofs of Supporting Lemmas}
\label{Appendix: Proofs of Supporting Lemmas}
\setcounter{equation}{0}
\renewcommand{\theequation}{\thesection.\arabic{equation}}

\begin{proof} {Proof of Lemma \ref{lem: RO_ChenS}.}
\label{proof:lem_RO_ChenS}
According to Definition \ref{def:Chen-S}, except the completely-$\calS$ condition, for fixed partition $(a,b)$ for $\mathcal{J}$, the Chen-$\calS$ condition is satisfied if the optimization model (\ref{RO:ChenS_opt1}) has strictly negative optimal value.
\begin{equation}
\label{RO:ChenS_opt1}
\begin{split}
\max_{h, v}\quad\quad &h_a^\top \theta_a +h_a^\top R_{ab}v \\
s.t.\quad &R_bv = -\theta_b\\
& h \in \mathbb{R}^{J}_{++},  v \in \mathbb{R}^{|b|}_{+}
\end{split}
\end{equation}

(\ref{RO:ChenS_opt1}) is equivalent to
\begin{equation}
\label{RO:ChenS_opt2}
\begin{split}
{\max_{h\geq e} \quad  h_a^\top \theta_a -\;\; \max_{\alpha^b \in \mathbb{R}^{|b|}}}\quad &\alpha^b\theta_b\\
s.t.\quad & {\alpha^{b}}^{\top}R_b \geq h_a^{\top}R_{ab}
\end{split},
\end{equation}
which is no larger than 
\begin{equation}
\label{RO:ChenS_opt3}
\begin{split}
\max_{h, \alpha^b}\quad\quad &h_a^\top \theta_a -\alpha_b^\top \theta_{b} \\
s.t.\quad &{\alpha^{b}}^{\top}R_b \geq h_a^{\top}R_{ab}\\
& \alpha^b \in \mathbb{R}^{|b|},  h \geq e 
\end{split}.
\end{equation}
Thus, we can rewrite \eqref{RO:ChenS_opt3} and add the completely-$\calS$ constraint, which gives an optimization model \eqref{RO:ChenS_opt4} that if \eqref{RO:ChenS_opt4} has strictly positive optimal value, then Chen-$\calS$ property holds for a fixed queue-ratio policy (so fixed $\Delta$), where $M$ is a large enough positive number.
\begin{equation}
\label{RO:ChenS_opt4}
\begin{split}
{\max_{h,\alpha^b,\epsilon, \eta^b}}\quad\epsilon\\
s.t.\quad &h_a^\top \theta_a - \alpha^{b\top} \theta_b \leq -\epsilon \qquad \forall (a,b)\\
& h_a^\top R_{ab} \leq \alpha^{b\top}R_b  \qquad \forall (a,b) \\
& \epsilon e\le R_b\eta^b \qquad \forall  b \subseteq J\\
& \alpha^b \in \mathbb{R}^{|b|}, h\geq e, \epsilon \leq M , \eta^b \geq e^b
\end{split}
\end{equation}

As a result, the SP is robustly queue-ratio stable if $\rho <e $ and the min-max problem
\begin{equation}
\label{RO:ChenS_final}
\begin{split}
{\min_{\Delta \in \Bbb{R}_{+}^{K\times J}: CM\Delta = I} \;\; \max_{h,\alpha^b,\epsilon, \eta^b}}\quad\epsilon\\
s.t.\quad &h_a^\top \theta_a - \alpha^{b\top} \theta_b \leq -\epsilon \qquad \forall (a,b)\\
& h_a^\top R_{ab} \leq \alpha^{b\top}R_b  \qquad \forall (a,b) \\
& \epsilon e\le R_b\eta^b \qquad \forall  b \subseteq J\\
& \alpha^b \in \mathbb{R}^{|b|}, h\geq e, \epsilon \leq M , \eta^b \geq e^b
\end{split}
\end{equation}
has a strictly positive solution. \eProof
\end{proof}

\begin{proof} {Proof of Lemma \ref{lem:Schur-S is stronger than Chen-S}.}
\label{proof:Schur-S is stronger than Chen-S}
This proof follows from that of Corollary 2.6 in \cite{chen1996sufficient}.
Since $R$ is \textit{Schur-S} all of its principal sub-matrices are non-singular, so that $R_b$ is invertible. Thus, $\left\{ v\in R_{+}^{|b|}: \theta_b + R_bv=0 \right\}$ contains the single vector $u = -R_b^{-1}\theta_b \in R_{+}^{|b|}$. As in the proof of Corollary 2.8 in \cite{chen1996sufficient}, if $R$ is both \textit{completely-$\calS$} and \textit{Schur-S}, there exists a positive vector $h$ such that $h_a^{'}\left[\theta_a-R_{ab}R_b^{-1}\theta_b\right]<0$ for any partition $(a,b)$ of $\mathcal{J}$. This same $h$ satisfies, in particular, 
$h_a^{'}\left[\theta_a + R_{ab}u\right] <0$
for any partition $(a,b)$ of $\mathcal{J}$ with $u = -R_b^{-1}\theta_b \in R_{+}^{|b|}$, and thus, all $u \in \left\{ v\in R_{+}^{|b|}: \theta_b + R_bv=0 \right\}$. The pair $(R, \theta)$ is therefore Chen-$\calS$. \eProof
\end{proof}

\vspace{0.2cm}
\begin{proof} {Proof of Lemma \ref{lem:decomposition with A}.}
\label{proof:decomposition with A}
First, we prove that
$R_{\Delta_1}$ and $R_{\Delta_2}$ are proportional in the $k^{th}$ row with the ratio
\begin{equation}
\label{equ:R1_R2_proportional_in_k_row}
(R_{\Delta_1})_{k.} = \frac{det(R_{\Delta_2}^{-1})}{det(R_{\Delta_1}^{-1})} (R_{\Delta_2})_{k.}.
\end{equation}
Let $(R_{\Delta_1}^{-1})^*$ and $(R_{\Delta_2}^{-1})^*$ be the adjugate matrices of $R_{\Delta_1}^{-1}$ and $R_{\Delta_2}^{-1}$ respectively, then with the known fact,
\begin{equation}
\label{equ:R_and_adjugent_matrix}
    R_{\Delta_1} = (R_{\Delta_1}^{-1})^{-1} = \frac{(R_{\Delta_1}^{-1})^*}{det(R_{\Delta_1}^{-1})} ,  \;   R_{\Delta_2} = (R_{\Delta_2}^{-1})^{-1} = \frac{(R_{\Delta_2}^{-1})^*}{det(R_{\Delta_2}^{-1})},
\end{equation}

In our settings, $R_{\Delta_1}^{-1}$ and $R_{\Delta_2}^{-1}$ differ only in the $k^{th}$ column, whose cofactor matrices $Co(R_{\Delta_1}^{-1})$ 
and $Co(R_{\Delta_2}^{-1})$ then have an identical $k^{th}$ column,
and $(R_{\Delta_1}^{-1})^* = Co^\top(R_{\Delta_1}^{-1})$ and $(R_{\Delta_2}^{-1})^* = Co^\top(R_{\Delta_2}^{-1})$ have the same $k^{th}$ row.

Thus, from \eqref{equ:R_and_adjugent_matrix}, we have 
\begin{equation*}
\begin{split}
     (R_{\Delta_1})_{k.} 
     & = \frac{(R_{\Delta_1}^{-1})^*_{k.}}{det(R_{\Delta_1}^{-1})}
     =\frac{(R_{\Delta_2}^{-1})^*_{k.}}{det(R_{\Delta_1}^{-1})}
      =\frac{det(R_{\Delta_2}^{-1})}{det(R_{\Delta_1}^{-1})}(R_{\Delta_2})_{k.},
\end{split}
\end{equation*} 
which implies that $R_{\Delta_1}$ and $R_{\Delta_2}$ are proportional in the $k^{th}$ row with the ratio in \eqref{equ:R1_R2_proportional_in_k_row}.

Second, from \eqref{equ:R3_with_uv}, using the 
Sherman-Morrison formula, we can write
\begin{equation*}
    R_{\Delta_3} = (R_{\Delta_1}^{-1} + (1-\lambda)uv^\top)^{-1} = R_{\Delta_1} - \frac{(1-\lambda)R_{\Delta_1}uv^\top R_{\Delta_1}}{1+(1-\lambda)v^\top R_{\Delta_1} u}.
\end{equation*}
Similarly, from \eqref{equ:R3_with_uv_2}, 
\begin{equation*}
R_{\Delta_3} = (R_{\Delta_2}^{-1} -\lambda uv^\top)^{-1} = R_{\Delta_2} + \frac{\lambda R_{\Delta_2}uv^\top R_{\Delta_2}}{1-\lambda v^\top R_{\Delta_2} u}.
\end{equation*}
 Suppose
 \begin{equation*}
      A_1 = \frac{uv^\top R_{\Delta_1}}{1+(1-\lambda)v^\top R_{\Delta_1} u} , \;  A_2 = \frac{uv^\top R_{\Delta_2}}{1-\lambda v^\top R_{\Delta_2} u}, 
 \end{equation*}
then we get $R_{\Delta_3} = R_{\Delta_1}(I-(1-\lambda)A_1) = R_{\Delta_2}(I+\lambda A_2)$.

Now, we prove that $A_1=A_2$. Let
\begin{equation*}
    E_1 = \frac{R_{\Delta_1}}{1+(1-\lambda)v^\top R_{\Delta_1} u} , E_2 = \frac{R_{\Delta_2}}{1-\lambda v^\top R_{\Delta_2} u}.
\end{equation*}
Since only the $k^{th}$ column of the matrix $uv^\top$ is non-zero, if we prove that $E_1$ and $E_2$ have the same $k^{th}$ row, then we can conclude that 
\begin{equation*}
    A_1=uv^\top E_1 = uv^\top E_2 = A_2.
\end{equation*}
Using \eqref{equ:R1_R2_proportional_in_k_row}, to show that $E_1$ and $E_2$ have an identical $k^{th}$ row, it suffices to show 
\begin{equation}
    \frac{1+(1-\lambda)v^\top R_{\Delta_1} u}{1-\lambda v^\top R_{\Delta_2} u} = \frac{det(R_{\Delta_2}^{-1})}{det(R_{\Delta_1}^{-1})}.\label{eq:1}
\end{equation}
Define 
\begin{equation*}
    B_1 = v^\top R_{\Delta_1} = (R_{\Delta_1})_k., \; B_2 = v^\top R_{\Delta_2} = (R_{\Delta_2})_k. , \; w = \frac{det(R_{\Delta_2}^{-1})}{det(R_{\Delta_1}^{-1})}.
\end{equation*}    
Then, we have $B_1 = wB_2$ and
\begin{equation*}
    \frac{1+(1-\lambda)v^\top R_{\Delta_1} u}{1-\lambda v^\top R_{\Delta_2} u} = \frac{1+(1-\lambda)B_1 u}{1-\lambda B_2 u} = \frac{1+(1-\lambda)wB_2 u}{1-\lambda B_2 u}.
\end{equation*}
Thus, to establish \eqref{eq:1}, it suffices to prove that
\begin{equation*}
    \frac{1+(1-\lambda)wB_2 u}{1-\lambda B_2 u} = w.
\end{equation*}
Namely, we have
\begin{equation*}
\begin{split}
    & 1+(1-\lambda)wB_2 u = w-w\lambda B_2 u\\ \Leftrightarrow \; & 1+wB_2 u - \lambda w B_2 u = w-w\lambda B_2 u \\
    \Leftrightarrow \; & 1+wB_2u = w\\ \Leftrightarrow \; & w(1-B_2u) = 1\\ \Leftrightarrow \; & w(1-B_2u)B_2 = B_2 \\
    \Leftrightarrow \; & B_1 - B_2uB_1 = B_2.
\end{split}    
\end{equation*} 
Since $uv^\top = CMQ(\Delta_2- \Delta_1) = R_{\Delta_2}^{-1}-R_{\Delta_1}^{-1}$, we have
\begin{equation*}
    u = (R_{\Delta_2}^{-1})_{.k} - (R_{\Delta_1}^{-1})_{.k}.
\end{equation*}
As
\begin{equation*}
    B_1(R_{\Delta_1}^{-1})_{.k} = 1, \; B_2(R_{\Delta_2}^{-1})_{.k} = 1, 
\end{equation*}
we have
\begin{equation*}
    B_2uB_1 = B_2((R_{\Delta_2}^{-1})_{.k} - (R_{\Delta_1}^{-1})_{.k})B_1 = (1-B_2(R_{\Delta_1}^{-1})_{.k}) B_1  = B_1 - B_2.
\end{equation*}
Overall, we have $B_1 - B_2uB_1 = B_1 - (B_1-B_2) = B_2$, and we can conclude that $A_1 = A_2$. Define $A = A_1 = A_2$, then $R_{\Delta_3} = R_{\Delta_1}(I-(1-\lambda)A) = R_{\Delta_2}(I+\lambda A)$, where $A = \displaystyle\frac{uv^\top R_{\Delta_1}}{1+(1-\lambda)v^\top R_{\Delta_1} u} = \frac{uv^\top R_{\Delta_2}}{1-\lambda v^\top R_{\Delta_2} u}$.\eProof
\end{proof}

\vspace{0.4cm}

\begin{proof}{Proof of Lemma \ref{lem: ensure R_Delta invertible}.}
\label{proof: ensure R_Delta invertible}

We suppose that $R_{\Delta_1}^{-1}$ and $R_{\Delta_2}^{-1}$ are invertible. Using the Sherman-Morrison formula, we have that $R_{\Delta_3}^{-1} = R_{\Delta_1}^{-1} + (1-\lambda)uv^\top$ is invertible if and only if
\begin{equation*}
    1+(1-\lambda)v^\top R_{\Delta_1}^{-1}u \neq 0.
\end{equation*}
Since $v^\top R_{\Delta_1} = (R_{\Delta_1})_{k.} = \displaystyle\frac{det(R_{\Delta_2}^{-1})}{det(R_{\Delta_1}^{-1})}(R_{\Delta_2})_{k.}$, $u =(R_{\Delta_2}^{-1})_{k.} - (R_{\Delta_1}^{-1})_{k.}$, 

\begin{equation*}
\begin{split}
     & 1+(1-\lambda)v^\top R_{\Delta_1}^{-1}u \neq 0\\
     \Leftrightarrow \;\;
     & 1+(1-\lambda) (R_{\Delta_1})_{k.} ((R_{\Delta_2}^{-1})_{k.} - (R_{\Delta_1}^{-1})_{k.})\neq 0\\
     \Leftrightarrow \;\;
     & 1+(1-\lambda)((R_{\Delta_1})_{k.}(R_{\Delta_2}^{-1})_{k.}-1) \neq 0\\
     \Leftrightarrow \;\;
     & 1+(1-\lambda)\left(\frac{det(R_{\Delta_2}^{-1})}{det(R_{\Delta_1}^{-1})}(R_{\Delta_2})_{k.}(R_{\Delta_2}^{-1})_{k.}-1\right)\neq 0\\
     \Leftrightarrow \;\;
     & 1+(1-\lambda)\left(\frac{det(R_{\Delta_2}^{-1})}{det(R_{\Delta_1}^{-1})}-1\right)\neq 0\\
     \Leftrightarrow \;\;
     & \frac{det(R_{\Delta_1}^{-1})}{det(R_{\Delta_2}^{-1})} \neq 1-\frac{1}{\lambda}.\\
\end{split}    
\end{equation*}


Since $\lambda \in (0,1)$, we have $1-\frac{1}{\lambda} \in (-\infty,0)$. Thus, if $det(R_{\Delta_1}^{-1})$ and $det(R_{\Delta_2}^{-1})$ have the same sign, then $R_{\Delta_3}^{-1}$ is invertible.
\eProof
\end{proof}

\vspace{0.2cm}

\begin{proof}
{Proof of Lemma \ref{lemma: LEGO paper (Lemma 15)}.}
\label{proof: LEGO paper (Lemma 15)}
Fix a regular time $t$, so that $\dot{\wbar{Z}}(t)$ and $\dot{\wbar{\epsilon}}(t)$ exist. Let $j \in \mathcal{J}$ be such that $\wbar{\epsilon}_k(t) = 0$ for all $k\in \mathcal{C}(j)$. Suppose that there exists $k$ such that $\dot{\epsilon}_k(t) <0$. For such $k$, it must be the case that $\wbar{\epsilon}_k(t-) >0$ and $\wbar{\epsilon}_k(t+) <0$. Let $B:= \{k \in \mathcal{C}(j):\wbar{\epsilon}_k(t-) \geq 0\}$. Then $\sum_{k \in B} u_k = 1$ at $t-$. However, since at $t+$, these are negative, it will be the case that $\sum_{k \in B} u_k = 0$ at $t+$ so that $\dot{T}(t)$ does not exist at $t$, which is a contradiction to its regularity. A similar argument can be applied if  $\dot{\epsilon}_k(t) >0$.\eProof
\end{proof}

\vspace{0.2cm}

\begin{proof}
{Proof of Lemma \ref{lemma: LEGO paper (Lemma 16)}.}
\label{proof: LEGO paper (Lemma 16)}

First, note that $F_k(u) = G_k(u) - \delta_k\sum_{l\in\mathcal{C}(j)}m_lG_l(u)$. Since $F_k(u) = 0$, we have $G_k(u)= \delta_k\sum_{l\in\mathcal{C}(j)}m_lG_l(u)$. If $G_k(u)\leq 0$ for some $k$ with $\delta_k >0$, then it must be the case that $\sum_{l\in \mathcal{C}(j)}m_lG_l(u) \leq 0$. If there exists $l$ with $G_l(u) > 0$, then by the above we would have a
contradiction because $0 < G_l(u) = \delta_l\sum_{l\in\mathcal{C}(j)}m_lG_l(u) \leq 0$.

Second, recall that $F(u) = (I -\Delta C M)G(u)$. We can consider the polyhedron of values $y\in \Bbb{R}^{|\mathcal{C}(j)|}$ such that $y_l-\delta_l\sum_{l\in\mathcal{C}(j)}m_iy_i =0$ for all $l\in C(j)$ and $y_l \leq 0$ for all $l\in \mathcal{C}(j)$. The extreme points of this polyhedron have $y_l = 0$ for all but one class. Any
point in this polyhedron is a convex combination of these extreme points. Take the true allocation $G_k(u)$ in this station. Then, the above says that $G_k(u)$ is a convex combination of these extreme points and since $G_k(u)$ is a linear function and $(I - P')$ is invertible, it must be the case that $u$ is a convex combination of $u_k$ that solve $G_k(u) = 0$ for all $k \in \mathcal{C}(j)$ but one.\eProof
\end{proof}

\section{Detailed Proofs for Balanced DHV Networks}

\setcounter{equation}{0}
\renewcommand{\theequation}{\thesection.\arabic{equation}}

\subsection{Chen-$\calS$ in Balanced DHV Networks
\label{Appendix: Chen-S in Balanced DHV Networks}}

In this appendix, we prove that the pairs $(R,\theta)$ for the $8$ static-priority policies in a balanced DHV network ($m_1+m_4=1$; $m_2+m_5=1$; $m_3+m_6=1$) with $\alpha_1<1$, are all Chen-$\calS$ with invertible reflection matrices and same-sign determinants. In turn, Theorem \ref{Thm: full statement for corner->inter} allows us to say that the refection matrix is invertible and Chen-$\calS$ for any ratio matrix $\Delta$.

For each of the 8 static priority policy, to prove that the corresponding $(R, \theta)$ is Chen-$\calS$, we must show that $R$ is completely-$\calS$, and consider all the 7 partitions of $\mathcal{J}$ and show that there exists a positive vector $h \in \mathbb{R}^3$ such that (\ref{equ:ChenS_def}) is satisfied for all the 7 partitions.\\

\begin{proof} {Proof for priority 2-4-6 (Chen-S).}
\label{proof: Proof of Chen-S for priority 2,4,6}

In the balanced DHV network under priority 2-4-6, we have $J = 3, K = 6, \rho = \left[\alpha_1 \ \alpha_1 \ \alpha_1\right]',$\\
$$R = \frac{1}{m_1m_3m_5+m_2m_4m_6}
 \left[
 \begin{matrix}
   m_1m_3m_5 & m_1m_4m_6 & -m_1m_4m_5 \\
   -m_2m_5 & m_1m_5 & m_5(m_2-m_1)\\
   m_3(m_2-m_3) & -m_3m_6 & m_3m_5
  \end{matrix}
  \right] .
$$

Remind that a matrix is completely-$\calS$ if all of its principal submatrices are $\calS$-matrices; and a square matrix $P$ is called an $\calS$-matrix if there exists a positive vector $u$ such that $Pu > 0$. Thus, to prove that the reflection matrix $R$ is completely-$\calS$, we need to prove that all its submatrices, say matrix $P$, are $\calS$-matrices.
\begin{itemize}
    \item [(i)] All the three $1\times 1$ principle submatrices are $\calS$-matrices because all of $R$'s diagonal elements are positive, and there obviously exists a positive vector $u$ such that $Pu > 0$;
    \item [(ii)] All the $2\times 2$ principle submatrices are $\calS$-matrices:

    \vspace{0.2cm}

    When $P = \displaystyle\frac{1}{m_1m_3m_5+m_2m_4m_6}
 \left[
 \begin{matrix}
   m_1m_3m_5 & m_1m_4m_6 \\
   -m_2m_5 & m_1m_5 \\
  \end{matrix}
  \right] 
$, it is easy to find a $u>0$ with $u_2$ being large enough such that  $Pu > 0$ as required.
\vspace{0.2cm}

   When $P = \displaystyle\frac{1}{m_1m_3m_5+m_2m_4m_6}
 \left[
 \begin{matrix}
   m_1m_3m_5 & -m_1m_4m_5 \\
   m_3(m_2-m_3) & m_3m_5 \\
  \end{matrix}
  \right] 
$, if $m_2-m_3\geq0$, then it is easy to find a $u>0$ with $u_1$ being large enough such that  $Pu > 0$ as required. If $m_2-m_3<0$, suppose $Pu = k>0$, $u = P^{-1}k$. Then $P^{-1}>0$ is sufficient to show that $P$ is a $\calS$-matrix.  $P^{-1}>0$ is equivalent to $m_1m_3^2m_5^2 + m_1m_4m_5m_3 (m_2-m_3)>0 \Leftrightarrow m_1m_3m_5(m_3m_5 + m_4(m_2-m_3))>0$. Since $m_3m_5 + m_4(m_2-m_3) \geq m_3m_5 +m_2-m_3  = m_2-m_3m_2 = m_2m_6>0$, we conclude that $P^{-1}>0$, and $P$ is a $\calS$-matrix.
\vspace{0.2cm}

   When $P = \displaystyle\frac{1}{m_1m_3m_5+m_2m_4m_6}
 \left[
 \begin{matrix}
   m_1m_5 & m_5(m_2-m_1) \\
   -m_3m_6 & m_3m_5 \\
  \end{matrix}
  \right] 
$, if $m_2-m_1\geq0$, then it is easy to find a $u>0$ with $u_1$ being large enough such that  $Pu > 0$ as required. If $m_2-m_1<0$, similar as above, $P^{-1}>0$ is sufficient to show that $P$ is a $\calS$-matrix.  $P^{-1}>0$ is equivalent to $m_1m_3m_5^2+m_3m_6m_5(m_2-m_1)>0 \Leftrightarrow m_3m_5(m_1m_5+m_6(m_2-m_1))>0$. Since $m_1m_5+m_6(m_2-m_1) \geq m_1m_5+m_2-m_1 = m_2-m_1m_2 = m_2m_4>0$, we conclude that $P^{-1}>0$, and $P$ is a $\calS$-matrix.

  \item [(iii)] The $3\times 3$ principle submatrix, which is $R$, is a $\calS$-matrices since $R^{-1} = CMQ\Delta >0$.
\end{itemize}
Thus, we can conclude that the reflection matrix $R$ is completely-$\calS$.\\

Let $ B = m_1m_3m_5+m_2m_4m_6>0$, then
$$
R = \frac{1}{B}
 \left[
 \begin{matrix}
   m_1m_3m_5 & m_1m_4m_6 & -m_1m_4m_5 \\
   -m_2m_5 & m_1m_5 & m_5(m_2-m_1)\\
   m_3(m_2-m_3) & -m_3m_6 & m_3m_5
  \end{matrix}
  \right] ,
$$
\begin{center}
\begin{equation}\nonumber
\begin{aligned}
\theta = R(\rho-e)&= \frac{\alpha_1-1}{B}
\left[
 \begin{matrix}
   m_1m_3m_5+m_1m_4m_6-m_1m_4m_5 \\
   m_1m_5-m_2m_5+m_5m_2-m_5m_1\\
   m_3m_2-m_3m_3-m_3m_6+m_3m_5
  \end{matrix}
  \right] \\
&= \frac{\alpha_1-1}{B}
\left[
 \begin{matrix}
   m_1m_3m_5+m_1m_4m_6-m_1m_4m_5 \\
   0\\
   0
  \end{matrix}
  \right] .
\end{aligned}
\end{equation}
\end{center}

Let $A = (\alpha_1-1)(m_1m_3m_5+m_1m_4m_6-m_1m_4m_5)$. Since
\begin{equation*}
\begin{split}
       m_1m_3m_5+m_1m_4m_6-m_1m_4m_5  &=  m_1(m_3m_5+m_4m_6-m_4m_5)\\
    &=  m_1(m_3m_5+m_4(m_6-m_5)),
\end{split}
\end{equation*}
$A = (\alpha_1-1)m_1(m_3m_5+m_4(m_6-m_5)).$


If $m_6-m_5\geq 0$, we have $$m_3m_5+m_4(m_6-m_5) >0;$$ 
if $m_6-m_5 < 0$, we have
\begin{equation*}
\begin{split}
      m_3m_5+m_4(m_6-m_5)>\; m_3m_5+(m_6-m_5) =\;  m_6-m_5(1-m_3)
     =\;  m_6-m_5m_6=m_2m_6>0.
\end{split}
\end{equation*}

Thus, we have $m_3m_5+m_4(m_6-m_5)>0$, and obviously, with $\alpha_1 < 1$, we have $A<0$,
so that
\begin{center}
\begin{equation}\nonumber
\begin{aligned}
\theta
= 
\left[
 \begin{matrix}
   \displaystyle\frac{A}{B} \\
   0\\
   0
  \end{matrix}
  \right] \leq 0 .
\end{aligned}
\end{equation}
\end{center}

\noindent \textbf{Partition 1: $a = 1; b= 2, 3$}\\
$$\theta_b = \theta_{23} = \left[
 \begin{matrix}
   0 \\
   0
  \end{matrix}
  \right] \Rightarrow u = - R_b^{-1}\theta_b = \left[
 \begin{matrix}
   0 \\
   0
  \end{matrix}
  \right].$$\\
  Thus, $h_a^{'}\left[\theta_a + R_{ab}u\right] = h_1\left[\theta_1+0\right] = h_1\theta_1 = h_1\displaystyle\frac{A}{B}$. Since $A<0, B>0$, taking $h_1 >0 $, we ensure $h_a^{'}\left[\theta_a + R_{ab}u\right] <0$.\\
~\\

\noindent \textbf{Partition 2: $a = 2; b= 1,3$}\\
$$\theta_b = \theta_{13} = \left[
 \begin{matrix}
   \displaystyle\frac{A}{B} \mcr
   0
  \end{matrix}
  \right] ,
  \;
  R_b^{-1} = R_{13}^{-1}= \left[
 \begin{matrix}
   \displaystyle\frac{1}{m_1} & \displaystyle\frac{m_4}{m_3} \mcr
   -\displaystyle\frac{m_2m_6-m_3m_5}{m_1m_5} & 1
  \end{matrix}
  \right] = \left[
 \begin{matrix}
   \displaystyle\frac{1}{m_1} & \displaystyle\frac{m_4}{m_3} \mcr
   \displaystyle\frac{m_3-m_2}{m_1m_5} & 1
  \end{matrix}
  \right]$$
$$\Rightarrow u = - R_b^{-1}\theta_b = \frac{A}{B}\left[
 \begin{matrix}
   -\displaystyle\frac{1}{m_1} \mcr
   \displaystyle\frac{m_2-m_3}{m_1m_5}
  \end{matrix}
  \right]; R_{ab} = \frac{1}{B}\left[
 \begin{matrix}
   -m_2m_5 & m_5(m_2-m_1) 
  \end{matrix}
  \right]; \theta_a= \theta_2 = 0.$$\\

\noindent Thus, 
\begin{center}
\begin{equation}\nonumber
\begin{aligned}
h_a^{'}\left[\theta_a + R_{ab}u\right] &= h_2\left[0+\frac{m_2m_5A}{B^2m_1}+\frac{Am_5(m_2-m_1)(m_2-m_3)}{B^2m_1m_5}\right]\\
&= h_2\frac{A}{B^2}\left[\frac{m_2m_5}{m_1}+\frac{(m_2-m_1)(m_2-m_3)}{m_1}\right]\\
&= h_2\frac{A}{B^2m_1}\left[m_2m_5+(m_2-m_1)(m_2-m_3)\right].\\
\end{aligned}
\end{equation}
\end{center}

\noindent We have 
\begin{equation*}
\begin{split}
    & m_2m_5+(m_2-m_1)(m_2-m_3)
    =\;  m_2m_5+m_2m_2-m_2m_3-m_2m_1+m_1m_3\\
    =\; & m_2-m_2m_3-m_2m_1+m_1m_3 
    =\;  m_2m_6+m_1m_3-m_2m_1
    =\;  m_2m_6+m_1(m_3-m_2),
\end{split}    
\end{equation*}
if $m_3-m_2 \geq 0$, then $$m_2m_6+m_1(m_3-m_2)>0; $$
if $m_3-m_2 <0$, then
\begin{equation*}
\begin{split}
    & m_2m_6+m_1(m_3-m_2) > \;  m_2m_6+(m_3-m_2)
    = \;  m_3-m_2(1-m_6)
    = \; m_3-m_2m_3=m_5m_3>0.
\end{split}    
\end{equation*}
Thus, we get $m_2m_5+(m_2-m_1)(m_2-m_3)>0$.
Since we also have $A<0, B>0,m_1>0$, taking $h_2 >0 $, we ensure $h_a^{'}\left[\theta_a + R_{ab}u\right] <0$.\\

\noindent \textbf{Partition 3: $a = 3; b= 1,2$}\\
$$\theta_b = \theta_{12} = \left[
 \begin{matrix}
   \displaystyle\frac{A}{B} \mcr
   0
  \end{matrix}
  \right] ,
  \;
  R_b^{-1} = R_{12}^{-1}= \left[
 \begin{matrix}
   1 & -\displaystyle\frac{m_4m_6}{m_3m_5+m_5m_6} \mcr
   \displaystyle\frac{m_2}{m_1} & m_3
  \end{matrix}
  \right] = \left[
 \begin{matrix}
   1 & -\displaystyle\frac{m_4m_6}{m_5} \mcr
   \displaystyle\frac{m_2}{m_1} & m_3
  \end{matrix}
  \right]$$
$$\Rightarrow u = - R_b^{-1}\theta_b = -\frac{A}{B}\left[
 \begin{matrix}
   1 \mcr
   \displaystyle\frac{m_2}{m_1}
  \end{matrix}
  \right]; R_{ab} = \frac{1}{B}\left[
 \begin{matrix}
   m_3(m_2-m_3) & -m_3m_6 
  \end{matrix}
  \right]; \theta_a= \theta_3 = 0.$$\\

\noindent Thus, 
\begin{center}
\begin{equation}\nonumber
\begin{aligned}
h_a^{'}\left[\theta_a + R_{ab}u\right] &= h_3\left[0+\frac{m_3(m_3-m_2)A}{B^2}+\frac{Am_3m_6m_2}{B^2m_1}\right]\\
&= h_3\frac{A}{B^2}m_3\left[m_3-m_2+\frac{m_6m_2}{m_1}\right]\\
&= h_3\frac{A}{B^2}\frac{m_3}{m_1}\left[m_1(m_3-m_2)+m_6m_2\right] .\\
\end{aligned}
\end{equation}
\end{center}
Since $m_1(m_3-m_2)+m_6m_2>0$ (we already get this when we check $A<0$),  $A<0, B>0,\displaystyle\frac{m_3}{m_1}>0$, taking $h_3 >0 $, we ensure $h_a^{'}\left[\theta_a + R_{ab}u\right] <0$.\\

\noindent \textbf{Partition 4: $a = 1,2; b= 3$}\\
$$\theta_b = \theta_{3} =0;R_b^{-1}=R_3^{-1} \;{\rm exists} \Rightarrow u = - R_b^{-1}\theta_b =0;$$
$$R_{ab} = \frac{1}{B}\left[
 \begin{matrix}
   -m_1m_4m_5 \mcr
   m_5(m_2-m_1) 
  \end{matrix}
  \right] \Rightarrow R_{ab}u =\left[
 \begin{matrix}
   0 \mcr
   0
  \end{matrix}
  \right];\theta_a= \theta_{12} = \left[
 \begin{matrix}
   \displaystyle\frac{A}{B} \mcr
   0 
  \end{matrix}
  \right]$$
Thus, 
\begin{center}
\begin{equation}\nonumber
\begin{aligned}
h_a^{'}\left[\theta_a + R_{ab}u\right] &= \left[
 \begin{matrix}
   h_1 & h_2
  \end{matrix}
  \right]
  \left[
 \begin{matrix}
   \displaystyle\frac{A}{B} \mcr
   0 
  \end{matrix}
  \right]
&=\frac{h_1A}{B}
\end{aligned}
\end{equation}
\end{center}
Since $A<0, B>0$, taking $h_1 >0 $, we ensure $h_a^{'}\left[\theta_a + R_{ab}u\right] <0$.\\

\noindent \textbf{Partition 5: $a = 1,3; b=2$}\\
$$\theta_b = \theta_{2} =0;R_b^{-1}=R_2^{-1}\; {\rm exists} \Rightarrow u = - R_b^{-1}\theta_b =0;$$
$$R_{ab} = \frac{1}{B}\left[
 \begin{matrix}
   m_1m_4m_6 \mcr
   -m_3m_6 
  \end{matrix}
  \right] \Rightarrow R_{ab}u =\left[
 \begin{matrix}
   0 \mcr
   0
  \end{matrix}
  \right];\theta_a= \theta_{13} = \left[
 \begin{matrix}
   \displaystyle\frac{A}{B} \mcr
   0 
  \end{matrix}
  \right]$$
Thus, 
\begin{center}
\begin{equation}\nonumber
\begin{aligned}
h_a^{'}\left[\theta_a + R_{ab}u\right] &= \left[
 \begin{matrix}
   h_1 & h_3
  \end{matrix}
  \right]
  \left[
 \begin{matrix}
   \displaystyle\frac{A}{B} \mcr
   0 
  \end{matrix}
  \right]
&= \frac{h_1A}{B}
\end{aligned}
\end{equation}
\end{center}
Since $A<0, B>0$, taking $h_1 >0 $, we ensure $h_a^{'}\left[\theta_a + R_{ab}u\right] <0$.\\

\noindent \textbf{Partition 6: $a = 2,3; b= 1$}\\
$$\theta_b = \theta_{1} =\frac{A}{B};\;
R_b^{-1}=R_1^{-1}=\frac{m_1m_3m_5+m_2m_4m_6}{m_1m_3m_5}=1+\frac{m_2m_4m_6}{m_1m_3m_5}$$
$$\Rightarrow u = - R_b^{-1}\theta_b =-(1+\frac{m_2m_4m_6}{m_1m_3m_5})\frac{A}{B};\;
R_{ab} = \frac{1}{B}\left[
 \begin{matrix}
   -m_2m_5 \mcr
   m_3(m_2-m_3) 
  \end{matrix}
  \right]; \;
\theta_a= \theta_{23} = \left[
 \begin{matrix}
   0 \mcr
   0
  \end{matrix}
  \right]$$
Thus, 
\begin{center}
\begin{equation}\nonumber
\begin{aligned}
h_a^{'}\left[\theta_a + R_{ab}u\right] &= \left[
 \begin{matrix}
   h_2 & h_3
  \end{matrix}
  \right]
  \left[
 \begin{matrix}
   \displaystyle\frac{A}{B^2}m_2m_5(1+\frac{m_2m_4m_6}{m_1m_3m_5}) \mcr
   \displaystyle\frac{A}{B^2}m_3(m_3-m_2)(1+\frac{m_2m_4m_6}{m_1m_3m_5}) 
  \end{matrix}
  \right]\mcr
&= \frac{A}{B^2}(1+\frac{m_2m_4m_6}{m_1m_3m_5})\left[m_2m_5h_2+m_3(m_3-m_2)h_3\right]
\end{aligned}
\end{equation}
\end{center}
Since $A<0, B>0, (1+\displaystyle\frac{m_2m_4m_6}{m_1m_3m_5})>0$, fixing $h_2>0$ large enough and $h_3>0$ small enough, we ensure $h_a^{'}\left[\theta_a + R_{ab}u\right] <0$.\\

\noindent \textbf{Partition 7: $a = 1,2,3; b= \varnothing$}\\
$$h_a^{'}\left[\theta_a + R_{ab}u\right] = h'\theta = h_1\frac{A}{B}$$
Since $A<0, B>0, h_1>0$, we ensure  $h_a^{'}\left[\theta_a + R_{ab}u\right] <0$.\\

We have shown, then, that exists a positive vector $h \in{R^3}$ with $h_2$ being a large enough positive number, and $h_3$ being a small enough positive number as required so that $(R,\theta)$ under the static-priority policy 2-4-6 is Chen-$\calS$.
\eProof\\

\end{proof}

\begin{proof} {Proof for priority 1-5-3 (Chen-S).}
\label{proof: Proof of Chen-S for priority 1,5,3}

\noindent In the balanced DHV network under priority 1-5-3, we have $J = 3, K = 6, \rho = \left[\alpha_1 \ \alpha_1 \ \alpha_1\right]',$\\
$$R = 
 \left[
 \begin{matrix}
   \displaystyle\frac{1}{m_2} & -\displaystyle\frac{m_4}{m_2} & 0 \mcr
   -\displaystyle\frac{m_5}{m_4} & 1 & 0\mcr
   \displaystyle\frac{m_3-m_2}{m_2m_4} & -\displaystyle\frac{m_3}{m_2} & 1
  \end{matrix}
  \right] ,
$$
\begin{center}
\begin{equation}\nonumber
\begin{aligned}
\theta = R(\rho-e)&= (\alpha_1-1)
\left[
 \begin{matrix}
   \displaystyle\frac{m_1}{m_2} \mcr
   1-\displaystyle\frac{m_5}{m_4}\mcr
   \displaystyle\frac{m_1}{m_2m_4}(m_3-m_2)
  \end{matrix}
  \right] \mcr
&= \frac{\alpha_1-1}{m_2m_4}
\left[
 \begin{matrix}
   m_1m_4 \mcr
   m_2m_4-m_2m_5\mcr
   m_1m_3-m_1m_2
  \end{matrix}
  \right] .
\end{aligned}
\end{equation}
\end{center}

Similar to the priority 2-4-6 case, it can be proved that the reflection matrix $R$ is completely-$\calS$, we omit the details here.\\

\noindent \textbf{Partition 1: $a = 1; b= 2,3$}\\
$$\theta_b = \theta_{23} = \frac{\alpha_1-1}{m_2m_4}\left[
 \begin{matrix}
   m_2m_4-m_2m_5 \mcr
   m_1m_3-m_1m_2
  \end{matrix}
  \right],
    \;
  R_b^{-1} = R_{23}^{-1}= \left[
 \begin{matrix}
   1 & 0 \mcr
   \displaystyle\frac{m_3}{m_2} & 1
  \end{matrix}
  \right] 
  \Rightarrow u = - R_b^{-1}\theta_b = \displaystyle\frac{1-\alpha_1}{m_4}\left[
 \begin{matrix}
   m_4-m_5 \mcr
   m_3-m_1
  \end{matrix}
  \right];$$
  
  $$
   R_{ab} = \left[
 \begin{matrix}
  -\displaystyle\frac{m_4}{m_2} & 0 
  \end{matrix}
  \right] \Rightarrow 
  R_{ab}u = (\alpha_1-1)(1-\frac{m_1}{m_2}).$$
  
  \noindent Thus, $h_a^{'}\left[\theta_a + R_{ab}u\right] = h_1(\alpha_1-1) $. Since $\alpha_1<1$, taking $h_1 >0 $, we ensure $h_a^{'}\left[\theta_a + R_{ab}u\right] <0$.\\

\noindent \textbf{Partition 2: $a = 2; b= 1,3$}\\
$$\theta_b = \theta_{13} = \frac{\alpha_1-1}{m_2m_4}\left[
 \begin{matrix}
   m_1m_4\mcr
   m_1m_3-m_1m_2
  \end{matrix}
  \right] ,
  \;
  R_b^{-1} = R_{13}^{-1}= \left[
 \begin{matrix}
   m_2 & 0 \mcr
   \displaystyle\frac{m_2-m_3}{m_4} & 1
  \end{matrix}
  \right] 
  \Rightarrow u = - R_b^{-1}\theta_b = \left[
 \begin{matrix}
   (1-\alpha_1)m_1 \mcr
  0
  \end{matrix}
  \right];$$
  $$R_{ab} = \left[
 \begin{matrix}
   -\displaystyle\frac{m_5}{m_4} & 0 
  \end{matrix}
  \right] \Rightarrow 
  R_{ab}u = (\alpha_1-1)\frac{m_1m_5}{m_4}.$$

\noindent Thus, $h_a^{'}\left[\theta_a + R_{ab}u\right] 
= (\alpha_1-1)h_2(\displaystyle\frac{m_1}{m_2}+\displaystyle\frac{m_1m_5}{m_4})$. Since $\alpha_1<1$, taking $h_2 >0 $, we ensure $h_a^{'}\left[\theta_a + R_{ab}u\right] <0$.\\

\noindent \textbf{Partition 3: $a = 3; b= 1,2$}\\
$$\theta_b = \theta_{12} = \frac{\alpha_1-1}{m_2m_4}\left[
 \begin{matrix}
   m_1m_4\mcr
   m_2m_4-m_2m_5
  \end{matrix}
  \right] ,
  \;
  R_b^{-1} = R_{12}^{-1}= \left[
 \begin{matrix}
   1 & \displaystyle\frac{m_4}{m_2} \mcr
   \displaystyle\frac{m_5}{m_4} & \displaystyle\frac{1}{m_2}
  \end{matrix}
  \right] 
  \Rightarrow u = - R_b^{-1}\theta_b = \left[
 \begin{matrix}
   1-\alpha_1 \mcr
   1-\alpha_1
  \end{matrix}
  \right];$$
  $$R_{ab} = \left[
 \begin{matrix}
   \displaystyle\frac{m_3-m_2}{m_2m_4} & -\displaystyle\frac{m_3}{m_2} 
  \end{matrix}
  \right] \Rightarrow 
  R_{ab}u = (1-\alpha_1)\frac{m_3m_1-m_2}{m_2m_4}.$$

\noindent Thus, $h_a^{'}\left[\theta_a + R_{ab}u\right] 
= h_3(\alpha_1-1)$. Since $\alpha_1<1$, taking $h_3 >0 $, we ensure $h_a^{'}\left[\theta_a + R_{ab}u\right] <0$.\\

\noindent \textbf{Partition 4: $a = 1,2; b= 3$}\\
$$\theta_b = \theta_{3} = \frac{\alpha_1-1}{m_2m_4}(m_1m_3-m_1m_2)\; , R_b^{-1}=R_3^{-1}=1 \; \Rightarrow u = - R_b^{-1}\theta_b = \frac{1-\alpha_1}{m_2m_4}(m_1m_3-m_1m_2);$$
$$R_{ab} = \left[
 \begin{matrix}
   0 \mcr
   0 
  \end{matrix}
  \right] \Rightarrow R_{ab}u =\left[
 \begin{matrix}
   0 \mcr
   0
  \end{matrix}
  \right];\theta_a= \theta_{12} = \frac{\alpha_1-1}{m_2m_4}\left[
 \begin{matrix}
   m_1m_4 \mcr
   m_2m_4-m_2m_5
  \end{matrix}
  \right]$$
Thus, 
\begin{center}
\begin{equation}\nonumber
\begin{aligned}
h_a^{'}\left[\theta_a + R_{ab}u\right] &= \frac{\alpha_1-1}{m_2m_4}\left[
 \begin{matrix}
   h_1 & h_2
  \end{matrix}
  \right]
  \left[
 \begin{matrix}
   m_1m_4 \mcr
   m_2m_4-m_2m_5
  \end{matrix}
  \right]
&=\frac{\alpha_1-1}{m_2m_4}\left[m_1m_4h_1+m_2(m_4-m_5)h_2\right]
\end{aligned}
\end{equation}
\end{center}
Since $\alpha_1<1$, fixing $h_1>0$ large enough and  $h_2 >0$ small enough, we ensure $h_a^{'}\left[\theta_a + R_{ab}u\right] <0$.\\

\noindent \textbf{Partition 5: $a = 1,3; b= 2$}\\
$$\theta_b = \theta_{2} = (\alpha_1-1)(1-\frac{m_5}{m_4})\; , R_b^{-1}=R_2^{-1}=1 \; \Rightarrow u = - R_b^{-1}\theta_b = (1-\alpha_1)(1-\frac{m_5}{m_4});$$
$$R_{ab} = \left[
 \begin{matrix}
   -\displaystyle\frac{m_4}{m_2} \mcr
   -\displaystyle\frac{m_3}{m_2} 
  \end{matrix}
  \right] \Rightarrow R_{ab}u =\frac{\alpha_1-1}{m_2m_4}\left[
 \begin{matrix}
   m_4(m_4-m_5) \mcr
   m_3(m_4-m_5)
  \end{matrix}
  \right];\theta_a= \theta_{13} = \frac{\alpha_1-1}{m_2m_4}\left[
 \begin{matrix}
   m_1m_4 \mcr
   m_1m_3-m_1m_2
  \end{matrix}
  \right]$$
Thus, 
\begin{center}
\begin{equation}\nonumber
\begin{aligned}
h_a^{'}\left[\theta_a + R_{ab}u\right] &= \frac{\alpha_1-1}{m_2m_4}\left[
 \begin{matrix}
   h_1 & h_3
  \end{matrix}
  \right]
  \left[
 \begin{matrix}
   m_1m_4+m_4^2-m_4m_5\mcr
   m_1m_3-m_1m_2+m_3m_4-m_3m_5
  \end{matrix}
  \right]
&=(\alpha_1-1)\left[h_1+\frac{m_3-m_1}{m_4}h_3\right]
\end{aligned}
\end{equation}
\end{center}
Since $\alpha_1<1$, fixing $h_1>0$ large enough and $h_3 >0$ small enough, we ensure $h_a^{'}\left[\theta_a + R_{ab}u\right] <0$.\\

\noindent \textbf{Partition 6: $a = 2,3; b= 1$}\\
$$\theta_b = \theta_{1} = (\alpha_1-1)\frac{m_1}{m_2}\; , R_b^{-1}=R_1^{-1}=m_2 \; \Rightarrow u = - R_b^{-1}\theta_b = (1-\alpha_1)m_1;$$
$$R_{ab} = \left[
 \begin{matrix}
   -\displaystyle\frac{m_5}{m_4} \mcr
   \displaystyle\frac{m_3-m_2}{m_2m_4}
  \end{matrix}
  \right] \Rightarrow R_{ab}u =\frac{\alpha_1-1}{m_2m_4}\left[
 \begin{matrix}
   m_1m_2m_5    \mcr
   m_2m_1-m_3m_1
  \end{matrix}
  \right];\theta_a= \theta_{23} = \frac{\alpha_1-1}{m_2m_4}\left[
 \begin{matrix}
   m_2m_4-m_2m_5 \mcr
   m_1m_3-m_1m_2
  \end{matrix}
  \right]$$
Thus, 
\begin{center}
\begin{equation}\nonumber
\begin{aligned}
h_a^{'}\left[\theta_a + R_{ab}u\right] &= \frac{\alpha_1-1}{m_4}\left[
 \begin{matrix}
   h_2 & h_3
  \end{matrix}
  \right]
  \left[
 \begin{matrix}
   m_4-m_5+m_1m_5 \mcr
   0
  \end{matrix}
  \right]
&=(\alpha_1-1)m_2h_2
\end{aligned}
\end{equation}
\end{center}
Since $\alpha_1<1$, taking $h_2>0$, we ensure $h_a^{'}\left[\theta_a + R_{ab}u\right] <0$.\\

\noindent \textbf{Partition 7: $a = 1,2,3; b= \varnothing$}\\
$$h_a^{'}\left[\theta_a + R_{ab}u\right] = h'\theta = \frac{\alpha_1-1}{m_2m_4}\left[m_1m_4h_1+m_2(m_4-m_5)h_2+m_1(m_3-m_2)h_3)\right]$$
Since $\alpha_1<1$, fixing $h_1>0$ large enough, $h_2>0$ and $h_3>0$ small enough, we ensure $h_a^{'}\left[\theta_a + R_{ab}u\right] <0$.
\\

\noindent We have shown, then, that exists a positive vector $h \in{R^3}$ with $h_1$ being a large enough positive number, $h_2$ and $h_3$ being small enough positive numbers as required so that $(R,\theta)$ under the static-priority policy 1-5-3  is Chen-$\calS$.
\eProof\\
\end{proof}

\begin{proof} {Proof for priority 1-5-6 (Chen-S).}
\label{proof: Proof of Chen-S for priority 1,5,6}

\noindent In the balanced DHV network under priority 1-5-6, we have $J = 3, K = 6, \rho = \left[\alpha_1 \ \alpha_1 \ \alpha_1\right]',$\\
$$R = 
 \left[
 \begin{matrix}
   \displaystyle\frac{1}{m_3} & 0& -\displaystyle\frac{m_4}{m_3} \mcr
   -\displaystyle\frac{m_5}{m_4} & 1 & 0\mcr
   \displaystyle\frac{m_3-m_2}{m_2m_4} & -\displaystyle\frac{m_3}{m_2} & 1
  \end{matrix}
  \right] ,
$$
\begin{center}
\begin{equation}\nonumber
\begin{aligned}
\theta = R(\rho-e)&= (\alpha_1-1)
\left[
 \begin{matrix}
   \displaystyle\frac{m_1}{m_3} \mcr
   \displaystyle\frac{m_4-m_5}{m_4}\mcr
   \displaystyle\frac{m_2m_4-m_3m_4+m_3-m_2}{m_2m_4}
  \end{matrix}
  \right] \mcr
&= (\alpha_1-1)
\left[
 \begin{matrix}
   \displaystyle\frac{m_1}{m_3} \mcr
   \displaystyle\frac{m_4-m_5}{m_4}\mcr
   \displaystyle\frac{m_1(m_3-m_2)}{m_2m_4}
  \end{matrix}
  \right] .
\end{aligned}
\end{equation}
\end{center}

Similar to the priority 2-4-6 case, it can be proved that the reflection matrix $R$ is completely-$\calS$, we omit the details here.\\

\noindent \textbf{Partition 1: $a = 1; b= 2,3$}\\
$$\theta_b = \theta_{23} = \frac{\alpha_1-1}{m_2m_4}\left[
 \begin{matrix}
   m_2m_4-m_2m_5 \mcr
   m_1m_3-m_1m_2
  \end{matrix}
  \right],
    \;
  R_b^{-1} = R_{23}^{-1}= \left[
 \begin{matrix}
   1 & 0 \mcr
   \displaystyle\frac{m_3}{m_2} & 1
  \end{matrix}
  \right] 
  \Rightarrow u = - R_b^{-1}\theta_b = \frac{1-\alpha_1}{m_4}\left[
 \begin{matrix}
   m_4-m_5 \mcr
   m_3-m_1
  \end{matrix}
  \right];$$
  
  $$
   R_{ab} = \left[
 \begin{matrix}
  0 & -\displaystyle\frac{m_4}{m_3} 
  \end{matrix}
  \right] \Rightarrow 
  R_{ab}u = (\alpha_1-1)(1-\displaystyle\frac{m_1}{m_3}).$$
  
  \noindent Thus, $h_a^{'}\left[\theta_a + R_{ab}u\right] = h_1(\alpha_1-1) $. Since $\alpha_1<1$, taking $h_1 >0 $, we ensure $h_a^{'}\left[\theta_a + R_{ab}u\right] <0$.\\

\noindent \textbf{Partition 2: $a = 2; b= 1,3$}\\
$$\theta_b = \theta_{13} = \frac{\alpha_1-1}{m_2m_3m_4}\left[
 \begin{matrix}
   m_1m_2m_4\mcr
   m_1m_3(m_3-m_2)
  \end{matrix}
  \right] ,
  \;
  R_b^{-1} = R_{13}^{-1}= \left[
 \begin{matrix}
   m_2 & \displaystyle\frac{m_2m_4}{m_3} \mcr
   \displaystyle\frac{m_2-m_3}{m_4} & \displaystyle\frac{m_2}{m_3}
  \end{matrix}
  \right] 
  \Rightarrow u = - R_b^{-1}\theta_b = \left[
 \begin{matrix}
   (1-\alpha_1)m_1 \mcr
  0
  \end{matrix}
  \right];$$
  $$R_{ab} = \left[
 \begin{matrix}
   -\displaystyle\frac{m_5}{m_4} & 0 
  \end{matrix}
  \right] \Rightarrow 
  R_{ab}u = (\alpha_1-1)\frac{m_1m_5}{m_4}.$$

\noindent Thus, $h_a^{'}\left[\theta_a + R_{ab}u\right] 
= (\alpha_1-1)m_2h_2$. Since $\alpha_1<1$, taking $h_2 >0 $, we ensure $h_a^{'}\left[\theta_a + R_{ab}u\right] <0$.\\

\noindent \textbf{Partition 3: $a = 3; b= 1,2$}\\
$$\theta_b = \theta_{12} = \frac{\alpha_1-1}{m_3m_4}\left[
 \begin{matrix}
   m_1m_4\mcr
   m_3m_4-m_3m_5
  \end{matrix}
  \right] ,
  \;
  R_b^{-1} = R_{12}^{-1}= \left[
 \begin{matrix}
   m_3 & 0 \mcr
   \displaystyle\frac{m_3m_5}{m_4} & 1
  \end{matrix}
  \right] 
  \Rightarrow u = - R_b^{-1}\theta_b =(1-\alpha_1) \left[
 \begin{matrix}
   m_1\mcr
   m_2
  \end{matrix}
  \right];$$
  $$R_{ab} = \left[
 \begin{matrix}
   \displaystyle\frac{m_3-m_2}{m_2m_4} & -\displaystyle\frac{m_3}{m_2} 
  \end{matrix}
  \right] \Rightarrow 
  R_{ab}u = \frac{1-\alpha_1}{m_2m_4}(m_1m_3-m_1m_2-m_2m_3m_4).$$

\noindent Thus, $h_a^{'}\left[\theta_a + R_{ab}u\right] 
= h_3(\alpha_1-1)m_3$. Since $\alpha_1<1$, taking $h_3 >0 $, we ensure $h_a^{'}\left[\theta_a + R_{ab}u\right] <0$.\\

\noindent \textbf{Partition 4: $a = 1,2; b= 3$}\\
$$\theta_b = \theta_{3} = \frac{\alpha_1-1}{m_2m_4}(m_1m_3-m_1m_2)\; , R_b^{-1}=R_3^{-1}=1 \; \Rightarrow u = - R_b^{-1}\theta_b = \frac{1-\alpha_1}{m_2m_4}(m_1m_3-m_1m_2);$$
$$R_{ab} = \left[
 \begin{matrix}
   -\displaystyle\frac{m_4}{m_3} \mcr
   0 
  \end{matrix}
  \right] \Rightarrow R_{ab}u =\left[
 \begin{matrix}
   \displaystyle\frac{\alpha_1-1}{m_2m_3}(m_1m_3-m_1m_2) \mcr
   0
  \end{matrix}
  \right];\theta_a= \theta_{12} = \frac{\alpha_1-1}{m_3m_4}\left[
 \begin{matrix}
   m_1m_4 \mcr
   m_3m_4-m_3m_5
  \end{matrix}
  \right]$$
Thus, 
\begin{center}
\begin{equation}\nonumber
\begin{aligned}
h_a^{'}\left[\theta_a + R_{ab}u\right] &= \frac{\alpha_1-1}{m_2m_4}\left[
 \begin{matrix}
   h_1 & h_2
  \end{matrix}
  \right]
  \left[
 \begin{matrix}
   m_1m_4 \mcr
   m_2m_4-m_2m_5
  \end{matrix}
  \right]
&=\frac{\alpha_1-1}{m_2m_4}\left[m_1m_4h_1+m_2(m_4-m_5)h_2\right]
\end{aligned}
\end{equation}
\end{center}
Since $\alpha_1<1$, fixing $h_1>0$ large enough and  $h_2>0 $ small enough, we ensure $h_a^{'}\left[\theta_a + R_{ab}u\right] <0$.\\

\noindent \textbf{Partition 5: $a = 1,3; b= 2$}\\
$$\theta_b = \theta_{2} = \frac{\alpha_1-1}{m_4}(m_4-m_5)\; , R_b^{-1}=R_2^{-1}=1 \; \Rightarrow u = - R_b^{-1}\theta_b = \frac{1-\alpha_1}{m_4}(m_4-m_5);$$
$$R_{ab} = \left[
 \begin{matrix}
   0 \mcr
   -\displaystyle\frac{m_3}{m_2} 
  \end{matrix}
  \right] \Rightarrow R_{ab}u =\displaystyle\frac{(\alpha_1-1)m_3}{m_2m_4}\left[
 \begin{matrix}
   0 \mcr
   m_4-m_5
  \end{matrix}
  \right];\theta_a= \theta_{13} = \frac{\alpha_1-1}{m_2m_3m_4}\left[
 \begin{matrix}
   m_1m_2m_4 \mcr
   m_1m_3^2-m_1m_2m_3
  \end{matrix}
  \right]$$
Thus, 
\begin{center}
\begin{equation}\nonumber
\begin{aligned}
h_a^{'}\left[\theta_a + R_{ab}u\right] &= \frac{\alpha_1-1}{m_2m_3m_4}\left[
 \begin{matrix}
   h_1 & h_3
  \end{matrix}
  \right]
  \left[
 \begin{matrix}
   m_1m_2m_4\mcr
   m_2m_3^2-m_1m_2m_3
  \end{matrix}
  \right]
&=(\alpha_1-1)(\frac{m_1}{m_3}h_1+\frac{m_3-m_1}{m_4}h_3)
\end{aligned}
\end{equation}
\end{center}
Since $\alpha_1<1$, fixing $h_1>0$ large enough and  $h_3>0$ small enough, we ensure $h_a^{'}\left[\theta_a + R_{ab}u\right] <0$.\\

\noindent \textbf{Partition 6: $a = 2,3; b= 1$}\\
$$\theta_b = \theta_{1} = (\alpha_1-1)\frac{m_1}{m_3}\; , R_b^{-1}=R_1^{-1}=m_3 \; \Rightarrow u = - R_b^{-1}\theta_b = (1-\alpha_1)m_1;$$
$$R_{ab} = \left[
 \begin{matrix}
   -\displaystyle\frac{m_5}{m_4} \mcr
   \displaystyle\frac{m_3-m_2}{m_2m_4}
  \end{matrix}
  \right] \Rightarrow R_{ab}u =\frac{\alpha_1-1}{m_2m_4}\left[
 \begin{matrix}
   m_1m_2m_5    \mcr
   m_2m_1-m_3m_1
  \end{matrix}
  \right];\theta_a= \theta_{23} = \frac{\alpha_1-1}{m_2m_4}\left[
 \begin{matrix}
   m_2m_4-m_2m_5 \mcr
   m_1m_3-m_1m_2
  \end{matrix}
  \right]$$
Thus, 
\begin{center}
\begin{equation}\nonumber
\begin{aligned}
h_a^{'}\left[\theta_a + R_{ab}u\right] &= \frac{\alpha_1-1}{m_2m_4}\left[
 \begin{matrix}
   h_2 & h_3
  \end{matrix}
  \right]
  \left[
 \begin{matrix}
   m_2m_4-m_2m_5+m_1m_5m_2\mcr
   0
  \end{matrix}
  \right]
&=(\alpha_1-1)m_2h_2
\end{aligned}
\end{equation}
\end{center}
Since $\alpha_1<1$, taking $h_2>0$, we ensure $h_a^{'}\left[\theta_a + R_{ab}u\right] <0$.\\

\noindent \textbf{Partition 7: $a = 1,2,3; b= \varnothing$}\\
$$h_a^{'}\left[\theta_a + R_{ab}u\right] = h'\theta = \frac{\alpha_1-1}{m_2m_3m_4}\left[m_1m_2m_4h_1+m_2m_3(m_4-m_5)h_2+m_1m_3(m_3-m_2)h_3)\right]$$
Since $\alpha_1<1$, fixing $h_1>0$ large enough, $h_2>0$ and $h_3>0$ small enough, we ensure $h_a^{'}\left[\theta_a + R_{ab}u\right] <0$.
\\

\noindent We have shown, then, that exists a positive vector $h \in{R^3}$ with $h_1$ being a large enough positive number, $h_2$ and $h_3$ being small enough positive numbers as required so that $(R,\theta)$ under the static-priority policy 1-5-6  is Chen-$\calS$.
\eProof\\

\end{proof}

\begin{proof} {Proof for priority 4-5-6 (Chen-S).}
\label{proof: Proof of Chen-S for priority 4,5,6}

\noindent In the balanced DHV network under priority 4-5-6, we have $J = 3, K = 6, \rho = \left[\alpha_1 \ \alpha_1 \ \alpha_1\right]',$\\
$$R = 
 \left[
 \begin{matrix}
   1 & 0 & -m_4 \mcr
   -\displaystyle\frac{m_2}{m_1} & 1 & \displaystyle\frac{m_2}{m_1}-1\mcr
   0 & -\displaystyle\frac{m_3}{m_2} & \displaystyle\frac{m_3}{m_2}
  \end{matrix}
  \right] ,
$$
\begin{center}
\begin{equation}\nonumber
\begin{aligned}
\theta = R(\rho-e)&= (\alpha_1-1)
\left[
 \begin{matrix}
   1-m_4 \mcr
   1-\displaystyle\frac{m_2}{m_1}+\displaystyle\frac{m_2}{m_1}-1\mcr
   \displaystyle\frac{m_3}{m_2}-\displaystyle\frac{m_3}{m_2}
  \end{matrix}
  \right] \mcr
&= (\alpha_1-1)
\left[
 \begin{matrix}
   m_1 \mcr
   0\mcr
   0
  \end{matrix}
  \right] .
\end{aligned}
\end{equation}
\end{center}

Similar to the priority 2-4-6 case, it can be proved that the reflection matrix $R$ is completely-$\calS$, we omit the details here.\\

\noindent \textbf{Partition 1: $a = 1; b= 2,3$}\\
$$\theta_b = \theta_{23} = \left[
 \begin{matrix}
   0\mcr
   0
  \end{matrix}
  \right],
    \;
  R_b^{-1} = R_{23}^{-1}= \left[
 \begin{matrix}
   \displaystyle\frac{m_1}{m_2} & \displaystyle\frac{m_1m_5-m_2m_4}{m_3} \mcr
   \displaystyle\frac{m_1}{m_2} & \displaystyle\frac{m_1}{m_3}
  \end{matrix}
  \right] 
  \Rightarrow u = - R_b^{-1}\theta_b = \left[
 \begin{matrix}
   0\mcr
   0
  \end{matrix}
  \right];$$
  
  $$
   R_{ab} = \left[
 \begin{matrix}
  0 & -m_4
  \end{matrix}
  \right] \Rightarrow 
  R_{ab}u = 0.$$
  
  \noindent Thus, $h_a^{'}\left[\theta_a + R_{ab}u\right] = h_1(\alpha_1-1)m_1 $. Since $\alpha_1<1$, taking $h_1 >0 $, we ensure $h_a^{'}\left[\theta_a + R_{ab}u\right] <0$.\\

\noindent \textbf{Partition 2: $a = 2; b= 1,3$}\\
$$\theta_b = \theta_{13} =(\alpha_1-1)\left[
 \begin{matrix}
   m_1\mcr
   0
  \end{matrix}
  \right] ,
  \;
  R_b^{-1} = R_{13}^{-1}= \left[
 \begin{matrix}
   1 & \displaystyle\frac{m_2m_4}{m_3} \mcr
   0 & \displaystyle\frac{m_2}{m_3}
  \end{matrix}
  \right] 
  \Rightarrow u = - R_b^{-1}\theta_b = \left[
 \begin{matrix}
   (1-\alpha_1)m_1 \mcr
  0
  \end{matrix}
  \right];$$
  $$R_{ab} = \left[
 \begin{matrix}
   -\displaystyle\frac{m_2}{m_1} & \displaystyle\frac{m_2}{m_1}-1
  \end{matrix}
  \right] \Rightarrow 
  R_{ab}u = (\alpha_1-1)m_2.$$

\noindent Thus, $h_a^{'}\left[\theta_a + R_{ab}u\right] 
= (\alpha_1-1)m_2h_2$. Since $\alpha_1<1$, taking $h_2 >0 $, we ensure $h_a^{'}\left[\theta_a + R_{ab}u\right] <0$.\\

\noindent \textbf{Partition 3: $a = 3; b= 1,2$}\\
$$\theta_b = \theta_{12} = (\alpha_1-1)\left[
 \begin{matrix}
   m_1\\
   0
  \end{matrix}
  \right] ,
  \;
  R_b^{-1} = R_{12}^{-1}= \left[
 \begin{matrix}
   1 & 0 \mcr
   \displaystyle\frac{m_2}{m_1} & 1
  \end{matrix}
  \right] 
  \Rightarrow u = - R_b^{-1}\theta_b =(1-\alpha_1) \left[
 \begin{matrix}
   m_1 \mcr
   m_2
  \end{matrix}
  \right];$$
  $$R_{ab} = \left[
 \begin{matrix}
   0 & -\displaystyle\frac{m_3}{m_2} 
  \end{matrix}
  \right] \Rightarrow 
  R_{ab}u = (\alpha_1-1)m_3.$$

\noindent Thus, $h_a^{'}\left[\theta_a + R_{ab}u\right] 
= h_3(\alpha_1-1)m_3$. Since $\alpha_1<1$, taking $h_3 >0 $, we ensure $h_a^{'}\left[\theta_a + R_{ab}u\right] <0$.\\

\noindent \textbf{Partition 4: $a = 1,2; b= 3$}\\
$$\theta_b = \theta_{3} = 0\; , R_b^{-1}=R_3^{-1}=\frac{m_2}{m_3} \; \Rightarrow u = - R_b^{-1}\theta_b = 0;$$
$$R_{ab} = \left[
 \begin{matrix}
   -m_4 \mcr
   \displaystyle\frac{m_2}{m_1}-1
  \end{matrix}
  \right] \Rightarrow R_{ab}u =\left[
 \begin{matrix}
   0 \mcr
   0
  \end{matrix}
  \right];\theta_a= \theta_{12} = (\alpha_1-1)\left[
 \begin{matrix}
   m_1 \mcr
   0
  \end{matrix}
  \right]$$
Thus, 
\begin{center}
\begin{equation}\nonumber
\begin{aligned}
h_a^{'}\left[\theta_a + R_{ab}u\right] &= (\alpha_1-1)\left[
 \begin{matrix}
   h_1 & h_2
  \end{matrix}
  \right]
  \left[
 \begin{matrix}
   m_1 \mcr
   0
  \end{matrix}
  \right]
&=(\alpha_1-1)m_1h_1
\end{aligned}
\end{equation}
\end{center}
Since $\alpha_1<1$, taking $h_1>0$, we ensure $h_a^{'}\left[\theta_a + R_{ab}u\right] <0$.\\

\noindent \textbf{Partition 5: $a = 1,3; b= 2$}\\
$$\theta_b = \theta_{2} = 0\; , R_b^{-1}=R_2^{-1}=1 \; \Rightarrow u = - R_b^{-1}\theta_b = 0;$$
$$R_{ab} = \left[
 \begin{matrix}
  0 \mcr
   -\displaystyle\frac{m_3}{m_2} 
  \end{matrix}
  \right] \Rightarrow R_{ab}u =\left[
 \begin{matrix}
   0 \mcr
   0
  \end{matrix}
  \right];\theta_a= \theta_{13} = (\alpha_1-1)\left[
 \begin{matrix}
   m_1 \mcr
   0
  \end{matrix}
  \right]$$
Thus, 
\begin{center}
\begin{equation}\nonumber
\begin{aligned}
h_a^{'}\left[\theta_a + R_{ab}u\right] &= (\alpha_1-1)\left[
 \begin{matrix}
   h_1 & h_3
  \end{matrix}
  \right]
  \left[
 \begin{matrix}
   m_1\mcr
   0
  \end{matrix}
  \right]
&=(\alpha_1-1)m_1h_1
\end{aligned}
\end{equation}
\end{center}
Since $\alpha_1<1$, taking $h_1>0$, we ensure $h_a^{'}\left[\theta_a + R_{ab}u\right] <0$.\\

\noindent \textbf{Partition 6: $a = 2,3; b= 1$}\\
$$\theta_b = \theta_{1} = (\alpha_1-1)m_1\; , R_b^{-1}=R_1^{-1}=1 \; \Rightarrow u = - R_b^{-1}\theta_b = (1-\alpha_1)m_1;$$
$$R_{ab} = \left[
 \begin{matrix}
   -\displaystyle\frac{m_2}{m_1} \mcr
   0
  \end{matrix}
  \right] \Rightarrow R_{ab}u =(\alpha_1-1)\left[
 \begin{matrix}
   m_2   \mcr
   0
  \end{matrix}
  \right];\theta_a= \theta_{23} = \left[
 \begin{matrix}
   0 \mcr
   0
  \end{matrix}
  \right]$$
Thus, 
\begin{center}
\begin{equation}\nonumber
\begin{aligned}
h_a^{'}\left[\theta_a + R_{ab}u\right] &= (\alpha_1-1)\left[
 \begin{matrix}
   h_2 & h_3
  \end{matrix}
  \right]
  \left[
 \begin{matrix}
   m_2 \mcr
   0
  \end{matrix}
  \right]
&=(\alpha_1-1)m_2h_2
\end{aligned}
\end{equation}
\end{center}
Since $\alpha_1<1$, taking $h_2>0$, we ensure $h_a^{'}\left[\theta_a + R_{ab}u\right] <0$.\\

\noindent \textbf{Partition 7: $a = 1,2,3; b= \varnothing$}\\
$$h_a^{'}\left[\theta_a + R_{ab}u\right] = h'\theta = (\alpha_1-1)m_1h_1$$

Since $\alpha_1<1$, taking $h_1>0$, we ensure $h_a^{'}\left[\theta_a + R_{ab}u\right] <0$.
\\

\noindent We have shown, then, that exists a positive vector $h \in{R^3}$ as required so that $(R,\theta)$ under the static-priority policy 4-5-6  is Chen-$\calS$.
\eProof\\
\end{proof}

\begin{proof} {Proof for priority 1-2-6 (Chen-S).}
\label{proof: Proof of Chen-S for priority 1,2,6}

\noindent In the balanced DHV network under priority 1-2-6, we have $J = 3, K = 6, \rho = \left[\alpha_1 \ \alpha_1 \ \alpha_1\right]',$\\
$$R = 
 \left[
 \begin{matrix}
   1 & \displaystyle\frac{m_4m_6}{m_3m_5} & -\displaystyle\frac{m_4}{m_3} \mcr
   -\displaystyle\frac{m_5}{m_4} & 1 & 0\mcr
   0 & -\displaystyle\frac{m_6}{m_5} & 1
  \end{matrix}
  \right] ,
$$
\begin{center}
\begin{equation}\nonumber
\begin{aligned}
\theta = R(\rho-e)&= (\alpha_1-1)
\left[
 \begin{matrix}
   1+\displaystyle\frac{m_4m_6}{m_3m_5}-\displaystyle\frac{m_4}{m_3} \mcr
   1-\displaystyle\frac{m_5}{m_4}\mcr
   1-\displaystyle\frac{m_6}{m_5}
  \end{matrix}
  \right] \mcr
&= \frac{\alpha_1-1}{m_3m_4m_5}
\left[
 \begin{matrix}
   m_3m_4m_5+m_4^2m_6-m_4^2m_5 \mcr
   m_3m_5(m_4-m_5)\mcr
   m_3m_4(m_5-m_6)
  \end{matrix}
  \right] .
\end{aligned}
\end{equation}
\end{center}

Similar to the priority 2-4-6 case, it can be proved that the reflection matrix $R$ is completely-$\calS$, we omit the details here.\\

\noindent \textbf{Partition 1: $a = 1; b= 2,3$}\\
$$\theta_b = \theta_{23} = \frac{\alpha_1-1}{m_4m_5}\left[
 \begin{matrix}
   m_5(m_4-m_5) \mcr
   m_4(m_5-m_6)
  \end{matrix}
  \right],
    \;
  R_b^{-1} = R_{23}^{-1}= \left[
 \begin{matrix}
   1 & 0 \mcr
   \displaystyle\frac{m_6}{m_5} & 1
  \end{matrix}
  \right] 
  \Rightarrow u = - R_b^{-1}\theta_b = \frac{1-\alpha_1}{m_4}\left[
 \begin{matrix}
   m_4-m_5 \mcr
   m_4-m_6
  \end{matrix}
  \right];$$
  
  $$
   R_{ab} = \left[
 \begin{matrix}
  \displaystyle\frac{m_4m_6}{m_3m_5} & -\displaystyle\frac{m_4}{m_3}
  \end{matrix}
  \right] \Rightarrow 
  R_{ab}u = \displaystyle\frac{1-\alpha_1}{m_3m_5}m_4(m_6-m_5).$$
  
  \noindent Thus, $h_a^{'}\left[\theta_a + R_{ab}u\right] = h_1(\alpha_1-1) $. Since $\alpha_1<1$, taking $h_1 >0 $, we ensure $h_a^{'}\left[\theta_a + R_{ab}u\right] <0$.\\

\noindent \textbf{Partition 2: $a = 2; b= 1,3$}\\
$$\theta_b = \theta_{13} = \frac{\alpha_1-1}{m_3m_5}\left[
 \begin{matrix}
   m_3m_5+m_4m_6-m_4m_5\mcr
   m_3(m_5-m_6)
  \end{matrix}
  \right] ,
  \;
  R_b^{-1} = R_{13}^{-1}= \left[
 \begin{matrix}
   1 & \displaystyle\frac{m_4}{m_3} \mcr
   0 & 1
  \end{matrix}
  \right] 
  \Rightarrow u = - R_b^{-1}\theta_b = (1-\alpha_1)\left[
 \begin{matrix}
   1 \mcr
  \displaystyle\frac{m_5-m_6}{m_5}
  \end{matrix}
  \right];$$
  $$R_{ab} = \left[
 \begin{matrix}
   -\displaystyle\frac{m_5}{m_4} & 0 
  \end{matrix}
  \right] \Rightarrow 
  R_{ab}u = (\alpha_1-1)\frac{m_5}{m_4}.$$

\noindent Thus, $h_a^{'}\left[\theta_a + R_{ab}u\right] 
= (\alpha_1-1)h_2$. Since $\alpha_1<1$, taking $h_2 >0 $, we ensure $h_a^{'}\left[\theta_a + R_{ab}u\right] <0$.\\

\noindent \textbf{Partition 3: $a = 3; b= 1,2$}\\
$$\theta_b = \theta_{12} = \frac{\alpha_1-1}{m_3m_4m_5}\left[
 \begin{matrix}
   m_3m_4m_5+m_4^2m_6-m_4^2m_5\mcr
   m_3m_4m_5-m_3m_5^2
  \end{matrix}
  \right] ,
  \;
  R_b^{-1} = R_{12}^{-1}= \left[
 \begin{matrix}
   m_3 & -\displaystyle\frac{m_4m_6}{m_5} \\
   \displaystyle\frac{m_3m_5}{m_4} & m_3
  \end{matrix}
  \right] $$
  $$
  \Rightarrow u = - R_b^{-1}\theta_b = (1-\alpha_1)\left[
 \begin{matrix}
   m_3-m_4+m_6 \mcr
   m_6-m_5+m_3
  \end{matrix}
  \right];$$
  $$R_{ab} = \left[
 \begin{matrix}
   0 & -\displaystyle\frac{m_6}{m_5} 
  \end{matrix}
  \right] \Rightarrow 
  R_{ab}u = (\alpha_1-1)\frac{m_6}{m_5}(1-m_5).$$

\noindent Thus, $h_a^{'}\left[\theta_a + R_{ab}u\right] 
= h_3(\alpha_1-1)m_3$. Since $\alpha_1<1$, taking $h_3 >0 $, we ensure $h_a^{'}\left[\theta_a + R_{ab}u\right] <0$.\\

\noindent \textbf{Partition 4: $a = 1,2; b= 3$}\\
$$\theta_b = \theta_{3} = (\alpha_1-1)(1-\frac{m_6}{m_5})\; , R_b^{-1}=R_3^{-1}=1 \; \Rightarrow u = - R_b^{-1}\theta_b = (1-\alpha_1)(1-\frac{m_6}{m_5});$$
$$R_{ab} = \left[
 \begin{matrix}
   -\displaystyle\frac{m_4}{m_3} \mcr
   0 
  \end{matrix}
  \right] \Rightarrow R_{ab}u =\left[
 \begin{matrix}
   \displaystyle\frac{\alpha_1-1}{m_3m_5}(m_4m_5-m_4m_6) \mcr
   0
  \end{matrix}
  \right];\theta_a= \theta_{12} = \frac{\alpha_1-1}{m_3m_4m_5}\left[
 \begin{matrix}
   m_3m_4m_5+m_4^2m_6-m_4^2m_5 \mcr
   m_3m_5(m_4-m_5)
  \end{matrix}
  \right]$$
Thus, 
\begin{center}
\begin{equation}\nonumber
\begin{aligned}
h_a^{'}\left[\theta_a + R_{ab}u\right] &= (\alpha_1-1)\left[
 \begin{matrix}
   h_1 & h_2
  \end{matrix}
  \right]
  \left[
 \begin{matrix}
   1\mcr
  1-\displaystyle\frac{m_5}{m_4}
  \end{matrix}
  \right]
&=(\alpha_1-1)\left[h_1+(1-\frac{m_5}{m_4})h_2\right]
\end{aligned}
\end{equation}
\end{center}
Since $\alpha_1<1$, fixing $h_1>0$ large enough and  $h_2>0$ small enough, we ensure $h_a^{'}\left[\theta_a + R_{ab}u\right] <0$.\\

\noindent \textbf{Partition 5: $a = 1,3; b= 2$}\\
$$\theta_b = \theta_{2} = (\alpha_1-1)(1-\frac{m_5}{m_4})\; , R_b^{-1}=R_2^{-1}=1 \; \Rightarrow u = - R_b^{-1}\theta_b = (1-\alpha_1)(1-\frac{m_5}{m_4});$$
$$R_{ab} = \left[
 \begin{matrix}
   \displaystyle\frac{m_4m_6}{m_3m_5} \mcr
   -\displaystyle\frac{m_6}{m_5} 
  \end{matrix}
  \right] \Rightarrow R_{ab}u =(\alpha_1-1)\left[
 \begin{matrix}
   \displaystyle\frac{m_6(m_5-m_4)}{m_3m_5} \mcr
   \displaystyle\frac{m_6(m_4-m_5)}{m_4m_5}
  \end{matrix}
  \right];\theta_a= \theta_{13} = \frac{\alpha_1-1}{m_3m_5}\left[
 \begin{matrix}
   m_3m_5+m_4m_6-m_4m_5 \mcr
   m_3m_5-m_3m_6
  \end{matrix}
  \right]$$
Thus, 
\begin{center}
\begin{equation}\nonumber
\begin{aligned}
h_a^{'}\left[\theta_a + R_{ab}u\right] &= \frac{\alpha_1-1}{m_3m_4}\left[
 \begin{matrix}
   h_1 & h_3
  \end{matrix}
  \right]
  \left[
 \begin{matrix}
   m_1m_4\mcr
  m_3(m_4-m_6)
  \end{matrix}
  \right]
&=(\alpha_1-1)\left[\frac{m_1}{m_3}h_1+\frac{m_4-m_6}{m_4}h_3\right]
\end{aligned}
\end{equation}
\end{center}
Since $\alpha_1<1$, fixing $h_1>0$ large enough and $h_3>0$ small enough, we ensure $h_a^{'}\left[\theta_a + R_{ab}u\right] <0$.\\

\noindent \textbf{Partition 6: $a = 2,3; b= 1$}\\
$$\theta_b = \theta_{1} = \frac{\alpha_1-1}{m_3m_5}(m_3m_5+m_4m_6-m_4m_5)\; , R_b^{-1}=R_1^{-1}=1 \; \Rightarrow u = - R_b^{-1}\theta_b = \frac{1-\alpha_1}{m_3m_5}(m_3m_5+m_4m_6-m_4m_5);$$
$$R_{ab} = \left[
 \begin{matrix}
   -\displaystyle\frac{m_5}{m_4} \mcr
  0
  \end{matrix}
  \right] \Rightarrow R_{ab}u =\frac{\alpha_1-1}{m_3m_4}\left[
 \begin{matrix}
   m_3m_5+m_4m_6-m_4m_5   \mcr
   0
  \end{matrix}
  \right];\theta_a= \theta_{23} = \frac{\alpha_1-1}{m_3m_4m_5}\left[
 \begin{matrix}
   m_3m_4m_5-m_3m_5^2\mcr
   m_3m_4m_5-m_3m_4m_6
  \end{matrix}
  \right]$$
Thus, 
\begin{center}
\begin{equation}\nonumber
\begin{aligned}
h_a^{'}\left[\theta_a + R_{ab}u\right] 
  &= \frac{\alpha_1-1}{m_3m_5}\left[
 \begin{matrix}
   h_2 & h_3
  \end{matrix}
  \right]
  \left[
 \begin{matrix}
   m_5m_2 \mcr
   m_3(m_5-m_6)
  \end{matrix}
  \right]
&=(\alpha_1-1)(\frac{m_2}{m_3}h_2+\frac{m_5-m_6}{m_5}h_3)
\end{aligned}
\end{equation}
\end{center}
Since $\alpha_1<1$, fixing $h_2>0$ large enough and  $h_3>0$ small enough, we ensure $h_a^{'}\left[\theta_a + R_{ab}u\right] <0$.\\

\noindent \textbf{Partition 7: $a = 1,2,3; b= \varnothing$}\\
$$h_a^{'}\left[\theta_a + R_{ab}u\right] = h'\theta = (\alpha_1-1)\left[(1+\frac{m_4m_6}{m_3m_5}-\frac{m_4}{m_3})h_1+(1-\frac{m_5}{m_4})h_2+(1-\frac{m_6}{m_5})h_3\right]$$
Let $K = 1+\displaystyle\frac{m_4m_6}{m_3m_5}-\frac{m_4}{m_3} = \frac{1}{m_3m_5}(m_3m_5+m_4m_6-m_4m_5)$.\\

\noindent If $m_6 - m_5\geq 0$, then $$K =\frac{1}{m_3m_5}(m_3m_5+m_4m_6-m_4m_5) = \frac{1}{m_3m_5}(m_3m_5+m_4(m_6-m_5))>0;$$
If $m_6 - m_5\leq 0$, then
\begin{equation*}
\begin{split}
  K = \; & \frac{1}{m_3m_5}(m_3m_5+m_4m_6-m_4m_5)
    = \;  \frac{1}{m_3m_5}((1-m_6)m_5+(1-m_1)m_6-m_4m_5)\mcr
    = \; & \frac{1}{m_3m_5}((m_5-m_5m_6+m_6-m_1m_6-m_4m_5)
    = \;  \frac{1}{m_3m_5}((1-m_4)m_5+(1-m_5)m_6-m_1m_6)\mcr
    = \; & \frac{1}{m_3m_5}(m_1m_5+m_2m_6-m_1m_6)
    = \; \frac{1}{m_3m_5}(m_1(m_5-m_6)+m_2m_6) >0.
\end{split}    
\end{equation*}

\noindent Thus, we have $K>0$. Since $\alpha_1<1$, fixing $h_1>0$ large enough, $h_2>0$ and $h_3>0$ small enough, we ensure $h_a^{'}\left[\theta_a + R_{ab}u\right] <0$.\\

\noindent We have shown, then, that exists a positive vector $h \in{R^3}$ with $h_1$ being a relatively large enough positive number comparing to $h_2$, and $h_2$ being a relatively large enough positive number comparing to $h_3$ as required so that $(R,\theta)$ under the static-priority policy 1-2-6  is Chen-$\calS$.
\eProof\\
\end{proof}

\begin{lemma}
\label{lem:lower triangular matrix, Schur-S}
(\cite{chen1996sufficient} Lemma 3.2) A lower triangular matrix whose inverse exists and is non-negative is both completely-$\calS$ and Schur-S.
\end{lemma}

\vspace{0.2cm}
 
\begin{lemma}
\label{lem:M-matrix, Schur-S}
(\cite{chen1996sufficient} Lemma 3.1) A nonsingular M-matrix is both completely-$\calS$ and Schur-S.
\end{lemma}
 
\vspace{0.2cm}
 
\begin{proof} {Proof for priority 1-2-3/4-2-3/4-5-3 (Chen-S).}
\label{proof: Proof of Chen-S for priority 1,2,3}
\noindent In the balanced DHV network under priority 1,2,3, we have
$$
R_{1,2,3}^{-1}(\Delta) =  \left[
 \begin{matrix}
    1 & 0 & 0 \mcr
    \displaystyle\frac{m_5}{m_4} & 1 & 0 \mcr
    \displaystyle\frac{m_6}{m_4} & \displaystyle\frac{m_6}{m_5} & 1 
  \end{matrix}
  \right],\;\;
R_{1,2,3}(\Delta) =  \left[
 \begin{matrix}
    1 & 0 & 0 \mcr
    -\displaystyle\frac{m_5}{m_4} & 1 & 0 \mcr
    0 & -\displaystyle\frac{m_6}{m_5} & 1 
  \end{matrix}
  \right],\\
  $$
where $R_{1,2,3}(\Delta)$ is a lower triangular matrix whose inverse exists and is non-negative. Thus, with Lemma \ref{lem:lower triangular matrix, Schur-S}, $R_{1,2,3}(\Delta)$ is both \textit{completely-$\calS$} and \textit{Schur-S}, and with $\alpha_1 <1$, it satisfies Chen-$\calS$ according to Lemma \ref{lem:Schur-S is stronger than Chen-S}.\\

\noindent Similarly, with

$$R_{4,2,3}(\Delta)^{-1} =  \left[
 \begin{matrix}
    \displaystyle\frac{1}{m_1} & 0  & 0  \mcr
    \displaystyle\frac{1}{m_1} & 1  & 0  \mcr
    \displaystyle\frac{1}{m_1} & \displaystyle\frac{m_6}{m_5}  & 1  
  \end{matrix}
  \right],\;\;
 R_{4,2,3}(\Delta) =  \left[
 \begin{matrix}
    m_1 &  0 & 0 \mcr
    -1 & 1 & 0  \mcr
    \displaystyle\frac{m_2-m_3}{m_5} & -\displaystyle\frac{m_6}{m_5} & 1  
  \end{matrix}
  \right];\\
  $$
  
  $$R_{4,5,3}(\Delta)^{-1} =  \left[
 \begin{matrix}
    \displaystyle\frac{1}{m_1} & \displaystyle\frac{m_4}{m_2} & 0 \mcr
    \displaystyle\frac{1}{m_1} & \displaystyle\frac{1}{m_2} & 0 \mcr
    \displaystyle\frac{1}{m_1} & \displaystyle\frac{1}{m_2} & 1 
  \end{matrix}
  \right],\;\;
 R_{4,5,3}(\Delta) =  \left[
 \begin{matrix}
    1 & m_1-1 & 0 \mcr
    -\displaystyle\frac{m_2}{m_1} & \displaystyle\frac{m_2}{m_1} & 0 \mcr
    0 & -1 & 1 
  \end{matrix}
  \right];\\
  $$
we know that  $R_{4,2,3}(\Delta)$ is lower triangular matrix whose inverse exists and is non-negative; $R_{4,5,3}(\Delta)$ is a non-singular \textit{M-matrix}. Thus, with Lemma \ref{lem:lower triangular matrix, Schur-S} and Lemma \ref{lem:M-matrix, Schur-S}, 
we can conclude that $R_{4,2,3}(\Delta)$ and $R_{4,5,3}(\Delta)$ satisfies Chen-$\calS$ with $\alpha_1<1$.
\eProof
\end{proof}

\vspace{0.2cm}

This completes the proof that the $(R,\theta)$ for the 8 static-priority cases are all Chen-$\calS$.

As stated above, the reflection matrices for the 8 balanced DHV networks following static-priority policies are all invertible. Now we check the sign of their determinants. We can calculate that $det(R_{4,2,3}) = m_1$, $det(R_{4,5,3}) = m_2$, $det(R_{4,5,6}) = m_3$, $det(R_{4,2,6}) = \displaystyle\frac{m_1m_3m_5}{m_1m_3m_5+m_2m_4m_6}$, $det(R_{1,2,3}) = det(R_{1,5,3}) =det(R_{1,2,6}) =det(R_{1,5,6}) =1$, where the subscript label indicates the classes with high priority in each station. The sign of the determinants are all positive. Thus, we can conclude that the reflection matrix for a balanced DHV network is invertible for any ratio matrix $\Delta$.

According to Theorem \ref{Thm: full statement for corner->inter}, in a balanced DHV network with $\alpha_1 <1$, the reflection matrix is invertible and $(R, \theta)$ is Chen-$\calS$ for any ratio matrix $\Delta$.

\subsection{The SSC Inequalities for Balanced DHV Networks
\label{Appendix: The SSC Inequalities for Balanced DHV Networks}}


The requirements in Theorem \ref{thm: LEGO paper (Theorem 2)} are satisfied if for all static priority policies, there exists a vector $h\in \mathbb{R}^K_{++}$ and a constant $r >0$ common to all static priority policies such that 
\begin{align}
\sum_{k \in \mathcal{H}_{\pi}:z_k>0} h_k\dot{\wbar{Z}}_k(t) \leq -r \text{ for any regular time } t \text{ with } \lVert \wbar{Z}_{\mathcal{H}_{\pi}}(t) \rVert >0.\label{equ: appendix_general_SSC_DHV}
\end{align}
We check the 8 static priority policies for (\ref{equ: appendix_general_SSC_DHV}) of balanced DHV networks in this appendix.

Given a policy $\pi$, let $\wbar{Z}^\pi_k(t)$ be the fluid level of class $k$ at time $t$ under the policy $\pi$, and $\dot{Y}^\pi(t) = \sum_{k \in \mathcal{H}_\pi:z_k>0} h_k\dot{\wbar{Z}^\pi_k}(t)$. Then, we need $\dot{Y}^\pi(t) <0$ with a common vector $h \in \Bbb{R}^K_{++}$ for any regular time $t$ with $\lVert \wbar{Z}_{\mathcal{H}_\pi}(t) \rVert >0$, i.e., at least one of the high-priority queues is not empty, under all the 8 static priority policies. We omit $\pi$ when fixed. Note that in balanced DHV networks, $m_k<1$ and thus $\mu_k>1$ for all $k \in \mathcal{K}$. Denote $I_k(t)$ as the input rate of class $k$ at time $t$; $O_k(t)$ as the output rate of class $k$ at time $t$. Then, $\dot{\wbar{Z}}_k(t) = I_k(t) - O_k(t) \text{ for } k\in \mathcal{K}; I_k(t) = O_{k-1}(t) \text{ for } k \in \mathcal{K}/\{1\};  I_1(t) = \alpha_1$. In this section, we use a set of high-priority classes to represent static priority policies, for example,
(4, 2, 6) stands for a static priority policy that class 4, 2 and 6 have high priorities in each station. \\

\noindent \textbf{Case 1:} (4, 2, 6)

\vspace{0.2cm}

\begin{enumerate}
    \item $\wbar{Z}_2(t)>0, \wbar{Z}_4(t)>0, \wbar{Z}_6(t)>0.$
    
    Since class 4,2,6 are all not empty and have priorities in their stations, servers focus on class 4,2,6 jobs, and $O_1(t) = O_3(t)= O_5(t) =0$, which means that $I_4(t) =I_2(t) = I_6(t) = 0$. Thus, $\dot{Y}(t) = h_2\dot{\wbar{Z}}_2(t)+h_4\dot{\wbar{Z}}_4(t)+h_6\dot{\wbar{Z}}_6(t) <0$, since $O_2(t), O_4(t), O_6(t) >0$.

    \vspace{0.2cm}
    
    \item $\wbar{Z}_2(t)>0, \wbar{Z}_6(t)>0, \wbar{Z}_4(t)=0.$
    
    Since class 2 and 6 are not empty and have priorities in their stations, $O_3(t) = O_5(t) = 0$. Thus, $\dot{\wbar{Z}}_6(t) = 0 - \mu_6 = -\mu_6$ and $O_2(t) = \mu_2$. Since class 4 have high priority in station 1 but queue 4 is empty and $I_4(t) = O_3(t) = 0$, station 1 focuses on class 1 jobs. If queue 1 is not empty, $O_1(t) = \mu_1$, and $\dot{\wbar{Z}}_2(t) = \mu_1-\mu_2$; if queue 1 is empty, $O_1(t) = \alpha_1$ (since $\alpha_1 <1<\mu_1$), and $\dot{\wbar{Z}}_2(t) = \alpha_1-\mu_2$. Thus, $\dot{Y}(t) = h_2\dot{\wbar{Z}}_2(t)+h_6\dot{\wbar{Z}}_6(t)= h_2(\mu_1-\mu_2) -h_6\mu_6$ or $h_2(\alpha_1-\mu_2) -h_6\mu_6$. Since $h_2 >0, \alpha_1 <\mu_1$, we have $h_2(\mu_1-\mu_2) -h_6\mu_6 > h_2(\alpha_1-\mu_2) -h_6\mu_6$. Thus, we need $h_2(\mu_1-\mu_2) -h_6\mu_6 <0$.
    
    \vspace{0.2cm}
    
    \item $\wbar{Z}_2(t)>0, \wbar{Z}_4(t)>0, \wbar{Z}_6(t)=0.$
    
    Similarly, since class 2 and 4 are not empty and have priorities in their stations, $O_1(t) = O_5(t)=0$. Thus, $\dot{\wbar{Z}}_2(t) = 0 - \mu_2 = -\mu_2$ and $O_4(t) = \mu_4$. Since class 6 have high priority in station 3 but queue 6 is empty and $I_6(t) =0$, station 3 focuses on class 3 jobs. If queue 3 is not empty, the $O_3(t) = \mu_3$, and $\dot{\wbar{Z}}_4(t) = \mu_3-\mu_4$; if queue 3 is empty, $O_3(t) = min\{\mu_2,\mu_3\}$, and $\dot{\wbar{Z}}_4(t) = min\{\mu_2,\mu_3\}-\mu_4 < \mu_3-\mu_4$. Thus, we need $\dot{Y}(t) = h_2\dot{\wbar{Z}}_2(t)+h_4\dot{\wbar{Z}}_4(t) \leq h_4(\mu_3-\mu_4) -h_2\mu_2 <0 $.

    \vspace{0.2cm}
    
    \item $\wbar{Z}_4(t)>0, \wbar{Z}_6(t)>0, \wbar{Z}_2(t)=0.$
    
    Since class 4 and 6 are not empty and have priorities in their stations, $O_1(t) =O_3(t) =0$. Thus, $\dot{\wbar{Z}}_4(t) = 0 - \mu_4 = -\mu_4$ and $O_6(t) = \mu_6$. Since class 2 have high priority in station 2 but queue 2 is empty and $I_2(t) = 0$, station 2 focuses on class 5 jobs. If queue 5 is not empty, $O_5(t) = \mu_5$, and $\dot{\wbar{Z}}_6(t) = \mu_5-\mu_6$; if queue 5 is empty, $O_5(t) = min\{\mu_4,\mu_5\}$, and $\dot{\wbar{Z}}_6(t) = min\{\mu_4,\mu_5\}-\mu_6 < \mu_5-\mu_6$. Thus, we need $\dot{Y}(t) = h_4\dot{\wbar{Z}}_4(t)+h_6\dot{\wbar{Z}}_6(t) \leq h_6(\mu_5-\mu_6) -h_4\mu_4 <0 $.

    \vspace{0.2cm}

    \item $\wbar{Z}_2(t)>0, \wbar{Z}_4(t)= \wbar{Z}_6(t)=0.$
    
    We need $\dot{Y}(t) = h_2 \dot{\wbar{Z}}_2(t) <0$, which is equivalent to $\dot{\wbar{Z}}_2(t) <0$. Since class 2 have priority in station 2 and queue 2 is not empty, $O_5(t)=0$ and $O_2(t) = \mu_2$. As queue 6 is empty, $I_6(t)= O_5(t) = 0$, we know that $O_6(t) = 0$ and station 3 focus on class 3 jobs. If queue 3 is not empty, $O_3(t) = \mu_3$; if queue 3 is empty, $O_3(t) = min\{\mu_2,\mu_3\} \leq\mu_3$. Thus, $\dot{\wbar{T}}_4(t) = \frac{O_3(t)}{\mu_4} \geq \frac{min\{\mu_2,\mu_3\}}{\mu_4}$.
    
    Since class 4 have priority in station 1 and queue 4 is empty, we have $\mu_4\dot{\wbar{T}}_4(t) = O_3(t)$. $\dot{\wbar{Z}}_2(t) <0$ is equivalent to $O_1(t) = \mu_1 \dot{\wbar{T}}_1(t) < O_2(t) = \mu_2$. Since $0\leq \dot{\wbar{T}}_1(t) \leq 1- \dot{\wbar{T}}_4(t)$, it is sufficient to consider $\mu_1(1-\frac{min\{\mu_2,\mu_3\}}{\mu_4}) < \mu_2$, which requires $\mu_1(1-\frac{\mu_2}{\mu_4}) < \mu_2$ and $\mu_1(1-\frac{\mu_3}{\mu_4}) < \mu_2$.

    $\mu_1(1-\frac{\mu_2}{\mu_4}) < \mu_2 \Leftrightarrow 1-\frac{\mu_2}{\mu_4} < \frac{\mu_2}{\mu_1}  \Leftrightarrow \frac{\mu_2}{\mu_4} + \frac{\mu_2}{\mu_1} >1 \Leftrightarrow  \frac{m_4}{m_2} + \frac{m_1}{m_2} >1\Leftrightarrow \frac{1}{m_2} > 1
    $ is true since $m_2<1$. 

    $\mu_1(1-\frac{\mu_3}{\mu_4}) < \mu_2 \Leftrightarrow 1-\frac{\mu_3}{\mu_4} < \frac{\mu_2}{\mu_1}  \Leftrightarrow \frac{\mu_2}{\mu_1} + \frac{\mu_3}{\mu_4} >1 \Leftrightarrow  \frac{m_1}{m_2} + \frac{m_4}{m_3} >1
    $ is true since $\frac{m_1}{m_2} + \frac{m_4}{m_3} >m_1+m_4=1$. 

    Thus, $\dot{Y}(t) = h_2 \dot{\wbar{Z}}_2(t) <0$ is satisfied.

    \vspace{0.2cm}
    
    \item $\wbar{Z}_4(t)>0, \wbar{Z}_2(t)=\wbar{Z}_6(t)=0.$
    
    We need $\dot{Y}(t) = h_4 \dot{\wbar{Z}}_4(t) <0$, which is equivalent to $\dot{\wbar{Z}}_4(t) <0$. Since class 4 have priority in station 1 and queue 4 is not empty, $O_1(t)=0$ and $O_4(t) = \mu_4$. As queue 2 is empty, $I_2(t)= O_1(t) = 0$, we know that $O_2(t) = 0$ and station 2 focus on class 5 jobs. If queue 5 is not empty, $O_5(t) = \mu_5$; if queue 5 is empty, $O_5(t) = min\{\mu_4,\mu_5\} \leq \mu_5$. Thus, $\dot{\wbar{T}}_6(t) = \frac{O_5(t)}{\mu_6} \geq \frac{min\{\mu_4,\mu_5\}}{\mu_6}$.
    
    Since class 6 have priority in station 3 and queue 6 is empty, we have $\mu_6\dot{\wbar{T}}_6(t) = O_5(t)$. $\dot{\wbar{Z}}_4(t) <0$ is equivalent to $O_3(t) = \mu_3\dot{\wbar{T}}_3(t) < O_4(t) = \mu_4$. Since $0\leq \dot{\wbar{T}}_3(t) \leq 1- \dot{\wbar{T}}_6(t)$, it is sufficient to consider $\mu_3(1-\frac{min\{\mu_4,\mu_5\}}{\mu_6}) < \mu_4$, which requires $\mu_3(1-\frac{\mu_4}{\mu_6}) < \mu_4$ and $\mu_3(1-\frac{\mu_5}{\mu_6}) < \mu_4$.

    $\mu_3(1-\frac{\mu_4}{\mu_6}) < \mu_4 \Leftrightarrow 1-\frac{\mu_4}{\mu_6} < \frac{\mu_4}{\mu_3}  \Leftrightarrow \frac{\mu_4}{\mu_3} + \frac{\mu_4}{\mu_6} >1 \Leftrightarrow  \frac{m_3}{m_4} + \frac{m_6}{m_4} >1\Leftrightarrow \frac{1}{m_4} > 1
    $ is true since $m_4<1$. 

    $\mu_3(1-\frac{\mu_5}{\mu_6}) < \mu_4 \Leftrightarrow 1-\frac{\mu_5}{\mu_6} < \frac{\mu_4}{\mu_3}  \Leftrightarrow \frac{\mu_4}{\mu_3} + \frac{\mu_5}{\mu_6} >1 \Leftrightarrow  \frac{m_3}{m_4} + \frac{m_6}{m_5} >1
    $ is true since $\frac{m_3}{m_4} + \frac{m_6}{m_5} >m_3+m_6=1$. 

    Thus, $\dot{Y}(t) = h_4 \dot{\wbar{Z}}_4(t) <0$ is satisfied.

    \vspace{0.2cm}

    \item $\wbar{Z}_6(t)>0, \wbar{Z}_2(t)= \wbar{Z}_4(t)=0.$

    We need $\dot{Y}(t) = h_6 \dot{\wbar{Z}}_6(t) <0$, which is equivalent to $\dot{\wbar{Z}}_6(t) <0$. Since class 6 have priority in station 3 and queue 6 is not empty, $O_3(t)=0$ and $O_6(t) = \mu_6$. As queue 4 is empty, $I_4(t)= O_3(t) = 0$, we know that $O_4(t) = 0$ and station 1 focus on class 1 jobs. If queue 1 is not empty, $O_1(t) = \mu_1$; if queue 1 is empty, $O_1(t) = min\{\alpha_1,\mu_1\} =\alpha_1$. Thus, $\dot{\wbar{T}}_2(t) = \frac{O_1(t)}{\mu_2} \geq \frac{\alpha_1}{\mu_2}$.
    


    Since class 2 have priority in station 5 and queue 2 is empty, we have $\mu_2\dot{\wbar{T}}_2(t) = O_1(t)$. $\dot{\wbar{Z}}_6(t) <0$ is equivalent to $O_5(t) = \mu_5 \dot{\wbar{T}}_5(t) < O_6(t) = \mu_6$. Since $0\leq \dot{\wbar{T}}_5(t) \leq 1- \dot{\wbar{T}}_6(t)$, it is sufficient to consider $\mu_5(1-\frac{\alpha_1}{\mu_2}) < \mu_6$.  $\mu_5(1-\frac{\alpha_1}{\mu_2}) < \mu_6 \Leftrightarrow 1-\frac{\alpha_1}{\mu_2} < \frac{\mu_6}{\mu_5}  \Leftrightarrow \frac{\mu_6}{\mu_5} + \frac{\alpha_1}{\mu_2} >1 \Leftrightarrow  \frac{m_5}{m_6} + \alpha_1 m_2 >1 \Leftrightarrow \alpha_1 > \frac{m_6-m_5}{m_2m_6}$. 

    Thus, $\dot{Y}(t) = h_6 \dot{\wbar{Z}}_6(t) <0$ is satisfied when $\alpha_1 > \frac{m_6-m_5}{m_2m_6}$.

    Note that $\frac{m_6-m_5}{m_2m_6}<1$, thus, the condition $\alpha_1 > \frac{m_6-m_5}{m_2m_6}$ does not conflict with the nominal condition($\alpha_1<1$). It is clear that $\frac{m_6-m_5}{m_2m_6}<1$ since $\frac{m_6-m_5}{m_2m_6}<1 \Leftrightarrow m_6-m_5 < m_2m_6 \Leftrightarrow  m_5 > m_6(1-m_2) \Leftrightarrow m_5 > m_6m_5 
    \Leftrightarrow m_5(1-m_6) >0 \Leftrightarrow m_5m_3 >0$, which is true obviously.

\end{enumerate}

\vspace{0.2cm}

To sum up, except for the positive condition for common vector $h$, we need the following conditions from case 1:
\begin{align*}
    & h_2(\mu_1-\mu_2)-h_6\mu_6 <0,  \\
    & h_4(\mu_3-\mu_4)-h_2\mu_2 <0,  \\
    & h_6(\mu_5-\mu_6)-h_4\mu_4 <0,\\
    & \alpha_1 > \frac{m_6-m_5}{m_2m_6}.
\end{align*}

\noindent \textbf{Case 2:} (1, 2, 3)

\vspace{0.2cm}

\begin{enumerate}
    \item $\wbar{Z}_1(t)>0, \wbar{Z}_2(t)>0, \wbar{Z}_3(t)>0.$

    It is clear that $\dot{\wbar{Z}}_1(t) = \alpha_1-\mu_1$, $\dot{\wbar{Z}}_2(t) = \mu_1-\mu_2$ and $\dot{\wbar{Z}}_3(t) = \mu_2-\mu_3$.  Thus, we need $\dot{Y}(t) = h_1\dot{\wbar{Z}}_1(t)+h_2\dot{\wbar{Z}}_2(t)+h_3\dot{\wbar{Z}}_3(t) = h_1(\alpha_1-\mu_1) + h_2(\mu_1-\mu_2) + h_3(\mu_2-\mu_3) <0$.

    \vspace{0.2cm}

    \item $\wbar{Z}_1(t)>0, \wbar{Z}_2(t)>0, \wbar{Z}_3(t)=0.$

    It is clear that $\dot{\wbar{Z}}_1(t) = \alpha_1-\mu_1$ and $\dot{\wbar{Z}}_2(t) = \mu_1-\mu_2$. Thus, we need $\dot{Y}(t) = h_1\dot{\wbar{Z}}_1(t)+h_2\dot{\wbar{Z}}_2(t)= h_1(\alpha_1-\mu_1) + h_2(\mu_1-\mu_2) <0$.

    \vspace{0.2cm}

    \item $\wbar{Z}_1(t)>0, \wbar{Z}_3(t)>0, \wbar{Z}_2(t)=0.$
    
    It is clear that $\dot{\wbar{Z}}_1(t) = \alpha_1-\mu_1$. Since class 2 has priority in station 2 and queue 2 is empty, we have  $O_2(t) = min\{\mu_1,\mu_2\}$. If $O_2(t) = \mu_2$, $t$ is not differentiable, since it will turn to the situation that $\wbar{Z}_2(t)>0$ immediately. Thus, $O_2(t) = \mu_1$. Then we can get $\dot{\wbar{Z}}_3(t) = \mu_1-\mu_3$ and we need $\dot{Y}(t) = h_1\dot{\wbar{Z}}_1(t)+h_3\dot{\wbar{Z}}_3(t)= h_1(\alpha_1-\mu_1) + h_3(\mu_1-\mu_3) <0$.

    \vspace{0.2cm}

    \item $\wbar{Z}_2(t)>0, \wbar{Z}_3(t)>0, \wbar{Z}_1(t)=0.$

    Since class 1 has priority in station 1 and queue 1 is empty, we have  $O_1(t) = min\{\mu_1,\alpha_1\} = \alpha_1$. Then we can get $\dot{\wbar{Z}}_2(t) = \alpha_1-\mu_2$, $\dot{\wbar{Z}}_3(t) = \mu_2-\mu_3$ and we need $\dot{Y}(t) = h_2\dot{\wbar{Z}}_2(t)+h_3\dot{\wbar{Z}}_3(t)=h_2(\alpha_1-\mu_2) + h_3(\mu_2-\mu_3) <0$.

    \vspace{0.2cm}

    \item $\wbar{Z}_1(t)>0, \wbar{Z}_2(t)= \wbar{Z}_3(t)=0.$

    It is clear that $\dot{\wbar{Z}}_1(t) = \alpha_1-\mu_1 <0$ and $\dot{Y}(t) = h_1 \dot{\wbar{Z}}_1(t) <0$.

    \vspace{0.2cm}

    \item $\wbar{Z}_2(t)>0, \wbar{Z}_1(t)=\wbar{Z}_3(t)=0.$
    
    Since class 1 has priority in station 1 and queue 1 is empty, we have  $O_1(t) = min\{\mu_1,\alpha_1\} = \alpha_1$. Thus, it is clear that $\dot{\wbar{Z}}_2(t) = \alpha_1-\mu_2 <0$ and $\dot{Y}(t) = h_2 \dot{\wbar{Z}}_2(t) <0$.

    \vspace{0.2cm}

    \item $\wbar{Z}_3(t)>0, \wbar{Z}_1(t)= \wbar{Z}_2(t)=0.$

    Similarly, it is clear that $\dot{\wbar{Z}}_3(t) = \alpha_1-\mu_3 <0$ and $\dot{Y}(t) = h_3 \dot{\wbar{Z}}_3(t) <0$.

\end{enumerate}

\vspace{0.2cm}

To sum up, except for the positive condition for common vector $h$, we need the following conditions from case 2:
\begin{align*}
    & h_1(\alpha_1-\mu_1) + h_2(\mu_1-\mu_2) + h_3(\mu_2-\mu_3) <0,  \\
    & h_1(\alpha_1-\mu_1) + h_2(\mu_1-\mu_2) <0,  \\
    &h_1(\alpha_1-\mu_1) + h_3(\mu_1-\mu_3) <0,\\
    &h_2(\alpha_1-\mu_2) + h_3(\mu_2-\mu_3) <0.
\end{align*}

\noindent \textbf{Case 3:} (4, 2, 3)
\vspace{0.2cm}

\begin{enumerate}
    \item $\wbar{Z}_2(t)>0, \wbar{Z}_3(t)>0, \wbar{Z}_4(t)>0.$
    
    Since class 4 is not empty and has high priority in station 1, thus $O_1(t) = 0$ and  $\dot{\wbar{Z}}_2(t) = 0-\mu_2 = -\mu_2$. It is easy to know that $\dot{\wbar{Z}}_3(t) = \mu_2-\mu_3$ and $\dot{\wbar{Z}}_4(t) = \mu_3-\mu_4$.  Thus, we need $\dot{Y}(t) = h_2\dot{\wbar{Z}}_2(t)+h_3\dot{\wbar{Z}}_3(t)+h_4\dot{\wbar{Z}}_4(t) = h_3(\mu_2-\mu_3) + h_4(\mu_3-\mu_4) -h_2\mu_2<0$.

    \vspace{0.2cm}

    \item $\wbar{Z}_2(t)>0, \wbar{Z}_3(t)>0, \wbar{Z}_4(t)=0.$
    
    Since class 2 has priority in station 2 and queue 2 is not empty, we have $O_2(t) = \mu_2$. Similarly, $O_3(t) = \mu_3$. Then, $\dot{\wbar{Z}}_3(t) = \mu_3-\mu_2$. To make $t$ differentiable, we have  $\dot{\wbar{Z}}_4(t) = \mu_3-\mu_4\dot{\wbar{T}}_4(t) = 0 \Rightarrow \dot{\wbar{T}}_4(t) = \frac{\mu_3}{\mu_4} \Rightarrow \dot{\wbar{T}}_1(t)_{max} = 1-\frac{\mu_3}{\mu_4} \Rightarrow O_1(t)_{max} = \mu_1(1-\frac{\mu_3}{\mu_4}) \Rightarrow \dot{\wbar{Z}}_2(t)_{max} = \mu_1(1-\frac{\mu_3}{\mu_4}) - \mu_2$. Thus, we need $\dot{Y}(t)_{max} = h_2\dot{\wbar{Z}}_2(t)_{max}+h_3\dot{\wbar{Z}}_3(t)=h_2\left[\mu_1(1-\frac{\mu_3}{\mu_4}) - \mu_2\right] + h_3(\mu_2-\mu_3) <0$.

    \vspace{0.2cm}

    \item $\wbar{Z}_2(t)>0, \wbar{Z}_4(t)>0, \wbar{Z}_3(t)=0.$

    Since class 4 has priority in station 1 and queue 4 is not empty, we have $O_1(t) = 0$ and $O_4(t) = \mu_4$. It is clear that $O_2(t) = \mu_2$. Then, $\dot{\wbar{Z}}_2(t) = 0-\mu_2=-\mu_2$. To make $t$ differentiable, we have  $\dot{\wbar{Z}}_3(t) = \mu_2-O_3(t) = 0 \Rightarrow O_3(t) = \mu_2$. Then $\dot{\wbar{Z}}_4(t) = \mu_2-\mu_4$. Thus, we need $\dot{Y}(t) = h_2\dot{\wbar{Z}}_2(t)+h_4\dot{\wbar{Z}}_4(t)=h_4(\mu_2-\mu_4) -h_2\mu_2<0$.

    \vspace{0.2cm}

    \item $\wbar{Z}_3(t)>0, \wbar{Z}_4(t)>0, \wbar{Z}_2(t)=0.$

    It is clear that $O_3(t) = \mu_3$, $O_4(t) = \mu_4$, then $\dot{\wbar{Z}}_4(t) = \mu_3 -\mu_4$. Since class 4 has priority in station 1 and queue 4 is not empty, $O_1(t) = 0$. Since queue 2 is empty, we have $O_2(t) = 0 \Rightarrow dot{\wbar{Z}}_3(t) = 0-\mu_3 =-\mu_3$. Thus, we need $\dot{Y}(t) = h_3\dot{\wbar{Z}}_3(t)+h_4\dot{\wbar{Z}}_4(t)=h_4(\mu_3-\mu_4) -h_3\mu_3<0$.

    \vspace{0.2cm}

    \item $\wbar{Z}_2(t)>0, \wbar{Z}_3(t)= \wbar{Z}_4(t)=0.$

    We need $\dot{Y}(t) = h_2 \dot{\wbar{Z}}_2(t) <0$, which is equivalent to $\dot{\wbar{Z}}_2(t) <0$. It is clear that $O_2(t) = \mu_2$. To make $t$ differentiable, we have  $\dot{\wbar{Z}}_3(t) = \mu_2-O_3(t) = 0 \Rightarrow O_3(t) = \mu_2$ and $\dot{\wbar{Z}}_4(t) = \mu_2-\mu_4\dot{\wbar{T}}_4(t) = 0 \Rightarrow \dot{\wbar{T}}_4(t)  = \frac{\mu_2}{\mu_4} \Rightarrow \dot{\wbar{T}}_1(t)_{max}  = 1-\frac{\mu_2}{\mu_4} \Rightarrow  O_1(t)_{max} = \mu_1(1-\frac{\mu_2}{\mu_4}) \Rightarrow \dot{\wbar{Z}}_2(t)_{max} = O_1(t)_{max} - \mu_2 = \mu_1(1-\frac{\mu_2}{\mu_4}) - \mu_2$, which is needed to be negative. $\mu_1(1-\frac{\mu_2}{\mu_4}) - \mu_2 < 0$ is true since $\mu_1(1-\frac{\mu_2}{\mu_4}) - \mu_2 < 0 \Leftrightarrow \frac{\mu_2}{\mu_4}+ \frac{\mu_2}{\mu_1} >1 \Leftrightarrow \frac{m_4}{m_2}+ \frac{m_1}{m_2} >1 \Leftrightarrow \frac{1}{m_2} >1$, which is satisfied for sure.

    \vspace{0.2cm}

    \item $\wbar{Z}_3(t)>0, \wbar{Z}_4(t)=\wbar{Z}_2(t)=0.$

    We need $\dot{Y}(t) = h_3 \dot{\wbar{Z}}_3(t) <0$, which is equivalent to $\dot{\wbar{Z}}_3(t) <0$. It is clear that $O_3(t) = \mu_3$. To make $t$ differentiable, we have  $\dot{\wbar{Z}}_4(t) = \mu_3-O_4(t) = 0 \Rightarrow O_4(t) = \mu_4\dot{\wbar{T}}_4(t) = \mu_3 \Rightarrow \dot{\wbar{T}}_4(t) = \frac{\mu_3}{\mu_4} \Rightarrow \dot{\wbar{T}}_1(t)_{max} = 1-\frac{\mu_3}{\mu_4} \Rightarrow O_1(t)_{max} = \mu_1\dot{\wbar{T}}_1(t)_{max} =\mu_1 (1-\frac{\mu_3}{\mu_4})  $ and $\dot{\wbar{Z}}_2(t) = O_1(t)- O_2(t) = 0 \Rightarrow O_2(t)_{max} = O_1(t)_{max} = \mu_1 (1-\frac{\mu_3}{\mu_4}) \Rightarrow \dot{\wbar{Z}}_3(t)_{max} = O_2(t)_{max} - \mu_3 = \mu_1 (1-\frac{\mu_3}{\mu_4})-\mu_3$, which is needed to be negative. $\mu_1 (1-\frac{\mu_3}{\mu_4})-\mu_3 < 0$ is true since $ \mu_1 (1-\frac{\mu_3}{\mu_4})-\mu_3< 0 \Leftrightarrow \frac{\mu_3}{\mu_1}+ \frac{\mu_3}{\mu_4} >1 \Leftrightarrow \frac{m_4}{m_3}+ \frac{m_1}{m_3} >1 \Leftrightarrow \frac{1}{m_3} >1$, which is satisfied for sure.

    \vspace{0.2cm}

    \item $\wbar{Z}_4(t)>0, \wbar{Z}_2(t)= \wbar{Z}_3(t)=0.$

    We need $\dot{Y}(t) = h_4 \dot{\wbar{Z}}_4(t) <0$, which is equivalent to $\dot{\wbar{Z}}_4(t) <0$. Since class 4 has priority in station 1 and queue 4 is not empty, we have $O_4(t) = \mu_4$ and $O_1(t) = 0$. Since queue 2 is empty, we have $O_2(t) = 0$. Again, since queue 3 is empty, we have $O_3(t) = 0$. Thus, $\dot{\wbar{Z}}_4(t) = 0-\mu_4 <0$ for sure.

\end{enumerate}

\vspace{0.2cm}

To sum up, except for the positive condition for common vector $h$, we need the following conditions from case 3:
\begin{align*}
    & h_3(\mu_2-\mu_3) + h_4(\mu_3-\mu_4) - h_2\mu_2<0,  \\
    & h_2\left[\mu_1(1-\frac{\mu_3}{\mu_4})-\mu_2\right] + h_3(\mu_2-\mu_3)<0,  \\
    &h_4(\mu_2-\mu_4) -h_2\mu_2 <0,\\
    &h_4(\mu_3-\mu_4)-h_3\mu_3 <0.
\end{align*}

\noindent \textbf{Case 4:} (4, 5, 3)

\vspace{0.2cm}

\begin{enumerate}
    \item $\wbar{Z}_4(t)>0, \wbar{Z}_5(t)>0, \wbar{Z}_3(t)>0.$

    It is clear that $O_3(t) = \mu_3$, $O_4(t) = \mu_4$, $O_5(t) = \mu_5$,
    $O_2(t) = 0$, then $\dot{\wbar{Z}}_3(t) = 0-\mu_3 = -\mu_3$,
    $\dot{\wbar{Z}}_4(t) = \mu_3-\mu_4$ and $\dot{\wbar{Z}}_5(t) = \mu_4-\mu_5$. Thus, we need $\dot{Y}(t) = h_3\dot{\wbar{Z}}_3(t)+h_4\dot{\wbar{Z}}_4(t)+h_5\dot{\wbar{Z}}_5(t) = h_4(\mu_3-\mu_4) +h_5(\mu_4-\mu_5) - h_3\mu_3<0$.

    \vspace{0.2cm}

    \item $\wbar{Z}_3(t)>0, \wbar{Z}_4(t)>0, \wbar{Z}_5(t)=0.$

    It is clear that $O_3(t) = \mu_3$, $O_4(t) = \mu_4$,
    $O_1(t) = 0$, then $\dot{\wbar{Z}}_4(t) = \mu_3-\mu_4$. To make $t$ differentiable, we have $\dot{\wbar{Z}}_5(t) = \mu_4-O_5(t) = 0 \Rightarrow O_5(t) = \mu_5\dot{\wbar{T}}_5(t) = \mu_4 \Rightarrow \dot{\wbar{T}}_5(t) = \frac{\mu_4}{\mu_5} \Rightarrow \dot{\wbar{T}}_2(t)_{max} = 1- \frac{\mu_4}{\mu_5} \Rightarrow O_2(t)_{max} = \mu_2 \dot{\wbar{T}}_2(t)_{max}  = \mu_2(1- \frac{\mu_4}{\mu_5} ) \Rightarrow \dot{\wbar{Z}}_3(t)_{max} = O_2(t)_{max} - \mu_3 = \mu_2(1- \frac{\mu_4}{\mu_5} ) - \mu_3$. Thus, we need $\dot{Y}(t)_{max} = h_3\dot{\wbar{Z}}_3(t)_{max}+h_4\dot{\wbar{Z}}_4(t) = h_3\left[\mu_2(1- \frac{\mu_4}{\mu_5} ) - \mu_3\right] + h_4(\mu_3-\mu_4) <0$.

    \vspace{0.2cm}

    \item $\wbar{Z}_3(t)>0, \wbar{Z}_5(t)>0, \wbar{Z}_4(t)=0.$

    It is clear that $O_3(t) = \mu_3$, $O_5(t) = \mu_5$,
    $O_2(t) = 0$, then $\dot{\wbar{Z}}_3(t) = 0-\mu_3 = -\mu_3$. To make $t$ differentiable, we have $\dot{\wbar{Z}}_4(t) = \mu_3-O_4(t) = 0 \Rightarrow O_4(t) =  = \mu_3 \Rightarrow \dot{\wbar{Z}}_5(t) = \mu_3-\mu_5$. Thus, we need $\dot{Y}(t) = h_3\dot{\wbar{Z}}_3(t)+h_5\dot{\wbar{Z}}_5(t) = h_5(\mu_3-\mu_5) - h_3\mu_3<0$.

    \vspace{0.2cm}

    \item $\wbar{Z}_4(t)>0, \wbar{Z}_5(t)>0, \wbar{Z}_3(t)=0.$

    It is clear that $O_4(t) = \mu_4$, $O_5(t) = \mu_5$,
    $O_2(t) = 0$, then $\dot{\wbar{Z}}_5(t) = \mu_4-\mu_5$. To make $t$ differentiable, we have $\dot{\wbar{Z}}_3(t) = 0-O_3(t) = 0 \Rightarrow O_3(t) = 0 \Rightarrow \dot{\wbar{Z}}_4(t) = 0-\mu_4 = -\mu_4$. Thus, we need $\dot{Y}(t) = h_4\dot{\wbar{Z}}_4(t)+h_5\dot{\wbar{Z}}_5(t) = h_5(\mu_4-\mu_5) - h_4\mu_4<0$.

    \vspace{0.2cm}

    \item $\wbar{Z}_3(t)>0, \wbar{Z}_4(t)= \wbar{Z}_5(t)=0.$

     We need $\dot{Y}(t) = h_3 \dot{\wbar{Z}}_3(t) <0$, which is equivalent to $\dot{\wbar{Z}}_3(t) <0$. It is clear that $O_3(t) = \mu_3$. To make $t$ differentiable, we have $\dot{\wbar{Z}}_4(t) = O_4(t)-\mu_3 = 0 \Rightarrow O_4(t) = \mu_3$ and $\dot{\wbar{Z}}_5(t) = O_4(t)-O_5(t) = 0 \Rightarrow O_5(t) = \mu_5 \dot{\wbar{T}}_5(t)= \mu_3 \Rightarrow \dot{\wbar{T}}_5(t) = \frac{\mu_3}{\mu_5} \Rightarrow \dot{\wbar{T}}_2(t)_{max} = 1- \frac{\mu_3}{\mu_5} \Rightarrow O_2(t)_{max} = \mu_2\dot{\wbar{T}}_2(t)_{max} = \mu_2(1- \frac{\mu_3}{\mu_5}) \Rightarrow  \dot{\wbar{Z}}_3(t)_{max} = O_2(t)_{max} - \mu_3 = \mu_2(1- \frac{\mu_3}{\mu_5}) - \mu_3$, which is needed to be negative. $ \mu_2(1- \frac{\mu_3}{\mu_5}) - \mu_3< 0$ is true since $  \mu_2(1- \frac{\mu_3}{\mu_5}) - \mu_3< 0 \Leftrightarrow \frac{\mu_3}{\mu_2}+ \frac{\mu_3}{\mu_5} >1 \Leftrightarrow \frac{m_2}{m_3}+ \frac{m_5}{m_3} >1 \Leftrightarrow \frac{1}{m_3} >1$, which is satisfied for sure.

     \vspace{0.2cm}

    \item $\wbar{Z}_4(t)>0, \wbar{Z}_3(t)=\wbar{Z}_5(t)=0.$

    We need $\dot{Y}(t) = h_4 \dot{\wbar{Z}}_4(t) <0$, which is equivalent to $\dot{\wbar{Z}}_4(t) <0$. It is clear that $O_4(t) = \mu_4$. To make $t$ differentiable, we have  $\dot{\wbar{Z}}_5(t) = \mu_4-O_5(t) = 0 \Rightarrow O_5(t) = \mu_5\dot{\wbar{T}}_5(t) = \mu_4 \Rightarrow \dot{\wbar{T}}_5(t) = \frac{\mu_4}{\mu_5} \Rightarrow \dot{\wbar{T}}_2(t)_{max} = 1-\frac{\mu_4}{\mu_5} \Rightarrow O_2(t)_{max} = \mu_2\dot{\wbar{T}}_2(t)_{max} =\mu_2 (1-\frac{\mu_4}{\mu_5})  $ and $\dot{\wbar{Z}}_3(t) = O_2(t)- O_3(t) = 0 \Rightarrow O_2(t)_{max} = O_3(t)_{max} = \mu_2 (1-\frac{\mu_4}{\mu_5}) \Rightarrow \dot{\wbar{Z}}_4(t)_{max} = O_3(t)_{max} - \mu_4 = \mu_2(1-\frac{\mu_4}{\mu_5})-\mu_4$, which is needed to be negative. $\mu_2 (1-\frac{\mu_4}{\mu_5})-\mu_4 < 0$ is true since $ \mu_2 (1-\frac{\mu_4}{\mu_5})-\mu_4< 0 \Leftrightarrow \frac{\mu_4}{\mu_2}+ \frac{\mu_4}{\mu_5} >1 \Leftrightarrow \frac{m_2}{m_4}+ \frac{m_5}{m_4} >1 \Leftrightarrow \frac{1}{m_4} >1$, which is satisfied for sure.

    \vspace{0.2cm}

    \item $\wbar{Z}_5(t)>0, \wbar{Z}_3(t)= \wbar{Z}_4(t)=0.$

    We need $\dot{Y}(t) = h_5 \dot{\wbar{Z}}_5(t) <0$, which is equivalent to $\dot{\wbar{Z}}_5(t) <0$. Since class 5 has priority in station 2 and queue 5 is not empty, we have $O_5(t) = \mu_5$ and $O_2(t) = 0$. Since queue 3 is empty, we have $O_3(t) = 0$. Again, since queue 4 is empty, we have $O_4(t) = 0$. Thus, $\dot{\wbar{Z}}_5(t) = 0-\mu_5 <0$ for sure.

\end{enumerate}

\vspace{0.2cm}

To sum up, except for the positive condition for common vector $h$, we need the following conditions from case 4:
\begin{align*}
    & h_4(\mu_3-\mu_4)+h_5(\mu_4-\mu_5)-h_3\mu_3 <0,  \\
    & h_3\left[\mu_2(1-\frac{\mu_4}{\mu_5})-\mu_3\right]+h_4(\mu_3-\mu_4) <0,  \\
    & h_5(\mu_3-\mu_5)-h_3\mu_3 <0,  \\
    & h_5(\mu_4-\mu_5)-h_4\mu_4 <0. 
\end{align*}

\noindent \textbf{Case 5:} (4, 5, 6)

\vspace{0.2cm}

\begin{enumerate}
    \item $\wbar{Z}_4(t)>0, \wbar{Z}_5(t)>0, \wbar{Z}_6(t)>0.$

    It is clear that $O_4(t) = \mu_4$, $O_5(t) = \mu_5$, $O_6(t) = \mu_6$,
    $O_3(t) = 0$, then $\dot{\wbar{Z}}_4(t) = 0-\mu_4 = -\mu_4$,
    $\dot{\wbar{Z}}_5(t) = \mu_4-\mu_5$ and $\dot{\wbar{Z}}_6(t) = \mu_5-\mu_6$. Thus, we need $\dot{Y}(t) = h_4\dot{\wbar{Z}}_4(t)+h_5\dot{\wbar{Z}}_5(t)+h_6\dot{\wbar{Z}}_6(t) = h_5(\mu_4-\mu_5) +h_6(\mu_5-\mu_6) - h_4\mu_4<0$.

    \vspace{0.2cm}

    \item $\wbar{Z}_4(t)>0, \wbar{Z}_5(t)>0, \wbar{Z}_6(t)=0.$

    Since class 4 has priority in station 1 and queue 4 is not empty, we have $O_4(t) = \mu_4$. Similarly, $O_5(t) = \mu_5$. Then, $\dot{\wbar{Z}}_5(t) = \mu_4-\mu_5$. To make $t$ differentiable, we have  $\dot{\wbar{Z}}_6(t) = \mu_5-\mu_6\dot{\wbar{T}}_6(t) = 0 \Rightarrow \dot{\wbar{T}}_6(t) = \frac{\mu_5}{\mu_6} \Rightarrow \dot{\wbar{T}}_3(t)_{max} = 1-\frac{\mu_5}{\mu_6} \Rightarrow O_3(t)_{max} = \mu_3(1-\frac{\mu_5}{\mu_6}) \Rightarrow \dot{\wbar{Z}}_4(t)_{max} = \mu_3(1-\frac{\mu_5}{\mu_6}) - \mu_4$. Thus, we need $\dot{Y}(t)_{max} = h_4\dot{\wbar{Z}}_4(t)_{max}+h_5\dot{\wbar{Z}}_5(t)=h_4\left[\mu_3(1-\frac{\mu_5}{\mu_6}) - \mu_4\right] + h_5(\mu_4-\mu_5) <0$.

    \vspace{0.2cm}

    \item $\wbar{Z}_4(t)>0, \wbar{Z}_6(t)>0, \wbar{Z}_5(t)=0.$

    Since class 6 has priority in station 3 and queue 6 is not empty, we have $O_3(t) = 0$ and $O_6(t) = \mu_6$. It is clear that $O_4(t) = \mu_4$. Then, $\dot{\wbar{Z}}_4(t) = 0-\mu_4=-\mu_4$. To make $t$ differentiable, we have  $\dot{\wbar{Z}}_5(t) = \mu_4-O_5(t) = 0 \Rightarrow O_5(t) = \mu_4$. Then $\dot{\wbar{Z}}_6(t) = \mu_4-\mu_6$. Thus, we need $\dot{Y}(t) = h_4\dot{\wbar{Z}}_4(t)+h_6\dot{\wbar{Z}}_6(t)=h_6(\mu_4-\mu_6) -h_4\mu_4<0$.

    \vspace{0.2cm}

    \item $\wbar{Z}_5(t)>0, \wbar{Z}_6(t)>0, \wbar{Z}_4(t)=0.$

    It is clear that $O_5(t) = \mu_5$, $O_6(t) = \mu_6$,
    $O_3(t) = 0$, then $\dot{\wbar{Z}}_6(t) = \mu_5-\mu_6$. To make $t$ differentiable, we have $\dot{\wbar{Z}}_4(t) = 0-O_4(t) = 0 \Rightarrow O_4(t)  = 0 \Rightarrow \dot{\wbar{Z}}_5(t) = 0-\mu_5 = -\mu_5$. Thus, we need $\dot{Y}(t) = h_5\dot{\wbar{Z}}_5(t)+h_6\dot{\wbar{Z}}_6(t) = h_6(\mu_5-\mu_6) - h_5\mu_5<0$.

    \vspace{0.2cm}

    \item $\wbar{Z}_4(t)>0, \wbar{Z}_5(t)= \wbar{Z}_6(t)=0.$

    We need $\dot{Y}(t) = h_4 \dot{\wbar{Z}}_4(t) <0$, which is equivalent to $\dot{\wbar{Z}}_4(t) <0$. It is clear that $O_4(t) = \mu_4$. To make $t$ differentiable, we have $\dot{\wbar{Z}}_5(t) = \mu_4-O_5(t) = 0 \Rightarrow O_5(t) = \mu_4$ and $\dot{\wbar{Z}}_6(t) = O_5(t)-O_6(t) = 0 \Rightarrow O_6(t) = \mu_6 \dot{\wbar{T}}_6(t)= \mu_4 \Rightarrow \dot{\wbar{T}}_6(t) = \frac{\mu_4}{\mu_6} \Rightarrow \dot{\wbar{T}}_3(t)_{max} = 1- \frac{\mu_4}{\mu_6} \Rightarrow O_3(t)_{max} = \mu_3\dot{\wbar{T}}_3(t)_{max} = \mu_3(1- \frac{\mu_4}{\mu_6}) \Rightarrow  \dot{\wbar{Z}}_4(t)_{max} = O_3(t)_{max} - \mu_4 = \mu_3(1- \frac{\mu_4}{\mu_6}) - \mu_4$, which is needed to be negative. $ \mu_3(1- \frac{\mu_4}{\mu_6}) - \mu_4< 0$ is true since $ \mu_3(1- \frac{\mu_4}{\mu_6}) - \mu_4< 0 \Leftrightarrow \frac{\mu_4}{\mu_3}+ \frac{\mu_4}{\mu_6} >1 \Leftrightarrow \frac{m_3}{m_4}+ \frac{m_6}{m_4} >1 \Leftrightarrow \frac{1}{m_4} >1$, which is satisfied for sure.

    \vspace{0.2cm}

    \item $\wbar{Z}_5(t)>0, \wbar{Z}_4(t)=\wbar{Z}_6(t)=0.$

    We need $\dot{Y}(t) = h_5 \dot{\wbar{Z}}_5(t) <0$, which is equivalent to $\dot{\wbar{Z}}_5(t) <0$. It is clear that $O_5(t) = \mu_5$. To make $t$ differentiable, we have  $\dot{\wbar{Z}}_6(t) = \mu_5-O_6(t) = 0 \Rightarrow O_6(t) = \mu_6\dot{\wbar{T}}_6(t) = \mu_5 \Rightarrow \dot{\wbar{T}}_6(t) = \frac{\mu_5}{\mu_6} \Rightarrow \dot{\wbar{T}}_3(t)_{max} = 1-\frac{\mu_5}{\mu_6} \Rightarrow O_3(t)_{max} = \mu_3\dot{\wbar{T}}_3(t)_{max} =\mu_3 (1-\frac{\mu_5}{\mu_6})  $ and $\dot{\wbar{Z}}_4(t) = O_3(t)- O_4(t) = 0 \Rightarrow O_4(t)_{max} = O_3(t)_{max} = \mu_3 (1-\frac{\mu_5}{\mu_6}) \Rightarrow \dot{\wbar{Z}}_5(t)_{max} = O_4(t)_{max} - \mu_5 = \mu_3(1-\frac{\mu_5}{\mu_6})-\mu_5$, which is needed to be negative. $\mu_3 (1-\frac{\mu_5}{\mu_6})-\mu_5 < 0$ is true since $ \mu_3 (1-\frac{\mu_5}{\mu_6})-\mu_5< 0 \Leftrightarrow \frac{\mu_5}{\mu_6}+ \frac{\mu_5}{\mu_3} >1 \Leftrightarrow \frac{m_3}{m_5}+ \frac{m_6}{m_5} >1 \Leftrightarrow \frac{1}{m_5} >1$, which is satisfied for sure.

    \vspace{0.2cm}

    \item $\wbar{Z}_6(t)>0, \wbar{Z}_4(t)= \wbar{Z}_5(t)=0.$

     We need $\dot{Y}(t) = h_6 \dot{\wbar{Z}}_6(t) <0$, which is equivalent to $\dot{\wbar{Z}}_6(t) <0$. Since class 6 has priority in station 3 and queue 6 is not empty, we have $O_6(t) = \mu_6$ and $O_3(t) = 0$. Since queue 4 is empty, we have $O_4(t) = 0$. Again, since queue 5 is empty, we have $O_5(t) = 0$. Thus, $\dot{\wbar{Z}}_6(t) = 0-\mu_6 <0$ for sure.
      
\end{enumerate}

\vspace{0.2cm}

To sum up, except for the positive condition for common vector $h$, we need the following conditions from case 5:
\begin{align*}
    & h_5(\mu_4-\mu_5)+h_6(\mu_5-\mu_6)-h_4\mu_4 <0,\\
    & h_4\left[\mu_3(1-\frac{\mu_5}{\mu_6})-\mu_4\right]+h_5(\mu_4-\mu_5) <0,\\
    & h_6(\mu_4-\mu_6)-h_4\mu_4 <0,\\
    & h_6(\mu_5-\mu_6)-h_5\mu_5 <0.
\end{align*}

\noindent \textbf{Case 6:} (1, 5, 6)

\vspace{0.2cm}

\begin{enumerate}
    \item $\wbar{Z}_1(t)>0, \wbar{Z}_5(t)>0, \wbar{Z}_6(t)>0.$

    It is clear that $O_1(t) = \mu_1$, $O_5(t) = \mu_5$, $O_6(t) = \mu_6$,
    $O_4(t) = 0$, then $\dot{\wbar{Z}}_1(t) = \alpha_1-\mu_1$,
    $\dot{\wbar{Z}}_5(t) = 0-\mu_5 = -\mu_5$ and $\dot{\wbar{Z}}_6(t) = \mu_5-\mu_6$. Thus, we need $\dot{Y}(t) = h_1\dot{\wbar{Z}}_1(t)+h_5\dot{\wbar{Z}}_5(t)+h_6\dot{\wbar{Z}}_6(t) = h_1(\alpha_1-\mu_1) +h_6(\mu_5-\mu_6) - h_5\mu_5<0$.

    \vspace{0.2cm}

    \item $\wbar{Z}_1(t)>0, \wbar{Z}_5(t)>0, \wbar{Z}_6(t)=0.$

    It is clear that $O_1(t) = \mu_1$, $O_5(t) = \mu_5$,
    $O_4(t) = 0$, then $\dot{\wbar{Z}}_1(t) = \alpha_1-\mu_1$, $\dot{\wbar{Z}}_5(t) = 0-\mu_5 = -\mu_5$. Thus, we need $\dot{Y}(t) = h_1\dot{\wbar{Z}}_1(t)+h_5\dot{\wbar{Z}}_5(t) = h_1(\alpha_1-\mu_1) - h_5\mu_5<0$, which is obviously satisfied since $\alpha_1 < 1< \mu_1$.

    \vspace{0.2cm}

    \item $\wbar{Z}_1(t)>0, \wbar{Z}_6(t)>0, \wbar{Z}_5(t)=0.$

    It is clear that $O_1(t) = \mu_1$, $O_6(t) = \mu_6$,
    $O_4(t) = 0$, then $\dot{\wbar{Z}}_1(t) = \alpha_1-\mu_1$. Since $O_4(t) = 0$ and queue 5 is empty, we have $O_5(t) = 0 \Rightarrow \dot{\wbar{Z}}_6(t)  = 0-\mu_6 = -\mu_6$. Thus, we need $\dot{Y}(t) = h_1\dot{\wbar{Z}}_1(t)+h_6\dot{\wbar{Z}}_6(t) = h_1(\alpha_1-\mu_1) - h_6\mu_6<0$, which is obviously satisfied since $\alpha_1 < 1< \mu_1$.

    \vspace{0.2cm}

    \item $\wbar{Z}_5(t)>0, \wbar{Z}_6(t)>0, \wbar{Z}_1(t)=0.$

    Since class 6 has priority in station 3 and queue 6 is not empty, we have $O_6(t) = \mu_6$ and $O_3(t) = 0$. Similarly, $O_5(t) = \mu_5$. Then, $\dot{\wbar{Z}}_6(t) = \mu_5-\mu_6$. To make $t$ differentiable, we have  $\dot{\wbar{Z}}_1(t) = \alpha_1-\mu_1\dot{\wbar{T}}_1(t) = 0 \Rightarrow \dot{\wbar{T}}_1(t) = \frac{\alpha_1}{\mu_1} \Rightarrow \dot{\wbar{T}}_4(t)_{max} = 1-\frac{\alpha_1}{\mu_1} \Rightarrow O_4(t)_{max} = \mu_4(1-\frac{\alpha_1}{\mu_1}) \Rightarrow \dot{\wbar{Z}}_5(t)_{max} = \mu_4(1-\frac{\alpha_1}{\mu_1}) - \mu_5$. Thus, we need $\dot{Y}(t)_{max} = h_5\dot{\wbar{Z}}_5(t)_{max}+h_6\dot{\wbar{Z}}_6(t)=h_5\left[\mu_4(1-\frac{\alpha_1}{\mu_1}) - \mu_5\right] + h_6(\mu_5-\mu_6) <0$.

    \vspace{0.2cm}

    \item $\wbar{Z}_1(t)>0, \wbar{Z}_5(t)= \wbar{Z}_6(t)=0.$

    We need $\dot{Y}(t) = h_1 \dot{\wbar{Z}}_1(t) <0$, which is equivalent to $\dot{\wbar{Z}}_1(t) <0$. It is clear that $\dot{\wbar{Z}}_1(t) = \alpha_1 - \mu_1 <0$.

    \vspace{0.2cm}

    \item $\wbar{Z}_5(t)>0, \wbar{Z}_1(t)=\wbar{Z}_6(t)=0.$



    We need $\dot{Y}(t) = h_5 \dot{\wbar{Z}}_5(t) <0$, which is equivalent to $\dot{\wbar{Z}}_5(t) <0$. It is clear that $O_5(t) = \mu_5$. To make $t$ differentiable, we have  $\dot{\wbar{Z}}_1(t) = \alpha_1-\mu_1\dot{\wbar{T}}_1(t) = 0 \Rightarrow \dot{\wbar{T}}_1(t) = \frac{\alpha_1}{\mu_1} \Rightarrow \dot{\wbar{T}}_4(t)_{max} = 1- \frac{\alpha_1}{\mu_1} \Rightarrow O_4(t)_{max} = \mu_4\dot{\wbar{T}}_4(t)_{max}  = \mu_4( 1- \frac{\alpha_1}{\mu_1} ) \Rightarrow \dot{\wbar{Z}}_5(t)_{max} = O_4(t)_{max} - \mu_5 = \mu_4( 1- \frac{\alpha_1}{\mu_1} ) - \mu_5$, which is needed to be negative. $\mu_4( 1- \frac{\alpha_1}{\mu_1} ) - \mu_5<0 \Leftrightarrow 1-m_1\alpha_1 < \frac{m_4}{m_5} \Leftrightarrow \alpha_1 > \frac{m_5-m_4}{m_1m_5}$. Thus, $\dot{Y}(t) = h_5 \dot{\wbar{Z}}_5(t) <0$ is satisfied when $\alpha_1 > \frac{m_5-m_4}{m_1m_5}$.

    Note that $\frac{m_5-m_4}{m_1m_5}<1$, thus, the condition $\alpha_1 > \frac{m_5-m_4}{m_1m_5}$ does not conflict with the nominal condition($\alpha_1<1$). It is clear that $\frac{m_5-m_4}{m_1m_5}<1$ since $\frac{m_5-m_4}{m_1m_5}<1 \Leftrightarrow m_5-m_4 < m_1m_5 \Leftrightarrow  m_4 > m_5(1-m_1) \Leftrightarrow m_4 > m_5m_4 
    \Leftrightarrow m_4(1-m_5) >0 \Leftrightarrow m_4m_2 >0$, which is true obviously. 

    \vspace{0.2cm}

    \item $\wbar{Z}_6(t)>0, \wbar{Z}_1(t)= \wbar{Z}_5(t)=0.$



    We need $\dot{Y}(t) = h_6 \dot{\wbar{Z}}_6(t) <0$, which is equivalent to $\dot{\wbar{Z}}_6(t) <0$. It is clear that $O_6(t) = \mu_6$. To make $t$ differentiable, we have $\dot{\wbar{Z}}_5(t) = O_4(t)- O_5(t) = 0 \Rightarrow O_5(t) = O_4(t)$ and  $\dot{\wbar{Z}}_1(t) = \alpha_1-\mu_1\dot{\wbar{T}}_1(t) = 0 \Rightarrow \dot{\wbar{T}}_1(t) = \frac{\alpha_1}{\mu_1} \Rightarrow \dot{\wbar{T}}_4(t)_{max} = 1- \frac{\alpha_1}{\mu_1} \Rightarrow O_4(t)_{max} = \mu_4\dot{\wbar{T}}_4(t)_{max}  = \mu_4( 1- \frac{\alpha_1}{\mu_1} ) \Rightarrow O_5(t)_{max} = O_4(t)_{max} = \mu_4( 1- \frac{\alpha_1}{\mu_1} ) \Rightarrow \dot{\wbar{Z}}_6(t)_{max} = O_5(t)_{max} - \mu_6 = \mu_4( 1- \frac{\alpha_1}{\mu_1} ) - \mu_6$, which is needed to be negative. $\mu_4( 1- \frac{\alpha_1}{\mu_1} ) - \mu_6<0 \Leftrightarrow 1-m_1\alpha_1 <\frac{m_4}{m_6} \Leftrightarrow \alpha_1 > \frac{m_6-m_4}{m_1m_6}$. Thus, $\dot{Y}(t) = h_6 \dot{\wbar{Z}}_6(t) <0$ is satisfied when $\alpha_1 > \frac{m_6-m_4}{m_1m_6}$.

    Note that $\frac{m_6-m_4}{m_1m_6}<1$, thus, the condition $\alpha_1 > \frac{m_6-m_4}{m_1m_6}$ does not conflict with the nominal condition($\alpha_1<1$). It is clear that $\frac{m_6-m_4}{m_1m_6}<1$ since $\frac{m_6-m_4}{m_1m_6}<1 \Leftrightarrow m_6-m_4 < m_1m_6 \Leftrightarrow  m_4 > m_6(1-m_1) \Leftrightarrow m_4 > m_6m_4 
    \Leftrightarrow m_4(1-m_6) >0 \Leftrightarrow m_4m_3 >0$, which is true obviously.

\end{enumerate}

\vspace{0.2cm}

To sum up, except for the positive condition for common vector $h$, we need the following conditions from case 6:
\begin{align*}
    & h_1(\alpha_1-\mu_1)+h_6(\mu_5-\mu_6)-h_5\mu_5 <0 ,\\
    & h_5\left[\mu_4(1-\frac{\alpha_1}{\mu_1})-\mu_5\right]+h_6(\mu_5-\mu_6) <0,\\
    & \alpha_1 > \frac{m_5-m_4}{m_1m_5},\\
    & \alpha_1 > \frac{m_6-m_4}{m_1m_6}.
\end{align*}

\noindent \textbf{Case 7:} (1, 2, 6)

\vspace{0.2cm}

\begin{enumerate}
    \item $\wbar{Z}_1(t)>0, \wbar{Z}_2(t)>0, \wbar{Z}_6(t)>0.$

    It is clear that $O_1(t) = \mu_1$, $O_2(t) = \mu_2$, $O_6(t) = \mu_6$,
    $O_5(t) = 0$, then $\dot{\wbar{Z}}_1(t) = \alpha_1-\mu_1$,
    $\dot{\wbar{Z}}_2(t) = \mu_1-\mu_2 $ and $\dot{\wbar{Z}}_6(t) = 0-\mu_6 = -\mu_6$. Thus, we need $\dot{Y}(t) = h_1\dot{\wbar{Z}}_1(t)+h_2\dot{\wbar{Z}}_2(t)+h_6\dot{\wbar{Z}}_6(t) = h_1(\alpha_1-\mu_1) +h_2(\mu_1-\mu_2) - h_6\mu_6<0$.

    \vspace{0.2cm}

    \item $\wbar{Z}_1(t)>0, \wbar{Z}_2(t)>0, \wbar{Z}_6(t)=0.$

    It is clear that $O_1(t) = \mu_1$, $O_2(t) = \mu_2$,
    then $\dot{\wbar{Z}}_1(t) = \alpha_1-\mu_1$, $\dot{\wbar{Z}}_2(t) = \mu_1-\mu_2 $. Thus, we need $\dot{Y}(t) = h_1\dot{\wbar{Z}}_1(t)+h_2\dot{\wbar{Z}}_2(t) = h_1(\alpha_1-\mu_1) - h_2(\mu_1-\mu_2)<0$.

    \vspace{0.2cm}

    \item $\wbar{Z}_1(t)>0, \wbar{Z}_6(t)>0, \wbar{Z}_2(t)=0.$

    It is clear that $O_1(t) = \mu_1$, $O_6(t) = \mu_6$,
    $O_4(t) = 0$, then $\dot{\wbar{Z}}_1(t) = \alpha_1-\mu_1$. To make $t$ differentiable, we have  $\dot{\wbar{Z}}_2(t) = \mu_1-O_2(t) = \mu_1-\mu_2\dot{\wbar{T}}_2(t) = 0  \Rightarrow \dot{\wbar{T}}_2(t) = \frac{\mu_1}{\mu_2} \Rightarrow \dot{\wbar{T}}_5(t)_{max} = 1-\frac{\mu_1}{\mu_2} \Rightarrow O_5(t)_{max} = \mu_5\dot{\wbar{T}}_5(t)_{max}  = \mu_5(1-\frac{\mu_1}{\mu_2}) \Rightarrow \dot{\wbar{Z}}_6(t)_{max} = O_5(t)_{max} - \mu_6 = \mu_5(1-\frac{\mu_1}{\mu_2}) - \mu_6$. Thus, we need $\dot{Y}(t)_{max} = h_1\dot{\wbar{Z}}_1(t)+h_6\dot{\wbar{Z}}_6(t)_{max} = h_1(\alpha_1-\mu_1) +h_6\left[\mu_5(1-\frac{\mu_1}{\mu_2}) - \mu_6\right] <0$, which is true since $\alpha_1-\mu_1$ is obviously negative, and we can prove that $\mu_5(1-\frac{\mu_1}{\mu_2}) - \mu_6 <0$ as the following. If $1-\frac{\mu_1}{\mu_2} \leq 0$, then  $\mu_5(1-\frac{\mu_1}{\mu_2}) - \mu_6 <0$ for sure; if $1-\frac{\mu_1}{\mu_2} >0$, then $\mu_5(1-\frac{\mu_1}{\mu_2}) - \mu_6 <0 \Leftrightarrow \frac{\mu_1}{\mu_2} + \frac{\mu_6}{\mu_5} > 1 \Leftrightarrow \frac{m_2}{m_1} + \frac{m_5}{m_6} > 1 $, which is satisfied since $\frac{m_2}{m_1} + \frac{m_5}{m_6} > m_5+m_2 = 1$.

    \vspace{0.2cm}

    \item $\wbar{Z}_2(t)>0, \wbar{Z}_6(t)>0, \wbar{Z}_1(t)=0.$

    It is clear that $O_2(t) = \mu_2$, $O_6(t) = \mu_6$,
    $O_5(t) = 0$, then $\dot{\wbar{Z}}_6(t) = 0-\mu_6 = -\mu_6$. To make $t$ differentiable, we have  $\dot{\wbar{Z}}_1(t) = \alpha_1-O_1(t) = 0  \Rightarrow O_1(t) =\alpha_1 \Rightarrow \dot{\wbar{Z}}_2(t) = \alpha_1 - \mu_2$. It obviously holds that $\dot{Y}(t) = h_2\dot{\wbar{Z}}_2(t)+h_6\dot{\wbar{Z}}_6(t) = h_2(\alpha_1-\mu_2) - h_6\mu_6<0$ as required.

    \vspace{0.2cm}

    \item $\wbar{Z}_1(t)>0, \wbar{Z}_2(t)= \wbar{Z}_6(t)=0.$

    We need $\dot{Y}(t) = h_1 \dot{\wbar{Z}}_1(t) <0$, which is equivalent to $\dot{\wbar{Z}}_1(t) <0$. It is clear that $\dot{\wbar{Z}}_1(t) = \alpha_1 - \mu_1 <0$.

    \vspace{0.2cm}

    \item $\wbar{Z}_2(t)>0, \wbar{Z}_1(t)=\wbar{Z}_6(t)=0.$

    We need $\dot{Y}(t) = h_2 \dot{\wbar{Z}}_2(t) <0$, which is equivalent to $\dot{\wbar{Z}}_2(t) <0$.  It is clear that $O_2(t) = \mu_2$. To make $t$ differentiable, we have  $\dot{\wbar{Z}}_1(t) = \alpha_1-O_1(t) = 0 \Rightarrow O_1(t) = \alpha_1$. Then, $\dot{\wbar{Z}}_2(t) = \alpha_1-\mu_2 <0$ for sure.

    \vspace{0.2cm}

    \item $\wbar{Z}_6(t)>0, \wbar{Z}_1(t)= \wbar{Z}_2(t)=0.$



     We need $\dot{Y}(t) = h_6 \dot{\wbar{Z}}_6(t) <0$, which is equivalent to $\dot{\wbar{Z}}_6(t) <0$. It is clear that $O_6(t) = \mu_6$. To make $t$ differentiable, we have  $\dot{\wbar{Z}}_1(t) = \alpha_1-O_1(t) = 0 \Rightarrow O_1(t) = \alpha_1$ and $\dot{\wbar{Z}}_2(t) = \alpha_1-O_2(t) = 0 \Rightarrow O_2(t) = \mu_2 \dot{\wbar{T}}_2(t) = \alpha_1 \Rightarrow \dot{\wbar{T}}_2(t) = \frac{\alpha_1}{\mu_2} \Rightarrow \dot{\wbar{T}}_5(t)_{max}= 1- \frac{\alpha_1}{\mu_2} \Rightarrow O_5(t)_{max} = \mu_5\dot{\wbar{T}}_5(t)_{max} = \mu_5(1- \frac{\alpha_1}{\mu_2})\Rightarrow \dot{\wbar{Z}}_6(t)_{max} = O_5(t)_{max} - \mu_6 =  \mu_5(1- \frac{\alpha_1}{\mu_2})-\mu_6$, which is needed to be negative. $\mu_5( 1- \frac{\alpha_1}{\mu_2} ) - \mu_6<0 \Leftrightarrow 1-m_2\alpha_1 < \frac{m_5}{m_6} \Leftrightarrow \alpha_1 > \frac{m_6-m_5}{m_2m_6}$. Thus, $\dot{Y}(t) = h_6 \dot{\wbar{Z}}_6(t) <0$ is satisfied when $\alpha_1 > \frac{m_6-m_5}{m_2m_6}$.

     We state in the analysis of case 1 ($\wbar{Z}_6(t)>0, \wbar{Z}_2(t)= \wbar{Z}_4(t)=0$) that  $\frac{m_6-m_5}{m_2m_6}<1$, thus, the condition $\alpha_1 > \frac{m_6-m_5}{m_2m_6}$ does not conflict with the nominal condition($\alpha_1<1$).
     
\end{enumerate}

\vspace{0.2cm}

To sum up, except for the positive condition for common vector $h$, we need the following conditions from case 7:
\begin{align*}
    & h_1(\alpha_1-\mu_1)+h_2(\mu_1-\mu_2)-h_6\mu_6 <0,\\
    & h_1(\alpha_1-\mu_1)+h_2(\mu_1-\mu_2) <0,\\
    & \alpha_1 > \frac{m_6-m_5}{m_2m_6}.
\end{align*}

\noindent \textbf{Case 8:} (1, 5, 3)

\vspace{0.2cm}

\begin{enumerate}
    \item $\wbar{Z}_1(t)>0, \wbar{Z}_3(t)>0, \wbar{Z}_5(t)>0.$

    It is clear that $O_1(t) = \mu_1$, $O_3(t) = \mu_3$, $O_5(t) = \mu_5$,
    $O_4(t) = 0$, $O_2(t) = 0$, then $\dot{\wbar{Z}}_1(t) = \alpha_1-\mu_1$,
    $\dot{\wbar{Z}}_3(t) = 0-\mu_3= -\mu_3$ and  $\dot{\wbar{Z}}_5(t) = 0-\mu_5= -\mu_5$. It obviously holds that $\dot{Y}(t) = h_1\dot{\wbar{Z}}_1(t)+h_3\dot{\wbar{Z}}_3(t)+h_5\dot{\wbar{Z}}_5(t) = h_1(\alpha_1-\mu_1) -h_3\mu_3 - h_5\mu_5<0$ as required.

    \vspace{0.2cm}

    \item $\wbar{Z}_1(t)>0, \wbar{Z}_3(t)>0, \wbar{Z}_5(t)=0.$

    It is clear that $O_1(t) = \mu_1$, $O_3(t) = \mu_3$,
    $O_4(t) = 0$, then $\dot{\wbar{Z}}_1(t) = \alpha_1-\mu_1$. Since queue 5 is empty and $O_4(t) = 0$, we have $O_5(t) = 0 \Rightarrow O_2(t)_{max} = \mu_2 \Rightarrow \dot{\wbar{Z}}_3(t)_{max} = O_2(t)_{max} - \mu_3  = \mu_2-\mu_3$. Thus, we need $\dot{Y}(t)_{max} = h_1\dot{\wbar{Z}}_1(t)+h_3\dot{\wbar{Z}}_3(t)_{max}= h_1(\alpha_1-\mu_1) +h_3(\mu_2-\mu_3)<0$.

    \vspace{0.2cm}

    \item $\wbar{Z}_1(t)>0, \wbar{Z}_5(t)>0, \wbar{Z}_3(t)=0.$

    It is clear that $O_1(t) = \mu_1$, $O_5(t) = \mu_5$,
    $O_4(t) = 0$, then $\dot{\wbar{Z}}_1(t) = \alpha_1-\mu_1$, $\dot{\wbar{Z}}_5(t) = 0-\mu_5 = -\mu_5$. It obviously holds that $\dot{Y}(t) = h_1\dot{\wbar{Z}}_1(t)+h_5\dot{\wbar{Z}}_5(t)= h_1(\alpha_1-\mu_1) -h_5\mu_5<0$ as required.

    \vspace{0.2cm}

    \item $\wbar{Z}_3(t)>0, \wbar{Z}_5(t)>0, \wbar{Z}_1(t)=0.$



    It is clear that $O_3(t) = \mu_3$, $O_5(t) = \mu_5$,
    $O_2(t) = 0$, then $\dot{\wbar{Z}}_3(t) = 0-\mu_3 = -\mu_3$. To make $t$ differentiable, we have  $\dot{\wbar{Z}}_1(t) = \alpha_1-O_1(t) =\alpha_1-\mu_1\dot{\wbar{T}}_1(t) = 0 \Rightarrow \dot{\wbar{T}}_1(t) = \frac{\alpha_1}{\mu_1} \Rightarrow \dot{\wbar{T}}_4(t)_{max} = 1- \frac{\alpha_1}{\mu_1} \Rightarrow O_4(t)_{max} = \mu_4\dot{\wbar{T}}_4(t)_{max}  = \mu_4(1- \frac{\alpha_1}{\mu_1}) \Rightarrow \dot{\wbar{Z}}_5(t)_{max} = O_4(t)_{max} - \mu_5 = \mu_4(1- \frac{\alpha_1}{\mu_1}) -\mu_5$. Thus, we need $\dot{Y}(t)_{max} = h_3\dot{\wbar{Z}}_3(t)+h_5\dot{\wbar{Z}}_5(t)_{max}= h_5\left[\mu_4(1- \frac{\alpha_1}{\mu_1}) -\mu_5\right]- h_3\mu_3<0$. 

    In the analysis of case 6 ($\wbar{Z}_5(t)>0, \wbar{Z}_1(t)=\wbar{Z}_6(t)=0$), we need the condition $\alpha_1 > \frac{m_5-m_4}{m_1m_5}$, with which 
    $\mu_4( 1- \frac{\alpha_1}{\mu_1} ) - \mu_5<0$. Thus, without adding more conditions, we have $\dot{Y}(t)_{max} = h_5\left[\mu_4(1- \frac{\alpha_1}{\mu_1}) -\mu_5\right]- h_3\mu_3<0$.

    \vspace{0.2cm}

    \item $\wbar{Z}_1(t)>0, \wbar{Z}_3(t)= \wbar{Z}_5(t)=0.$

    We need $\dot{Y}(t) = h_1 \dot{\wbar{Z}}_1(t) <0$, which is equivalent to $\dot{\wbar{Z}}_1(t) <0$. It is clear that $\dot{\wbar{Z}}_1(t)  = \alpha_1 - \mu_1<0$.

    \vspace{0.2cm}

    \item $\wbar{Z}_3(t)>0, \wbar{Z}_1(t)=\wbar{Z}_5(t)=0.$

    We need $\dot{Y}(t) = h_3 \dot{\wbar{Z}}_3(t) <0$, which is equivalent to $\dot{\wbar{Z}}_3(t) <0$. It is clear that $O_3(t) = \mu_3$. To make $t$ differentiable, we have  $\dot{\wbar{Z}}_1(t) = \alpha_1-O_1(t) =0 \Rightarrow O_1(t) =\mu_1\dot{\wbar{T}}_1(t)= \alpha_1 \Rightarrow \dot{\wbar{T}}_1(t) = \frac{\alpha_1}{\mu_1}$ and $\dot{\wbar{Z}}_5(t) = O_4(t)-O_5(t) =0 \Rightarrow O_5(t) = O_4(t) \Rightarrow  \dot{\wbar{T}}_5(t) = \frac{O_5(t)}{\mu_5} = \frac{O_4(t)}{\mu_5} \Rightarrow  \dot{\wbar{T}}_5(t)_{min} = \frac{O_4(t)_{min}}{\mu_5}\Rightarrow \dot{\wbar{T}}_2(t)_{max} = 1- \dot{\wbar{T}}_5(t)_{min} = 1-\frac{O_4(t)_{min}}{\mu_5}$. Then $O_2(t)_{max} = \mu_2\dot{\wbar{T}}_2(t)_{max}  = \mu_2(1-\frac{O_4(t)_{min}}{\mu_5})$. If queue 4 is not empty, $O_4(t) = \mu_4(1-\dot{\wbar{T}}_1(t)) = \mu_4(1-\frac{\alpha_1}{\mu_1})$; if queue 4 is empty, $O_4(t) = min\{\mu_3, \mu_4(1-\frac{\alpha_1}{\mu_1})\}$. Thus, $O_4(t)_{min} = \mu_3$ or $\mu_4(1-\frac{\alpha_1}{\mu_1})$, and $O_2(t)_{max} = \mu_2(1-\frac{\mu_3}{\mu_5})$ or $\mu_2(1-\frac{\mu_4(1-\frac{\alpha_1}{\mu_1})}{\mu_5})$. Then, $\dot{\wbar{Z}}_3(t) <0$ is equivalent to $\mu_2(1-\frac{\mu_3}{\mu_5})-\mu_3<0$ and $\mu_2\left[1-\frac{\mu_4}{\mu_5}(1-\frac{\alpha_1}{\mu_1})\right] - \mu_3 <0$. 

    $\mu_2(1-\frac{\mu_3}{\mu_5})-\mu_3<0$ is satisfied since $\mu_2(1-\frac{\mu_3}{\mu_5})-\mu_3<0 \Leftrightarrow \frac{\mu_3}{\mu_2} + \frac{\mu_3}{\mu_5} >1 \Leftrightarrow \frac{m_2}{m_3} + \frac{m_5}{m_3} >1 \Leftrightarrow \frac{1}{m_3}>1$, which is true obviously.

    Let $f(\alpha_1) = \mu_2\left[1-\frac{\mu_4}{\mu_5}(1-\frac{\alpha_1}{\mu_1})\right] - \mu_3$. Then $\dot{f}(\alpha_1) = \frac{\mu_2\mu_4}{\mu_5\mu_1} > 0$. We know that $f(1) = \mu_2\left[1-\frac{\mu_4}{\mu_5}(1-\frac{1}{\mu_1})\right] - \mu_3=  \mu_2\left[1-\frac{\mu_4}{\mu_5}(1-m_1)\right] - \mu_3 = \mu_2(1-\frac{m_5}{m_4}\times m_4)-\mu_3 = \mu_2(1-m_5)-\mu_3= 1-\mu_3 <0 $. So, for any $\alpha_1 <1$, we have $f(\alpha_1) < f(1)<0$, which means $\mu_2\left[1-\frac{\mu_4}{\mu_5}(1-\frac{\alpha_1}{\mu_1})\right] - \mu_3 <0$. 

    \vspace{0.2cm}

    \item $\wbar{Z}_5(t)>0, \wbar{Z}_1(t)= \wbar{Z}_3(t)=0.$



    We need $\dot{Y}(t) = h_5 \dot{\wbar{Z}}_5(t) <0$, which is equivalent to $\dot{\wbar{Z}}_5(t) <0$. It is clear that $O_5(t) = \mu_5$. To make $t$ differentiable, we have $\dot{\wbar{Z}}_1(t)= \alpha_1- O_1(t)= 0 \Rightarrow O_1(t)= \mu_1\dot{\wbar{T}}_1(t)= \alpha_1 \Rightarrow \dot{\wbar{T}}_1(t) = \frac{\alpha_1}{\mu_1} \Rightarrow \dot{\wbar{T}}_4(t)_{max} = 1- \frac{\alpha_1}{\mu_1} \Rightarrow O_4(t)_{max} = \mu_4\dot{\wbar{T}}_4(t)_{max}  = \mu_4 (1- \frac{\alpha_1}{\mu_1} ) \Rightarrow \dot{\wbar{Z}}_5(t)_{max} = O_4(t)_{max} - \mu_5 =\mu_4 (1- \frac{\alpha_1}{\mu_1})-\mu_5 $, which is needed to be negative. $\mu_4 (1- \frac{\alpha_1}{\mu_1})-\mu_5 <0 \Leftrightarrow 1-m_1\alpha_1 < \frac{m_4}{m_5} \Leftrightarrow \alpha_1 > \frac{m_5-m_4}{m_1m_5}$. Thus, $\dot{Y}(t) = h_5 \dot{\wbar{Z}}_5(t) <0$ is satisfied when $\alpha_1 > \frac{m_5-m_4}{m_1m_5}$.

    We state in the analysis of case 6 ($\wbar{Z}_5(t)>0, \wbar{Z}_1(t)= \wbar{Z}_6(t)=0$) that  $\frac{m_5-m_4}{m_1m_5}<1$, thus, the condition $\alpha_1 > \frac{m_5-m_4}{m_1m_5}$ does not conflict with the nominal condition($\alpha_1<1$).

\end{enumerate}

\vspace{0.2cm}

To sum up, except for the positive condition for common vector $h$, we need the following condition from case 8:
\begin{align*}
    &h_1(\alpha_1-\mu_1)+h_3(\mu_2-\mu_3) <0,\\
    &\alpha_1 > \frac{m_5-m_4}{m_1m_5}.
\end{align*}

\subsection{Sufficient Conditions for SSC of Balanced DHV Networks
\label{Appendix: Sufficient Conditions for SSC of Balanced DHV Networks}}

The question of whether there exists such a vector $h$ satisfying the 23 inequalities that we mentioned in \S \ref{Balanced DHV Network Example (SSC)} can be turned into a LP problem $(P_h)$.
If $(P_h)$ has a strictly positive solution, then a vector $h$ that satisfies all 23 inequalities exists. The inequalities below are reordered by their algebraic form.
\begin{alignat*}{3}
(P_h)= \;
& \text{max} \; \;&& p\\
& \text{s.t.}
&& h_1, h_2, h_3, h_4, h_5, h_6 \geq p \\
&&& h_2(\mu_1-\mu_2)-h_6\mu_6 \leq -p\\
&&& h_4(\mu_3-\mu_4)-h_2\mu_2 \leq -p\\
&&& h_6(\mu_5-\mu_6)-h_4\mu_4 \leq -p\\
&&& h_1(\alpha_1-\mu_1)+h_3(\mu_2-\mu_3) \leq -p\\
&&& h_3(\mu_1-\mu_3)+h_1(\alpha_1-\mu_1) \leq -p\\
&&& h_4(\mu_2-\mu_4)-h_2\mu_2 \leq -p\\
&&& h_5(\mu_3-\mu_5)-h_3\mu_3 \leq -p\\
&&& h_6(\mu_4-\mu_6)-h_4\mu_4 \leq -p\\
&&& h_3(\mu_2-\mu_3)+h_2(\alpha_1-\mu_2) \leq -p\\
&&& h_4(\mu_3-\mu_4)-h_3\mu_3 \leq -p\\
&&& h_5(\mu_4-\mu_5)-h_4\mu_4 \leq -p\\
&&& h_6(\mu_5-\mu_6)-h_5\mu_5 \leq -p\\
&&& h_1(\alpha_1-\mu_1)+h_2(\mu_1-\mu_2) \leq -p\\
&&& h_2\left[\mu_1(1-\frac{\mu_3}{\mu_4})-\mu_2\right]+h_3(\mu_2-\mu_3) \leq -p\\
&&& h_3\left[\mu_2(1-\frac{\mu_4}{\mu_5})-\mu_3\right]+h_4(\mu_3-\mu_4) \leq -p\\
&&& h_4\left[\mu_3(1-\frac{\mu_5}{\mu_6})-\mu_4\right]+h_5(\mu_4-\mu_5) \leq -p\\
&&& h_5\left[\mu_4(1-\frac{\alpha_1}{\mu_1})-\mu_5\right]+h_6(\mu_5-\mu_6) \leq -p\\
&&& h_2(\mu_1-\mu_2)+h_3(\mu_2-\mu_3)+h_1(\alpha_1-\mu_1) \leq -p\\
&&& h_3(\mu_2-\mu_3)+h_4(\mu_3-\mu_4)-h_2\mu_2 \leq -p\\
&&& h_4(\mu_3-\mu_4)+h_5(\mu_4-\mu_5)-h_3\mu_3 \leq -p\\
&&& h_5(\mu_4-\mu_5)+h_6(\mu_5-\mu_6)-h_4\mu_4 \leq -p\\
&&& h_6(\mu_5-\mu_6)+h_1(\alpha_1-\mu_1)-h_5\mu_5 \leq -p\\
&&& p \leq 1
\end{alignat*}


The constraint $p\leq 1$ guarantees that the objective value is bounded from above, without which, if ($P_h$) has a positive solution, then by scaling both $h$ and $p$, this solution is infinite. We care only about whether a positive solution to ($P_h$) exists.


One obvious feasible solution to $(P_h)$ is $h_1= h_2 = h_3 = h_4 = h_5 = h_6 = p = 0$, with a corresponding objective value of 0. Since $(P_h)$ is a maximization problem, the optimal value of $(P_h)$ is non-negative. If the optimal value is strictly positive, then we can scale $h$ and $p$ to make it larger. Thus, the optimal value of $(P_h)$ is always $0$ or $1$. According to strong duality, the dual optimal of $(P_h)$ is also $0$ or $1$; call the dual problem as $(D_h)$.
\begin{footnotesize}
\begin{alignat}{2}
(D_h)= \max\quad & x_{29} \nonumber\\
\mbox{s.t.}\quad 
& -x_1 +(\alpha_1-\mu_1)(x_{10}+x_{11}+x_{19}+x_{24}+x_{28}) = 0 \label{SSC_DHV_appendix_eq1}\\
& -x_2 + (\mu_1-\mu_2)(x_7+x_{19}+x_{24})-\mu_2(x_8+x_{12}+x_{25})+(\alpha_1-\mu_2)x_{15}+\left[\mu_1(1-\frac{\mu_3}{\mu_4})-\mu_2\right]x_{20} = 0 \label{SSC_DHV_appendix_eq2}\\
& -x_3 +(\mu_1-\mu_3)x_{11}+(\mu_2-\mu_3)(x_{10}+x_{15}+x_{20}+x_{24}+x_{25})-\mu_3(x_{13}+x_{16}+x_{26})\\
& + \left[\mu_2(1-\frac{\mu_4}{\mu_5})-\mu_3\right]x_{21} = 0 \label{SSC_DHV_appendix_eq3}\\
& -x_4 +(\mu_3-\mu_4)(x_8+x_{16}+x_{21}+x_{25}+x_{26})-\mu_4(x_9+x_{14}+x_{17}+x_{27})+(\mu_2-\mu_4)x_{12}\\
& + \left[\mu_3(1-\frac{\mu_5}{\mu_6})-\mu_4\right]x_{22} = 0 \label{SSC_DHV_appendix_eq4}\\
& -x_5 +(\mu_3-\mu_5)x_{13}+(\mu_4-\mu_5)(x_{17}+x_{22}+x_{26}+x_{27})-\mu_5(x_{18}+x_{28}) + \left[\mu_4(1-\frac{\alpha_1}{\mu_1})-\mu_5\right]x_{23} = 0 \label{SSC_DHV_appendix_eq5}\\
& -x_6+(\mu_4-\mu_6)x_{14}+(\mu_5-\mu_6)(x_9+x_{18}+x_{23}+x_{27}+x_{28})-\mu_6x_7 =0  \label{SSC_DHV_appendix_eq6}\\
& \sum_{i=1}^{29} x_i = 1 \label{SSC_DHV_appendix_eq7}\\
& x_1, x_2, ..., x_{29} \geq 0 \label{SSC_DHV_appendix_eq8}
\end{alignat}
\end{footnotesize}
 If the optimal value of of $(D_h)$ is $0$, it is clear that $x_{29} = 0$. If the optimal value is $1$, with constraints (\ref{SSC_DHV_appendix_eq7}) and (\ref{SSC_DHV_appendix_eq8}), $x_{29} = 1$, $x_1 = x_2 = ... = x_{28} = 0$. This means that if $(D_h)$ has $x_{29} = 1$, $x_1 = x_2 = ... = x_{28} = 0$ as its optimal solution, then the vector $h$ that satisfies all the 23 inequalities exists.


Let us now consider the signs of coefficients of the variables in $(D_h)$. 
We have
\begin{align*}
    & \frac{m_1}{m_2}+\frac{m_4}{m_3} > m_1+ m_4 = 1 \\
    \Rightarrow \;& \frac{\mu_2}{\mu_1}+\frac{\mu_3}{\mu_4} >1\\
    \Rightarrow \;& 1-\frac{\mu_3}{\mu_4} < \frac{\mu_2}{\mu_1}\\
    \Rightarrow \;& \mu_1 ( 1-\frac{\mu_3}{\mu_4})-\mu_2<0.
\end{align*}
Similarly, we have
\begin{align*}
   \mu_2 ( 1-\frac{\mu_4}{\mu_5})-\mu_3<0,\\
   \mu_3 ( 1-\frac{\mu_5}{\mu_6})-\mu_4<0,\\
   \mu_4 ( 1-\frac{\alpha_1}{\mu_1})-\mu_5<0.
\end{align*}

Form constraints (\ref{SSC_DHV_appendix_eq1}) and (\ref{SSC_DHV_appendix_eq8}), since $\alpha_1 - \mu_1 <0$, $x_1, x_{10}, x_{11}, x_{19}, x_{24}$ and $x_{28}$ should always be $0$. $(D_h)$ can be simplified as:
\begin{footnotesize}
\begin{alignat}{2}
(D_{h1})= \max\quad & x_{29} \nonumber\\
\mbox{s.t.}\quad
& -x_2 + (\mu_1-\mu_2)x_7-\mu_2(x_8+x_{12}+x_{25})+(\alpha_1-\mu_2)x_{15} +\left[\mu_1(1-\frac{\mu_3}{\mu_4})-\mu_2\right])x_{20} = 0 \label{SSC_DHV_appendix_eq9}\\
& -x_3 +(\mu_2-\mu_3)(x_{15}+x_{20}+x_{25})-\mu_3(x_{13}+x_{16}+x_{26}) + \left[\mu_2(1-\frac{\mu_4}{\mu_5})-\mu_3\right]x_{21} = 0 \label{SSC_DHV_appendix_eq10}\\
& -x_4 +(\mu_3-\mu_4)(x_8+x_{16}+x_{21}+x_{25}+x_{26})-\mu_4(x_9+x_{14}+x_{17}+x_{27})+(\mu_2-\mu_4)x_{12}\\
& + \left[\mu_3(1-\frac{\mu_5}{\mu_6})-\mu_4\right]x_{22} = 0 \label{SSC_DHV_appendix_eq11}\\
& -x_5 +(\mu_3-\mu_5)x_{13}+(\mu_4-\mu_5)(x_{17}+x_{22}+x_{26}+x_{27})-\mu_5x_{18} + \left[\mu_4(1-\frac{\alpha_1}{\mu_1})-\mu_5\right]x_{23} = 0 \label{SSC_DHV_appendix_eq12}\\
& -x_6+(\mu_4-\mu_6)x_{14}+(\mu_5-\mu_6)(x_9+x_{18}+x_{23}+x_{27})-\mu_6x_7 =0  \label{SSC_DHV_appendix_eq13}\\
& \sum_{i=1}^{29} x_i = 1 \label{SSC_DHV_appendix_eq14}\\
& x_1, x_2, ..., x_{29} \geq 0 \label{SSC_DHV_appendix_eq15}
\end{alignat}
\end{footnotesize}
We divide our analysis into 4 cases. The universe
\begin{align*}
U 
&= \{\mu_1 \leq \mu_2\} \cup \{(\mu_4 \leq \mu_6 )\cap (\mu_5 \leq \mu_6 )\} \cup \{(\mu_1 >\mu_2) \cap  ((\mu_4>\mu_6) \cup (\mu_5 > \mu_6)) \}\\
&= \{\mu_1 \leq \mu_2\} \cup \{(\mu_4 \leq \mu_6 )\cap (\mu_5 \leq \mu_6 )\} \cup \{(\mu_1 >\mu_2) \cap  ((\mu_1<\mu_3) \cup (\mu_2< \mu_3)) \}\\
&= \{\mu_1 \leq \mu_2\} \cup \{(\mu_4 \leq \mu_6 )\cap (\mu_5 \leq \mu_6 )\} \cup \{(\mu_1 >\mu_2) \cap (\mu_3 > \mu_2) \}\\
&= \{\mu_1 \leq \mu_2\} \cup \{(\mu_4 \leq \mu_6 )\cap (\mu_5 \leq \mu_6 )\} \cup \{(\mu_1 >\mu_2) \cap (\mu_3 > \mu_2) \cap (\mu_3 \leq \mu_4) \cap (\mu_2 \leq \mu_4) \} \\
& \;\;\;\;\; \cup \{(\mu_1 >\mu_2) \cap (\mu_3 > \mu_2) \cap ((\mu_2>\mu_4) \cup (\mu_3 >\mu4))\}\\
&= \{\mu_1 \leq \mu_2\} \cup \{(\mu_4 \leq \mu_6 )\cap (\mu_5 \leq \mu_6 )\} \cup \{(\mu_1 > \mu_2 )\cap( \mu_4 \geq \mu_3 > \mu_2)\} \cup \{(\mu_1 > \mu_2) \cap (\mu_3 >\mu_2) \cap (\mu_3 > \mu_4)\}.
\end{align*}

\textbf{Case 1:} $\mu_1-\mu_2 \leq 0$. 
\vspace{0.2cm}

From constraint (\ref{SSC_DHV_appendix_eq9}), since $\mu_1-\mu_2 \leq 0$, $\alpha_1 - \mu_2 <0$ and $\mu_1(1-\frac{\mu_3}{\mu_4})-\mu_2 <0$, we have $x_2= x_7= x_8= x_{12}=x_{25}=x_{15}=x_{20} = 0$. $(D_{h1})$ can be further simplified as:
\begin{footnotesize}
\begin{alignat}{2}
(D_{h2})= \max\quad & x_{29} \nonumber\\
\mbox{s.t.}\quad
& -x_3 -\mu_3(x_{13}+x_{16}+x_{26}) + \left[\mu_2(1-\frac{\mu_4}{\mu_5})-\mu_3\right]x_{21} = 0 \label{SSC_DHV_appendix_eq16}\\
& -x_4 +(\mu_3-\mu_4)(x_{16}+x_{21}+x_{26})-\mu_4(x_9+x_{14}+x_{17}+x_{27}) + \left[\mu_3(1-\frac{\mu_5}{\mu_6})-\mu_4\right]x_{22} = 0 \label{SSC_DHV_appendix_eq17}\\
& -x_5 +(\mu_3-\mu_5)x_{13}+(\mu_4-\mu_5)(x_{17}+x_{22}+x_{26}+x_{27})-\mu_5x_{18} + \left[\mu_4(1-\frac{\alpha_1}{\mu_1})-\mu_5\right]x_{23} = 0 \label{SSC_DHV_appendix_eq18}\\
& -x_6+(\mu_4-\mu_6)x_{14}+(\mu_5-\mu_6)(x_9+x_{18}+x_{23}+x_{27})=0  \label{SSC_DHV_appendix_eq19}\\
& \sum_{i=1}^{29} x_i = 1 \label{SSC_DHV_appendix_eq20}\\
& x_1, x_2, ..., x_{29} \geq 0 \label{SSC_DHV_appendix_eq21}
\end{alignat}
\end{footnotesize}
Similarly, from constraint (\ref{SSC_DHV_appendix_eq16}), since $\mu_2 ( 1-\displaystyle\frac{\mu_4}{\mu_5})-\mu_3<0$, we have $x_3 = x_{13} =x_{16} = x_{26} = x_{21} = 0$. Repeat the above steps, from constraint (\ref{SSC_DHV_appendix_eq17}), we have $x_4 = x_9 = x_{14} = x_{17} = x_{27} = x_{22} = 0$; from constraint (\ref{SSC_DHV_appendix_eq18}), we have $x_5 = x_{18} = x_{23} =0$; from constraint (\ref{SSC_DHV_appendix_eq19}), we have $x_6 = 0$. So far, we have claimed that $x_1 = x_2 = ... = x_{28} = 0$, thus, $x_{29} = 1$, which is the optimal solution of $(D_h)$. 

\vspace{0.2cm}

\textbf{Case 2:} $(\mu_4 \leq \mu_6 )\cap (\mu_5 \leq \mu_6 )$.
\vspace{0.2cm}

From constraint (\ref{SSC_DHV_appendix_eq13}) of $(D_{h1})$, since $\mu_4 - \mu_6 \leq 0$ and $\mu_5 - \mu_6 \leq 0$, we have $x_7 = 0$. Once we have $x_7 = 0$, we follow the same argument as in Case 1. Thus, $x_1 = x_2 = ... = x_{28} = 0$, $x_{29} = 1$ is the optimal solution of $(D_h)$ under Case 2.

Before we analyze Case 3 and Case 4 separately, we find that Case 3 and Case 4 both satisfy $\mu_2- \mu_3 <0$ and $\mu_1-\mu_2 >0$. Thus, from constraint (\ref{SSC_DHV_appendix_eq10}), since $\mu_2- \mu_3 <0$ and $\mu_2(1-\displaystyle\frac{\mu_4}{\mu_5})-\mu_3 <0$, we have $x_3 = x_{15} = x_{20} = x_{25} = x_{13} = x_{16} = x_{26} = x_{21} = 0$. Also, from constraint (\ref{SSC_DHV_appendix_eq12}) in $(D_{h1})$, since $\mu_1-\mu_2 >0$ implies $\mu_4-\mu_5 <0$, along with $\mu_4(1-\displaystyle\frac{\alpha_1}{\mu_1})-\mu_5 <0$, the only coefficient that could be positive is $(\mu_3-\mu_5)$. However, $(\mu_3-\mu_5)$ is the coefficient of $x_{13}$ and we already know that $x_{13} = 0$. Thus, we have $x_5 = x_{17} = x_{22} = x_{26} = x_{27} = x_{18} = x_{23} = 0$. $(D_{h1})$ can be simplified as: 
\begin{alignat}{2}
(D_{h3})= \max\quad & x_{29} \nonumber\\
\mbox{s.t.}\quad
& -x_2 + (\mu_1-\mu_2)x_7-\mu_2(x_8+x_{12}) = 0 \label{SSC_DHV_appendix_eq22}\\
& -x_4 +(\mu_3-\mu_4)x_8-\mu_4(x_9+x_{14})+(\mu_2-\mu_4)x_{12} = 0 \label{SSC_DHV_appendix_eq23}\\
& -x_6+(\mu_4-\mu_6)x_{14}+(\mu_5-\mu_6)x_9-\mu_6x_7 =0  \label{SSC_DHV_appendix_eq24}\\
& x_2 + x_4 +x_6+ x_7 +x_8 +x_9+x_{12} +x_{14} = 1 \label{SSC_DHV_appendix_eq25}\\
& x_2,x_4 ..., x_{29} \geq 0 \label{SSC_DHV_appendix_eq26}
\end{alignat}

\textbf{Case 3:} $(\mu_1 > \mu_2 )\cap( \mu_4 \geq \mu_3 > \mu_2)$.
\vspace{0.2cm}

From constraint (\ref{SSC_DHV_appendix_eq23}) in $(D_{h3})$, since $\mu_3-\mu_4 \leq 0$ and $\mu_2-\mu_4 <0$, we have $x_4 = x_9 =x_{14} = 0$. Note that we cannot decide if $x_8 = x_{12} = 0$ here since 
$\mu_3-\mu_4$ could be $0$. Then, from constraint (\ref{SSC_DHV_appendix_eq24}) in $(D_{h3})$, since $x_{14} = x_9 = 0$, we have $x_6 = x_7 = 0$. Consequently, from constraint (\ref{SSC_DHV_appendix_eq22}) in $(D_{h3})$, since $x_7 = 0$, we have $x_2 = x_8 = x_{12} = 0$. Thus, $x_1 = x_2 = ... = x_{28} = 0$, $x_{29} = 1$ is the optimal solution of $(D_h)$ under Case 3. 

\vspace{0.2cm}

\textbf{Case 4:} $(\mu_1 > \mu_2) \cap (\mu_3 >\mu_2) \cap (\mu_3 > \mu_4)$.
\vspace{0.2cm}

In $(D_{h3})$, we see $x_2, x_4$ and $x_6$ as slackness variables, and constraints (\ref{SSC_DHV_appendix_eq22}) - (\ref{SSC_DHV_appendix_eq24}) can be rewritten as:
\begin{equation*}
\begin{split}
     &\left\{\begin{array}{lr}
         (\mu_1-\mu_2)x_7 \geq \mu_2(x_8+x_{12})\\
         (\mu_3-\mu_4)x_8 \geq \mu_4(x_9+x_{14})+(\mu_4-\mu_2)x_{12}\\
         (\mu_5-\mu_6)x_9 \geq \mu_6x_7 + (\mu_6-\mu_4)x_{14}
        \end{array}
    \right.\\
    \Leftrightarrow
       & \left\{\begin{array}{lr}
         (\mu_1-\mu_2)x_7 \geq \mu_2(x_8+x_{12})\\
         (\mu_3-\mu_4)(x_8+x_{12}) \geq \mu_4(x_9+x_{14})+(\mu_3-\mu_2)x_{12}\\
         (\mu_5-\mu_6)(x_9+x_{14}) \geq \mu_6x_7 + (\mu_5-\mu_4)x_{14}
        \end{array}
    \right.\\
    \Rightarrow
    & \;(\mu_1-\mu_2)x_7 \geq \mu_2 (x_8+x_{12})\\
    \Rightarrow
    & \;(\mu_1-\mu_2)x_7 \geq \frac{\mu_2\mu_4}{\mu_3-\mu_4}(x_9+x_{14})+ \frac{\mu_2(\mu_3-\mu_2)}{\mu_3-\mu_4}x_{12}\\
    \Rightarrow
    & \;(\mu_1-\mu_2)x_7 \geq \frac{\mu_2\mu_4}{\mu_3-\mu_4}\left[\frac{\mu_6}{\mu_5-\mu_6}x_7+\frac{\mu_5-\mu_4}{\mu_5-\mu_6}x_{14}\right]+ \frac{\mu_2(\mu_3-\mu_2)}{\mu_3-\mu_4}x_{12}\\
    \Leftrightarrow
    & \;(\mu_1-\mu_2)(\mu_3-\mu_4)(\mu_5-\mu_6)x_7 \geq \mu_2\mu_4\mu_6x_7 + \mu_2\mu_4(\mu_5-\mu_4)x_{14}+\mu_2(\mu_3-\mu_2)(\mu_5-\mu_6)x_{12}\\
    \Leftrightarrow
    & \;\left[(\mu_1-\mu_2)(\mu_3-\mu_4)(\mu_5-\mu_6)-\mu_2\mu_4\mu_6\right]x_7 \geq \mu_2\mu_4(\mu_5-\mu_4)x_{14}+\mu_2(\mu_3-\mu_2)(\mu_5-\mu_6)x_{12}
\end{split}    
\end{equation*}

Note that $(\mu_1 > \mu_2) \cap (\mu_3 >\mu_2) \cap (\mu_3 > \mu_4)$ implies $(\mu_5 > \mu_4) \cap (\mu_5 >\mu_6) \cap (\mu_1 > \mu_6)$. So we have $\mu_5-\mu_4 >0$, $\mu_3-\mu_2 >0$ and $\mu_5-\mu_6 >0$. Thus, the right hand side of the above inequality is larger than or equal to $0$. As Chen and Ye proved in \cite{chen2001existence}, under $m_2+m_4+m_6 <2$, we have $\left[(\mu_1-\mu_2)(\mu_3-\mu_4)(\mu_5-\mu_6)-\mu_2\mu_4\mu_6\right]<0$ since this is equivalent to $m_1m_3m_5-(m_2-m_1)(m_4-m_3)(m_6-m_5) >0$ and we can verify that  $m_1m_3m_5-(m_2-m_1)(m_4-m_3)(m_6-m_5) = (m_1+m_3-m_2)\left[(m_2-m_1)m_6+m_1m_5\right] = (2-m_4-m_6-m_2)\left[(m_2-m_1)m_6+m_1m_5\right]>0 $. 

\vspace{0.2cm}

Consequently, under $m_2+m_4+m_6 <2$, we have $x_7 = x_{14} = x_{12} = 0$. From constraint (\ref{SSC_DHV_appendix_eq22}), since $x_7 = 0$,we have $x_2=x_8 = 0$. From constraint (\ref{SSC_DHV_appendix_eq23}), since $x_8 = x_{12} = 0$, we have $x_4 = x_9 = 0$. From constraint (\ref{SSC_DHV_appendix_eq24}), since $x_{14} = x_9 =0$, we have $x_6=0$. So far, we have argued that under $m_2+m_4+m_6 <2$, $x_1 = x_2 = ... = x_{28} = 0$, thus $x_{29} = 1$, and this is the optimal solution of $(D_h)$ under Case 4 when $m_2+m_4+m_6 <2$.

To sum up, if we have $m_2+m_4+m_6<2$, then the optimal value of $(D_h)$ is $1$, which means that there exists a vector $h$ that satisfies all the 23 inequalities we derived for a balanced DHV network. According to Theorem \ref{thm: LEGO paper (Theorem 2)}, we can conclude that under $m_2+m_4+m_6<2$, $\alpha_1 \in (\max\{\frac{m_6-m_5}{m_2m_6}, \frac{m_5-m_4}{m_1m_5}, \frac{m_6-m_4}{m_1m_6},0\},1)$, a balanced DHV network satisfies SSC under any queue-ratio policies.

\section{Detailed Proofs for Balanced Push Started Lu-Kumar Networks}

\setcounter{equation}{0}
\renewcommand{\theequation}{\thesection.\arabic{equation}}

\subsection{Chen-$\calS$ in Balanced push started Lu-Kumar Networks
\label{Appendix: Chen-S in Balanced DHV2 Networks}}

In this appendix, we prove that the reflection matrices for all the static-priority cases of a balanced push started Lu-Kumar network ($m_1+m_3+m_4=1$; $m_2+m_5=1$) with $\alpha_1<1$ are Chen-$\calS$, and invertible and have same-sign determinants if $m_2m_4 > m_3m_5$. In turn, Theorem \ref{Thm: full statement for corner->inter} allows us to say that the refection matrix is invertible and Chen-$\calS$ for any ratio matrix $\Delta$ of a balanced push started Lu-Kumar network if $m_2m_4 > m_3m_5$.

For each of the 6 static priority policies, to prove that the corresponding $(R, \theta)$ is Chen-$\calS$, we must show that $R$ is completely-$\calS$, and consider all the 3 partitions of $\mathcal{J}$ and show that there exists a positive vector $h \in \mathbb{R}^2$ such that (\ref{equ:ChenS_def}) satisfies for all the 3 partitions.

\vspace{0.2cm}

\begin{proof} {Proof for lowest-priority 1, 2 (Chen-S).}
\label{proof: Proof of Chen-S for lowest-priority 1,2 (DHV2}

\noindent In the balanced push started Lu-Kumar network with class 1 having the lowest priority in station 1 and class 2 having the lowest priority in station 2, we have $J = 2$, $K = 5$, $\rho = \left[\alpha_1 \ \alpha_1\right]',$\\
$$R = 
 \left[
 \begin{matrix}
   1 & m_1-1 \mcr
   -\displaystyle\frac{m_2}{m_1} & \displaystyle\frac{m_2}{m_1} 
  \end{matrix}
  \right] ,
$$
\begin{center}
\begin{equation}\nonumber
\begin{aligned}
\theta = R(\rho-e)
= (\alpha_1-1)
\left[
 \begin{matrix}
   m_1 \mcr
   0
  \end{matrix}
  \right] .
\end{aligned}
\end{equation}
\end{center}

It is clear that the reflection matrix $R$ is completely-$\calS$, since $R$ is $2\times2$, $R^{-1} = CMQ\Delta>0$ and all its diagonal elements are positive.\\

\noindent \textbf{Partition 1: $a = 1; b= 2$}\\
$$\theta_b = \theta_{2} = 0,
    \;
  R_b^{-1} = R_{2}^{-1}= \frac{m_1}{m_2}
  \Rightarrow u = - R_b^{-1}\theta_b = 0 \Rightarrow 
  R_{ab}u = 0.$$
  
  \noindent Thus, $h_a^{'}\left[\theta_a + R_{ab}u\right] = h_1(\alpha_1-1)m_1 $. Since $\alpha_1<1$, taking $h_1 >0 $, we ensure $h_a^{'}\left[\theta_a + R_{ab}u\right] <0$.\\

\noindent \textbf{Partition 2: $a = 2; b= 1$}\\
$$\theta_b = \theta_{1} = (\alpha_1-1)m_1,
    \;
  R_b^{-1} = R_{1}^{-1}= 1
  \Rightarrow u = - R_b^{-1}\theta_b = (1-\alpha_1)m_1;$$
  
  $$
   R_{ab} = -\frac{m_2}{m_1} \Rightarrow 
  R_{ab}u = (\alpha_1-1)m_2.$$
  
  \noindent Thus, $h_a^{'}\left[\theta_a + R_{ab}u\right] = h_2(\alpha_1-1)m_2 $. Since $\alpha_1<1$, taking $h_2 >0 $, we ensure $h_a^{'}\left[\theta_a + R_{ab}u\right] <0$.\\

\noindent \textbf{Partition 3: $a = 1, 2; b= \varnothing$}\\
$$h_a^{'}\left[\theta_a + R_{ab}u\right] = h'\theta = (\alpha_1-1)h_1m_1$$
Since $\alpha_1<1$, taking $h_1>0$, we ensure $h_a^{'}\left[\theta_a + R_{ab}u\right] <0$.
\\

\noindent We have shown, then, that exists a positive vector $h \in{R^2}$  as required so that $(R,\theta)$ under the static-priority policy 1, 2 is Chen-$\calS$.
\eProof\\
\end{proof}

\begin{proof} {Proof for lowest-priority 1, 5 (Chen-S).}
\label{proof: Proof of Chen-S for lowest-priority 1,5 (DHV2)}

\noindent In the balanced push started Lu-Kumar network with class 1 having the lowest priority in station 1 and class 5 having the lowest priority in station 2, we have $J = 2, K = 5, \rho = \left[\alpha_1 \ \alpha_1\right]',$\\
$$R = 
 \left[
 \begin{matrix}
   m_1 & 0 \mcr
  -1 & 1
  \end{matrix}
  \right] ,
$$
\begin{center}
\begin{equation}\nonumber
\begin{aligned}
\theta = R(\rho-e)
= (\alpha_1-1)
\left[
 \begin{matrix}
   m_1 \mcr
   0
  \end{matrix}
  \right] .
\end{aligned}
\end{equation}
\end{center}

It is clear that the reflection matrix $R$ is completely-$\calS$, since $R$ is $2\times2$, $R^{-1} = CMQ\Delta>0$ and all its diagonal elements are positive.\\

\noindent \textbf{Partition 1: $a = 1; b= 2$}\\
$$\theta_b = \theta_{2} = 0,
    \;
  R_b^{-1} = R_{2}^{-1}= 1
  \Rightarrow u = - R_b^{-1}\theta_b = 0 \Rightarrow 
  R_{ab}u = 0.$$
  
  \noindent Thus, $h_a^{'}\left[\theta_a + R_{ab}u\right] = h_1(\alpha_1-1)m_1 $. Since $\alpha_1<1$, taking $h_1 >0 $, we ensure $h_a^{'}\left[\theta_a + R_{ab}u\right] <0$.\\

\noindent \textbf{Partition 2: $a = 2; b= 1$}\\
$$\theta_b = \theta_{1} = (\alpha_1-1)m_1,
    \;
  R_b^{-1} = R_{1}^{-1}= \frac{1}{m_1}
  \Rightarrow u = - R_b^{-1}\theta_b = 1-\alpha_1;$$
  
  $$
   R_{ab} = -1\Rightarrow 
  R_{ab}u = \alpha_1-1.$$
  
  \noindent Thus, $h_a^{'}\left[\theta_a + R_{ab}u\right] = h_2(\alpha_1-1)$. Since $\alpha_1<1$, taking $h_2 >0 $, we ensure $h_a^{'}\left[\theta_a + R_{ab}u\right] <0$.\\

\noindent \textbf{Partition 3: $a = 1,2; b= \varnothing$}\\
$$h_a^{'}\left[\theta_a + R_{ab}u\right] = h'\theta = (\alpha_1-1)h_1m_1$$
Since $\alpha_1<1$, taking $h_1>0$, we ensure $h_a^{'}\left[\theta_a + R_{ab}u\right] <0$.
\\

\noindent We have shown, then, that exists a positive vector $h \in{R^2}$  as required so that $(R,\theta)$ under the static-priority policy 1, 5 is Chen-$\calS$.
\eProof\\
\end{proof}

\begin{proof} {Proof for lowest-priority 3,2 (Chen-S).}
\label{proof: Proof of Chen-S for lowest-priority 3,2 (DHV2)}

\noindent In the balanced push started Lu-Kumar network with class 3 having the lowest priority in station 1 and class 2 having the lowest priority in station 2, we have $J = 2, K = 5, \rho = \left[\alpha_1 \ \alpha_1\right]',$\\
$$R = 
 \left[
 \begin{matrix}
   \displaystyle\frac{m_3}{m_2(1-m_1)} & -\displaystyle\frac{m_3}{m_2} \mcr
   -\displaystyle\frac{m_5}{1-m_1} & 1 
  \end{matrix}
  \right] ,
$$
\begin{center}
\begin{equation}\nonumber
\begin{aligned}
\theta = R(\rho-e)
= (\alpha_1-1)
\left[
 \begin{matrix}
   \displaystyle\frac{m_3m_1}{m_2(1-m_1)} \mcr
   -\displaystyle\frac{m_2-m_1}{1-m_1}
  \end{matrix}
  \right] .
\end{aligned}
\end{equation}
\end{center}

It is clear that the reflection matrix $R$ is completely-$\calS$, since $R$ is $2\times2$, $R^{-1} = CMQ\Delta>0$ and all its diagonal elements are positive.\\

\noindent \textbf{Partition 1: $a = 1; b= 2$}\\
$$\theta_b = \theta_{2} = (\alpha_1-1)\frac{m_2-m_1}{1-m_1},
    \;
  R_b^{-1} = R_{2}^{-1}= 1
  \Rightarrow u = - R_b^{-1}\theta_b = (1-\alpha_1)\frac{m_2-m_1}{1-m_1};$$
  
  $$
   R_{ab} = -\frac{m_3}{m_2} \Rightarrow 
  R_{ab}u = (\alpha_1-1)\frac{m_3(m_2-m_1)}{m_2(1-m_1)}.$$
  
  \noindent Thus, $h_a^{'}\left[\theta_a + R_{ab}u\right] = h_1(\alpha_1-1)\displaystyle\frac{m_1m_3+m_3(m_2-m_1)}{m_2(1-m_1)} = h_1(\alpha_1-1)\frac{m_3}{1-m_1} $. Since $\alpha_1<1$, taking $h_1 >0 $, we ensure $h_a^{'}\left[\theta_a + R_{ab}u\right] <0$.\\

\noindent \textbf{Partition 2: $a = 2; b= 1$}\\
$$\theta_b = \theta_{1} = (\alpha_1-1)\frac{m_3m_1}{m_2(1-m_1)},
    \;
  R_b^{-1} = R_{1}^{-1}= \frac{m_2(1-m_1)}{m_3}
  \Rightarrow u = - R_b^{-1}\theta_b = (1-\alpha_1)m_1;$$
  
  $$
   R_{ab} = -\frac{m_5}{1-m_1} \Rightarrow 
  R_{ab}u = (\alpha_1-1)\frac{m_1m_5}{1-m_1}.$$
  
  \noindent Thus, $h_a^{'}\left[\theta_a + R_{ab}u\right] = h_2(\alpha_1-1)\displaystyle\frac{m_2-m_1+m_1m_5}{1-m_1} = h_2(\alpha_1-1)\frac{m_2-m_1(1-m_5)}{1-m_1} = h_2(\alpha_1-1)\frac{m_2-m_1m_2}{1-m_1} = (\alpha_1-1)h_2m_2$. Since $\alpha_1<1$, taking $h_2 >0 $, we ensure $h_a^{'}\left[\theta_a + R_{ab}u\right] <0$.\\

\noindent \textbf{Partition 3: $a = 1,2; b= \varnothing$}\\
$$h_a^{'}\left[\theta_a + R_{ab}u\right] = h'\theta = (\alpha_1-1)\left[\frac{m_3m_1}{m_2(1-m_1)}h_1+\frac{m_2-m_1}{1-m_1}h_2\right]$$
Since $\alpha_1<1$, $\displaystyle\frac{m_3m_1}{m_2(1-m_1)}>0$, fixing $h_1>0$ large enough, we ensure $h_a^{'}\left[\theta_a + R_{ab}u\right] <0$.
\\

\noindent We have shown, then, that exists a positive vector $h \in{R^2}$ with $h_1$ being large enough as required so that $(R,\theta)$ under the static-priority policy 3, 2 is Chen-$\calS$.
\eProof\\
\end{proof}

\begin{proof} {Proof for lowest-priority 3, 5 (Chen-S).}
\label{proof: Proof of Chen-S for lowest-priority 3,5 (DHV2)}

\noindent In the balanced push started Lu-Kumar network with class 3 having lowest priority in station 1 and class 5 having lowest priority in station 2, we have $J = 2$, $K = 5$, $\rho = \left[\alpha_1 \ \alpha_1\right]',$\\
$$R = 
 \left[
 \begin{matrix}
   \displaystyle\frac{m_3}{1-m_1} & 0 \mcr
   -\displaystyle\frac{m_5}{1-m_1} & 1
  \end{matrix}
  \right] ,
$$
\begin{center}
\begin{equation}\nonumber
\begin{aligned}
\theta = R(\rho-e)
= \frac{\alpha_1-1}{1-m_1}
\left[
 \begin{matrix}
   m_3 \mcr
   m_2-m_1
  \end{matrix}
  \right] .
\end{aligned}
\end{equation}
\end{center}

It is clear that the reflection matrix $R$ is completely-$\calS$, since $R$ is $2\times2$, $R^{-1} = CMQ\Delta>0$ and all its diagonal elements are positive.\\

\noindent \textbf{Partition 1: $a = 1; b= 2$}\\
$$\theta_b = \theta_{2} = (\alpha_1-1)\frac{m_2-m_1}{1-m_1},
    \;
  R_b^{-1} = R_{2}^{-1}= 1
  \Rightarrow u = - R_b^{-1}\theta_b = (\alpha_1-1)\frac{m_1-m_2}{1-m_1};$$
  
  $$
   R_{ab} = 0 \Rightarrow 
  R_{ab}u = 0.$$
  
  \noindent Thus, $h_a^{'}\left[\theta_a + R_{ab}u\right] = h_1(\alpha_1-1)\displaystyle\frac{m_3}{1-m_1}$. Since $\alpha_1<1$, taking $h_1 >0 $, we ensure $h_a^{'}\left[\theta_a + R_{ab}u\right] <0$.\\

\noindent \textbf{Partition 2: $a = 2; b= 1$}\\
$$\theta_b = \theta_{1} = (\alpha_1-1)\frac{m_3}{1-m_1},
    \;
  R_b^{-1} = R_{1}^{-1}= \frac{1-m_1}{m_3}
  \Rightarrow u = - R_b^{-1}\theta_b = 1-\alpha_1;$$
  
  $$
   R_{ab} = -\frac{m_5}{1-m_1} \Rightarrow 
  R_{ab}u = (\alpha_1-1)\frac{m_5}{1-m_1}.$$
  
  \noindent Thus, $h_a^{'}\left[\theta_a + R_{ab}u\right] = h_2(\alpha_1-1)\displaystyle\frac{m_2-m_1+m_5}{1-m_1} = (\alpha_1-1)h_2$. Since $\alpha_1<1$, taking $h_2 >0 $, we ensure $h_a^{'}\left[\theta_a + R_{ab}u\right] <0$.\\

\noindent \textbf{Partition 3: $a = 1,2; b= \varnothing$}\\
$$h_a^{'}\left[\theta_a + R_{ab}u\right] = h'\theta = \frac{\alpha_1-1}{1-m_1}\left[m_3h_1+(m_2-m_1)h_2\right]$$ Since $\alpha_1<1$, $1-m_1 >0$, fixing $h_1>0$ large enough, we ensure $h_a^{'}\left[\theta_a + R_{ab}u\right] <0$. 
\\

\noindent We have shown, then, that exists a positive vector $h \in{R^2}$ with $h_1$ being large enough as required so that $(R,\theta)$ under the static-priority policy 3, 5 is Chen-$\calS$.
\eProof\\
\end{proof}

\begin{proof} {Proof for lowest-priority 4, 2 (Chen-S).}
\label{proof: Proof of Chen-S for lowest-priority 4, 2 (DHV2)}

\noindent In the balanced push started Lu-Kumar network with class 4 having lowest priority in station 1 and class 2 having lowest priority in station 2, we have $J = 2$, $K = 5$, $\rho = \left[\alpha_1 \ \alpha_1\right]',$\\
$$R = 
 \left[
 \begin{matrix}
   \displaystyle\frac{m_4}{m_2m_4-m_3m_5} & -\displaystyle\frac{m_4(1-m_1)}{m_2m_4-m_3m_5}  \mcr
   -\displaystyle\frac{m_2m_5}{m_2m_4-m_3m_5}  & \displaystyle\frac{m_2m_4}{m_2m_4-m_3m_5} 
  \end{matrix}
  \right] ,
$$
\begin{center}
\begin{equation}\nonumber
\begin{aligned}
\theta = R(\rho-e)
= \frac{\alpha_1-1}{m_2m_4-m_3m_5}
\left[
 \begin{matrix}
   m_1m_4\mcr
   m_2(m_4-m_5)
  \end{matrix}
  \right] .
\end{aligned}
\end{equation}
\end{center}

It is clear that the reflection matrix $R$ is completely-$\calS$ if $m_2m_4>m_3m_5$, since $R$ is $2\times2$, $R^{-1} = CMQ\Delta>0$ and all its diagonal elements are positive when $m_2m_4>m_3m_5$.\\

\noindent \textbf{Partition 1: $a = 1; b= 2$}\\
$$\theta_b = \theta_{2} = \frac{\alpha_1-1}{m_2m_4-m_3m_5}m_2(m_4-m_5),
    \;
  R_b^{-1} = R_{2}^{-1}= \frac{m_2m_4-m_3m_5}{m_2m_4}
  \Rightarrow u = - R_b^{-1}\theta_b = (1-\alpha_1)\frac{m_4-m_5}{m_4};$$
  
  $$
   R_{ab} = -\frac{m_4(1-m_1)}{m_2m_4-m_3m_5} \Rightarrow 
  R_{ab}u = \frac{\alpha_1-1}{m_2m_4-m_3m_5}(m_4-m_5)(1-m_1).$$
  
  \noindent Thus, $h_a^{'}\left[\theta_a + R_{ab}u\right] = \displaystyle\frac{\alpha_1-1}{m_2m_4-m_3m_5}h_1(m_4-m_5)(1-m_1)(m_1m_4+m_4-m_1m_4-m_5+m_1m_5)= \frac{\alpha_1-1}{m_2m_4-m_3m_5}h_1(m_4-m_5+m_1m_5)$. Since $m_4-m_5+m_1m_5 = m_4-m_5(1-m_1)=m_4-m_5(m_3+m_4)=m_4(1-m_5)-m_3m_5=m_4m_2-m_3m_5$, $\left[\theta_a + R_{ab}u\right] = (\alpha_1-1)h_1$. As $\alpha_1<1$, taking $h_1 >0 $, we ensure $h_a^{'}\left[\theta_a + R_{ab}u\right] <0$.\\

\noindent \textbf{Partition 2: $a = 2; b= 1$}\\
$$\theta_b = \theta_{1} = \frac{\alpha_1-1}{m_2m_4-m_3m_5}m_1m_4,
    \;
  R_b^{-1} = R_{1}^{-1}= \frac{m_2m_4-m_3m_5}{m_4}
  \Rightarrow u = - R_b^{-1}\theta_b = (1-\alpha_1)m_1;$$
  
  $$
   R_{ab} = -\frac{m_2m_5}{m_2m_4-m_3m_5} \Rightarrow 
  R_{ab}u = (\alpha_1-1)\frac{m_1m_2m_5}{m_2m_4-m_3m_5}.$$
  
  \noindent Thus, $h_a^{'}\left[\theta_a + R_{ab}u\right] = \displaystyle\frac{\alpha_1-1}{m_2m_4-m_3m_5}h_2\left[m_2(m_4-m_5)+m_1m_2m_5\right]$. Since $m_2(m_4-m_5)+m_1m_2m_5=m_2m_4+m_2m_5(m_1-1) = m_2m_4+m_2m_5(-m_3-m_4)=m_2m_5(1-m_5)-m_2m_3m_5=m_2^2m_4-m_2m_3m_5=m_2(m_2m_4-m_3m_5)$, we have $\left[\theta_a + R_{ab}u\right] = (\alpha_1-1)h_2m_2$. As $\alpha_1<1$, taking $h_2 >0 $, we ensure $h_a^{'}\left[\theta_a + R_{ab}u\right] <0$.\\

\noindent \textbf{Partition 3: $a = 1,2; b= \varnothing$}\\
$$h_a^{'}\left[\theta_a + R_{ab}u\right] = h'\theta = \frac{\alpha_1-1}{m_2m_4-m_3m_5}\left[m_1m_4h_1+m_2(m_4-m_5)h_2\right]$$
Since $\alpha_1<1$, $m_1m_4>0$, and we are not sure the sigh of $m_2(m_4-m_5)$, if $m_2m_4>m_3m_5$ is satisfied, fixing $h_1 > 0$ large enough, and $h_2 >0$, we ensure $h_a^{'}\left[\theta_a + R_{ab}u\right] <0$.
\\

\noindent We have shown, if $m_2m_4>m_3m_5$, then, that exists a positive vector $h \in{R^2}$ with $h_1$ being large enough as required so that $(R,\theta)$ under the static-priority policy 4, 2 is Chen-$\calS$.
\eProof\\
\end{proof}

\begin{proof} {Proof for lowest-priority 4, 5 (Chen-S).}
\label{proof: Proof of Chen-S for lowest-priority 4,5 (DHV2)}

\noindent In the balanced push started Lu-Kumar network with class 4 having lowest priority in station 1 and class 5 having lowest priority in station 2, we have $J = 2$, $K = 5$, $\rho = \left[\alpha_1 \ \alpha_1\right]',$\\
$$R = 
 \left[
 \begin{matrix}
   1 & 0 \mcr
   -\displaystyle\frac{m_5}{m_4} & 1 
  \end{matrix}
  \right] ,
$$
\begin{center}
\begin{equation}\nonumber
\begin{aligned}
\theta = R(\rho-e)
= (\alpha_1-1)
\left[
 \begin{matrix}
   1 \mcr
   \displaystyle\frac{m_4-m_5}{m_4}
  \end{matrix}
  \right] .
\end{aligned}
\end{equation}
\end{center}

It is clear that the reflection matrix $R$ is completely-$\calS$, since $R$ is $2\times2$, $R^{-1} = CMQ\Delta>0$ and all its diagonal elements are positive.\\

\noindent \textbf{Partition 1: $a = 1; b= 2$}\\
$$\theta_b = \theta_{2} = (\alpha_1-1)\frac{m_4-m_5}{m_4},
    \;
  R_b^{-1} = R_{2}^{-1}= 1
  \Rightarrow u = - R_b^{-1}\theta_b = (1-\alpha_1)\frac{m_4-m_5}{m_4};$$
  
  $$
   R_{ab} = 0 \Rightarrow 
  R_{ab}u = 0.$$
  
  \noindent Thus, $h_a^{'}\left[\theta_a + R_{ab}u\right] = h_1(\alpha_1-1) $. Since $\alpha_1<1$, taking $h_1 >0 $, we ensure $h_a^{'}\left[\theta_a + R_{ab}u\right] <0$.\\

\noindent \textbf{Partition 2: $a = 2; b= 1$}\\
$$\theta_b = \theta_{1} = \alpha_1-1,
    \;
  R_b^{-1} = R_{1}^{-1}= 1
  \Rightarrow u = - R_b^{-1}\theta_b = 1-\alpha_1;$$
  
  $$
   R_{ab} = -\frac{m_5}{m_4} \Rightarrow 
  R_{ab}u = (\alpha_1-1)\frac{m_5}{m_4}.$$
  
  \noindent Thus, $h_a^{'}\left[\theta_a + R_{ab}u\right] = (\alpha_1-1)h_2\displaystyle\frac{m_4-m_5+m_5}{m_4} = (\alpha_1-1)h_2 $. Since $\alpha_1<1$, taking $h_2 >0 $, we ensure $h_a^{'}\left[\theta_a + R_{ab}u\right] <0$.\\

\noindent \textbf{Partition 3: $a = 1,2; b= \varnothing$}\\
$$h_a^{'}\left[\theta_a + R_{ab}u\right] = h'\theta = (\alpha_1-1)(h_1+\frac{m_4-m_5}{m_4}h_2)$$
Since $\alpha_1<1$, fixing $h_1>0$ large enough, we ensure $h_a^{'}\left[\theta_a + R_{ab}u\right] <0$.
\\

\noindent We have shown, then, that exists a positive vector $h \in{R^2}$ with $h_1$ being large enough as required so that $(R,\theta)$ under the static-priority policy 4, 5 is Chen-$\calS$.
\eProof
\end{proof}

\vspace{0.2cm}

\noindent Now, we complete the proof that the 6 static-priority cases are all Chen-$\calS$ if $m_2m_4 > m_3m_5$.\\

\noindent The reflection matrices for the 6 balanced push started Lu-Kumar networks following static-priority policies are all invertible. Now we check the sign of determinants of them. We can calculate that $det(R_{1,2}) = m_2>0$, $det(R_{1,5}) = m_1 >0$, $det(R_{3,2}) = \displaystyle\frac{m_3}{1-m_1} >0 $, $det(R_{3,5}) = \displaystyle\frac{m_3}{1-m_1} >0$, $det(R_{4,2}) = \displaystyle\frac{m_2m_4}{m_2m_4-m_3m_5}$, $det(R_{4,5}) = 1 >0$, where the subscript label indicates the classes with the lowest priority in each station. If we want the sign of the determinants to be the same, we need $m_2m_4>m_3m_5$. Thus, we can conclude that the reflection matrix for a balanced push started Lu-Kumar network is invertible for any ratio matrix $\Delta$ if $m_2m_4 > m_3m_5$.

According to Theorem \ref{Thm: full statement for corner->inter}, in a balanced push started Lu-Kumar network with $\alpha_1 <1$, the reflection matrix is invertible and $(R,\theta)$ is Chen-$\calS$ for any ratio matrix $\Delta$ if $m_2m_4 > m_3m_5$.\\

\subsection{The SSC Inequalities for Balanced Push Started Lu-Kumar Networks}
\label{Appendix: The SSC Inequalities for Balanced DHV2 Networks}

\noindent Similarly as the analysis in Appendix \ref{Appendix: The SSC Inequalities for Balanced DHV Networks}, we get the SSC conditions for balanced push started Lu-Kumar networks in this appendix. Note that in balanced push started Lu-Kumar networks, we have $m_1+m_3+m_4 = 1$ and $m_2+m_5 = 1$. 
We check the 6 static priority policies of balanced push started Lu-Kumar networks for (\ref{equ: appendix_general_SSC_DHV}) in the following analysis. As mentioned in \S \ref{Balanced DHV2 Network Example (SSC)}, we use a fixed priority order as $1>3>4$ (in station 1) within the set $\{\ell \in \mathcal{C}(j):Z_l(t)-\delta_lW_j(t) >0 \}$ as the ``tie-breaking" rule.

Same as Appendix \ref{Appendix: The SSC Inequalities for Balanced DHV Networks}, given a policy $\pi$, let $\dot{Y}^\pi(t) = \sum_{k \in \mathcal{H}_\pi:z_k>0} h_k\dot{\wbar{Z}_k^\pi}(t)$. Then, we want $\dot{Y}^\pi(t) <0$ with a common vector $h \in \Bbb{R}^K_{++}$ for any regular time $t$ with $\lVert \wbar{Z}^\pi_{\mathcal{H}_\pi}(t) \rVert >0$, i.e., at least one of the high-priority queues is not empty, under all the 6 static priority policies. We omit $\pi$ when it is fixed. In balanced push started Lu-Kumar networks, $m_k<1$ and thus $\mu_k>1$ for all $k \in \mathcal{K}$. Denote $I_k(t)$ as the input rate of class $k$ at time $t$; $O_k(t)$ as the output rate of class $k$ at time $t$. Then, $\dot{\wbar{Z}}_k(t) = I_k(t) - O_k(t) \text{ for } k\in \mathcal{K}; I_k(t) = O_{k-1}(t) \text{ for } k \in \mathcal{K}/\{1\};  I_1(t) = \alpha_1$. In this section, we use a set of lowest-priority classes to represent static priority policies, for example,
(1, 2) stands for a static priority policy that class 1 and 2 have the lowest priorities in each station.

\vspace{0.2cm}

\noindent \textbf{Case 1:} (1, 2)
    
\vspace{0.2cm}

\begin{enumerate}
    \item $\wbar{Z}_3(t)>0, \wbar{Z}_5(t)>0.$

    It is clear that server 1 serves queue 3, server 2 serves queue 5, then $O_3(t) = \mu_3$, $O_5(t) = \mu_5$, $O_1(t) = O_2(t) = O_4(t) = 0$. So $\dot{\wbar{Z}}_3(t) = 0-\mu_3 = -\mu_3$,
    $\dot{\wbar{Z}}_4(t) = \mu_3-0 = \mu_3 $ and $\dot{\wbar{Z}}_5(t) = 0-\mu_5 = -\mu_5$. Thus, we need $\dot{Y}(t) =  -h_3\mu_3 +h_4\mu_3 - h_5\mu_5<0$.

    \vspace{0.2cm}
    
    \item $\wbar{Z}_3(t)=0, \wbar{Z}_4(t)>0, \wbar{Z}_5(t)>0.$ 

    It is clear that $O_1(t) = 0$, server 2 serves queue 5, then $O_5(t) = \mu_5$, $O_2(t) = 0$. Since $I_3(t)=O_2(t) = 0$ and $\wbar{Z}_3(t)=0$, we have $O_3(t) = 0$ and serve 1 serves queue 4, which means $O_4(t) = \mu_4$. So $\dot{\wbar{Z}}_3(t) = 0-0 = 0$,
    $\dot{\wbar{Z}}_4(t) = 0-\mu_4  = -\mu_4 $ and $\dot{\wbar{Z}}_5(t) = \mu_4-\mu_5$. Thus, we need $\dot{Y}(t) = -h_4\mu_4 +h_5(\mu_4-\mu_5)<0$.
    
    \vspace{0.2cm}
    
    \item $\wbar{Z}_3(t)>0, \wbar{Z}_5(t)=0.$

    It is clear that server 1 serves queue 3, then $O_3(t) = \mu_3$, $O_1(t) = O_4(t) = 0$. Since $O_4(t) = 0$ and $\wbar{Z}_5(t)=0$, we have $O_5(t) =0$. 
    So $\dot{\wbar{Z}}_4(t) = \mu_3- 0  = \mu_3$ and $\dot{\wbar{Z}}_5(t) = 0-0 = 0$.

    If $\wbar{Z}_2(t) =0$, since $O_1(t) = 0$, we have $O_2(t) = 0$, then $\dot{\wbar{Z}}_3(t) = 0-\mu_3  = -\mu_3$. Thus, we need $\dot{Y}(t) =  -h_3\mu_3 +h_4\mu_3<0$. 
    
    If $\wbar{Z}_2(t)> 0$, since $O_5(t) = 0$, server 2 serves queue 2, we have $O_2(t) = \mu_2$, then $\dot{\wbar{Z}}_3(t) = \mu_2-\mu_3$. Thus, we need $\dot{Y}(t) = h_3(\mu_2-\mu_3) +h_4\mu_3<0$. 
    
    Note that $h_3(\mu_2-\mu_3) +h_4\mu_3 \geq -h_3\mu_3 +h_4\mu_3$, then it is enough to only consider $\dot{Y}(t) = h_3(\mu_2-\mu_3) +h_4\mu_3 <0$ in this situation. 
    
    \vspace{0.2cm}

    \item $\wbar{Z}_3(t)=0, \wbar{Z}_4(t)=0, \wbar{Z}_5(t)>0.$ 

    It is clear that server 2 serves queue 5, then $O_5(t) = \mu_5$, $O_2(t) = 0$. Since $I_3(t)=O_2(t) = 0$ and $\wbar{Z}_3(t)=0$, we have $O_3(t) =0$. Since $I_4(t)=O_3(t) = 0$ and $\wbar{Z}_4(t)=0$, we have $O_4(t) =0$. So $\dot{\wbar{Z}}_3(t) = 0-0 = 0$,
    $\dot{\wbar{Z}}_4(t) =0-0 = 0 $ and $\dot{\wbar{Z}}_5(t) = 0-\mu_5 = -\mu_5$. It obviously holds that $\dot{Y}(t) = - h_5\mu_5<0$ as required. 

    \vspace{0.2cm}

    \item $\wbar{Z}_3(t)=0, \wbar{Z}_4(t)>0, \wbar{Z}_5(t)=0.$ 

     It is clear that server 1 serves queue 4, then $O_4(t) = \mu_4$, $O_1(t)=O_3(t) = 0$. So $\dot{\wbar{Z}}_4(t) = 0-\mu_4 =-\mu_4$. If $\mu_4 > \mu_5$,  $\dot{\wbar{Z}}_5(t) = \mu_4-O_5(t) > 0$, then $t$ is not differentiable; if $\mu_4 \leq \mu_5$, $O_5(t) = \mu_5$, $\dot{\wbar{Z}}_5(t) = \mu_4-\mu_4 = 0$, $\dot{\wbar{T}}_5(t) = \frac{\mu_4}{\mu_5}$. If $\wbar{Z}_2(t)>0$, we have $O_2(t) = \mu_2(1-\frac{\mu_4}{\mu_5})$, and $\dot{\wbar{Z}}_3(t) =  \mu_2(1-\frac{\mu_4}{\mu_5})-0 > 0$, then $t$ is not differentiable; if $\wbar{Z}_2(t)=0$, since $I_2(t) = O_1(t) = 0$, we have $O_2(t) =0$, so $\dot{\wbar{Z}}_3(t) = 0-0 = 0$. It obviously holds that $\dot{Y}(t) = - h_4\mu_4<0$ as required.

    \vspace{0.2cm}
    \end{enumerate}  

To sum up, except for the positive condition for common vector $h$, we need the following conditions from case 1:
\begin{align*}
  &-h_3\mu_3 +h_4\mu_3 - h_5\mu_5<0,\\
  &-h_4\mu_4 +h_5(\mu_4-\mu_5)<0,\\
  &h_3(\mu_2-\mu_3) +h_4\mu_3 <0  .
\end{align*}
 Note that $-h_3\mu_3 +h_4\mu_3 - h_5\mu_5 \leq h_3(\mu_2-\mu_3) +h_4\mu_3$, thus, it is enough to consider only the following two conditions from case 1:
\begin{align*}
  &-h_4\mu_4 +h_5(\mu_4-\mu_5)<0,\\
  &h_3(\mu_2-\mu_3) +h_4\mu_3 <0  .
\end{align*}

\noindent \textbf{Case 2:} (1, 5)
    
\vspace{0.2cm}

\begin{enumerate}
    \item $\wbar{Z}_2(t)>0, \wbar{Z}_3(t)>0.$

     It is clear that server 1 serves queue 3, server 2 serves queue 2, then $O_3(t) = \mu_3$, $O_2(t) = \mu_2$, $O_1(t) = O_4(t) = O_5(t) = 0$. So $\dot{\wbar{Z}}_2(t) = 0-\mu_2 = -\mu_2$,
    $\dot{\wbar{Z}}_3(t) = \mu_2-\mu_3$ and $\dot{\wbar{Z}}_4(t) = \mu_3-0 = \mu_3$. Thus, we need $\dot{Y}(t) = -h_2\mu_2 +h_3(\mu_2-\mu_3) + h_4\mu_3<0$.
    
    \vspace{0.2cm}
    
    \item $\wbar{Z}_2(t)=0, \wbar{Z}_3(t)>0$ 

    It is clear that server 1 serves queue 3, then $O_3(t) = \mu_3$, $O_1(t) = O_4(t) = 0$. Since $I_2(t)=O_1(t)  = 0$ and $\wbar{Z}_2(t)=0$, we have $O_2(t) =0$. 
    So $\dot{\wbar{Z}}_2(t) = 0-0 = 0$,
    $\dot{\wbar{Z}}_3(t) = 0-\mu_3 = -\mu_3$ and $\dot{\wbar{Z}}_4(t) = \mu_3-0 = \mu_3$. Thus, we need $\dot{Y}(t) = - h_3\mu_3 + h_4\mu_3<0$.

    \vspace{0.2cm}
    
    \item $\wbar{Z}_2(t)>0, \wbar{Z}_3(t)=0, \wbar{Z}_4(t)>0.$

    It is clear that server 1 serves queue 4, and server 2 serves queue 2, then $O_4(t) = \mu_4$, $O_2(t) =\mu_2$,  $O_1(t) = O_3(t)  = O_5(t)= 0$. Thus,  $\dot{\wbar{Z}}_3(t) = \mu_2- 0  = \mu_2 >0$, which means that $t$ is not differentiable.
    
    \vspace{0.2cm}

    \item $\wbar{Z}_2(t)=0, \wbar{Z}_3(t)=0, \wbar{Z}_4(t)>0.$ 

     It is clear that server 1 serves queue 4, then $O_4(t) = \mu_4$, $O_1(t) = O_3(t)= 0$. Since $I_2(t)=O_1(t) = 0$ and $\wbar{Z}_2(t)=0$, we have $O_2(t) =0$. Since $I_3(t)=O_2(t) = 0$ and $\wbar{Z}_3(t)=0$, we have $O_3(t) =0$. So $\dot{\wbar{Z}}_2(t) = 0-0 = 0$,
    $\dot{\wbar{Z}}_3(t) =0-0 = 0 $ and $\dot{\wbar{Z}}_4(t) = 0-\mu_4 = -\mu_4$. It obviously holds that$\dot{Y}(t) = - h_4\mu_4<0$ as required.

    \vspace{0.2cm}

    \item $\wbar{Z}_2(t)>0, \wbar{Z}_3(t)=0, \wbar{Z}_4(t)=0.$ 

    It is clear that server 2 serves queue 2, then $O_2(t) = \mu_2$, $O_5(t)= 0$. If $\mu_2 > \mu_3$, we have  $O_3(t)= \mu_3$, and then $\dot{\wbar{Z}}_3(t) = \mu_2-\mu_3 > 0$, then $t$ is not differentiable; if $\mu_2 \leq \mu_3$, $O_3(t) = \mu_2$, $\dot{\wbar{Z}}_3(t) = \mu_2-\mu_2 = 0$, $\dot{\wbar{T}}_3(t) = \frac{\mu_2}{\mu_3}$. If $\mu_2 > \mu_4(1-\frac{\mu_2}{\mu_3})$, we have $O_4(t) = \mu_4(1-\frac{\mu_2}{\mu_3})$, and $\dot{\wbar{Z}}_4(t) =  \mu_2 - \mu_4(1-\frac{\mu_2}{\mu_3}) > 0$, then $t$ is not differentiable; if $\mu_2 \leq \mu_4(1-\frac{\mu_2}{\mu_3})$, then $O_4(t) = \mu_2$, $\dot{\wbar{Z}}_4(t) = \mu_2-\mu_2 = 0$ and $\dot{\wbar{T}}_1(t) = 1-\frac{\mu_2}{\mu_3} - \frac{\mu_2}{\mu_4}$. So, $O_1(t) = \alpha_1$ or $O_1(t)= \mu_1(1-\frac{\mu_2}{\mu_3} - \frac{\mu_2}{\mu_4})$. Thus, $\dot{\wbar{Z}}_2(t) =\alpha_1 - \mu_2 <0$ or $\dot{\wbar{Z}}_2(t) =\mu_1(1-\frac{\mu_2}{\mu_3} - \frac{\mu_2}{\mu_4}) - \mu_2$. We have $\mu_1(1-\frac{\mu_2}{\mu_3} - \frac{\mu_2}{\mu_4}) - \mu_2 <0$ since $\mu_1(1-\frac{\mu_2}{\mu_3} - \frac{\mu_2}{\mu_4}) - \mu_2 <0 \Leftrightarrow \frac{\mu_2}{\mu_1} + \frac{\mu_2}{\mu_3} + \frac{\mu_2}{\mu_4}>1 \Leftrightarrow \frac{m_1+m_3+m_4}{m_2}>1 \Leftrightarrow \frac{1}{m_2}>1$, which is true obviously. So, we have $\dot{\wbar{Z}}_2(t) <0$,  $\dot{\wbar{Z}}_3(t) = 0$ and $\dot{\wbar{Z}}_4(t) = 0$. Thus, $\dot{Y}(t) <0$ is satisfied obviously, as required.

    \vspace{0.2cm}
    \end{enumerate}  

To sum up, except for the positive condition for common vector $h$, we need the following conditions from case 2:
\begin{align*}
  &-h_2\mu_2 +h_3(\mu_2-\mu_3) + h_4\mu_3<0\\
  &- h_3\mu_3 + h_4\mu_3<0.
\end{align*}

\noindent \textbf{Case 3:} (3, 2)
    
\vspace{0.2cm}

\begin{enumerate}
    \item $\wbar{Z}_1(t)>0, \wbar{Z}_5(t)>0.$

    It is clear that $I_1(t) = \alpha_1$, server 1 serves queue 1, server 2 serves queue 5, then $O_1(t) = \mu_1$, $O_5(t) = \mu_5$, $O_2(t) = O_3(t) = O_4(t) = 0$. So $\dot{\wbar{Z}}_1(t) = \alpha_1-\mu_1$,
    $\dot{\wbar{Z}}_4(t) = 0-0 = 0$ and $\dot{\wbar{Z}}_5(t) = 0-\mu_5 = -\mu_5$. Thus, we need $\dot{Y}(t) =  h_1(\alpha_1-\mu_1) - h_5\mu_5<0$, which is satisfied obviously since $\alpha_1-\mu_1 <0$.
    
    \vspace{0.2cm}
    
    \item $\wbar{Z}_1(t)=0, \wbar{Z}_4(t)>0, \wbar{Z}_5(t)>0$ 

    It is clear that server 1 serves queue 4, then $O_4(t) = \mu_4$, $O_1(t)  = O_3(t)= 0$. Thus,  $\dot{\wbar{Z}}_1(t) = \alpha_1- 0 >0$, which means that $t$ is not differentiable.

    \vspace{0.2cm}
    
    \item $\wbar{Z}_1(t)>0, \wbar{Z}_5(t)=0.$

    It is clear that server 1 serves queue 1, then $O_1(t) = \mu_1$, $O_3(t) = O_4(t) = 0$, so  $\dot{\wbar{Z}}_1(t) = \alpha_1- \mu_1$, $\dot{\wbar{Z}}_4(t) = 0- 0 = 0$. Since $I_5(t)=O_4(t)  = 0$ and $\wbar{Z}_5(t)=0$, we have $O_5(t) =0$. 
    So $\dot{\wbar{Z}}_5(t) = 0-0 = 0$. It obviously holds that $\dot{Y}(t)  = h_1(\alpha_1-\mu_1) <0$ as required since $\alpha_1- \mu_1<0$.
    
    \vspace{0.2cm}

    \item $\wbar{Z}_1(t)=0, \wbar{Z}_4(t)=0, \wbar{Z}_5(t)>0.$ 

    It is clear that server 2 serves queue 5, then $O_5(t) = \mu_5$, $O_2(t) = 0$. To make $t$ differentiable, we have $\dot{\wbar{Z}}_1(t) = \alpha_1-O_1(t) = 0$, then  $O_1(t) = \alpha_1$ and $\dot{\wbar{T}}_1(t) = \frac{\alpha_1}{\mu_1}$.

    If $\wbar{Z}_3(t)=0$, since $I_3(t) = O_2(t) = 0$, we have $O_3(t) = 0$. Since $I_4(t) = O_3(t) = 0$, $\wbar{Z}_4(t)=0$, we have $O_4(t) = 0$. Then $\dot{\wbar{Z}}_4(t) = 0-0 = 0$, $\dot{\wbar{Z}}_5(t) = 0-\mu_5 = -\mu_5$. It obviously holds that $\dot{Y}(t) = -h_5\mu_5 <0$ as required.

    If $\wbar{Z}_3(t)>0$, to keep $t$ differentiable, we have $\dot{\wbar{Z}}_4(t) = O_3(t)-O_4(t) = 0$, which means $\mu_3(1-\frac{\alpha_1}{\mu_1}-\dot{\wbar{T}}_4(t)) = \mu_4\dot{\wbar{T}}_4(t) \Leftrightarrow (\mu_3+\mu_4)\dot{\wbar{T}}_4(t) = \mu_3(1-\frac{\alpha_1}{\mu_1}) \Leftrightarrow \dot{\wbar{T}}_4(t) = \frac{\mu_3(1-\frac{\alpha_1}{\mu_1})}{\mu_3+\mu_4}$. Thus, if we want $\dot{Y}(t) = h_5\dot{\wbar{Z}}_5(t) <0$, we need $\dot{\wbar{Z}}_5(t) = O_4(t)-\mu_5 = \mu_4\dot{\wbar{T}}_4(t)-\mu_5 = \frac{\mu_3\mu_4(1-\frac{\alpha_1}{\mu_1})}{\mu_3+\mu_4}-\mu_5 <0$. $\frac{\mu_3\mu_4(1-\frac{\alpha_1}{\mu_1})}{\mu_3+\mu_4}-\mu_5 <0 \Leftrightarrow 1-\frac{\alpha_1}{\mu_1} < \frac{\mu_5(\mu_3+\mu_4)}{\mu_3\mu_4} \Leftrightarrow 1-m_1\alpha_1 < \frac{m_3+m_4}{m_5} \Leftrightarrow m_1\alpha_1 > \frac{m_5-m_3-m_4}{m_5} \Leftrightarrow \alpha_1 > \frac{m_5-(1-m_1)}{m_1m_5} \Leftrightarrow \alpha_1 > \frac{m_1-m_2}{m_1m_5}$. Thus, $\dot{Y}(t) = h_5\dot{\wbar{Z}}_5(t) <0$ is satisfied if $\alpha_1 > \frac{m_1-m_2}{m_1m_5}$.

    Note that $\frac{m_1-m_2}{m_1m_5}<1$, thus, the condition $\alpha_1 > \frac{m_1-m_2}{m_1m_5}$ does not conflict with the nominal condition($\alpha_1<1$). It is clear that $\frac{m_1-m_2}{m_1m_5}<1$ since $\frac{m_1-m_2}{m_1m_5}<1 \Leftrightarrow m_1-m_2 < m_1m_5 \Leftrightarrow  m_2 > m_1(1-m_5) \Leftrightarrow m_2 > m_1m_2 
    \Leftrightarrow m_2(1-m_1) >0 \Leftrightarrow m_2(m_3+m_4) >0$, which is true obviously. 
    

    \vspace{0.2cm}

    \item $\wbar{Z}_1(t)=0, \wbar{Z}_4(t)>0, \wbar{Z}_5(t)=0.$ 

    It is clear that server 1 serves queue 4, then $O_4(t) = \mu_4$, $O_1(t) = O_3(t)= 0$. Thus,  $\dot{\wbar{Z}}_1(t) = \alpha_1- 0 >0$, which means that $t$ is not differentiable.

    \vspace{0.2cm}
    \end{enumerate}     

To sum up, except for the positive condition for common vector $h$, we need the following conditions from case 3:
\begin{align*}
\alpha_1 > \frac{m_1-m_2}{m_1m_5}.
\end{align*}

\vspace{0.2cm}

\noindent \textbf{Case 4:} (3, 5)
    
\vspace{0.2cm}

\begin{enumerate}
    \item $\wbar{Z}_1(t)>0, \wbar{Z}_2(t)>0.$

    It is clear that server 1 serves queue 1, server 2 serves queue 2, then $O_1(t) = \mu_1$, $O_2(t) = \mu_2$, $O_3(t) = O_4(t)= O_5(t) = 0$. So $\dot{\wbar{Z}}_1(t) = \alpha_1-\mu_1$,
    $\dot{\wbar{Z}}_2(t) = \mu_1-\mu_2 $ and $\dot{\wbar{Z}}_4(t) = 0-0 = 0$. Thus, we need $\dot{Y}(t) = h_1(\alpha_1-\mu_1)+ h_2(\mu_1-\mu_2) <0$.

    \vspace{0.2cm}
    
    \item $\wbar{Z}_1(t)=0, \wbar{Z}_2(t)>0, \wbar{Z}_4(t)>0$ 

    It is clear that server 1 serves queue 4, then $O_4(t) = \mu_4$, $O_1(t) = O_3(t)= 0$. Thus,  $\dot{\wbar{Z}}_1(t) = \alpha_1- 0 >0$, which means that $t$ is not differentiable.

    \vspace{0.2cm}
    
    \item $\wbar{Z}_1(t)>0, \wbar{Z}_2(t)=0.$

    It is clear that server 1 serves queue 1, then $O_1(t) = \mu_1$, $O_4(t) = O_3(t)= 0$. Thus,  $\dot{\wbar{Z}}_4(t) = 0- 0 =0$, $\dot{\wbar{Z}}_1(t) = \alpha_1- \mu_1 =0$. To make $t$ differentiable, we have $\dot{\wbar{Z}}_2(t) = \mu_1- O_2(t) =0$. It obviously holds that $\dot{Y}(t) = h_1(\alpha_1-\mu_1)<0$ as required since $\alpha_1 <\mu_1$.
    
    \vspace{0.2cm}

    \item $\wbar{Z}_1(t)=0, \wbar{Z}_2(t)=0, \wbar{Z}_4(t)>0.$ 

    It is clear that server 1 serves queue 4, then $O_4(t) = \mu_4$, $O_1(t) = O_3(t)= 0$. Thus,  $\dot{\wbar{Z}}_1(t) = \alpha_1- 0 >0$, which means that $t$ is not differentiable.

    \vspace{0.2cm}

    \item $\wbar{Z}_1(t)=0, \wbar{Z}_2(t)>0, \wbar{Z}_4(t)=0.$ 

    It is clear that server 2 serves queue 2, then $O_2(t) = \mu_2$, $O_5(t)= 0$. To make $t$ differentiable, we have  $\dot{\wbar{Z}}_1(t) = \alpha_1- O_1(t) =0$,  $\dot{\wbar{Z}}_4(t) = O_3(t)- O_4(t) =0$, then $O_1(t) = \alpha_1$ $\dot{\wbar{T}}_1(t) = \frac{\alpha_1}{\mu_1}$, $\dot{\wbar{Z}}_2(t) = \alpha_1- \mu_2<0$, $O_3(t) = O_4(t)$. To have $O_3(t) = O_4(t)$, we need $\mu_3(1-\frac{\alpha_1}{\mu_1}-\dot{\wbar{T}}_4(t)) = \mu_4\dot{\wbar{T}}_4(t) \Leftrightarrow \dot{\wbar{T}}_4(t) = \frac{1-\frac{\alpha_1}{\mu_1}}{\mu_3+\mu_4}$, which is possible since $\dot{\wbar{T}}_4(t) >0$. Thus, we have $\dot{Y}(t)<0$ for sure.

    \vspace{0.2cm}
    \end{enumerate} 

To sum up, except for the positive condition for common vector $h$, we need the following condition from case 4:
\begin{align*}
  h_1(\alpha_1-\mu_1)+ h_2(\mu_1-\mu_2) <0.
\end{align*}

\noindent \textbf{Case 5:} (4, 2)
    
\vspace{0.2cm}

\begin{enumerate}
    \item $\wbar{Z}_1(t)>0, \wbar{Z}_5(t)>0.$

    It is clear that server 1 serves queue 1, server 2 serves queue 5, then $O_1(t) = \mu_1$, $O_5(t) = \mu_5$, $O_2(t) = O_3(t)= O_4(t) = 0$. So $\dot{\wbar{Z}}_1(t) = \alpha_1-\mu_1$,
    $\dot{\wbar{Z}}_3(t) = 0-0=0 $ and $\dot{\wbar{Z}}_5(t) = 0-\mu_5 = -\mu_5$. It obviously holds that $\dot{Y}(t) = h_1(\alpha_1-\mu_1)-h_5\mu_5 <0$ as required since $\alpha_1 < \mu_1$.
    
    \vspace{0.2cm}
    
    \item $\wbar{Z}_1(t)=0, \wbar{Z}_3(t)>0, \wbar{Z}_5(t)>0$ 

     It is clear that server 2 serves queue 5, then $O_5(t) = \mu_5$, $O_2(t) = 0$. Thus,  $\dot{\wbar{Z}}_1(t) = \alpha_1- 0 >0$, which means that $t$ is not differentiable.

    \vspace{0.2cm}
    
    \item $\wbar{Z}_1(t)>0, \wbar{Z}_5(t)=0.$

    It is clear that server 1 serves queue 1, then $O_1(t) = \mu_1$, $O_3(t)= O_4(t) = 0$. So $\dot{\wbar{Z}}_1(t) = \alpha_1-\mu_1$. Since $I_5(t) = O_4(t) = 0$, $\wbar{Z}_5(t)=0$, we have $O_5(t) = 0$ and then $\dot{\wbar{Z}}_5(t) = 0-0 = 0$. 

    If $\wbar{Z}_2(t)>0$, then server 2 serves queue 2, we have $O_2(t) = \mu_2$ and $\dot{\wbar{Z}}_3(t) = \mu_2-0 = \mu_2$. Thus, we need $\dot{Y}(t) = h_1(\alpha_1-\mu_1)+h_3\mu_2 <0$.

    If $\wbar{Z}_2(t)=0$, when $\mu_1 > \mu_2$, $\dot{\wbar{Z}}_2(t) = \mu_1 - \mu_2 >0$ and $t$ is not differentiable; when $\mu_1 \leq \mu_2$, $O_2(t) = \mu_1$, $\dot{\wbar{Z}}_3(t) = \mu_1-0 = \mu_1$, and we need $\dot{Y}(t) = h_1(\alpha_1-\mu_1)+h_3\mu_1<0$.

    It is clear that when $\mu_1 \leq \mu_2$, $h_1(\alpha_1-\mu_1)+h_3\mu_2 \geq h_1(\alpha_1-\mu_1)+h_3\mu_1$. So it is sufficient to consider only $h_1(\alpha_1-\mu_1)+h_3\mu_2 <0$ in this situation.

    \vspace{0.2cm}

    \item $\wbar{Z}_1(t)=0, \wbar{Z}_3(t)=0, \wbar{Z}_5(t)>0.$ 

    It is clear that server 2 serves queue 5, then $O_5(t) = \mu_5$, $O_2(t)  = 0$. Since $I_3(t) = O_2(t) = 0$, $ \wbar{Z}_3(t)=0$, we have $O_3(t) = 0$ and $\dot{\wbar{Z}}_3(t) = 0-0 =0$. To keep $t$ differentiable, we have $\dot{\wbar{Z}}_1(t) = \alpha_1-O_1(t) = 0$, then $O_1(t) = \alpha_1$ and $\dot{\wbar{T}}_1 = \frac{\alpha_1}{\mu_1}$.

    If $\wbar{Z}_4(t)=0$, since $I_4(t) = O_3(t) = 0$, we have $O_4(t) = 0$. Then $\dot{\wbar{Z}}_5(t) = 0-\mu_5 = -\mu_5$, It obviously holds that $\dot{Y}(t) = -h_5\mu_5 <0$ as required.

    If $\wbar{Z}_4(t)>0$, then $O_4(t) = \mu_4(1-\frac{\alpha_1}{\mu_1})$, so $\dot{\wbar{Z}}_5(t) = \mu_4(1-\frac{\alpha_1}{\mu_1}) - \mu_5$. Thus, if we want $\dot{Y}(t) = h_5\dot{\wbar{Z}}_5(t) <0$, we need $\mu_4(1-\frac{\alpha_1}{\mu_1}) - \mu_5 <0 \Leftrightarrow 1-m_1\alpha_1 < \frac{m_4}{m_5} \Leftrightarrow \alpha_1 > \frac{m_5-m_4}{m_1m_5}$. Thus, $\dot{Y}(t) = h_5\dot{\wbar{Z}}_5(t) <0$ is satisfied if $\alpha_1 > \frac{m_5-m_4}{m_1m_5}$.

    Note that to make the condition $\alpha_1 > \frac{m_5-m_4}{m_1m_5}$ does not conflict with the nominal condition($\alpha_1<1$), we need $ \frac{m_5-m_4}{m_1m_5} <1 \Leftrightarrow m_5-m_4 <m_1m_5
    \Leftrightarrow m_4 > m_5(1-m_1)
    \Leftrightarrow m_4 > m_5(m_3+m_4)
    \Leftrightarrow m_4(1-m_5) > m_3m_5
    \Leftrightarrow m_2m_4 > m_3m_5$.


    \vspace{0.2cm}

    \item $\wbar{Z}_1(t)=0, \wbar{Z}_3(t)>0, \wbar{Z}_5(t)=0.$

     It is clear that server 1 serves queue 3, then $O_3(t) = \mu_3$, $O_1(t) = O_4(t) = 0$. Thus,  $\dot{\wbar{Z}}_1(t) = \alpha_1- 0 >0$, which means that $t$ is not differentiable.

    \vspace{0.2cm}
    \end{enumerate}

To sum up, except for the positive condition for common vector $h$, we need the following condition from case 5:
\begin{align*}
    & h_1(\alpha_1-\mu_1)+h_3\mu_2 <0 \\
    & m_2m_4>m_4m_5,\\
    & \alpha_1 > \frac{m_5-m_4}{m_1m_5}.
\end{align*}

\noindent \textbf{Case 6:} (4, 5)
    
\vspace{0.2cm}

\begin{enumerate}
    \item $\wbar{Z}_1(t)>0, \wbar{Z}_2(t)>0.$

    It is clear that server 1 serves queue 1, server 2 serves queue 2, then $O_1(t) = \mu_1$, $O_2(t) = \mu_2$, $O_3(t) = O_4(t)= O_5(t) = 0$. So $\dot{\wbar{Z}}_1(t) = \alpha_1-\mu_1$,
    $\dot{\wbar{Z}}_2(t) = \mu_1-\mu_2$ and $\dot{\wbar{Z}}_3(t) = \mu_2-0 = \mu_2$. Thus, we need $\dot{Y}(t) = h_1(\alpha_1-\mu_1)+h_2(\mu_1-\mu_2)+h_3\mu_2 <0$.
    
    \vspace{0.2cm}

    \item $\wbar{Z}_1(t)>0, \wbar{Z}_2(t)=0.$

    It is clear that server 1 serves queue 1, then $O_1(t) = \mu_1$, $O_3(t)= O_4(t) = 0$. So $\dot{\wbar{Z}}_1(t) = \alpha_1-\mu_1$. Since $O_1(t) = \mu_1$, $\wbar{Z}_2(t)=0$,  when $\mu_1 > \mu_2$,  we have $O_2(t) = \mu_2$, and then $\dot{\wbar{Z}}_2(t) = \mu_1 - \mu_2 >0$, so $t$ is not differentiable; when $\mu_1 \leq \mu_2$, we have $O_2(t) = \mu_1$, then $\dot{\wbar{Z}}_2(t) = \mu_1 - \mu_1 =0$, $\dot{\wbar{Z}}_3(t) = \mu_1-0 = \mu_1$, and we need $\dot{Y}(t) = h_1(\alpha_1-\mu_1)+h_3\mu_1<0$ (if $\mu_1 \leq \mu_2$).

    \vspace{0.2cm}

    \item $\wbar{Z}_1(t)=0,  \wbar{Z}_3(t)>0$ 

    It is clear that server 1 serves queue 3, then $O_3(t) = \mu_3$, $O_1(t) = O_4(t) = 0$. Thus,  $\dot{\wbar{Z}}_1(t) = \alpha_1- 0 >0$, which means that $t$ is not differentiable.

    \vspace{0.2cm}

    \item $\wbar{Z}_1(t)=0, \wbar{Z}_2(t)>0, \wbar{Z}_3(t)=0.$ 

    It is clear that server 2 serves queue 2, then $O_2(t) = \mu_2$, $O_5(t) = 0$. To make $t$ differentiable, we have $\dot{\wbar{Z}}_1(t) = \alpha_1 - O_1(t) =0$, then $O_1(t) = \alpha_1$, $\dot{\wbar{T}}_1(t) = \frac{\alpha_1}{\mu_1}$ and $\dot{\wbar{Z}}_2(t) = \alpha_1- \mu_2$. Thus, the maximum of $O_3(t) = \mu_3(1-\frac{\alpha_1}{\mu_1})$. When $\mu_2 > \mu_3(1-\frac{\alpha_1}{\mu_1})$, we have $O_3(t) = \mu_3(1-\frac{\alpha_1}{\mu_1})$, and then $\dot{\wbar{Z}}_3(t) = \mu_2 - \mu_3(1-\frac{\alpha_1}{\mu_1}) >0$, so $t$ is not differentiable; when $\mu_2 \leq \mu_3(1-\frac{\alpha_1}{\mu_1})$, we have $O_3(t) = \mu_2$, then $\dot{\wbar{Z}}_3(t) = \mu_2 - \mu_2=0$, and it obviously holds that $\dot{Y}(t) = h_2(\alpha_1-\mu_2)<0$ as required since $\alpha_1<\mu_2$.

    \vspace{0.2cm}
    \end{enumerate}    

To sum up, except for the positive condition for common vector $h$, we need the following condition from case 6:
\begin{align*}
    & h_1(\alpha_1-\mu_1) + h_2(\mu_1-\mu_2)+h_3\mu_2<0 \\
    & h_1(\alpha_1-\mu_1)+h_3\mu_1 <0 \;\;\text{if}\;\; \mu_1 \leq \mu_2.
\end{align*}




\section{Detailed Proofs for Balanced Lu-Kumar Networks}

\setcounter{equation}{0}
\renewcommand{\theequation}{\thesection.\arabic{equation}}

\subsection{Chen-$\calS$ in Balanced Lu-Kumar Networks}
\label{Appendix: Chen-S in Balanced LuKumar Networks}

In this appendix, we prove that the reflection matrices for the four static-priority polices of a balanced Lu-Kumar network ($m_1+m_4=1$; $m_2+m_3=1$) with $\alpha_1<1$ are Chen-$\calS$, invertible and have same-sign determinants if $m_2+m_4 < 1$. Theorem \ref{Thm: full statement for corner->inter} then guarantees that the refection matrix is invertible and Chen-$\calS$ for any ratio matrix $\Delta$ of a balanced Lu-Kumar network if $m_2+m_4 < 1$.

Given a static priority ratio matrix $\Delta$, to prove that the corresponding $(R_{\Delta}, \theta_{\Delta})$ is Chen-$\calS$, we must show that (i) $R$ is completely-$\calS$, and (ii) that there exists a positive vector $h \in \mathbb{R}^2$ such that (\ref{equ:ChenS_def}) satisfies for each of the 3 partitions of $\mathcal{J}$. 

\vspace{0.2cm}

\begin{proof} {Proof for priority 1, 2 (Chen-S).}
\label{proof: Proof of Chen-S for priority 1,2 (LuKumar)}

\noindent In the balanced Lu-Kumar network with class 1 having the priority in station 1 and class 2 having the priority in station 2, we have $J = 2$, $K = 4$, $\rho = \left[\alpha_1 \ \alpha_1\right]',$\\
$$R = 
 \left[
 \begin{matrix}
   1 & -\displaystyle\frac{m_4}{m_3}\mcr
   0 & 1
  \end{matrix}
  \right] ,
$$
\begin{center}
\begin{equation}\nonumber
\begin{aligned}
\theta = R(\rho-e)
= (\alpha_1-1)
\left[
 \begin{matrix}
   1-\displaystyle\frac{m_4}{m_3} \mcr
   1
  \end{matrix}
  \right] .
\end{aligned}
\end{equation}
\end{center}

It is clear that the reflection matrix $R$ is completely-$\calS$, since $R$ is $2\times2$, $R^{-1} = CMQ\Delta>0$ and all its diagonal elements are positive.\\

\noindent \textbf{Partition 1: $a = 1; b= 2$}
$$\theta_b = \theta_{2} = \alpha_1 -1,
    \;
  R_b^{-1} = R_{2}^{-1}= 1
  \Rightarrow u = - R_b^{-1}\theta_b = 1-\alpha_1 \Rightarrow 
  R_{ab}u = (\alpha_1-1)\displaystyle\frac{m_4}{m_3}.$$ Thus, $h_a^{'}\left[\theta_a + R_{ab}u\right] = h_1(\alpha_1-1) $. Since $\alpha_1<1$, taking $h_1 >0 $, we ensure $h_a^{'}\left[\theta_a + R_{ab}u\right] <0$.\\

\noindent \textbf{Partition 2: $a = 2; b= 1$}
$$\theta_b = \theta_{1} = (\alpha_1-1)(1-\displaystyle\frac{m_4}{m_3}),
    \;
  R_b^{-1} = R_{1}^{-1}= 1
  \Rightarrow u = - R_b^{-1}\theta_b = (1-\alpha_1)\left(1-\displaystyle\frac{m_4}{m_3}\right);$$  $$
   R_{ab} = 0 \Rightarrow 
  R_{ab}u = 0.$$ 
  \noindent Thus, $h_a^{'}\left[\theta_a + R_{ab}u\right] = h_2(\alpha_1-1) $. Since $\alpha_1<1$, taking $h_2 >0 $, we ensure $h_a^{'}\left[\theta_a + R_{ab}u\right] <0$.\\

\noindent \textbf{Partition 3: $a = 1, 2; b= \varnothing$}\\
$$h_a^{'}\left[\theta_a + R_{ab}u\right] = h'\theta = (\alpha_1-1)\left[\left(1-\displaystyle\frac{m_4}{m_3}\right)h_1 + h_2\right].$$
Since $\alpha_1<1$, fixing $h_2>0$ large enough, we guarantee that $h_a^{'}\left[\theta_a + R_{ab}u\right] <0$.
\\

\noindent We have shown, then, that exists a positive vector $h \in{R^2}$  as required so that $(R,\theta)$ under the static-priority policy 1, 2 is Chen-$\calS$.
\eProof\\
\end{proof}

\begin{proof} {Proof for priority 1, 3 (Chen-S).}
\label{proof: Proof of Chen-S for priority 1,3 (LuKumar)}

\noindent In the balanced Lu-Kumar network with class 1 having the priority in station 1 and class 3 having the priority in station 2, we have $J = 2, K = 4, \rho = \left[\alpha_1 \ \alpha_1\right]',$\\
$$R = 
 \left[
 \begin{matrix}
   1 & -m_4 \mcr
  0 & m_2
  \end{matrix}
  \right] ,
$$
\begin{center}
\begin{equation}\nonumber
\begin{aligned}
\theta = R(\rho-e)
= (\alpha_1-1)
\left[
 \begin{matrix}
   m_1 \mcr
   m_2
  \end{matrix}
  \right] .
\end{aligned}
\end{equation}
\end{center}

It is clear that the reflection matrix $R$ is completely-$\calS$, since $R$ is $2\times2$, $R^{-1} = CMQ\Delta>0$ and all its diagonal elements are positive.\\

\noindent \textbf{Partition 1: $a = 1; b= 2$}\\
$$\theta_b = \theta_{2} = (\alpha_1 -1)m_2,
    \;
  R_b^{-1} = R_{2}^{-1}= \displaystyle\frac{1}{m_2}
  \Rightarrow u = - R_b^{-1}\theta_b = 1-\alpha_1 \Rightarrow 
  R_{ab}u = (\alpha_1-1)m_4.$$
  
  \noindent Thus, $h_a^{'}\left[\theta_a + R_{ab}u\right] = h_1(\alpha_1-1) $. Since $\alpha_1<1$, taking $h_1 >0 $, we ensure $h_a^{'}\left[\theta_a + R_{ab}u\right] <0$.\\

\noindent \textbf{Partition 2: $a = 2; b= 1$}\\
$$\theta_b = \theta_{1} = (\alpha_1-1)m_1,
    \;
  R_b^{-1} = R_{1}^{-1}= 1
  \Rightarrow u = - R_b^{-1}\theta_b = (1-\alpha_1)m_1;$$
  $$
   R_{ab} = 0\Rightarrow 
  R_{ab}u = 0.$$
  
  \noindent Thus, $h_a^{'}\left[\theta_a + R_{ab}u\right] = h_2(\alpha_1-1)m_2$. Since $\alpha_1<1$, taking $h_2 >0 $, we ensure $h_a^{'}\left[\theta_a + R_{ab}u\right] <0$.\\

\noindent \textbf{Partition 3: $a = 1,2; b= \varnothing$}\\
$$h_a^{'}\left[\theta_a + R_{ab}u\right] = h'\theta = (\alpha_1-1)[m_1h_1+m_2h_2]$$
Since $\alpha_1<1$, taking $h_1>0$, $h_2>0$, we ensure $h_a^{'}\left[\theta_a + R_{ab}u\right] <0$.
\\

\noindent We have shown, then, that exists a positive vector $h \in{R^2}$  as required so that $(R,\theta)$ under the static-priority policy 1, 3 is Chen-$\calS$.
\eProof\\
\end{proof}

\begin{proof} {Proof for priority 4,2 (Chen-S).}
\label{proof: Proof of Chen-S for priority 4,2 (LuKumar)}

\noindent In the balanced Lu-Kumar network with class 4 having the priority in station 1 and class 2 having the priority in station 2, we have $J = 2, K = 4, \rho = \left[\alpha_1 \ \alpha_1\right]',$\\
$$R = \displaystyle\frac{1}{m_1-m_2}
 \left[
 \begin{matrix}
   m_1m_3 & -m_1m_4 \mcr
   -m_3 & m_3
  \end{matrix}
  \right] ,
$$
\begin{center}
\begin{equation}\nonumber
\begin{aligned}
\theta = R(\rho-e)
= (\alpha_1-1)
\left[
 \begin{matrix}
   m_1\mcr
   0
  \end{matrix}
  \right] .
\end{aligned}
\end{equation}
\end{center}

It is clear that the reflection matrix $R$ is completely-$\calS$ if $m_1-m_2>0$, equivalently $m_2+m_4 <1$, since $R$ is $2\times2$, $R^{-1} = CMQ\Delta>0$ and all its diagonal elements are positive when $m_1-m_2>0$.\\

\noindent \textbf{Partition 1: $a = 1; b= 2$}\\
$$\theta_b = \theta_{2} = 0,
    \;
  R_b^{-1} = R_{2}^{-1}= \displaystyle\frac{m_1-m_2}{m_3}
  \Rightarrow u = - R_b^{-1}\theta_b = 0;$$
  
  $$
   R_{ab} = \frac{m_1m_4}{m_2-m_1} \Rightarrow 
  R_{ab}u = 0.$$
  
  \noindent Thus, $h_a^{'}\left[\theta_a + R_{ab}u\right] = h_1(\alpha_1-1)m_1 $. Since $\alpha_1<1$, taking $h_1 >0 $, we ensure $h_a^{'}\left[\theta_a + R_{ab}u\right] <0$.\\

\noindent \textbf{Partition 2: $a = 2; b= 1$}\\
$$\theta_b = \theta_{1} = (\alpha_1-1)m_1,
    \;
  R_b^{-1} = R_{1}^{-1}= \frac{m_1-m_2}{m_1m_3}
  \Rightarrow u = - R_b^{-1}\theta_b = (\alpha_1-1)\frac{m_2-m_1}{m_3};$$
  
  $$
   R_{ab} = \frac{m_3}{m_2-m_1} \Rightarrow 
  R_{ab}u = \alpha_1-1.$$
  
  \noindent Thus, $h_a^{'}\left[\theta_a + R_{ab}u\right] = h_2(\alpha_1-1)$. Since $\alpha_1<1$, taking $h_2 >0 $, we ensure $h_a^{'}\left[\theta_a + R_{ab}u\right] <0$.\\

\noindent \textbf{Partition 3: $a = 1,2; b= \varnothing$}\\
$$h_a^{'}\left[\theta_a + R_{ab}u\right] = h'\theta = (\alpha_1-1)m_1h_1$$
Since $\alpha_1<1$, taking $h_1>0$, we ensure $h_a^{'}\left[\theta_a + R_{ab}u\right] <0$.
\\

\noindent We have shown, then, that exists a positive vector $h \in{R^2}$  as required so that $(R,\theta)$ under the static-priority policy 4, 2 is Chen-$\calS$.
\eProof\\
\end{proof}

\begin{proof} {Proof for priority 4, 3 (Chen-S).}
\label{proof: Proof of Chen-S for priority 4,3 (LuKumar)}

\noindent In the balanced Lu-Kumar network with class 4 having priority in station 1 and class 3 having priority in station 2, we have $J = 2$, $K = 4$, $\rho = \left[\alpha_1 \ \alpha_1\right]',$\\
$$R = 
 \left[
 \begin{matrix}
   1 & -m_4 \mcr
   -\displaystyle\frac{m_2}{m_1} & \displaystyle\frac{m_2}{m_1}
  \end{matrix}
  \right] ,
$$
\begin{center}
\begin{equation}\nonumber
\begin{aligned}
\theta = R(\rho-e)
= (\alpha_1-1)
\left[
 \begin{matrix}
   m_1 \mcr
   0
  \end{matrix}
  \right] .
\end{aligned}
\end{equation}
\end{center}

It is clear that the reflection matrix $R$ is completely-$\calS$, since $R$ is $2\times2$, $R^{-1} = CMQ\Delta>0$ and all its diagonal elements are positive.\\

\noindent \textbf{Partition 1: $a = 1; b= 2$}\\
$$\theta_b = \theta_{2} = 0,
    \;
  R_b^{-1} = R_{2}^{-1}= \displaystyle\frac{m_1}{m_2}
  \Rightarrow u = 0;$$
  $$
   R_{ab} = -m_4 \Rightarrow 
  R_{ab}u = 0.$$
  
  \noindent Thus, $h_a^{'}\left[\theta_a + R_{ab}u\right] = h_1(\alpha_1-1)m_1$. Since $\alpha_1<1$, taking $h_1 >0 $, we ensure $h_a^{'}\left[\theta_a + R_{ab}u\right] <0$.\\

\noindent \textbf{Partition 2: $a = 2; b= 1$}\\
$$\theta_b = \theta_{1} = (\alpha_1-1)m_1,
    \;
  R_b^{-1} = R_{1}^{-1}= 1
  \Rightarrow u = - R_b^{-1}\theta_b = (1-\alpha_1)m_1;$$
  
  $$
   R_{ab} = -\frac{m_2}{m_1} \Rightarrow 
  R_{ab}u = (\alpha_1-1)m_2.$$
  
  \noindent Thus, $h_a^{'}\left[\theta_a + R_{ab}u\right] = h_2(\alpha_1-1)m_2$. Since $\alpha_1<1$, taking $h_2 >0 $, we ensure $h_a^{'}\left[\theta_a + R_{ab}u\right] <0$.\\

\noindent \textbf{Partition 3: $a = 1,2; b= \varnothing$}\\
$$h_a^{'}\left[\theta_a + R_{ab}u\right] = h'\theta = h_1(\alpha_1-1)m_1$$ Since $\alpha_1<1$, taking $h_1>0$, we ensure $h_a^{'}\left[\theta_a + R_{ab}u\right] <0$. 
\\

\noindent We have shown, then, that exists a positive vector $h \in{R^2}$  as required so that $(R,\theta)$ under the static-priority policy 4, 3 is Chen-$\calS$.
\eProof\\
\end{proof}

\noindent We can now complete the proof that the 4 static-priority cases are all Chen-$\calS$ if $m_2+m_4 <1$.\\

\noindent The reflection matrices for the 4 static-priorty policies are invertible. We can calculate that $det(R_{1,2}) = 1>0$, $det(R_{1,3}) = m_2 >0$, $det(R_{4,2}) = \displaystyle\frac{m_1m_3}{1-m_4-m_2}$, $det(R_{4,3}) = m_2 >0$, where the subscript label indicates the classes with the priority in each station. These determinants have the same sign if $m_2+m_4<1$. We can conclude that the reflection matrix for a balanced push started Lu-Kumar network is invertible for any ratio matrix $\Delta$ provided that $m_2+m_4<1$.

According to Theorem \ref{Thm: full statement for corner->inter}, in a balanced Lu-Kumar network with $\alpha_1 <1$, the reflection matrix is invertible and $(R,\theta)$ is Chen-$\calS$ for any ratio matrix $\Delta$ if $m_2+m_4<1$.\\

\subsection{The SSC Inequalities for Balanced Lu-Kumar Networks}
\label{Appendix: The SSC Inequalities for Balanced LuKumar Networks}

\noindent In balanced Lu-Kumar networks, we have $m_1+m_4 = m_2+m_3 = 1$. We consider the 4 static priority policies of balanced Lu-Kumar networks for \eqref{equ: appendix_general_SSC_DHV} in the following analysis.


Given a policy $\pi$, we let $\dot{Y}^\pi(t) = \sum_{k \in \mathcal{H}_\pi:z_k>0} h_k\dot{\wbar{Z}^\pi_k}(t)$. We need $\dot{Y}^\pi(t) <0$ with a common vector $h \in \Bbb{R}^K_{++}$ for any regular time $t$ with $\lVert \wbar{Z}^\pi_{\mathcal{H}_\pi}(t) \rVert >0$, i.e., with at least one non-empty high-priority queue. We omit $\pi$ when it is fixed. In the balanced Lu-Kumar network, $m_k<1$ and thus $\mu_k>1$ for all $k \in \mathcal{K}$. Denote $I_k(t)$ as the input rate of class $k$ at time $t$; $O_k(t)$ as the output rate of class $k$ at time $t$. Then, $\dot{\wbar{Z}}_k(t) = I_k(t) - O_k(t) \text{ for } k\in \mathcal{K}; I_k(t) = O_{k-1}(t) \text{ for } k \in \mathcal{K}/\{1\};  I_1(t) = \alpha_1$. 

We use the set of high priority classes to represent static priority policies; (1,2) for example stands for a static priority policy that class 1 and 2 have the priorities in their respective stations. 

\vspace{0.2cm}

\noindent \textbf{Case 1:} (4, 2)
    
\vspace{0.2cm}

\begin{enumerate}
    \item $\wbar{Z}_2(t)>0, \wbar{Z}_4(t)>0.$


In this case, by the priority policy, at time $t$ server 1 serves queue 4 and server 2 serves queue 2 so that $O_4(t) = \mu_4$, $O_2(t) = \mu_2$, $O_1(t) = O_3(t) = 0$. Then, $\dot{\wbar{Z}}_2(t) = 0-\mu_2 = -\mu_2$ and
    $\dot{\wbar{Z}}_4(t) = 0-\mu_4 = -\mu_4$. It obviously holds that $\dot{Y}(t) = -h_2\mu_2 - h_4\mu_4<0$ as required.
    
    \vspace{0.2cm}
    
    \item $\wbar{Z}_2(t)>0, \wbar{Z}_4(t)=0$ 

In this scenario, server 2 serves queue 2 so that $O_2(t) = \mu_2$, $O_3(t) = 0$. Since $I_4(t)=O_3(t)  = 0$ and $\wbar{Z}_4(t)=0$, we have $O_4(t) =0$. 
    We have then that $\dot{\wbar{Z}}_4(t) = 0-0 = 0$. If $\wbar{Z}_1(t)=0$, $O_1(t) = \alpha_1$, 
    $\dot{\wbar{Z}}_2(t) = \alpha_1-\mu_2 <0$, thus, we have $\dot{Y}(t) = h_2(\alpha_1-\mu_2)<0$. If $\wbar{Z}_1(t)>0$, $O_1(t) = \mu_1$, 
    $\dot{\wbar{Z}}_2(t) = \mu_1-\mu_2$ and $\dot{Y}(t) = h_2(\mu_1-\mu_2)<0$, if $\mu_1-\mu_2 <0$, equivalently, $m_2< m_1 \Leftrightarrow m_2+m_4<1$. 

    \vspace{0.2cm}
    
    \item $\wbar{Z}_2(t)=0, \wbar{Z}_4(t)>0.$

    It is clear that server 1 serves queue 4, then $O_4(t) = \mu_4$ and $O_1(t) =0$. Since $\wbar{Z}_2(t)=0$, $O_2(t) = 0$, $\dot{\wbar{Z}}_2(t) = 0- 0  = 0$. If $\wbar{Z}_3(t)=0$, then since $O_2(t) = 0$, we have $O_3(t) = 0$ and $\dot{\wbar{Z}}_4(t) = 0- \mu_4  = -\mu_4 <0$, thus, it obviously holds that $\dot{Y}(t) = -\mu_4h_4<0$ as required. If $\wbar{Z}_3(t)>0$, then $O_3(t) = \mu_3$, and $\dot{\wbar{Z}}_4(t) = \mu_3- \mu_4$, thus, we need $\dot{Y}(t) = (\mu_3-\mu_4)h_4<0$, which is satisfied only when $\mu_3-\mu_4 <0$, equivalently, $m_4< m_3 \Leftrightarrow m_2+m_4<1$. 
   
    \vspace{0.2cm}
    \end{enumerate}  

We conclude that, in addition to the positivity of $h$, the condition $m_2+m_4<1$ arises from case 1. 

\vspace{0.2cm}

\noindent \textbf{Case 2:} (1, 2)
    
\vspace{0.2cm}

\begin{enumerate}
    \item $\wbar{Z}_1(t)>0, \wbar{Z}_2(t)>0.$

    In this scenario, server 1 serves queue 1, server 2 serves queue 2, so that $O_1(t) = \mu_1$, $O_2(t) = \mu_2$, $I_1(t) = \alpha_1$ and, in turn, $\dot{\wbar{Z}}_1(t) = \alpha_1-\mu_1 <0 $ and
    $\dot{\wbar{Z}}_2(t) = \mu_1-\mu_2$. Thus, we need $\dot{Y}(t) =  h_1(\alpha_1-\mu_1) +h_2(\mu_1-\mu_2) <0$. 

    \vspace{0.2cm}


    
\item $\wbar{Z}_1(t)>0, \wbar{Z}_2(t)=0.$ 

In this scenario, server 1 serves queue 1 which it prioritizes so that then $O_1(t) = \mu_1$, $\dot{\wbar{Z}}_1(t) = \alpha_1-\mu_1<0$.  If $\mu_1 > \mu_2$, then since $\wbar{Z}_2(t)=0$, $O_2(t) =\mu_2$,  $\dot{\wbar{Z}}_2(t) = \mu_1-\mu_2 >0$, which means that $t$ is not differentiable. If $\mu_1 \leq \mu_2$, then since $\wbar{Z}_2(t)=0$, $O_2(t) =\mu_1$,  $\dot{\wbar{Z}}_2(t) = \mu_1-\mu_1 =0$. $\dot{\wbar{Z}}_1(t) = \alpha_1-\mu_1 <0$. Thus, it obviously holds that $\dot{Y}(t) =  h_1(\alpha_1-\mu_1) <0$ as required.

    \vspace{0.2cm}
    
    \item $\wbar{Z}_1(t)=0, \wbar{Z}_2(t)>0.$

    It is clear that $\alpha_1 < \mu_1$, since $\wbar{Z}_1(t)=0$, we have $O_1(t) = \alpha_1$, then $\dot{\wbar{Z}}_1(t) = \alpha_1- \alpha_1  = 0$. Since $\wbar{Z}_2(t)>0$, server 2 serves queue 2, we have $O_2(t) = \mu_2$, then $\dot{\wbar{Z}}_2(t) = \alpha_1- \mu_2 < 0$. It obviously holds that $\dot{Y}(t) =  h_2(\alpha_1-\mu_2) < 0$ as required.

    \vspace{0.2cm}
    \end{enumerate}  

To sum up, in addition to the positivity of $h$,  we need the following conditions from case 2:
\begin{align*}
  &h_1(\alpha_1-\mu_1) +h_2(\mu_1-\mu_2) <0.
\end{align*}
 Note that when $m_2+m_4<1 \Leftrightarrow \mu_1-\mu_2<0$ (from Case 1), we have $h_1(\alpha_1-\mu_1) +h_2(\mu_1-\mu_2) <0$, thus, we do not need to consider any additional conditions from Case 2.

\vspace{0.2cm}

\noindent \textbf{Case 3:} (1, 3)
    
\vspace{0.2cm}

\begin{enumerate}
    \item $\wbar{Z}_1(t)>0, \wbar{Z}_3(t)>0.$

    In this case, $I_1(t) = \alpha_1$, server 1 serves queue 1, server 2 serves queue 3, then $O_1(t) = \mu_1$, $O_3(t) = \mu_3$, $O_2(t) = O_4(t) = 0$. So $\dot{\wbar{Z}}_1(t) = \alpha_1-\mu_1 <0$ and
    $\dot{\wbar{Z}}_3(t) = 0-\mu_3 = -\mu_3$. Thus, it obviously holds that $\dot{Y}(t) =  h_1(\alpha_1-\mu_1) - h_3\mu_3<0$ as required.
    
    \vspace{0.2cm}
    
    \item $\wbar{Z}_1(t)>0, \wbar{Z}_3(t)=0.$ 

    It is clear that server 1 serves queue 1, then $O_1(t) = \mu_1$,  $\dot{\wbar{Z}}_1(t) = \alpha_1- \mu_1 <0$. If $O_2(t) \geq \mu_3$, we have $O_3(t) = \mu_3$ since queue 3 has priority in station 2, then $\dot{\wbar{T}}_3(t) = 1 \Rightarrow \dot{\wbar{T}}_2(t) = 0 \Rightarrow O_2(t) = 0$, which is conflict with $O_2(t) \geq \mu_3$. As a result, $O_2(t) < \mu_3$, and then $O_3(t) = O_2(t)$, $\dot{\wbar{Z}}_3(t) = O_2(t)- O_3(t) =0$. It obviously holds that $\dot{Y}(t) =  h_1(\alpha_1-\mu_1)<0$ as required.

    \vspace{0.2cm}
    
    \item $\wbar{Z}_1(t)=0, \wbar{Z}_3(t)>0.$

    In this scenario, server 2 serves queue 3, then $O_3(t) = \mu_3$, $O_2(t) =0$, so $\dot{\wbar{Z}}_3(t) = 0- \mu_3 = -\mu_3$. Since $I_1(t)=\alpha_1 < \mu_1$ and $\wbar{Z}_1(t)=0$, we have $O_1(t) =\alpha_1$. 
    So $\dot{\wbar{Z}}_1(t) = \alpha_1-\alpha_1 = 0$. It obviously holds that $\dot{Y}(t)  = -\mu_3h_3<0$ as required.
    
    \vspace{0.2cm}
    \end{enumerate}     

To sum up, except for the positive condition for common vector $h$, we do not need any other conditions from Case 3.

\vspace{0.2cm}

\noindent \textbf{Case 4:} (4, 3)
    
\vspace{0.2cm}

\begin{enumerate}
    \item $\wbar{Z}_3(t)>0, \wbar{Z}_4(t)>0.$

    It is clear that server 1 serves queue 4, server 2 serves queue 3, then $O_4(t) = \mu_4$, $O_3(t) = \mu_3$, $O_1(t) = O_2(t) = 0$. So $\dot{\wbar{Z}}_3(t) = 0-\mu_3 = \mu_3$ and 
    $\dot{\wbar{Z}}_4(t) = \mu_3-\mu_4 $. Thus, we need $\dot{Y}(t) = -\mu_3h_3+ h_4(\mu_3-\mu_4) <0$.

    \vspace{0.2cm}
    
    \item $\wbar{Z}_3(t)=0, \wbar{Z}_4(t)>0.$ 

    In this case, server 1 serves queue 4, then $O_4(t) = \mu_4$, $O_1(t) = 0$. If $\wbar{Z}_2(t)=0$, since  $O_1(t) = 0$, we have $O_2(t) = 0$, then $O_3(t) = 0$ since $\wbar{Z}_3(t)=0$, thus, $\dot{\wbar{Z}}_3(t) = 0-0 = 0$ and 
    $\dot{\wbar{Z}}_4(t) = 0-\mu_4 =\mu_4$, and it obviously holds that $\dot{Y}(t) =-h_4\mu_4 <0$ as required. If $\wbar{Z}_2(t)>0$, to make $t$ differentiable, we have $\dot{\wbar{Z}}_3(t) = O_2(t)- O_3(t) =0 \Leftrightarrow \mu_2\dot{\wbar{T}}_2(t) -\mu_3\dot{\wbar{T}}_3(t) =0 \Leftrightarrow \mu_2\dot{\wbar{T}}_2(t) -\mu_3(1-\dot{\wbar{T}}_2(t)) =0 \Leftrightarrow (\mu_2+\mu_3)\dot{\wbar{T}}_2(t) =\mu_3 \Leftrightarrow \dot{\wbar{T}}_2(t) = \displaystyle\frac{\mu_3}{\mu_2+\mu_3} = \displaystyle\frac{\frac{1}{m_3}}{\frac{1}{m_2}+\frac{1}{m_3}} =\displaystyle\frac{m_2}{m_3+m_2} = m_2$. Thus, $O_3(t) = O_2(t) = \mu_2\dot{\wbar{T}}_2(t) = 1 $, $\dot{\wbar{Z}}_4(t) = 1-\mu_4 < 0$, and it obviously holds that $\dot{Y}(t) =h_4(1-\mu_4) <0$ as required. 

    \vspace{0.2cm}
    
    \item $\wbar{Z}_3(t)>0, \wbar{Z}_4(t)=0.$

    In this scenario, server 2 serves queue 3, then $O_3(t) = \mu_3$, $O_2(t) =  0$. Thus,  $\dot{\wbar{Z}}_3(t) = 0- \mu_3 =\mu_3$. To make $t$ differentiable, we have $\dot{\wbar{Z}}_4(t) = \mu_3- O_4(t) =0$. It obviously holds that $\dot{Y}(t) = -h_3\mu_3<0$ as required.
    
    \vspace{0.2cm}
    \end{enumerate} 

To sum up, except for the positive condition for common vector $h$, we need the following condition from Case 4:
\begin{align*}
 -\mu_3h_3+ h_4(\mu_3-\mu_4) <0.
\end{align*}

Note that when $m_2+m_4<1 \Leftrightarrow m_4<m_3 \Leftrightarrow \mu_3-\mu_4<0$ (from Case 1), we have $-\mu_3h_3+ h_4(\mu_3-\mu_4) <0$, thus, we do not need to consider any additional conditions from Case 4.\\

\end{appendices}
\end{document}